\documentclass[a4paper,12pt,frenchb]{article}

\usepackage{amsmath,amsbsy,amsfonts,amssymb,amsthm}
\usepackage[french]{babel}
\newtheorem*{lem}{Lemme}

\begin{document}

 \title{ Espaces $FC(\mathfrak{g}(F))$ et endoscopie }
\author{J.-L. Waldspurger}
\date{26 janvier 2020}
 \maketitle
 
 {\bf Abstract} Let $F$ be a $p$-adic field and let $G$ be a connected reductive group defined over $F$. We assume $p$ is big. Denote $\mathfrak{g}$ the Lie algebra of $G$.   We normalize suitably a Fourier-transform
  $f\mapsto \hat{f}$ on $C_{c}^{\infty}(\mathfrak{g}(F))$.  In a preceeding paper, we have defined the space $FC(\mathfrak{g}(F))$ of functions $f\in C_{c}^{\infty}(\mathfrak{g}(F))$ such that the orbital integrals  of $f$ and of $\hat{f}$  are $0$ for each element of $ \mathfrak{g}(F)$ which  is not topologically nilpotent. These  spaces are compatible with endoscopic transfer. We assume here that $G$ is absolutely quasi-simple and simply connected. We define a decomposition of  the space $FC(\mathfrak{g}(F))$ in a direct sum of subspaces such that the endoscopic transfer becomes (more or less) clear on each subspace. In particular, if $G$ is quasi-split, we describe the subspace $FC^{st}(\mathfrak{g}(F))$ of "stable" elements in $FC(\mathfrak{g}(F))$.

 \bigskip 

{\bf Introduction}

\bigskip

Soient $F$ une extension finie d'un corps ${\mathbb Q}_{p}$ et $G$ un groupe r\'eductif connexe d\'efini sur $F$. On impose que $p$ est grand relativement \`a $G$. On note $\mathfrak{g}$ l'alg\`ebre de Lie de $G$.  Introduisons l'espace $I(\mathfrak{g}(F))$ qui est le quotient de $C_{c}^{\infty}(\mathfrak{g}(F))$ par le sous-espace des fonctions $f$ dont  les int\'egrales orbitales $I^G(X,f)$ sont nulles pour tout $X\in \mathfrak{g}(F)$. Dans l'article \cite{W1}, on a introduit un sous-espace de dimension finie $FC(\mathfrak{g}(F))\subset I(\mathfrak{g}(F))$. En d\'efinissant une transformation de Fourier $f\mapsto \hat{f}$ dans $C_{c}^{\infty}(\mathfrak{g}(F))$ convenablement normalis\'ee, on peut le caract\'eriser de la fa\c{c}on  suivante. C'est l'image dans $I(\mathfrak{g}(F))$ du sous-espace des fonctions $f\in C_{c}^{\infty}(\mathfrak{g}(F))$ telles que $I^G(X,f)=I^G(X,\hat{f})=0$ pour tout $X\in \mathfrak{g}(F)$ qui n'est pas topologiquement nilpotent. On peut aussi le d\'ecrire \`a l'aide de l'immeuble de Bruhat-Tits et la th\'eorie des faisceaux-caract\`eres de Lusztig. 
Introduisons l'immeuble du groupe adjoint $G_{AD}$. Il est d\'ecompos\'e en facettes et on note $S(G)$ l'ensemble des sommets. A $s\in S(G)$, on associe un sous-groupe parahorique $K_{s}^0\subset G(F)$ et son plus grand sous-groupe distingu\'e pro-$p$-unipotent $K_{s}^+$. D'apr\`es Bruhat et Tits, il existe un groupe r\'eductif connexe $G_{s}$ d\'efini sur le corps r\'esiduel ${\mathbb F}_{q}$ tel que $K_{s}^0/K_{s}^+=G_{s}({\mathbb F}_{q})$.   Dans $\mathfrak{g}(F)$, on a de fa\c{c}on similaire une sous-$\mathfrak{o}_{F}$-alg\`ebre $\mathfrak{k}_{s}$ (o\`u $\mathfrak{o}_{F}$ est l'anneau des entiers de $F$) et une sous-$\mathfrak{o}_{F}$-alg\`ebre $\mathfrak{k}_{s}^+$, de sorte que $\mathfrak{k}_{s}/\mathfrak{k}_{s}^+=\mathfrak{g}_{s}({\mathbb F}_{q})$ (o\`u $\mathfrak{g}_{s}$ est l'alg\`ebre de Lie de $G_{s}$). Toute fonction d\'efinie sur $\mathfrak{g}_{s}({\mathbb F}_{q})$ se rel\`eve en une fonction sur $\mathfrak{k}_{s}$ et s'\'etend par $0$ hors de $\mathfrak{k}_{s}$ en une fonction d\'efinie sur $\mathfrak{g}(F)$.   Notons $FC(\mathfrak{g}_{s}({\mathbb F}_{q}))$ l'espace engendr\'e par les fonctions caract\'eristiques des faisceaux-caract\`eres d\'efinis sur $\mathfrak{g}_{s}$ qui sont cuspidaux, \`a support nilpotent et invariants par l'action de Frobenius. Par le proc\'ed\'e ci-dessus, on consid\`ere  $FC(\mathfrak{g}_{s}({\mathbb F}_{q}))$ comme un sous-espace de $C_{c}^{\infty}(\mathfrak{g}(F))$. Alors $FC(\mathfrak{g}(F))$ est l'image dans $I(\mathfrak{g}(F))$ de la somme de ces sous-espaces $FC(\mathfrak{g}_{s}({\mathbb F}_{q}))$ sur tous les sommets $s\in S(G)$. L'int\'er\^et de l'espace $FC(\mathfrak{g}(F))$ est qu'il nous semble jouer un r\^ole crucial dans deux probl\`emes concernant la th\'eorie de l'endoscopie de langlands: d'une part, les propri\'et\'es endoscopiques de l'espace des distributions invariantes \`a support nilpotent, cf. \cite{W2} et \cite{DK}; d'autre part, la d\'etermination des paquets stables de repr\'esentations de niveau $0$, cf. \cite{MW2}.

Dans cet article, nous supposons que $G$ est absolument quasi-simple et simplement connexe. Nous nous proposons d'\'eclaicir partiellement les propri\'et\'es de l'espace $FC(\mathfrak{g}(F))$ relatives \`a l'endoscopie. 

Il est facile de d\'ecrire une base de l'espace $FC(\mathfrak{g}(F))$ (ou plus exactement d'un sous-espace de $C_{c}^{\infty}(\mathfrak{g}(F))$ qui s'envoie bijectivement sur $FC(\mathfrak{g}(F))$). On fixe un sous-ensemble $\underline{S}(G)\subset S(G)$ de repr\'esentants des orbites pour l'action de $G(F)$ dans $S(G)$. Pour tout $s\in \underline{S}(G)$, on d\'ecrit les faisceaux-caract\`eres d\'efinis sur $\mathfrak{g}_{s}$ qui sont cuspidaux, \`a support nilpotent et invariants par l'action de Frobenius. Alors, quand $s$ d\'ecrit $FC(\mathfrak{g}(F))$, les fonctions caract\'eristiques de ces faisceaux  (identifi\'ees comme ci-dessus \`a des \'el\'ements de $C_{c}^{\infty}(\mathfrak{g}(F))$) forment une base de $FC(\mathfrak{g}(F))$. Remarquons que l'on tombe tout de suite sur une difficult\'e calculatoire: il n'y a pas de normalisation canonique des fonctions caract\'eristiques en question, notre base n'est bien d\'efinie qu'\`a des scalaires pr\`es. En fait, il est plus opportun de modifier la base ainsi d\'ecrite en tenant compte de l'action de $G_{AD}(F)$ sur l'immeuble, c'est-\`a-dire en regroupant les sommets qui sont reli\'es par cette action. On d\'ecrira une telle base sur laquelle l'action de $G_{AD}(F)$ se lit bien. Pr\'ecis\'ement, pour chaque groupe $G$, on d\'efinira un ensemble fini ${\cal X}$ de nature combinatoire. Pour chaque $x\in {\cal X}$, on d\'efinira un sous-espace $FC_{x}\subset FC(\mathfrak{g}(F))$ qui poss\`ede une base form\'ee de combinaisons lin\'eaires simples des \'el\'ements de base d\'ecrits ci-dessus. On prouvera que 
$$(1) \qquad FC(\mathfrak{g}(F))=\oplus_{x\in {\cal X}}FC_{x}.$$
 Dans de nombreux cas, $FC_{x}$ sera une droite mais il y a aussi des cas o\`u cet espace sera de dimension strictement sup\'erieure \`a $1$.

Supposons un instant que $G$ soit quasi-d\'eploy\'e. On d\'efinit l'espace $SI(\mathfrak{g}(F))$ comme le quotient de $C_{c}^{\infty}(\mathfrak{g}(F))$ par le sous-espace des fonctions $f$ dont les orbitales int\'egrales stables $S^G(X,f)$ sont nulles pour tout $X\in \mathfrak{g}(F)$ semi-simple r\'egulier. Levons l'hypoth\`ese de quasi-d\'eploiement. 
Introduisons un ensemble $Endo_{ell}(G)$ de repr\'esentants des classes d'\'equivalence de donn\'ees endoscopiques elliptiques de $G$. A une donn\'ee endoscopique ${\bf G}'$  est associ\'ee  un groupe endoscopique $G'$ qui est quasi-d\'eploy\'e sur $F$ (soulignons qu'une donn\'ee endoscopique est plus riche que la simple donn\'ee de ce groupe $G'$, deux donn\'ees non \'equivalentes pouvant avoir le m\^eme groupe associ\'e).  Soit ${\bf G}'\in Endo_{ell}(G)$. Modulo le choix d'un facteur de transfert, on d\'efinit l'application de transfert $transfert^{{\bf G}'}:I(\mathfrak{g}(F))\to SI(\mathfrak{g}'(F))$. On note $FC(\mathfrak{g}(F),{\bf G}')$ le sous-espace des $f\in FC(\mathfrak{g}(F))$ telles que $transfert^{{\bf G}''}(f)=0$ pour tout ${\bf G}''\in Endo_{ell}(G)$, ${\bf G}''\not={\bf G}'$. Il y a toujours une unique donn\'ee endoscopique elliptique dont le groupe associ\'e est une forme int\'erieure quasi-d\'eploy\'ee de $G$, on la note ${\bf G}$. Dans le cas o\`u $G$ est quasi-d\'eploy\'e, on pose $FC^{st}(\mathfrak{g}(F))=FC(\mathfrak{g}(F),{\bf G})$. Posons 
$$(2) \qquad FC^{{\cal E}}(\mathfrak{g}(F))=\oplus_{{\bf G}'\in Endo_{ell}(G)}FC^{st}(\mathfrak{g}'(F))^{Out({\bf G}')}.$$
Le groupe $Out({\bf G}')$ est le groupe fini des automorphismes "ext\'erieurs" de la donn\'ee ${\bf G}'$, il agit naturellement sur $FC^{st}(\mathfrak{g}'(F))$ et $FC^{st}(\mathfrak{g}'(F))^{Out({\bf G}')}$ est le sous-espace des invariants. 

On a prouv\'e dans \cite{W1} que
$$FC(\mathfrak{g}(F))=\oplus_{{\bf G}'\in Endo_{ell}(G)}FC(\mathfrak{g}(F),{\bf G}'),$$
et que l'application $transfert=\oplus_{{\bf G}'\in Endo_{ell}(G)}transfert^{{\bf G}'}$ \'etablissait un isomorphisme
$$(3)\qquad FC(\mathfrak{g}(F))\to FC^{{\cal E}}(\mathfrak{g}(F)).$$ 
On  d\'efinira un ensemble fini ${\cal Y}$ de nature combinatoire. Pour chaque $y\in {\cal Y}$, on d\'efinira un sous-espace $FC_{y}^{{\cal E}}\subset FC^{{\cal E}}(\mathfrak{g}(F))$. On d\'efinira une bijection $\varphi:{\cal X}\to {\cal Y}$ et  notre r\'esultat principal sera que, pour tout $x\in {\cal X}$, on a l'\'egalit\'e 
$$(4) \qquad transfert(FC_{x})=FC^{{\cal E}}_{\varphi(x)}.$$ 

On est conscient que la pr\'esentation ci-dessus est tr\`es incompl\`ete: pour obtenir ces r\'esultats, il suffirait de poser ${\cal X}={\cal Y}=\{1\}$, $FC_{1}=FC(\mathfrak{g}(F))$ et $FC_{1}^{{\cal E}}=FC^{{\cal E}}(\mathfrak{g}(F))$. Mais le lecteur constatera que nos d\'ecompositions des espaces $FC(\mathfrak{g}(F))$ et $FC^{{\cal E}}(\mathfrak{g}(F))$ sont tout-\`a-fait non triviales. Dans presque tous les cas, chaque sous-espace $FC^{st}(\mathfrak{g}'(F))^{Out({\bf G}')}$ de $FC^{{\cal E}}(\mathfrak{g}(F))$ est une somme de sous-espace $FC^{{\cal E}}_{y}$ (il y a une exception dans le cas o\`u $G$ est quasi--d\'eploy\'e et non d\'eploy\'e de type $D_{n}$ avec $n$ pair, cf. \ref{Dnpairpadiquenonram}; dans ce cas il y a des paires de donn\'ees endoscopiques que nous n'avons pas r\'eussi \`a distinguer).  En particulier, quand $G$ est quasi-d\'eploy\'e, nos r\'esultats d\'ecrivent dans tous les cas le sous-espace $FC^{st}(\mathfrak{g}(F))$ comme la somme des $FC_{x}$ pour $x$ parcourant un sous-ensemble ${\cal X}^{st}$ de ${\cal X}$. A titre d'exemple, donnons une cons\'equence dans le cas o\`u $G$ est de type $E_{8}$. Dans ce cas, on a l'\'egalit\'e $FC(\mathfrak{g}(F))=FC^{st}(\mathfrak{g}(F))$ et
$$dim(FC(\mathfrak{g}(F)))=3+4\delta_{3}(q-1)+2\delta_{4}(q-1)+4\delta_{5}(q-1),$$
o\`u, pour tout entier $n\geq1$, on  a not\'e $\delta_{n}$ la fonction caract\'eristique de $n{\mathbb Z}$. Remarquons que, si $q-1$ est divisible par $60$, la dimension ci-dessus est le nombre de repr\'esentations unipotentes cuspidales du groupe sur ${\mathbb F}_{q}$ de type $E_{8}$. 

Indiquons les grandes lignes de la d\'emonstration en nous restreignant au cas quasi-d\'eploy\'e. On d\'ecrit facilement l'ensemble $Endo_{ell}(G)$, cf. \ref{description}.  En raisonnant par r\'ecurrence sur le rang de $G$, on peut d\'ecrire les espaces $FC^{st}(\mathfrak{g}'(F))$ pour toute donn\'ee ${\bf G}'\in Endo_{ell}(G)$, ${\bf G}'\not={\bf G}$. L'action du groupe $Out({\bf G}')$ se lit bien dans nos descritions et on en d\'eduit l'espace $FC^{st}(\mathfrak{g}'(F))^{Out({\bf G}')}$. Puisque l'on a d\'ej\`a d\'ecrit l'espace $FC(\mathfrak{g}(F))$ et que l'on conna\^{\i}t donc sa dimension, l'isomorphisme (3) est alors suffisant pour d\'eterminer la dimension du sous-espace $FC^{st}(\mathfrak{g}(F))$, qui est le seul sous-espace de $FC^{{\cal E}}(\mathfrak{g}(F))$ que la r\'ecurrence ne permet pas de conna\^{\i}tre. Pour d\'eterminer la bijection $\varphi$ et pour prouver l'\'egalit\'e principale  (4), on utilise les trois arguments ci-dessous. 

D'abord l'action du groupe adjoint $G_{AD}(F)$. Toute donn\'ee endoscopique ${\bf G}'\in Endo_{ell}(G)$ d\'etermine un caract\`ere $\xi_{{\bf G}'}$  de ce groupe, trivial sur l'image naturelle de $G(F)$, et les \'el\'ements de $FC(\mathfrak{g}(F),{\bf G}')$ se transforment par $G_{AD}(F)$ selon ce caract\`ere. Dans beaucoup de cas, on pourra aussi associer \`a tout $x\in {\cal X}$ un  tel caract\`ere $\xi_{x}$ de sorte que les \'el\'ements de $FC_{x}$ se transforment par $G_{AD}(F)$ selon ce caract\`ere. Dans ce cas, on aura forc\'ement $FC_{x}\subset_{{\bf G}'\in Endo_{ell}(G),\xi_{{\bf G}'}=\xi_{x}}FC(\mathfrak{g}(F),{\bf G}')$. Remarquons que cet argument est suffisant dans le cas d'un groupe d\'eploy\'e de type $A_{n-1}$, les caract\`eres de $G_{AD}(F)$ \'etant alors suffisants pour distinguer des diff\'erentes droites $FC_{x}$ et les diff\'erentes donn\'ees ${\bf G}'$. 

Moy et Prasad ont d\'efini des $\mathfrak{o}_{F}$-r\'eseaux $\mathfrak{k}_{x,r}\subset \mathfrak{g}(F)$ o\`u $x$ est un point de l'immeuble de $G_{AD}$ et $r$ est un  nombre r\'eel. Ils sont d\'ecroissants en $r$. Pour un \'el\'ement $X\in \mathfrak{g}(F)$, on appelle la profondeur de $X$  la borne sup\'erieure (\'eventuellement infinie) de l'ensemble des r\'eels $r$ tels qu'il existe un point $x$ de sorte que $X\in \mathfrak{k}_{x,r}$. Si $X$ est semi-simple r\'egulier, cette borne est finie et est atteinte. On note $r(X)$ cette profondeur. Soit $y\in {\cal Y}$. Supposons pour simplifier que $FC^{{\cal E}}_{y}$ soit une droite \'egale \`a $FC^{st}(\mathfrak{g}'(F))^{Out({\bf G}')}$ pour une donn\'ee ${\bf G}'\in Endo_{ell}(G)$ (c'est souvent le cas). Fixons un \'el\'ement non nul $f'_{y}\in FC^{{\cal E}}_{y}$. Supposons que l'on connaisse un \'el\'ement $Y_{y}\in \mathfrak{g}'(F)$ tel que $S^{G'}(Y_{y},f'_{y})\not=0$ et que l'on connaisse la profondeur $r(Y_{y})$. Soit $x\in {\cal X}$. Supposons que l'on connaisse un r\'eel $r_{x}$ tel  que le support de tout \'el\'ement de $FC_{x}$ soit  form\'e d'\'el\'ements  $X\in \mathfrak{g}(F)$ tels que $r(X)\geq r_{x}$. Si $r_{x}> r_{Y_{y}}$, alors $transfert^{{\bf G}'}(f)=0$ pour tout $f\in FC_{x}$. En effet, \'ecrivons $transfert^{{\bf G}'}(f)=c f'_{y}$ avec $c\in {\mathbb C}$. On a alors $S^{G'}(Y_{y},transfert^{{\bf G}'}(f))=cS^{G'}(Y_{y},f'_{y})$. Par d\'efinition du transfert endoscopique, le membre de gauche est combinaison lin\'eaire d'int\'egrales orbitales $I^G(X,f)$ pour des \'el\'ements $X\in \mathfrak{g}(F)$ correspondant \`a $Y_{y}$. Un tel \'el\'ement a la m\^eme profondeur que $Y_{y}$, c'est-\`a-dire $r(Y_{y})$. Puisque $r_{x}> r(Y_{y})$, la classe de conjugaison de $X$ ne coupe pas le support de $f$ donc $I^G(X,f)=0$. Donc $0=S^{G'}(Y_{y},transfert^{{\bf G}'}(f))=cS^{G'}(Y_{y},f'_{y})$, d'o\`u $c=0$. 

Le troisi\`eme argument est que l'espace $FC(\mathfrak{g}(F))$ est contenu dans celui des fonctions cuspidales $I_{cusp}(\mathfrak{g}(F))\subset I(\mathfrak{g}(F))$ et est muni du produit elliptique, qui est un produit hermitien d\'efini positif. La d\'ecomposition (1) est orthogonale pour ce produit. Via l'isomorphisme (3), l'espace $FC^{{\cal E}}(\mathfrak{g}(F))$ est aussi muni d'un tel produit et la d\'ecomposition (2) est orthogonale. Cela a de multiples cons\'equences. Par exemple, si on r\'esussit \`a d\'eterminer l'image de $FC^{st}(\mathfrak{g}'(F))^{Out({\bf G}')}$ dans $FC(\mathfrak{g}(F))$ pour tout ${\bf G}'\not={\bf G}$, alors l'espace $FC^{st}(\mathfrak{g}(F))$ est lui-aussi d\'etermin\'e: c'est l'orthogonal de la somme des espaces pr\'ec\'edents pour ${\bf G}'\not={\bf G}$. Ou encore si le deuxi\`eme argument est valide pour tout couple $(x,y)\in {\cal X}\times {\cal Y}$, celui-ci  nous dit plus ou moins que l'isomorphisme (3) est triangulaire pour des ordres convenables sur ${\cal X}$ et ${\cal Y}$. Puisqu'il est hermitien, il est forc\'ement diagonal. On renvoie \`a \ref{ingredients} pour une \'elaboration plus pr\'ecise de cet argument. 

Le premier paragraphe est consacr\'e \`a divers pr\'eliminaires sur les immeubles et les donn\'ees endoscopiques. On y d\'emontre aussi quelques lemmes galoisiens \'el\'ementaires. Le deuxi\`eme paragraphe rappelle la description due \`a Lusztig des faisceaux-caract\`eres cuspidaux \`a support unipotent pour les groupes finis. Au troisi\`eme paragraphe, on pr\'esente la fa\c{c}on dont nous d\'ecrirons ensuite nos r\'esultats dans les diff\'erents cas. Les paragraphes 4 et 5 sont consacr\'es aux groupes de type $A_{n-1}$ (on a pr\'ef\'er\'e $A_{n-1}$ \`a $A_{n}$ car on utilise beaucoup d'alg\`ebre lin\'eaire et on pr\'ef\`ere travailler avec $GL(n)$ plut\^ot qu'avec $GL(n+1)$). En particulier, on traite en d\'etail le cas des groupes unitaires associ\'es \`a une extension quadratique $E/F$ ramifi\'ee. Pour les groupes classiques, c'est le r\'esultat le plus nouveau de notre \'etude, les autres cas \'etant largement contenus dans \cite{W2}. Les groupes de type $B_{n}$, $C_{n}$ ou $D_{n}$ sont consid\'er\'es dans les chapitres 6 et 7. Pour la raison que l'on vient de donner, on s'est plusieurs fois content\'e d'indiquer bri\`evement les constructions n\'ecessaires, en laissant les d\'etails au lecteur. Les groupes exceptionnels (y compris celui de $D_{4}$ trialitaire qui est de fait exceptionnel) sont trait\'es au chapitre 8. Ils ne sont d'ailleurs pas les plus difficiles. 

Nos r\'esultats sont complets quant \`a la d\'etermination des espaces $FC^{st}(\mathfrak{g}(F))$ pour les groupes $G$ quasi-d\'eploy\'es. Les r\'esultats de transfert sont beaucoup plus impr\'ecis. Pour obtenir une description exacte, il faudrait am\'eliorer trois points. D'abord, comme on l'a dit ci-dessus, d\'efinir pr\'ecis\'ement les normalisations des fonctions caract\'eristiques des faisceaux-caract\`eres. Cette question nous para\^{\i}t li\'ee \`a la conjecture de Lusztig reliant ces fonctions aux caract\`eres de repr\'esentations de groupes finis. D'autre part, calculer plus pr\'ecis\'ement diverses int\'egrales orbitales dont nous nous contentons de prouver la non-nullit\'e. Enfin calculer explicitement divers facteurs de transfert. Pour les groupes classiques, c'est une simple question de patience mais l'auteur avoue se sentir un peu d\'emuni dans le cas des groupes exceptionnels.

 \section{Pr\'eliminaires}
 
 \subsection{Notations}\label{notations}
 
 Soit $F$ un corps local non-archim\'edien de caract\'eristique nulle. On note  ${\mathbb F}_{q}$ le corps r\'esiduel de $F$, $q$ \'etant son nombre d'\'el\'ements, $p$ la caract\'eristique  de ${\mathbb F}_{q}$, $\mathfrak{o}_{F}$ l'anneau d'entiers de $F$, $\mathfrak{p}_{F}$ son id\'eal maximal.  On fixe une uniformisante $\varpi_{F}$ de $\mathfrak{p}_{F}$.   
 
 Pour $k=F$ ou $k={\mathbb F}_{q}$, on fixe une
 cl\^oture alg\'ebrique $\bar{k}$ de $k$.   Toutes les extensions de $k$ que l'on consid\'erera seront suppos\'ees incluses dans $\bar{k}$.   On note $\Gamma_{k}$ le groupe de Galois de $\bar{k}/k$. Pour toute extension galoisienne $k'$ de $k$, on note $\Gamma_{k'/k}$ le groupe de Galois de $k'/k$.  On note $F^{nr}$ l'extension non ramifi\'ee maximale de $F$ dont le corps r\'esiduel s'identifie \`a $\bar{{\mathbb F}}_{q}$. On note $I_{F}$ le sous-groupe d'inertie de $\Gamma_{F}$, c'est-\`a-dire le groupe $\Gamma_{F^{nr}}$. On pose simplement $\Gamma_{F}^{nr}=\Gamma_{F^{nr}/F}\simeq \Gamma_{{\mathbb F}_{q}}$. 
 On note $Fr$ l'\'el\'ement de Frobenius de  l'un ou l'autre de ces deux  groupes. On note $W_{F}$ le groupe de Weil de $F$, c'est-\`a-dire le sous-groupe des \'el\'ements de $\Gamma_{F}$ dont l'image dans $\Gamma_{F}^{nr}$ est une puissance enti\`ere du Frobenius.
 
 On note $val_{F}$ la valuation usuelle de $F$ et on la prolonge \`a $\bar{F}$ en une valuation \`a valeurs dans ${\mathbb Q}$.
 
 Pour tout groupe ab\'elien $A$, on note $A^{\vee}=Hom(A,{\mathbb C}^{\times})$ son groupe des caract\`eres. Consid\'erons un groupe $G$ agissant sur un ensemble $U$. Pour $u\in U$, on note $Z_{G}(u)$ le fixateur de $u$ dans $G$. Pour un sous-ensemble $V\subset U$, on note $Norm_{G}(V)$ le stabilisateur de $V$ dans $G$.

 Soit  encore $k=F$ ou $k={\mathbb F}_{q}$. Pour tout groupe alg\'ebrique $H$ d\'efini sur $k$, on note $H^0$ sa composante neutre et $Z(H)$ son centre. Pour tout tore $T$ d\'efini sur $k$, on note $X_{*}(T)$, resp. $X^*(T)$, son groupe de cocaract\`eres alg\'ebriques, resp. caract\`eres alg\'ebriques, d\'efinis sur $\bar{k}$. 
 
 Soit $G$ un groupe r\'eductif connexe d\'efini sur $k$. On note   $G_{AD}$ le groupe adjoint et $G_{SC}$ le rev\^etement simplement connexe de $G_{AD}$. On note $\pi:G\to G_{AD}$ la projection naturelle. Pour un sous-groupe alg\'ebrique $H$ de $G$, on note $H_{ad}=\pi(H)$ son image dans $G_{AD}$.
 On identifie $G$ \`a $G(\bar{k})$. Ce groupe et beaucoup d'objets qui lui sont reli\'es sont alors munis d'une action de $\Gamma_{k}$ que l'on note $\sigma\mapsto \sigma_{G}$. On note $\mathfrak{g}$ l'alg\`ebre de Lie de $G$ et on appelle conjugaison par $G$ l'action adjointe de $G$ sur $\mathfrak{g}$. On note $\mathfrak{g}_{reg}$ l'ensemble des \'el\'ements semi-simples r\'eguliers de $\mathfrak{g}$. On note $\mathfrak{g}_{ell}(F)$ le sous-ensemble des \'el\'ements elliptiques de $\mathfrak{g}_{reg}(F)$. Pour $X\in \mathfrak{g}$, on note  $G_{X}$ la composante neutre de $Z_{G}(X)$. On note $h(G)$ le plus grand des nombres de Coxeter des composantes irr\'eductibles sur $\bar{k}$ du groupe $G_{AD}$ et on note $rg(G)$ le rang de $G$ sur $\bar{k}$. 
 
 Dans chaque section de l'article, on consid\'erera un tel groupe $G$ d\'efini sur $F$ ou ${\mathbb F}_{q}$. On supposera toujours que
 
 $p\geq sup(2h(G)+1,rg(G)+2)$.
 
 On note ${\mathbb N}_{>0}={\mathbb N}-\{0\}$. Pour tout $n\in {\mathbb N}_{>0}$, on note $\delta_{n}$ la fonction caract\'eristique du sous-groupe $n{\mathbb Z}$ de ${\mathbb Z}$. Pour tout corps commutatif $k$ et tout $n\in {\mathbb N}_{>0}$, on note $\boldsymbol{\zeta}_{n}(k)$ le groupe des racines $n$-i\`emes de l'unit\'e dans $k^{\times}$ et $\boldsymbol{\zeta}_{n,prim}(k)$ le sous-ensemble des racines primitives d'ordre $n$.

 \subsection{Groupes sur ${\mathbb F}_{q}$}\label{groupessurFq}
 
 Soit $G$ un groupe r\'eductif connexe d\'efini sur ${\mathbb F}_{q}$. On note ${\cal C}(\mathfrak{g}({\mathbb F}_{q}))$ l'espace des fonctions sur $\mathfrak{g}({\mathbb F}_{q})$, \`a valeurs complexes.  Il est muni du produit hermitien non d\'eg\'en\'er\'e
$$(f',f)=\vert G({\mathbb F}_{q})\vert ^{-1}\sum_{X\in \mathfrak{g}({\mathbb F}_{q})}\bar{f}'(X)f(X).$$
On note  $C(\mathfrak{g}({\mathbb F}_{q}))$ le sous-espace des fonctions qui sont invariantes par conjugaison par $G({\mathbb F}_{q})$. 
 On note $C_{cusp}(\mathfrak{g}({\mathbb F}_{q}))$ le sous-espace  des fonctions cuspidales, c'est-\`a-dire des fonctions $f\in C(\mathfrak{g}({\mathbb F}_{q}))$ qui v\'erifient la condition suivante. Soient $P$ un sous-groupe parabolique propre de $G$ et $M$ une composante de Levi de $P$, tous deux d\'efinis sur ${\mathbb F}_{q}$; on note $U_{P}$ le radical unipotent de $P$.  Alors, pour tout $X\in \mathfrak{m}({\mathbb F}_{q})$, on a l'\'egalit\'e
$$\sum_{Y\in \mathfrak{u}_{P}({\mathbb F}_{q})}f(X+Y)=0.$$
Fixons une forme bilin\'eaire sym\'etrique non d\'eg\'en\'er\'ee $<.,.>$ sur $\mathfrak{g}({\mathbb F}_{q})$ invariante par conjugaison par $G({\mathbb F}_{q})$. Fixons aussi un caract\`ere non trivial $\psi:{\mathbb F}_{q}\to {\mathbb C}^{\times}$. On d\'efinit la transformation de Fourier $f\mapsto \hat{f}$ dans ${\cal C}(\mathfrak{g}({\mathbb F}_{q}))$ par
$$\hat{f}(X)=q^{-dim(\mathfrak{g})/2}\sum_{Y\in \mathfrak{g}({\mathbb F}_{q})}f(Y)\psi(<X,Y>)$$
pour tout $X\in \mathfrak{g}({\mathbb F}_{q})$. 

On note $C_{nil}(\mathfrak{g}({\mathbb F}_{q}))$ le sous-espace de $C(\mathfrak{g}({\mathbb F}_{q}))$ form\'e des fonctions \`a support nilpotent. On pose $C_{nil,cusp}(\mathfrak{g}({\mathbb F}_{q}))= C_{nil}(\mathfrak{g}({\mathbb F}_{q}))\cap C_{cusp}(\mathfrak{g}({\mathbb F}_{q}))$. On note $FC(\mathfrak{g}({\mathbb F}_{q}))$ le sous-espace des fonctions $f\in C_{nil}(\mathfrak{g}({\mathbb F}_{q}))$ telles que $\hat{f}$ appartient elle-aussi \`a $C_{nil}(\mathfrak{g}({\mathbb F}_{q}))$. 
On a vu en \cite{W1} 2(2)  que

(1) si $Z(G)^0\not=\{1\}$,  $FC(\mathfrak{g}({\mathbb F}_{q}))=0$.

Supposons donc $Z(G)^{0}=\{1\}$, c'est-\`a-dire $G$ semi-simple. Lusztig a prouv\'e en \cite{L} paragraphe 11 que $FC(\mathfrak{g}({\mathbb F}_{q}))$ avait pour base les fonctions caract\'eristiques des faisceaux-caract\`eres cuspidaux qui sont invariants par $\Gamma_{{\mathbb F}_{q}}$, consid\'er\'ees \`a homoth\'etie pr\`es.    Notons $fc(\mathfrak{g}({\mathbb F}_{q}))$  cet ensemble de fonctions. Ce r\'esultat implique que   $FC(\mathfrak{g}({\mathbb F}_{q}))\subset C_{nil, cusp}(\mathfrak{g}({\mathbb F}_{q}))$.

Rappelons quelques propri\'et\'es d'une telle fonction caract\'eristique $f\in fc(\mathfrak{g}({\mathbb F}_{q}))$. Consid\'erons un \'el\'ement $N\in \mathfrak{g}_{nil}({\mathbb F}_{q})$ et appelons orbite g\'eom\'etrique de $N$ l'ensemble des $N'\in \mathfrak{g}({\mathbb F}_{q})$ tels qu'il existe $g\in G$ de sorte que $N'=g^{-1}Ng$. Le groupe $\Gamma_{{\mathbb F}_{q}}$ agit sur $Z_{G}(N)/Z_{G}(N)^0$. Soit $\epsilon$ une repr\'esentation irr\'eductible de $Z_{G}(N)/Z_{G}(N)^0$  dont la classe est fixe par l'action galoisienne. On fixe un automorphisme $Fr_{\epsilon}$ de l'espace de cette repr\'esentation de sorte que $\epsilon(Fr(g))=Fr_{\epsilon}\epsilon(g)Fr_{\epsilon}^{-1}$ pour tout $g\in Z_{G}(N)/Z_{G}(N)^0$. On d\'efinit une fonction $f_{N,\epsilon,Fr_{\epsilon}}$ \`a support dans l'orbite g\'eom\'etrique de $N$ par la formule $f_{N,\epsilon,Fr_{\epsilon}}(N')=trace(\epsilon(gFr(g)^{-1})Fr_{\epsilon})$ pour tous $N'$ et $g$ comme ci-dessus. Alors, pour toute fonction   $f\in fc(\mathfrak{g}({\mathbb F}_{q}))$, il existe $N$, $\epsilon$ et $Fr_{\epsilon}$ de sorte que $f=f_{N,\epsilon,Fr_{\epsilon}}$. Le choix de l'automorphisme $Fr_{\epsilon}$ ne change la fonction que par multiplication par un scalaire. On supposera que $Fr_{\epsilon}$ est d'ordre fini et que c'est l'identit\'e quand $\epsilon$ est de degr\'e $1$, ce qui est le cas le plus fr\'equent. On oubliera le $Fr_{\epsilon}$ de la notation.   Consid\'erons une fonction caract\'eristique $f_{N,\epsilon}$ comme ci-dessus. L'orbite de $N$ est distingu\'ee. Fixons un $\mathfrak{sl}_{2}$-triplet $(N^-,H,N)$. On a donc $[H,N]=2N$, $[H,N^-]=-2N^-$ et $[N,N^-]=H$. Pour tout $i\in {\mathbb Z}$, posons $\mathfrak{g}_{i}=\{X\in \mathfrak{g}; [H,X]=iX\}$ et $\mathfrak{g}_{\geq i}=\oplus_{j\geq i}\mathfrak{g}_{j}$. Alors $\mathfrak{g}_{i}=\{0\}$ pour $i$ impair et on a $dim_{\bar{{\mathbb F}}_{q}}(\mathfrak{g}_{0})=dim_{\bar{{\mathbb F}}_{q}}(\mathfrak{g}_{2})$. Notons $P$ le sous-groupe parabolique de $G$ d'alg\`ebre de Lie $ \mathfrak{g}_{\geq 0}$ et $M$ sa composante de Levi d'alg\`ebre $\mathfrak{g}_{0}$. On a $N\in \mathfrak{g}_{2}({\mathbb F}_{q})$ et l'orbite de $N$ sous l'action de $M$ est un ouvert de Zariski de $\mathfrak{g}_{2}$ que l'on note $\tilde{\mathfrak{g}}_{2}$. Pour $N'\in \tilde{\mathfrak{g}}_{2}({\mathbb F}_{q})$, on a
$$\{u^{-1}N'u; u\in U_{P}({\mathbb F}_{q})\}=N'+ \mathfrak{g}_{\geq 3}({\mathbb F}_{q})$$
et
$$\{g\in G; g^{-1}N'g\in \mathfrak{u}_{P}\}=P.$$
Notons $\tilde{f}_{N,\epsilon}$ la fonction sur $\mathfrak{g}({\mathbb F}_{q})$, \`a support dans $\tilde{\mathfrak{g}}_{2}({\mathbb F}_{q})+\mathfrak{g}_{\geq 3}({\mathbb F}_{q})$, telle que, pour $N'\in \tilde{\mathfrak{g}}_{2}({\mathbb F}_{q})$ et $X\in \mathfrak{g}_{\geq 3}({\mathbb F}_{q})$, $\tilde{f}_{N,\epsilon}(N'+X)=f_{N,\epsilon}(N')$. A l'aide des propri\'et\'es pr\'ec\'edentes, on v\'erifie que
$$f_{N,\epsilon}(X)=\vert P({\mathbb F}_{q})\vert ^{-1}\sum_{g\in G({\mathbb F}_{q})}\tilde{f}_{N,\epsilon}(g^{-1}Xg)$$
pour tout $X\in \mathfrak{g}({\mathbb F}_{q})$.

 \subsection{Groupes $p$-adiques}\label{groupespadiques}
 Pour la suite de la section, $G$ est un groupe r\'eductif connexe d\'efini sur $F$. On note $A_{G}$ le plus grand sous-tore de $Z(G)$ qui est d\'eploy\'e sur $F$. 
 
  On note $\hat{G}$ le dual de Langlands de $G$, qui est un groupe complexe. Il est muni d'une action de $\Gamma_{F}$ que l'on note simplement $\sigma\mapsto \sigma_{G}$. Le $L$-groupe est le produit semi-direct $\hat{G}\rtimes W_{F}$.   On note $\hat{G}_{AD}$ le groupe adjoint de $\hat{G}$ et $\hat{G}_{SC}$ son rev\^etement simplement connexe (ce sont les groupes duaux de $G_{SC}$, resp. $G_{AD}$).

On note $Imm(G_{AD})$ l'immeuble de Bruhat-Tits du groupe adjoint $G_{AD}$ sur $F$. Le groupe $G(F)$, et plus g\'en\'eralement le groupe $G_{AD}(F)$, agissent sur cet immeuble. D'autre part, celui-ci est d\'ecompos\'e en facettes. On note $Fac(G)$ l'ensemble des facettes et  $S(G)\subset Fac(G)$ l'ensemble des sommets. On consid\'erera aussi l'immeuble de $G_{AD}$ sur des extensions $K$ de $F$ contenues dans $\bar{F}$. On ajoutera alors la lettre $K$ dans la notation: $Imm_{K}(G_{AD})$, $S_{K}(G)$ etc... Si $K/F$ est galoisienne, le groupe $\Gamma_{K/F}$ agit sur $Imm_{K}(G_{AD})$  et, si $K/F$ est mod\'er\'ement ramifi\'ee, $Imm(G_{AD})$ s'identifie \`a l'ensemble des points fixes par $\Gamma_{K/F}$ dans $Imm_{K}(G_{AD})$.  

 Pour tout ${\cal F}\in Fac(G)$, on note $K_{{\cal F}}^{\dag}$ le stabilisateur de ${\cal F}$ dans $G(F)$,  $K_{{\cal F}}^0\subset K_{{\cal F}}^{\dag}$ le sous-groupe parahorique, $K_{{\cal F}}^+$ son plus grand sous-groupe distigu\'e pro-$p$-unipotent et $G_{{\cal F}}$ le groupe r\'eductif connexe sur ${\mathbb F}_{q}$ d\'efini par Bruhat et Tits tel que $K_{{\cal F}}^0/K_{{\cal F}}^+\simeq G_{{\cal F}}({\mathbb F}_{q})$. Il y a des objets correspondants dans l'alg\`ebre de Lie: $\mathfrak{k}_{{\cal F}}^+\subset \mathfrak{k}_{{\cal F}}\subset \mathfrak{g}(F)$, avec $\mathfrak{k}_{{\cal F}}/\mathfrak{k}_{{\cal F}}^+\simeq \mathfrak{g}_{{\cal F}}({\mathbb F}_{q})$. On ajoute des indices $AD$ pour les m\^emes objets relatifs au groupe $G_{AD}$: $K_{{\cal F},AD}^{\dag}$ etc...
 
 \subsection{L'espace $FC(\mathfrak{g}(F))$}\label{lespaceFC}
 
 Fixons un caract\`ere $\psi$ de $F$ de conducteur $\mathfrak{o}_{F}$ et une forme bilin\'eaire sym\'etrique et non d\'eg\'en\'er\'ee $<.,.>$ sur $\mathfrak{g}(F)$, invariante par l'action par conjugaison de $G(F)$. On d\'efinit la transformation de Fourier $f\mapsto \hat{f}$ dans $C_{c}^{\infty}(\mathfrak{g}(F))$ par la formule usuelle
 $$\hat{f}(X)=\int_{\mathfrak{g}(F)}f(Y)\psi(<X,Y>)\,dY,$$
 o\`u $dY$ est la mesure autoduale. 
 On sait que l'on peut choisir la forme $<.,.>$ de sorte que, pour tout ${\cal F}\in Fac(G)$, on ait l'\'egalit\'e $\hat{{\bf 1}}_{\mathfrak{k}_{{\cal F}}}=\vert \mathfrak{g}_{{\cal F}}({\mathbb F}_{q})\vert ^{1/2}{\bf 1}_{\mathfrak{k}_{{\cal F}}^+}$, o\`u l'on a not\'e ${\bf 1}_{\mathfrak{k}_{{\cal F}}}$ et ${\bf 1}_{\mathfrak{k}_{{\cal F}}^+}$ les fonctions caract\'eristiques de $\mathfrak{k}_{{\cal F}}$ et $\mathfrak{k}_{{\cal F}}^+$. On suppose la forme ainsi choisie. La mesure sur $\mathfrak{g}(F)$ v\'erifie alors l'\'egalit\'e $mes(\mathfrak{k}_{{\cal F}}^+)=\vert \mathfrak{g}_{{\cal F}}({\mathbb F}_{q})\vert ^{-1/2}$ pour tout ${\cal F}\in Fac(G)$. On peut relever cette mesure en une mesure de Haar sur $G(F)$ telle que  $mes(K_{{\cal F}}^+)=mes(\mathfrak{k}_{{\cal F}}^+)$ pour tout ${\cal F}$.  Des mesures analogues seront choisies sur tout autre groupe r\'eductif connexe. 

Pour $f\in C_{c}^{\infty}(\mathfrak{g}(F))$ et $X\in \mathfrak{g}_{reg}(F)$, on d\'efinit l'int\'egrale orbitale
$$I^G(X,f)=D_{G}(X)^{1/2}\int_{G_{X}(F)\backslash G(F)}f(g^{-1}Xg)\,dg,$$
o\`u $D_{G}$ est le discriminant de Weyl. 
On note $I(\mathfrak{g}(F))$ le quotient de $C_{c}^{\infty}(\mathfrak{g}(F))$ par le sous-espace des fonctions $f$ telles que $I^G(X,f)=0$ pour tout $X\in \mathfrak{g}_{reg}(F)$. Les int\'egrales orbitales se descendent en des formes lin\'eaires sur $I(\mathfrak{g}(F))$. On note $I_{cusp}(\mathfrak{g}(F))$ le sous-espace de $I(\mathfrak{g}(F))$ form\'e des $f\in I(\mathfrak{g}(F))$ telles que $I^G(X,f)=0$ pour tout $X\in \mathfrak{g}_{reg}(F)$ tel que $X\not\in\mathfrak{g}_{ell}(F)$. 

Notons ${\cal T}_{ell}$ un  ensemble de repr\'esentants des classes de conjugaison par $G(F)$ dans l'ensemble des sous-tores maximaux elliptiques de $G$. Pour un tel tore $T$, posons $W^G(T)=Norm_{G(F)}(T)/T(F)$.
L'espace $I_{cusp}(\mathfrak{g}(F))$ est muni du produit scalaire elliptique
$$(f,f')_{ell}=\sum_{T\in {\cal T}_{ell}} \vert W^G(T)\vert ^{-1}mes(A_{G}(F)\backslash T(F))\int_{A_{G}(F)\backslash T(F)}I^G(X,\bar{f})I^G(X,f')\,dX.$$
C'est un produit hermitien d\'efini positif. 

Le groupe $G_{AD}(F)$ agit naturellement par conjugaison sur $I(\mathfrak{g}(F))$ et $I_{cusp}(\mathfrak{g}(F))$. Cette action se factorise en une action du groupe $G_{AD}(F)/\pi(G(F))$. Ce dernier est ab\'elien fini. On pose $\Xi=(G_{AD}(F)/\pi(G(F)))^{\vee}$. On notera ${\bf 1}$ l'\'el\'ement neutre de $\Xi$.  Pour $\xi\in \Xi$, on note $I_{cusp,\xi}(\mathfrak{g}(F))$ le sous-espace des $f\in I_{cusp}(\mathfrak{g}(F))$ telles que $I^G(g^{-1}Xg,f)=\xi(g)I^G(X,f)$ pour tous $g\in G_{AD}(F)$ et $X\in \mathfrak{g}_{reg}(F)$. On a la d\'ecomposition
$$(1) \qquad I_{cusp}(\mathfrak{g}(F))=\oplus_{\xi\in \Xi}I_{cusp,\xi}(\mathfrak{g}(F)).$$
Elle est orthogonale pour le produit scalaire elliptique.

 {\bf Remarque.}  L'homomorphisme de Langlands fournit un isomorphisme $H^1(W_{F},Z(\hat{G}_{SC}))\simeq \Xi$.   \bigskip

  Pour une facette ${\cal F}\in Fac(G)$, l'espace ${\cal C}(\mathfrak{g}_{{\cal F}}({\mathbb F}_{q}))$  s'identifie \`a un sous-espace de $C_{c}^{\infty}(\mathfrak{g}(F))$: une fonction sur $\mathfrak{g}_{{\cal F}}({\mathbb F}_{q})\simeq \mathfrak{k}_{{\cal F}}/\mathfrak{k}_{{\cal F}}^+$ se rel\`eve en une fonction sur $\mathfrak{k}_{{\cal F}}$ puis s'\'etend \`a $\mathfrak{g}(F)$ par $0$ hors de $\mathfrak{k}_{{\cal F}}$.  On note $FC(\mathfrak{g}(F))$ l'image  dans $I(\mathfrak{g}(F))$ de l'espace
  $$\sum_{s\in S(G)}FC(\mathfrak{g}_{s}({\mathbb F}_{q})).$$
   On sait qu'en fait, cette image est contenue dans le sous-espace $I_{cusp}(\mathfrak{g}(F))$. 
   
   \subsection{Orbites dans l'ensemble des facettes}\label{orbites}
   
    Pour ce paragraphe et jusqu'en \ref{donneesendoscopiques}, on  suppose $G$ simplement connexe et absolument quasi-simple.

Supposons d'abord $G$ quasi-d\'eploy\'e. On fixe un sous-groupe de Borel $B$ de $G$ et un sous-tore maximal $T$ de $B$ tous deux d\'efinis sur $F$. On note $T^{nr}$ le plus grand sous-tore de $T$ d\'eploy\'e sur $F^{nr}$. On note $\Sigma$  l'ensemble des racines de $T$ dans $G$ et $\Delta$ celui des racines simples d\'etermin\'e par $B$. On fixe un \'epinglage $(E_{\alpha})_{\alpha\in \Delta}$ conserv\'e par l'action galoisienne et on prolonge cet \'epinglage en une base de Chevalley sur $\bar{F}$. Cette base d\'etermine un unique point $s^{nr}\in S_{F^{nr}}(G_{AD})$ de sorte que $\mathfrak{k}_{s^{nr},F^{nr}}$ soit l'ensemble des points fixes par $I_{F}$ dans le $\mathfrak{o}_{\bar{F}}$-r\'eseau engendr\'e par  les \'el\'ements de la base.  On note $C^{nr}$ la chambre de $Imm_{F^{nr}}(G_{AD})$ qui v\'erifie les conditions suivantes: elle appartient \`a l'appartement associ\'e \`a $T^{nr}$, $s^{nr}$ appartient \`a la cl\^oture de $C^{nr}$ et $\mathfrak{k}_{C^{nr},F^{nr}}\cap B(F^{nr})=\mathfrak{k}_{s^{nr},F^{nr}}\cap B(F^{nr})$. Le point $s^{nr}$ est fix\'e par l'action galoisienne de $\Gamma_{F}^{nr}$ et la chambre $C^{nr}$ est conserv\'ee par cette action.  A la  chambre $C^{nr}$ est associ\'e un ensemble de racines affines (dont les annulateurs sont les murs de la chambre et qui sont positives sur celle-ci). Notons $\Delta_{a}^{nr}$ l'ensemble de racines sous-jacent \`a cet ensemble de racines affines. Ce sont des formes lin\'eaires sur $X_{*}(T)^{I_{F}}=X_{*}(T^{nr})$.  Ces formes lin\'eaires s\'eparent les \'el\'ements de $X_{*}(T)^{I_{F}}$ et l'espace des relations lin\'eaires entre \'el\'ements de $\Delta_{a}^{nr}$ est une droite (il y a une unique relation lin\'eaire, \`a un scalaire pr\`es).  
Cet ensemble de racines $\Delta_{a}^{nr}$ plus les relations de produit scalaire entre elles forment un diagramme ${\cal D}_{a}^{nr}$. Dans chaque cas, ce diagramme ${\cal D}_{a}^{nr}$ est d\'ecrit dans les tables de Tits sous le nom de "local index", cf. \cite{T}

Posons $N=K_{C^{nr},AD,F^{nr}}^{\dag}/K_{C^{nr},AD,F^{nr}}^0$.
 D'apr\`es \cite{HR} lemme 14, on a $N\simeq (Z(\hat{G}_{SC})^{I_{F}})^{\vee}$, en particulier, c'est un groupe ab\'elien. On sait que $H^1(\Gamma_{F},G_{AD})\simeq H^1(\Gamma_{F}^{nr},N)$ et ce dernier groupe n'est autre que le quotient $N_{\Gamma_{F}^{nr}}$ de $N$ par le sous-groupe des $Fr(n)n^{-1}$ pour $n\in N$.  Remarquons que $N_{\Gamma^{nr}_{F}}$ est aussi \'egal \`a $(Z(\hat{G}_{SC})^{\Gamma_{F}})^{\vee}$ et on retrouve l'isomorphisme de Kottwitz $H^1(\Gamma_{F},G_{AD})\simeq (Z(\hat{G}_{SC})^{\Gamma_{F}})^{\vee}$, cf. \cite{K} th\'eor\`eme 1.2.

Supprimons l'hypoth\`ese que $G$ est quasi-d\'eploy\'e. On peut fixer une forme quasi-d\'eploy\'ee $G^*$ de $G$ et un  torseur int\'erieur $\varphi:G\to G^*$. On effectue  les constructions pr\'ec\'edentes pour ce groupe $G^*$ (on ajoute si besoin est  des exposants $*$). Au torseur int\'erieur est associ\'e un \'el\'ement de $H^1(\Gamma_{F},G^*_{AD})$ que l'on note $u_{G}$ puis un \'el\'ement de $N_{\Gamma_{F}^{nr}}^*$. Relevons ce dernier en un \'el\'ement de $N^*$. On peut ensuite relever celui-ci en un \'el\'ement $n_{G}\in G^*_{AD}(F^{nr})$ qui se prolonge en un cocycle de $\Gamma_{F}^{nr}$ dans $G_{AD}^*(F^{nr})$, c'est-\`a-dire un cocycle $\underline{n}_{G}$ tel que $\underline{n}_{G}(Fr)=n_{G}$.  On note encore $\underline{n}_{G}$ le cocycle de $\Gamma_{F}$ dans $G^*(F^{nr})$ obtenu par inflation. On peut alors supposer que $G(F^{nr})=G^*(F^{nr})$, que $\varphi$ est l'identit\'e et que l'on a les \'egalit\'es d'actions galoisiennes $\sigma_{G}(g)=\underline{n}_{G}(\sigma)\sigma_{G^*}(g)\underline{n}_{G}(\sigma)^{-1}$ pour tous $g\in G(F^{nr})$ et $\sigma\in \Gamma_{F}$. 

Notons $\underline{S}(G)$ l'intersection de $S(G)$ et de la cl\^oture de $C^{nr}$.  C'est un ensemble de repr\'esentants des orbites de l'action de $G(F)$ dans $S(G)$.  
Le groupe $\Gamma_{F}$ agit   sur l'ensemble $\Delta_{a}^{nr}$ par $(\sigma,\alpha)\mapsto \underline{n}_{G}(\sigma)\circ \sigma_{G^*}(\alpha)$ pour $\sigma\in \Gamma_{F}$ et $\alpha\in \Delta_{a}^{nr}$. Les \'el\'ements de $\underline{S}(G)$ sont en bijection avec les orbites non vides de l'action de $\Gamma_{F}^{nr}$. Pour un sommet $s$ correspondant \`a une orbite ${\cal O}_{s}$, le diagramme de Dynkin associ\'e au groupe $G_{s}$, muni de son action de $\Gamma_{F}^{nr}$, s'obtient en supprimant les \'el\'ements de ${\cal O}_{s}$ du diagramme ${\cal D}_{a}^{nr}$. En particulier, un sous-tore maximal de $G_{s}$ a pour groupe de cocaract\`eres $X_{*}(T^*)^{I_{F}}$. On a vu en \cite{W1} 9.2 que l'espace
$$(1) \qquad \sum_{s\in \underline{S}(G)}FC(\mathfrak{g}_{s}({\mathbb F}_{q}))$$
s'envoyait bijectivement sur $FC(\mathfrak{g}(F))$. Nous identifierons ces deux espaces.

 \begin{lem}{Soit $s\in S(G,C^{nr})$ et ${\cal O}_{s}\subset \Delta_{a}^{nr}$ l'orbite galoisienne associ\'ee. Si ${\cal O}_{s}$ a au moins $2$ \'el\'ements, on a $FC(\mathfrak{g}_{s}({\mathbb F}_{q}))=\{0\}$.}\end{lem}
 
 Preuve. Puisqu'il n'y a, \`a un scalaire pr\`es, qu'une seule relation lin\'eaire entre les \'el\'ements de $\Delta_{a}^{nr}$, on voit que le nombre de racines simples de $G_{s}$ est \'egal \`a la dimension de $X_{*}(T^*)^{I_{F}}$ plus $1$ moins le nombre d'\'el\'ements de ${\cal O}_{s}$. Si ${\cal O}_{s}$ a au moins $2$ \'el\'ements, le nombre de racines simples de $G_{s}$ est strictement inf\'erieur \`a la dimension d'un sous-tore maximal, donc $Z(G_{s})^0\not=\{1\}$. La conclusion r\'esulte de \ref{groupessurFq} (1). $\square$ 
 
 \subsection{Alc\^oves}\label{alcoves}
 On suppose que $G$   n'est pas de type $A_{n-1}$. On se propose de d\'ecrire l'ensemble $\Delta_{a}^{nr}$ et l'alc\^ove $C^{nr}$. Puisque le corps de base est ici $F^{nr}$, on ne perd rien \`a supposer $G$ quasi-d\'eploy\'e. 
 
 Le premier cas est celui o\`u $G$ est d\'eploy\'e sur $F^{nr}$.   On pose ${\cal A}^{nr}=X_{*}(T)\otimes_{{\mathbb Z}}{\mathbb R}$. L'appartement associ\'e \`a $T$ est isomorphe \`a ${\cal A}^{nr}$, le point $s^{nr}$ d'identifiant \`a $0\in {\cal A}$. On pose $\Sigma^{nr}=\Sigma$. Pour $\alpha\in \Sigma^{nr}$ et $b\in {\mathbb Z}$, notons $\alpha[b]$ la forme affine sur ${\cal A}^{nr}$ d\'efinie par $\alpha[b](H)=\alpha(H)+b$. Dans le cas particulier $b=0$, on continue \`a noter   simplement $\alpha=\alpha[0]$. On note $\Sigma^{aff}$ l'ensemble des  racines affines, c'est-\`a-dire  des applications $\alpha[b]$ pour
  $\alpha\in \Sigma^{nr}$ et $b\in {\mathbb Z}$.
 Notons $\alpha_{0}$ l'oppos\'ee de la  plus grande racine dans $\Sigma^{nr}$. On a $\Delta^{nr}_{a}=\Delta\cup\{\alpha_{0}\}$. En \'ecrivant $-\alpha_{0}$ dans la base $\Delta$, on obtient une relation
$$(1) \qquad \sum_{\alpha\in \Delta_{a}}d(\alpha)\alpha=0$$
o\`u $d(\alpha_{0})=1$ et, pour $\alpha\in \Delta$, $d(\alpha)$ appartient \`a l'ensemble ${\mathbb N}_{>0}$ des entiers strictement positifs.   Les murs de $C^{nr}$ sont les annulateurs de $\alpha_{0}[1]$ et des $\alpha\in \Delta$. Plus pr\'ecis\'ement, $C^{nr}$ est l'ensemble des $H\in {\cal A}^{nr}$ tels que $\alpha_{0}[1](H)>0$ et $\alpha(H)>0$ pour tout $\alpha\in \Delta$. Il est commode de poser $\alpha^{aff}=\alpha$ pour $\alpha\in \Delta$ et $\alpha_{0}^{aff}=\alpha_{0}[1]$. La relation (1) devient la relation affine
$$(2) \qquad \sum_{\alpha\in \Delta_{a}}d(\alpha)\alpha^{aff}=1.$$

Supposons maintenant que $G$ n'est pas d\'eploy\'e sur $F^{nr}$. Notre hypoth\`ese sur $p$ implique que $G$ est d\'eploy\'e sur une extension mod\'er\'ement ramifi\'ee. Puisque le groupe de ramification mod\'er\'ee est ab\'elien, l'action alg\'ebrique  de $I_{F}$ sur $G$ se factorise par un homomorphisme surjectif $I_{F}\to {\mathbb Z}/e{\mathbb Z}$. Un \'el\'ement de $I_{F}$ d'image $1$ dans ce dernier groupe agit par un automorphisme $\theta$ de $G$ d'ordre $e$ qui pr\'eserve  $B$, $T$ et l'\'epinglage. Notre hypoth\`ese que $G$ n'est pas d\'eploy\'e sur $F^{nr}$ signifie que $e\geq2$. En inspectant tous les syst\`emes de racines, on voit que $e=2$ ou $e=3$. De plus, l'existence d'un tel $\theta$ implique que toutes les racines de $G$ ont m\^eme longueur.  Posons $\bar{{\cal A}}= X_{*}(T)\otimes_{{\mathbb Z}}{\mathbb R}$. C'est l'appartement associ\'e \`a $T$ dans l'immeuble de $G_{AD}$ sur une extension de $F$ qui d\'eploie $G$. On  fixe  sur cet espace un produit scalaire invariant par l'action du groupe de Weyl $W$ de $G$ de sorte que $(\check{\beta},\check{\beta})=2$ pour tout $\beta\in \Sigma$. Posons ${\cal A}^{nr}=X_{*}(T)^{\theta}\otimes_{{\mathbb Z}}{\mathbb R}$, que l'on munit de la restriction du produit scalaire pr\'ec\'edent. C'est l'appartement associ\'e \`a $T^{nr}$ dans $Imm_{F^{nr}}(G_{AD})$.  
Notons $\beta^{res}$ la restriction de $\beta$ \`a ${\cal A}^{nr}$. On note $\Sigma^{nr}$ l'ensemble des $\beta^{res}$ pour $\beta\in \Sigma$. Pour $\alpha\in \Sigma^{nr}$, choisissons $\beta\in \Sigma$ telle que $\alpha=\beta^{res}$.  Notons $e(\alpha)$ le plus petit entier $m\geq1$ tel que $\theta^{m}(\beta)=\beta$. Il ne d\'epend pas du choix de $\beta$.  On a forc\'ement $e(\alpha)=1$ ou $e$. 
Si  $e(\alpha)=1$, c'est-\`a-dire $\theta(\beta)=\beta$, on pose $\check{\alpha}=\check{\beta}$, qui appartient \`a ${\cal A}^{nr}$. Si $ e(\alpha)=e$, les $e$ racines $\beta,\theta(\beta),...,\theta^{e-1}(\beta)$ sont distinctes et orthogonales (car on a suppos\'e que $G$ n'\'etait pas de type $A_{n-1}$). On pose $\check{\alpha}=\check{\beta}+\theta(\check{\beta})+...+\theta^{e-1}(\check{\beta})$. C'est un \'el\'ement de ${\cal A}^{nr}$. L'ensemble  $\Sigma^{nr}$ est un syst\`eme de racines irr\'eductible et l'ensemble $\check{\Sigma}^{nr}=\{\check{\alpha};\alpha\in \Sigma^{nr}\}$ en est le syst\`eme de coracines.  Remarquons que $(\check{\alpha},\check{\alpha})$ vaut $2e(\alpha)$ pour tout $\alpha\in \Sigma^{nr}$.   On note $\Delta^{nr}$  l'ensemble des $\beta^{res}$  pour $\beta\in \Delta$. C'est une  base de $\Sigma^{res}$. Notons $\alpha_{0}$ l'\'el\'ement de $\Sigma^{nr}$ tel que $\check{\alpha}_{0}$ soit l'oppos\'e de la plus grande coracine du syst\`eme $\check{\Sigma}^{nr}$ (pour l'ordre d\'efini par la base $\check{\Delta}^{nr}=\{\check{\alpha}; \alpha\in \Delta^{nr}\}$).

{\bf Remarque.} On prendra garde \`a cette notation: en g\'en\'eral, $\alpha_{0}$ n'est pas l'oppos\'ee de la  plus grande racine de $\Sigma^{nr}$. 
\bigskip

 Posons  $\Delta_{a}^{nr}=\{\alpha_{0}\}\cup \Delta^{nr}$. En  \'ecrivant $-\check{ \alpha}_{0}$ dans la base $\check{\Delta}^{nr}$, on obtient une relation 
$$\sum_{\alpha\in \Delta_{a}^{nr}}d(\check{\alpha})\check{\alpha}=0,$$
avec $d(\check{ \alpha}_{0})=1$ et $d(\check{\alpha})\in {\mathbb N}_{>0}$ pour tout $\alpha\in \Delta^{nr}$. Puisque le produit scalaire permet d'identifier $\check{\alpha}$ \`a $e(\alpha)\alpha$, on en d\'eduit la relation
$$(3) \qquad \sum_{\alpha\in \Delta_{a}^{nr}}d(\alpha)\alpha=0,$$
o\`u $d(\alpha)=d(\check{\alpha})e(\alpha)$ pour tout $\alpha\in \Delta_{a}^{nr}$.

 A $\beta\in \Sigma$ est associ\'ee une sym\'etrie $s_{\beta}$ de l'espace $ \bar{{\cal A}}$ qui appartient au  groupe de Weyl $W$. A $\alpha\in \Sigma^{nr}$ est associ\'ee une sym\'etrie $s_{\alpha}$ de l'espace ${\cal A}^{nr}$. Supposons $\alpha=\beta^{res}$. 
 Si $e(\alpha)=1$, $s_{\alpha}$ est la restriction de $s_{\beta}$ \`a ${\cal A}^{nr}$. Si $e(\alpha)=e$, $s_{\alpha}$ est la restriction \`a ${\cal A}^{nr}$ de $s_{\beta}s_{\theta(\beta)}...s_{\theta^{e-1}(\beta)}$. Introduisons l'homomorphisme
$v:T(\bar{F})\to \bar{{\cal A}}$ tel que $v(x_{*}\otimes t)=x_{*}\otimes val_{F}(t)$ pour tout $x_{*}\in X_{*}(T)$ et $t\in \bar{F}^{\times}$. Le groupe de Weyl  affine $W_{a}^{nr}$  est l'ensemble de transformations de ${\cal A}^{nr}$ engendr\'e par les sym\'etries $s_{\alpha}$ pour $\alpha\in \Sigma^{nr}$ et l'ensemble de translations par des \'el\'ements de $v(T(F^{nr}))$. Pour $\alpha\in \Sigma^{nr}$ et $b\in {\mathbb R}$, notons $\alpha[b]$ la forme affine sur ${\cal A}^{nr}$ d\'efinie par $\alpha[b](H)=\alpha(H)+b$. 
Les racines affines sont les formes affines $ \alpha[b]$ telles qu'il existe un \'el\'ement $n\in W^{nr}_{a}$ de sorte que l'on ait l'\'egalit\'e $\{H\in {\cal A}^{nr};n(H)=H\}=\{H\in {\cal A}^{nr};\alpha[b](H)=0\}$. On voit que ce sont les $\alpha[b]$ o\`u $b\in \frac{1}{e(\alpha)}{\mathbb Z}$. On note $\Sigma^{aff}$ l'ensemble des racines affines. 
Montrons que

(4) les murs de $C^{nr}$ sont les annulateurs de $\alpha_{0}[\frac{1}{e(\alpha_{0})}]$ et des $\alpha\in \Delta^{nr}$; l'ensemble  $C^{nr}$ est celui des $H\in {\cal A}^{nr}$ tels que $\alpha_{0}[\frac{1}{e(\alpha_{0})}](H)>0$ et $\alpha(H)>0$ pour tout $\alpha\in \Delta^{nr}$. 

Il suffit de prouver que les annulateurs des racines affines ne coupent pas l'ensemble ainsi d\'efini. Consid\'erons une telle racine affine $\alpha[\frac{m}{e(\alpha)}]$ avec $m\in {\mathbb Z}$ et un \'el\'ement $H$ de l'ensemble ci-dessus.  Supposons $\alpha[\frac{m}{e(\alpha)}](H)=0$, c'est-\`a-dire $e(\alpha)\alpha(H)= -m$. Quitte \`a remplacer $\alpha$ par $-\alpha$ et $m$ par $-m$, on peut supposer $\alpha>0$ (pour l'ordre d\'efini par la base $\Delta^{nr}$). Ecrivons $\check{\alpha}$ dans la base $\check{\Delta}^{nr}$:
$$\check{\alpha}=\sum_{\alpha'\in \Delta^{nr}}m(\check{\alpha}')\check{\alpha}'$$
avec des $m(\alpha')\in {\mathbb N}$.
 En identifiant toute coracine $\check{\alpha}''$ \`a $e(\alpha'')\alpha''$, la relation  ci-dessus devient
$$e(\alpha)\alpha=\sum_{\alpha'\in \Delta^{nr}}m(\check{\alpha}')e(\alpha')\alpha'.$$
On  a $\alpha'(H)>0$ pour tout $\alpha'\in \Delta^{nr}$ et au moins un des coefficients $m(\check{\alpha}')$ est strictement positif. Donc $e(\alpha)\alpha(H)>0$. Puisque $-\check{\alpha}_{0}$ est l'oppos\'ee de la plus grande racine de $\check{\Sigma}^{nr}$, on sait que $m(\check{\alpha}')\leq d(\check{\alpha}')$ pour tout $\alpha'\in \Delta^{nr}$. D'o\`u
$$e(\alpha)\alpha(H)\leq \sum_{\alpha'\in \Delta^{nr}}d(\check{\alpha}')e(\alpha')\alpha'(H).$$
En utilisant (3) et l'\'egalit\'e $d(\check{\alpha}_{0})=1$, le membre de droite ci-dessus est $-e(\alpha_{0})\alpha_{0}(H)$. On a $\alpha_{0}[\frac{1}{e(\alpha_{0})}](H)>0$ donc ce membre de droite est strictement inf\'erieur \`a $1$. On a ainsi  d\'emontr\'e que $0< -m=e(\alpha)\alpha(H)<1$. Mais $m\in {\mathbb Z}$, contradiction. Cela d\'emontre (4). 

Il  r\'esulte de (4) que $\Delta_{a}^{nr} =\Delta^{nr}\cup\{\alpha_{0}\}$. 
  On pose   $\alpha^{aff}=\alpha$ pour $\alpha\in \Delta^{nr}$ et $\alpha_{0}^{aff}=\alpha_{0}[\frac{1}{e(\alpha_{0})}]$. La relation (3) devient la relation affine (2) qui est donc v\'erifi\'ee en tout cas.
  
 \subsection{Action de $G_{AD}(F)$ sur $FC(\mathfrak{g}(F))$}\label{actionsurFC}
 Introduisons l'homomorphisme de Kottwitz $w_{G_{AD}}:G_{AD}(F)\to (Z(\hat{G}_{SC})^{I_{F}})^{\vee,\Gamma_{{\mathbb F}_{q}}}$. On note $G_{AD}(F)_{0}$ son noyau. Introduisons le groupe $G(F^{nr})^{\sharp}=\{g\in G(F^{nr}); gFr(g)^{-1}\in Z(G)^{I_{F}}\}$. C'est l'image r\'eciproque de $G_{AD}(F)$ par l'application $\pi^{nr}:G(F^{nr})\to G_{AD}(F^{nr})$. Montrons que
  
  (1) l'homomorphisme $\pi^{nr}$ envoie surjectivement $G(F^{nr})^{\sharp}$ sur $G_{AD}(F)_{0}$.
   
   Posons $W_{{\mathbb F}_{q}}=W_{F}/I_{F}$. C'est le sous-groupe de $\Gamma_{{\mathbb F}_{q}}$ form\'e des puissances enti\`eres du Frobenius. On a la suite exacte
  $$(2) \qquad 1\to H^1(W_{{\mathbb F}_{q}},Z(\hat{G}_{SC})^{I_{F}})\stackrel{\iota}{\to} H^1(W_{F},Z(\hat{G}_{SC}))\to H^1(I_{F},Z(\hat{G}_{SC}))$$
 Remarquons que $H^1(W_{{\mathbb F}_{q}},Z(\hat{G}_{SC})^{I_{F}})$ n'est autre que l'espace des coinvariants $(Z(\hat{G}_{SC})^{I_{F}})_{\Gamma_{{\mathbb F}_{q}}}$ dont le dual est $(Z(\hat{G}_{SC})^{I_{F}})^{\vee,\Gamma_{{\mathbb F}_{q}}}$.  Dualement au premier homomorphisme de la suite ci-dessus, on a donc un homomorphisme
 $$  H^1(W_{F},Z(\hat{G}_{SC}))^{\vee}\stackrel{\iota^{\vee}}{\to} (Z(\hat{G}_{SC})^{I_{F}})^{\vee,\Gamma_{{\mathbb F}_{q}}}.$$
 L'homomorphisme $w_{G_{AD}}$ est le compos\'e de l'homomorphisme de Langlands $a_{G_{AD}}:G_{AD}(F)\to H^1(W_{F},Z(\hat{G}_{SC}))^{\vee}$ et de l'homomorphisme pr\'ec\'edent.

  Soit $m\geq1$ un entier, notons $F_{m}$  l'extension  non ramifi\'ee de $F$ de degr\'e $m$. On dispose de l'homomorphisme de corestriction $H^1(W_{F_{m}},Z(\hat{G}_{SC}))\stackrel{Cores_{m}}{\to} H^1(W_{F},Z(\hat{G}_{SC}))$ et la suite (2) s'ins\`ere dans un diagramme commutatif
  $$\begin{array}{ccccccc} 1&\to& H^1(W_{{\mathbb F}_{q^m}},Z(\hat{G}_{SC})^{I_{F}})&\to& H^1(W_{F_{m}},Z(\hat{G}_{SC}))&\to& H^1(I_{F},Z(\hat{G}_{SC}))\\&&\,\,\,\,\,\downarrow cores_{m}&&\,\,\,\,\,\downarrow Cores_{m}&&\,\,\downarrow C_{m}\\ 1&\to& H^1(W_{{\mathbb F}_{q}},Z(\hat{G}_{SC})^{I_{F}})&\stackrel{\iota}{\to}& H^1(W_{F},Z(\hat{G}_{SC}))&\to& H^1(I_{F},Z(\hat{G}_{SC}))\\ \end{array}$$
  L'homomorphisme $C_{m}$ se d\'ecrit de la fa\c{c}on suivante: pour un cocycle $u:I_{F}\to Z(\hat{G}_{SC})$, $C_{m}(u)$ est le cocycle $\sigma\mapsto \prod_{j=0,...,m-1}Fr^{-j}(u(Fr^j\sigma Fr^{-j}))$. On voit sur cette formule que $C_{m}=0$ si $m$ est assez divisible (on entend par l\`a qu'il existe un entier $N\geq1$ tel que la propri\'et\'e soit v\'erifi\'ee pour $m$ divisible par $N$). D'autre part, l'homomorphisme $cores_{m}$ est surjectif: il s'identifie \`a l'homomorphisme naturel  $(Z(\hat{G}_{SC})^{I_{F}})_{\Gamma_{{\mathbb F}_{q^m}}}\to (Z(\hat{G}_{SC})^{I_{F}})_{\Gamma_{{\mathbb F}_{q}}}$. De ces deux propri\'et\'es r\'esulte l'assertion suivante, o\`u on note $Cores_{m}^{\vee}$ et $\iota^{\vee}$ les homomorphismes duaux de $Cores_{m}$ et $\iota$: 
  
  (3) $ Ker(Cores_{m}^{\vee})=Ker(\iota^{\vee})$ si $m$ est assez divisible.

    Consid\'erons le diagramme commutatif
  $$\begin{array}{ccccccc}G(F)&\stackrel{\pi}{\to}& G_{AD}(F)&\stackrel{a_{G_{AD}}}{\to}& H^1(W_{F},Z(\hat{G}_{SC}))^{\vee}&\stackrel{\iota^{\vee}}{\to}& H^1(W_{{\mathbb F}_{q}},Z(\hat{G}_{SC})^{I_{F}})^{\vee}\\ \downarrow&& \downarrow&&\,\,\,\,\,\downarrow Cores_{m}^{\vee}&&\,\,\,\,\,\downarrow cores_{m}^{\vee}\\ 
  G(F_{m})&\stackrel{\pi_{m}}{\to}& G_{AD}(F_{m})&\stackrel{a_{G_{AD},m}}{\to}& H^1(W_{F_{m}},Z(\hat{G}_{SC}))^{\vee}&\stackrel{\iota_{m}^{\vee}}{\to}& H^1(W_{{\mathbb F}_{q^m}},Z(\hat{G}_{SC})^{I_{F}})^{\vee}
  
  \end{array}$$
  Notons $H_{m}$ l'image dans $G_{AD}(F)$ de $G(F_{m})\cap G(F^{nr})^{\sharp}$, autrement dit le sous-groupe des $g\in G_{AD}(F)$ dont l'image dans $G_{AD}(F_{m})$ appartient \`a l'image de $\pi_{m}:G(F_{m})\to G_{AD}(F_{m})$.  On sait que l'image de $\pi^m$ est le noyau de $a_{G_{AD},m}$. Donc
 $H_{m}$ est le noyau de $Cores_{m}^{\vee}\circ a_{G_{AD}}$. Si $m$ est assez divisible, en utilisant (3), c'est le noyau de $\iota^{\vee}\circ a_{G_{AD}}$, c'est-\`a-dire de $w_{G_{AD}}$. Cela prouve que $H_{m}=G_{AD}(F)_{0}$ si $m$ est assez divisible. Puisque $\pi^{nr}(G(F^{nr})^{\sharp})$ est la r\'eunion des $H_{m}$, cela d\'emontre (1). 

   L'assertion (1)   permet de d\'efinir un homomorphisme $\delta:G_{AD}(F)_{0}\to (Z(G)^{I_{F}})_{\Gamma_{{\mathbb F}_{q}}}$: pour $g\in  G_{AD}(F)_{0}$, on choisit $g'\in  G(F^{nr})^{\sharp}$ tel que $\pi^{nr}(g')=g$ et $\delta(g)$ est l'image de $g'Fr(g')^{-1}$ dans $ (Z(G)^{I_{F}})_{\Gamma_{{\mathbb F}_{q}}}$. Il est imm\'ediat que le noyau de $\delta$ est $\pi(G(F))$. 
    Soit $s\in S(G)$.   D'apr\`es \cite{HR} proposition 3,  le groupe $K_{AD,s}^0=K_{AD,s}^{\dag}\cap G_{AD}(F)_{0}$, a fortiori $K_{AD,s}^0$  est contenu dans $G_{AD}(F)_{0}$. Montrons que 
    
    (4) la restriction de $\delta$ \`a $K_{AD,s}^{0}$ est surjective, a fortiori $\delta$ est surjective. 
    
 Rappelons que, puisque $G$ est simplement connexe, $K_{s}^0=K_{s}^{\dag}$.  Introduisons le groupe $K_{s,F^{nr}}^{0}$, analogue de $K_{s}^0$ sur le corps de base $F^{nr}$.   Soit   $z\in (Z(G)^{I_{F}}) $. On a $(Z(G)^{I_{F}}) \subset  K_{s,F^{nr}}^{0}$.  On voit, essentiellement gr\^ace au th\'eor\`eme de Lang,  qu'il  existe $g\in K_{s,F^{nr}}^{0}$ tel que $gFr(g)^{-1}=z$. On a alors $g\in K_{s,F^{nr}}^{0}\cap G(F^{nr})^{\sharp}$. L'image $\pi^{nr}(g)$ appartient \`a $K_{AD,s}^{0}$ et son image par $\delta$ est l'image de $z$ dans $(Z(G)^{I_{F}})_{\Gamma_{{\mathbb F}_{q}}}$. Cela prouve (4).

Soit $s\in S(G)$. L'espace $FC(\mathfrak{g}_{s}({\mathbb F}_{q}))$ est somme des droites port\'ees par  les \'el\'ements de $fc(\mathfrak{g}_{s}({\mathbb F}_{q}))$, c'est-\`a-dire les fonctions caract\'eristiques des faisceaux-caract\`eres cuspidaux \`a support nilpotent et invariants par l'action de $\Gamma_{{\mathbb F}_{q}}$.    Le groupe $K_{AD,s}^{\dag}$ agit naturellement sur $C(\mathfrak{g}_{s}({\mathbb F}_{q}))$. La d\'efinition g\'eom\'etrique des faisceaux-caract\`eres entra\^{\i}ne que cette action conserve l'espace $FC(\mathfrak{g}_{s}({\mathbb F}_{q}))$ en permutant les droites pr\'ec\'edentes. Fixons une telle droite, port\'ee par une fonction caract\'eristique $f$.  Fixons un \'el\'ement $N\in \mathfrak{g}_{s,nil}({\mathbb F}_{q})$ dans le support de $f$ et une repr\'esentation irr\'eductible $\epsilon$ de $Z_{G_{s}}(N)/Z_{G_{s}}(N)^0$ d\'eterminant cette fonction \`a une constante pr\`es, cf. \ref{groupessurFq}.   Le groupe $Z(G)^{I_{F}}$ se plonge  naturellement dans $Z(G_{s})\subset Z_{G_{s}}(N)$. Notons $\epsilon_{G}$ le caract\`ere de $Z(G)^{I_{F}}$ par lequel ce groupe agit sur $\epsilon$, via ce plongement. Parce que la classe de $\epsilon$ est conserv\'ee par l'action galoisienne, $\epsilon_{G}$ se quotiente en un caract\`ere de  $(Z(G)^{I_{F}})_{\Gamma_{{\mathbb F}_{q}}}$. Notons $K_{AD,s}^{\dag}(f)$ le sous-groupe de $K_{AD,s}^{\dag}$ form\'e des \'el\'ements qui conservent la droite ${\mathbb C}f$. Ce groupe agit sur cette droite par un caract\`ere $\xi_{f}$. Montrons que

(5) $K_{AD,s}^0\subset K_{AD,s}^{\dag}(f)$ et la restriction de $\xi_{f}$ \`a $K_{AD,s}^0$ co\"{\i}ncide avec $\epsilon_{G}\circ \delta$. 
  
  Soit $g\in K_{AD,s}^{0}$. Gr\^ace \`a (1), on peut relever $g$ en un \'el\'ement $g'\in K_{s}^{nr}\cap G(F^{nr})^{\sharp}$. L'action de $g$ sur $FC(\mathfrak{g}_{s}({\mathbb F}_{q}))$ est la m\^eme que celle de $g'$. Il r\'esulte imm\'ediatement des d\'efinitions que cette derni\`ere conserve notre faisceau-caract\`ere et agit sur sa fonction-caract\'eristique par multiplication par $\epsilon_{G}(g'Fr(g')^{-1})$, c'est-\`a-dire par $\epsilon_{G}\circ \delta(g)$. D'o\`u (5).

 Le groupe $G_{AD}(F)$ agit naturellement dans $I(\mathfrak{g}(F))$ en conservant le sous-espace $FC(\mathfrak{g}(F))$. Quand on identifie ce dernier espace \`a l'espace (1) du paragraphe \ref{orbites}, l'action se d\'ecrit de la fa\c{c}on suivante. Soient $g\in G_{AD}(F)$ et $s\in \underline{S}(G)$. On choisit $h\in G(F)$ de sorte que $hgs\in \underline{S}(G)$.  Posons $g\star s=hgs$. Pour $f\in FC(\mathfrak{g}_{s}({\mathbb F}_{q}))$, la fonction $f\circ ad(hg)^{-1}$ appartient \`a $FC(\mathfrak{g}_{g\star s}({\mathbb F}_{q}))$, on la note $g\star f$.  Ces constructions ne d\'ependent pas du choix de $h$. 
On obtient une action $(g,s)\mapsto g\star s$ de $G_{AD}(F)$ sur $\underline{S}(G)$ et une action  
 $f\mapsto g\star f$ de $G_{AD}(F)$ sur $FC(\mathfrak{g}(F))$, qui est l'action cherch\'ee. Les restrictions de ces actions \`a $\pi(G(F))$ sont bien s\^ur triviales. 
 
 Les consid\'erations qui pr\'ec\`edent conduisent \`a la description suivante de  l'action de $G_{AD}(F)$ sur $FC(\mathfrak{g}(F))$.   Fixons un sous-ensemble $\underline{S}(G_{AD})\subset \underline{S}(G)$ de repr\'esentants des orbites de l'action $(g,s)\mapsto g\star s$ de $G_{AD}(F)$ dans $\underline{S}$. Pour tout $s\in \underline{S}(G_{AD})$, fixons un ensemble ${\cal B}_{s}\subset fc(\mathfrak{g}_{s}({\mathbb F}_{q}))$ de sorte que les droites ${\mathbb C}f$ pour $f\in {\cal B}_{s}$ forment un ensemble de repr\'esentants des orbites de l'action de $K_{AD,s}^{\dag}$ dans l'ensemble des droites port\'ees par des \'el\'ements de $fc(\mathfrak{g}_{s}({\mathbb F}_{q}))$. Pour $f\in {\cal B}_{s}$, notons $\Xi_{f}$ l'ensemble des caract\`eres de $G_{AD}(F)/\pi(G(F))$ dont la restriction \`a $K_{AD,s}^{\dag}(f)$ co\"{\i}ncide avec $\xi_{f}$. Pour tout $\xi\in \Xi_{f}$, on d\'efinit la fonction
 $$f_{\xi}=\sum_{g\in \pi(G(F))K_{AD,s}^{\dag}(f)\backslash G_{AD}(F)}\xi^{-1}(g) g\star f.$$
 Alors $FC(\mathfrak{g}(F))$ a pour base les $f_{\xi}$ quand $s$ d\'ecrit $\underline{S}(G_{AD})$, $f$ d\'ecrit ${\cal B}_{s}$ et $\xi$ d\'ecrit $\Xi_{f}$. Cette base est orthogonale pour le produit scalaire elliptique.
 
 Remarquons que la fonction $f_{\xi}$ appartient \`a $I_{cusp,\xi}(\mathfrak{g}(F))$. En utilisant les notations de (4) pour la fonction $f$ et en utilisant (4) et (5), on voit aussi que  l'\'el\'ement $\xi={\bf 1}$ ne peut appartenir \`a $\Xi_{f}$ que si $\epsilon$ vaut $1$ sur  l'image de $Z(G)^{I_{F}}$ dans $Z_{G_{s}}(N)$.

\subsection{Donn\'ees endoscopiques}\label{donneesendoscopiques}

Pour ce paragraphe, on l\`eve les hypoth\`eses sur $G$ en supposant seulement que c'est un groupe r\'eductif connexe d\'efini sur $F$. On fixe un sous-groupe de Borel $\hat{B}$ de $\hat{G}$ et un sous-tore maximal $\hat{T}$ de $\hat{B}$, tous deux conserv\'es par l'action galoisienne $\sigma\mapsto \sigma_{G}$. On note $Endo_{ell}(G)$ l'ensemble des classes d'\'equivalence de donn\'ees endoscopiques elliptiques de $G$ (on renvoie \`a \cite{MW} I.1.5 pour les d\'efinitions). Soit ${\bf G}'=(G',s,{\cal G}')$ une telle donn\'ee. Le groupe $\hat{G}'$ s'identifie \`a $Z_{\hat{G}}(s)^0$. Il est muni d'une action galoisienne $\sigma\mapsto \sigma_{G'}$. A \'equivalence pr\`es, on peut supposer et on suppose que $s\in \hat{T}$ et qu'il existe un sous-groupe de Borel $\hat{B}'$ de $\hat{G}'$ contenant $\hat{T}$ et un \'epinglage $\hat{\cal E}'$ de la paire $(\hat{B}',\hat{T})$ de sorte que $\hat{B}'$, $\hat{T}$ et $\hat{{\cal E}}'$  soient conserv\'ees par cette action galoisienne. On suppose fix\'e un facteur de transfert $\Delta^{{\bf G}'}$ unitaire sur $\mathfrak{g}'(F)\times \mathfrak{g}(F)$.

{\bf Remarque.} Plus pr\'ecis\'ement, ce facteur est d\'efini sur le sous-ensemble des couples $(X',X)$  d'\'el\'ements semi-simples tels que $X'$ soit $G$-r\'egulier et $X$ soit r\'egulier (rappelons que $X'$ est dit $G$-r\'egulier si et seulement si la classe de conjugaison stable dans $\mathfrak{g}(F)$ correspondant \`a celle de $X'$ est form\'ee d'\'el\'ements r\'eguliers). Pour ce que nous allons faire, on peut consid\'erer que $\Delta^{{\bf G}'}$ est nul hors  de ce sous-ensemble.
\bigskip

Pour $f'\in C_{c}^{\infty}(\mathfrak{g}'(F))$ et $X'\in \mathfrak{g}'_{reg}(F)$, on d\'efinit l'int\'egrale orbitale stable $S^{G'}(X',f')=\sum_{Y'}I^{G'}(Y',f')$ o\`u $Y'$ parcourt les \'el\'ements stablement conjugu\'es de $X'$ \`a conjugaison pr\`es par $G'(F)$. On note $SI(\mathfrak{g}'(F))$ le quotient de $C_{c}^{\infty}(\mathfrak{g}'(F))$ par le sous-espace des $f'$ telles que $S^{G'}(X',f')=0$ pour tout $X'\in \mathfrak{g}'_{reg}(F)$. On note $SI_{cusp}(\mathfrak{g}'(F))$ l'image de $I_{cusp}(\mathfrak{g}'(F))$ dans $SI(\mathfrak{g}'(F))$. 

On note $Aut({\bf G}')$ le groupe des $x\in \hat{G}$ tels que $xsx^{-1}\in Z(\hat{G})s$ et $x{\cal G}'x^{-1}={\cal G}'$ (un tel $x$ est appel\'e un automorphisme de ${\bf G}'$). Sa composante neutre est \'egale \`a $\hat{G}'$ et on pose $Out({\bf G}')=Aut({\bf G}')/\hat{G}'$ (ce groupe est appel\'e celui des automorphismes ext\'erieurs de ${\bf G}'$). Soit $x\in Aut({\bf G}')$. Quitte \`a multiplier $x$ par un \'el\'ement de $\hat{G}'$, on peut supposer que $Ad(x)$ conserve $\hat{{\cal E}}'$. Alors $Ad(x)$ commute \`a l'action $\sigma\mapsto \sigma_{G'}$. Fixons une paire de Borel \'epingl\'ee ${\cal E}'$ de $G'$ conserv\'ee par l'action galoisienne sur $G'$. C'est possible puisque $G'$ est quasi-d\'eploy\'e. Il correspond alors \`a $Ad(x)$ un  automorphisme $\alpha_{x}$ de $G'$ d\'efini sur $F$ et conservant ${\cal E}'$.   Fixons un \'el\'ement $s_{sc}\in \hat{G}_{SC}$ qui a m\^eme image dans $\hat{G}_{AD}$ que $s$. Puisque $xsx^{-1}\in Z(\hat{G})s$, il existe $z_{x}\in Z(\hat{G}_{SC})$ tel que $xs_{sc}x^{-1}=z_{x}s_{sc}$. On v\'erifie que $z_{x}$ est fixe par l'action galoisienne (toutes les actions galoisiennes co\"{\i}ncident sur le groupe $Z(\hat{G}_{SC})$). Rappelons que les deux groupes $Z(\hat{G}_{SC})^{\Gamma_{F}}$ et $H^1(\Gamma_{F};G^*_{AD})$ sont en dualit\'e, on note $(.,.)$ l'accouplement en question.
 Soient $X'\in \mathfrak{g}'(F)$ et $X\in \mathfrak{g}(F)$ deux \'el\'ements semi-simples qui se correspondent, avec $X\in \mathfrak{g}_{reg}(F)$. Alors $\alpha_{x}(X')$ et $X$ se correspondent aussi. Rappelons (cf. \ref{orbites}) qu'\`a $G$ est associ\'e un \'el\'ement $u_{G}\in H^1(\Gamma_{F};G^*_{AD})$. 

\begin{lem}{Sous ces hypoth\`eses, on a l'\'egalit\'e
$$\Delta^{{\bf G}'}(\alpha_{x}(X'),X)=(z_{x},u_{G})^{-1}\Delta^{{\bf G}'}(X',X).$$}\end{lem}

C'est un exercice sur les facteurs de transfert.

On d\'efinit une action de $\underline{Aut}({\bf G}')$ sur $SI(\mathfrak{g}'(F))$ par la formule $x(f)(X')=(z_{x},u_{G})^{-1}f(\alpha_{x}^{-1}(X'))$. Elle se quotiente en une action de $Out({\bf G}')$. 

Le transfert est l'homomorphisme $transfert^{{\bf G}'}:I(\mathfrak{g}(F))\to SI(\mathfrak{g}'(F))$ d\'efini par les \'egalit\'es 
$$SI^{G'}(X',transfert^{{\bf G}'}(f))=\sum_{X}\Delta^{{\bf G}'}(X',X)I^G(X,f)$$
pour tout $X'\in \mathfrak{g}'(F)$ $G$-r\'egulier, o\`u $X$ parcourt $\mathfrak{g}(F)$ modulo conjugaison par $G(F)$. Notons $I_{cusp}(\mathfrak{g}(F),{\bf G}')$ le sous-espace des $f\in I_{cusp}(\mathfrak{g}(F))$ telles que $transfert^{{\bf G}''}(f)=0$ pour tout ${\bf G}''\in Endo_{ell}(G)$, ${\bf G}''\not={\bf G}'$. Alors le transfert se restreint en un isomorphisme
$$(1) \qquad transfert^{{\bf G}'}:I_{cusp}(\mathfrak{g}(F),{\bf G}')\simeq SI_{cusp}(\mathfrak{g}'(F))^{Out({\bf G}')},$$
o\`u l'exposant $Out({\bf G}')$ signifie, comme c'est l'usage, que l'on prend les invariants par l'action de ce groupe. Cet isomorphisme est une similitude pour les produits scalaires elliptiques.  Pr\'ecis\'ement, on a l'\'egalit\'e
$$(f,f')_{ell}=c(G,{\bf G}')(transfert^{{\bf G}'}(f),transfert^{{\bf G}'}(f'))_{ell}$$
pour tous $f,f'\in I_{cusp}(\mathfrak{g}(F),{\bf G}')$, o\`u $c(G,{\bf G}')$ est la constante d\'efinie en \cite{MW} I.4.17.

 On sait que la donn\'ee ${\bf G}'$ d\'etermine un caract\`ere $\xi_{{\bf G}'}\in \Xi$  tel que, pour deux \'el\'ements semi-simples $X'\in \mathfrak{g}'(F)$ et $X\in \mathfrak{g}(F)$ qui se correspondent, avec $X\in \mathfrak{g}_{reg}(F)$, on  ait l'\'egalit\'e $\Delta^{{\bf G}'}(X',g^{-1}Xg)=\xi_{{\bf G}'}(g)^{-1}\Delta^{{\bf G}'}(X',X)$ pour tout $g\in G_{AD}(F)$.  On a rappel\'e le calcul de ce caract\`ere en \cite{MW} I.2.7 (on a $\xi_{{\bf G}'}=\omega_{\sharp}^{-1}$ avec les notations de cette r\'ef\'erence). On fixe $s_{sc}\in \hat{G}_{SC}$ d'image $s$ dans $\hat{G}=\hat{G}_{AD}$. Pour $\sigma\in \Gamma_{F}$, l'automorphisme $\sigma_{G'}$ de $\hat{T}$ se rel\`eve en un automorphisme encore not\'e $\sigma_{G'}$ de l'image r\'eciproque de $\hat{T}$ dans $\hat{G}_{SC}$. L'application $\sigma\mapsto s_{sc}^{-1}\sigma_{G'}(s_{sc})$ prend ses valeurs dans $Z(\hat{G}_{SC})$ et d\'efinit un \'el\'ement de $H^1(\Gamma_{F},Z(\hat{G}_{SC}))$ que l'on restreint en   un \'el\'ement de $H^1(W_{F},Z(\hat{G}_{SC}))$. Ce dernier groupe s'identifie \`a  $\Xi$ par l'homomorphisme de Langlands. Le caract\`ere $\xi_{{\bf G}'}$ est l'\'el\'ement de $\Xi$ auquel s'identifie le cocycle pr\'ec\'edent. Remarquons que l'on a l'inclusion
$I_{cusp}(\mathfrak{g}(F),{\bf G}')\subset I_{cusp,\xi_{{\bf G}'}}(\mathfrak{g}(F))$.

On a l'\'egalit\'e 
$$(2)\qquad I_{cusp}(\mathfrak{g}(F))=\oplus_{{\bf G}'\in Endo_{ell}(G)}I_{cusp}(\mathfrak{g}(F),{\bf G}').$$
Cette d\'ecomposition est orthogonale pour le produit scalaire elliptique. Le lien entre cette d\'ecomposition et l'\'egalit\'e (1) du paragraphe \ref{lespaceFC} est l'\'egalit\'e

(3) $I_{cusp,\xi}(\mathfrak{g}(F))=\oplus_{{\bf G}'\in Endo_{ell}(G); \xi_{{\bf G}'}=\xi}I_{cusp}(\mathfrak{g}(F),{\bf G}')$,

\noindent pour tout $\xi\in \Xi$.

Pour ${\bf G}'\in Endo_{ell}(G)$, on pose $FC(\mathfrak{g}(F),{\bf G}')=FC(\mathfrak{g}(F))\cap I_{cusp}(\mathfrak{g}(F),{\bf G}')$. L'ensemble $Endo_{ell}(G)$ contient une donn\'ee "principale" ${\bf G}=(G^*,1,{^LG})$. Dans le cas o\`u $G$ est quasi-d\'eploy\'e, on pose simplement $I_{cusp}^{st}(\mathfrak{g}(F))=I_{cusp}(\mathfrak{g}(F),{\bf G})$ et $FC^{st}(\mathfrak{g}(F))=FC(\mathfrak{g}(F),{\bf G})$. L'espace $I_{cusp}^{st}(\mathfrak{g}(F))$ s'identifie \`a $SI_{cusp}(\mathfrak{g}(F))$. On a d\'emontr\'e dans \cite{W1} les analogues de (1) et (2) pour nos espaces $FC(\mathfrak{g}(F))$, \`a savoir

- pour tout ${\bf G}'\in Endo_{ell}(G)$,  le transfert se restreint en un isomorphisme
$$transfert^{{\bf G}'}:FC(\mathfrak{g}(F),{\bf G}')\simeq FC^{st}(\mathfrak{g}'(F))^{Out({\bf G}')},$$

- on a l'\'egalit\'e
 $$FC(\mathfrak{g}(F))=\oplus_{{\bf G}'\in Endo_{ell}(G)}FC(\mathfrak{g}(F),{\bf G}').$$
 
 \noindent Il est commode de poser
 $$FC^{{\cal E}}(\mathfrak{g}(F))=\oplus_{{\bf G}'\in Endo_{ell}(G)}FC^{st}(\mathfrak{g}'(F))^{{\bf G}'}.$$
 On munit cet espace d'un produit scalaire elliptique en posant
 $$(\oplus_{{\bf G}'}f_{{\bf G}'},\oplus_{{\bf G}'}f'_{{\bf G}'})_{ell}=\sum_{{\bf G}'}c(G,{\bf G}')(f_{{\bf G}'},f'_{{\bf G}'})_{ell},$$
 pour tous $\oplus_{{\bf G}'}f_{{\bf G}'},\oplus_{{\bf G}'}f'_{{\bf G}'}\in FC^{{\cal E}}(\mathfrak{g}(F))$. On note $transfert:FC(\mathfrak{g}(F))\to FC^{{\cal E}}(\mathfrak{g}(F))$ la somme sur ${\bf G}'$ des applications $transfert^{{\bf G}'}$. Alors
 
 (4) l'application $transfert$ est une isom\'etrie de $FC(\mathfrak{g}(F))$ sur $FC^{{\cal E}}(\mathfrak{g}(F))$.
 
 \subsection{Description des donn\'ees endoscopiques elliptiques}\label{description}
 On impose de nouveau dans ce paragraphe et le suivant que $G$ est simplement connexe et absolument quasi-simple. On note $\hat{\Sigma}$ l'ensemble des racines de $\hat{T}$ dans $\hat{G}$ et $\hat{\Sigma}^+$, resp.    $\hat{\Delta}$, le sous-ensemble de racines positives, resp. simples, d\'etermin\'e par $\hat{B}$.  Notons $\hat{\cal D}$ le diagramme de Dynkin dont l'ensemble de sommets est $\hat{\Delta}$. Les hypoth\`eses sur $G$ impliquent que $\hat{G}$ est adjoint et que $\hat{\cal D}$ est  connexe.  Notons $\hat{\alpha}_{0}$ l'oppos\'ee de la plus grande racine dans $\hat{\Sigma}^+$ et $\hat{\Delta}_{a}=\{\hat{\alpha}_{0}\}\cup \hat{\Delta}$. En \'ecrivant $-\hat{\alpha}_{0}$ dans la base $\hat{\Delta}$, on obtient une relation
$$(1) \qquad \sum_{\hat{\alpha}\in \hat{\Delta}_{a}}d(\hat{\alpha})\hat{\alpha}=0$$
o\`u $d(\hat{\alpha}_{0})=1$ et, pour $\hat{\alpha}\in \hat{\Delta}$, $d(\hat{\alpha})$ appartient \`a l'ensemble ${\mathbb N}_{>0}$ des entiers strictement positifs. On sait que   l'espace des relations entre les \'el\'ements de $\hat{\Delta}_{a}$ est la droite port\'ee par la relation (1). Cela implique que la relation (1) est la seule relation lin\'eaire entre les \'el\'ements de $\hat{\Delta}_{a}$ dont les coefficients sont entiers relatifs et dont au moins un coefficient vaut $1$.

Notons $\hat{\cal D}_{a}$ le diagramme de Dynkin compl\'et\'e dont l'ensemble de sommets est $\hat{\Delta}_{a}$. Notons $Aut(\hat{\cal D})\subset Aut(\hat{\cal D}_{a})$ les groupes d'automorphismes de $\hat{\cal D}$ et $\hat{\cal D}_{a}$.  Remarquons que, puisque $\hat{G}$ est semi-simple,  le groupe $Aut(\hat{\cal D}_{a})$ se plonge naturellement dans le groupe d'automorphismes de $\hat{T}$. Puisque $(\hat{B},\hat{T})$ est conserv\'ee par l'action galoisienne, l'action galoisienne se restreint en une action sur $\hat{\Delta}$ et $\hat{\Delta}_{a}$ et peut \^etre consid\'er\'ee comme un homomorphisme $\Gamma_{F}\to Aut(\hat{\cal D})$. Comme on le sait, la donn\'ee de l'action galoisienne est d'ailleurs \'equivalente \`a celle de cet homomorphisme. On note $E_{G}$ l'extension galoisienne finie de $F$ telle que $\Gamma_{E_{G}}$ soit le noyau de cet homomorphisme.  Notons $\hat{W}$ le groupe de Weyl de $\hat{G}$ relatif \`a $\hat{T}$ et $\hat{\Omega}$ le sous-groupe des \'el\'ements de $\hat{W}$ qui conservent $\hat{\Delta}_{a}$. L'application qui \`a un \'el\'ement de $\hat{\Omega}$ associe son action sur $\hat{\cal D}_{a}$ identifie $\hat{\Omega}$ \`a un sous-groupe de $Aut(\hat{\cal D}_{a})$. On sait que c'est un sous-groupe ab\'elien distingu\'e de $Aut(\hat{\cal D}_{a})$ et que $Aut(\hat{\cal D}_{a})$ est le produit semi-direct de $\hat{\Omega}$ par $Aut(\hat{\cal D})$, cf. \cite{B} VI.4.3. 
Signalons que, pour $\tau\in Aut(\hat{\cal D}_{a})$, 
on a $d(\tau(\hat{\alpha}))=d(\hat{\alpha})$ pour tout $\hat{\alpha}\in \hat{\Delta}_{a}$:  en appliquant $\tau^{-1}$ \`a (1), on obtient
$\sum_{\hat{\alpha}\in \hat{\Delta}_{a}}d(\tau(\hat{\alpha}))\hat{\alpha}=0$; c'est une relation entre les \'el\'ements de $\hat{\Delta}_{a}$  \`a coefficients entiers dont au moins l'un d'eux vaut $1$ donc c'est la relation (1).

Notons ${\cal E}_{ell}(G)$ l'ensemble des couples $(\sigma\mapsto\sigma_{G'}, {\cal O})$ o\`u 

$\sigma\mapsto \sigma_{G'}$ est un homomorphisme de $\Gamma_{F}$ dans $Aut(\hat{\cal D}_{a})$ tel qu'il existe une application
$\omega_{G'}:\Gamma_{F}\to \hat{\Omega}$ de sorte que $ \sigma_{G'}=\omega_{G'}(\sigma)\sigma_{G}$  pour tout $\sigma\in \Gamma_{F}$; 

${\cal O}$ est un sous-ensemble non vide de $\hat{\Delta}_{a}$ qui est conserv\'e par l'action $\sigma\mapsto \sigma_{G'}$ et qui forme une unique orbite pour cette action.

Disons que deux tels couples $(\sigma\mapsto \sigma_{G'_{1}},{\cal O}_{1})$ et $(\sigma\mapsto \sigma_{G'_{2}},{\cal O}_{2})$ sont \'equivalents si et seulement s'il existe $\omega\in \hat{\Omega}$ tel que ${\cal O}_{2}=\omega({\cal O}_{1})$ et $\sigma_{G'_{2}}=\omega\sigma_{G'_{1}}\omega^{-1}$ pour tout $\sigma\in \Gamma_{F}$. Notons $E_{ell}(G)$ l'ensemble des classes d'\'equivalence dans ${\cal E}_{ell}(G)$.

Nous allons associer \`a tout \'el\'ement de ${\cal E}_{ell}(G)$ une donn\'ee endoscopique de $G$. La construction d\'epend du choix de racines de l'unit\'e dans ${\mathbb C}^{\times}$. Pr\'ecis\'ement, pour tout entier $d\geq1$, on fixe un \'el\'ement $\zeta_{d}\in {\mathbb C}^{\times}$ qui est une racine de l'unit\'e d'ordre $d$. 

 Consid\'erons un couple $(\sigma\mapsto \sigma_{G'},{\cal O})\in {\cal E}_{ell}(G)$.  Posons $d=\sum_{\hat{\alpha}\in {\cal O}}d(\hat{\alpha})$. Soit $s$ l'unique \'el\'ement de $\hat{T}$ tel que $\hat{\alpha}(s)=1$ pour $\hat{\alpha}\in \hat{\Delta}-({\cal O}\cap \hat{\Delta})$ et $\hat{\alpha}(s)=\zeta_{d}$ pour $\hat{\alpha}\in {\cal O}\cap \hat{\Delta}$. La relation   (1) implique que ces \'egalit\'es s'\'etendent \`a $\hat{\Delta}_{a}$, c'est-\`a-dire $\hat{\alpha}(s)=1$ pour $\hat{\alpha}\in \hat{\Delta}_{a}-{\cal O}$ et $\hat{\alpha}(s)=\zeta_{d}$ pour $\hat{\alpha}\in {\cal O}$. 
  Puisque ${\cal O}$ est stable par l'action galoisienne $\sigma\mapsto \sigma_{G'}$, le point $s$ est fixe par cette action. 
Posons $\hat{G'}=Z_{\hat{G}}(s)^0$. Ce groupe contient $\hat{T}$. Notons ${\cal G}'$ le sous-ensemble de $^LG$ form\'e des \'el\'ements $(g\dot{\omega}_{G'}(\sigma),\sigma)$ pour $g\in \hat{G'}$ et $\sigma\in \Gamma_{F}$, o\`u $\dot{\omega}_{G'}(\sigma)$ est un repr\'esentant quelconque de $\omega_{G'}(\sigma)$ dans $Norm_{\hat{G}}(\hat{T})$.   L'ensemble ${\cal G}'$ un  groupe qui normalise $\hat{G}'$. On en d\'eduit de la fa\c{c}on habituelle une $L$-action de $\Gamma_{F}$ sur $\hat{G}'$ et on introduit le groupe r\'eductif connexe $G'$ d\'efini sur $F$ et quasi-d\'eploy\'e dont $\hat{G}'$, muni de cette action galoisienne, est le groupe dual. On v\'erifie que le triplet ${\bf G}'=(G',s,{\cal G}')$ est une donn\'ee endoscopique  elliptique de $G$.  Il existe un sous-groupe de Borel $\hat{B}'$ de $\hat{G}'$ contenant $\hat{T}$ de sorte que $\hat{\Delta}_{a}-{\cal O}$ soit l'ensemble de racines simples de $\hat{T}$ dans $\hat{G}'$ associ\'e \`a $\hat{B}'$. 

Le groupe d'automorphismes ext\'erieurs $Out({\bf G}')$ est \'egal \`a celui des $\omega\in \hat{\Omega}$ tels que ${\cal O}=\omega({\cal O})$ et $\sigma_{G'}=\omega\sigma_{G'}\omega^{-1}$ pour tout $\sigma\in \Gamma_{F}$. 

Notons $E_{G'}$ l'extension galoisienne de $F$ telle que $\Gamma_{E_{G'}}$ soit le noyau de l'action $\sigma\mapsto \sigma_{G'}$. On a prouv\'e en \cite{LW} 1.4(4) que

(2) $E_{G}\subset E_{G'}$ et la restriction de $\omega_{G'}$ \`a $\Gamma_{E_{G}}$ se quotiente en un homomorphisme injectif de $\Gamma_{E_{G'}/E_{G}}$ dans $\hat{\Omega}$.

On a prouv\'e en \cite{LW}  2.3 que l'application qui, \`a $(\sigma\mapsto \sigma_{G'},{\cal O})\in {\cal E}_{ell}(G)$, associe ${\bf G}'$, se quotiente en une bijection de $E_{ell}(G)$ sur $Endo_{ell}(G)$. 

 \subsection{A propos du centre de $G'$}\label{centre}
 Consid\'erons un \'el\'ement $(\sigma\mapsto \sigma_{G'},{\cal O})$ de ${\cal E}_{ell}(G)$, notons ${\bf G}'$ la donn\'ee endoscopique elliptique qui lui correspond. Fixons $\hat{\alpha}\in {\cal O}$ et notons $E_{{\cal O}}$ l'extension de $F$ telle que $\Gamma_{E_{{\cal O}}}$ soit l'ensemble des $\sigma\in \Gamma_{F}$ tels que $\sigma_{G'}(\hat{\alpha})=\hat{\alpha}$.  Puisque ${\cal O}$ forme une unique orbite pour l'action galoisienne, la classe de conjugaison de $E_{{\cal O}}$ par $\Gamma_{F}$ ne d\'epend pas du choix de $\hat{\alpha}$.
 
 \begin{lem}{Si $E_{{\cal O}}/F$ n'est pas totalement ramifi\'ee, $FC^{st}(\mathfrak{g}'(F))=\{0\}$. Si $E_{{\cal O}}/F$ est totalement ramifi\'ee, $FC^{st}(\mathfrak{g}'(F))\simeq FC^{st}(\mathfrak{g}'_{SC}(F))$.}\end{lem}
 
 Preuve. Notons $A_{G'}^{nr}$ le plus grand sous-tore central de $G'$ d\'eploy\'e sur $F^{nr}$. D'apr\`es \cite{W1}, remarque du paragraphe 11, il suffit de prouver que  $A_{G'}^{nr}=\{1\}$ si et seulement si $E_{{\cal O}}/F$ est totalement ramifi\'ee. Pour tout ensemble $V$, notons ${\mathbb Q}[V]$ l'espace vectoriel sur ${\mathbb Q}$ de base $V$. Posons $X^*(\hat{T})_{{\mathbb Q}}=X^*(\hat{T})\otimes_{{\mathbb Z}}{\mathbb Q}$. On a un homomorphisme naturel $P:{\mathbb Q}[\hat{\Delta}_{a}]\to X^*(\hat{T})_{{\mathbb Q}}$ dont le noyau est la droite engendr\'ee par la relation \ref{description} (1).  On a l'\'egalit\'e $P({\mathbb Q}[\hat{\Delta}_{a}]^{I_{F}})=X^*(\hat{T})_{{\mathbb Q}}^{I_{F}}$. Puisque $\hat{\Delta}_{a}-{\cal O}$ est une base de l'ensemble de racines de $\hat{T}$ dans $\hat{G}'$, la
   condition $A_{G'}^{nr}=\{1\}$ \'equivaut \`a ce que $X^*(\hat{T})_{{\mathbb Q}}^{I_{F}}$ soit contenu dans  $P({\mathbb Q}[\hat{\Delta}_{a}-{\cal O}]^{I_{F}})$. Cela \'equivaut \`a ce que ${\mathbb Q}[\hat{\Delta}_{a}-{\cal O}]^{I_{F}}$ soit de codimension $1$ dans ${\mathbb Q}[\hat{\Delta}_{a}]^{I_{F}}$. Ou encore \`a ce que ${\mathbb Q}[{\cal O}]^{I_{F}}$ soit de dimension $1$. Or la repr\'esentation de $\Gamma_{F}$ dans ${\mathbb Q}[{\cal O}]$ est l'induite de $\Gamma_{E_{{\cal O}}}$ \`a $\Gamma_{F}$ de la repr\'esentation triviale de $\Gamma_{E_{{\cal O}}}$. La dimension de ${\mathbb Q}[{\cal O}]^{I_{F}}$ est le nombre d'\'el\'ements de l'ensemble de doubles classes $\Gamma_{E_{{\cal O}}}\backslash \Gamma_{F}/I_{F}$. Ce nombre vaut $1$ si et seulement si $E_{{\cal O}}/F$ est totalement ramifi\'ee. $\square$
   
   \subsection{Extensions di\'edrales}\label{extensionsdiedrales}
   
    Soit $K/F$ une  extension galoisienne de degr\'e $n$ premier \`a $p$. Rappelons que, si  $K/F$ est non ramifi\'ee ou totalement ramifi\'ee, $\Gamma_{K/F}$ est cyclique.  Il existe une telle extension totalement ramifi\'ee  si et seulement si $F$ contient les racines $n$-i\`emes de l'unit\'e, c'est-\`a-dire  $\delta_{n}(q-1)=1$. Dans ce cas, il y a $n$ telles extensions $K=F(\alpha)$ o\`u $\alpha^n=\varpi_{F}u$, $u$ parcourt $\mathfrak{o}_{F}^{\times}/\mathfrak{o}_{F}^{\times,n}$. On note $E_{0}$ l'extension quadratique non ramifi\'ee de $F$ et $Q$ l'extension biquadratique de $F$. 

Fixons une extension quadratique $E/F$. On consid\`ere une extension galoisienne $K/F$ contenant $E$ et telle que $K/E$ soit cyclique de degr\'e $n>1$. Fixons un g\'en\'erateur $\rho$ de $\Gamma_{K/E}$ et un \'el\'ement $\delta\in \Gamma_{K/F}-\Gamma_{K/E}$. On impose que ces \'el\'ements v\'erifient les \'egalit\'es $\delta^2=1$ et $\delta\rho=\rho^{-1}\delta$. Remarquons que

(1) L'extension $K/F$ n'est pas cyclique: 

En effet, si elle est cyclique,  elle est commutative. L'\'egalit\'e $\delta\rho=\rho^{-1}\delta$ implique  $\rho=\rho^{-1}$ et $n=2$. Il y a deux \'el\'ements distincts d'ordre $2$, $\rho$ et $\delta$. C'est impossible dans le cas cyclique.

On a aussi

(2) la plus grande extension non ramifi\'ee de $F$ contenue dans $K$ est au plus d'ordre $2$.

Preuve. Soit $H$ cette extension et $p:\Gamma_{K/F}\to \Gamma_{H/F}$ la projection. Puisque ce dernier groupe  est commutatif, l'\'egalit\'e $p(\delta\rho)=p(\rho^{-1}\delta)$ entra\^{\i}ne $p(\rho)^2=1$. On a aussi $p(\delta)^2=1$. 
  Puisque $\Gamma_{H/F}$ est cyclique et est engendr\'e par $p(\rho)$ et $p(\delta)$, il est au plus d'ordre $2$.

  Supposons $E=E_{0}$. Prouvons que

(3) une extension $K/F$ comme ci-dessus existe si et seulement si $\delta_{n}(q+1)=1$; dans ce cas $K$ est unique  et $K/E_{0}$ est totalement ramifi\'ee; on a K=$E_{0}(\alpha)$ o\`u $\alpha^n=\varpi_{F}$

Preuve. D'apr\`es (2), $K/E_{0}$ est totalement ramifi\'ee.  Donc $K=E_{0}(\alpha)$ avec $\alpha^n=u\varpi_{F}$ pour un  $u\in \mathfrak{o}_{E}^{\times}$ et $E_{0}$ contient les racines $n$-i\`emes de l'unit\'e. On a $\rho(\alpha)=\zeta \alpha$ pour  une racine primitive $\zeta$ de $1$ d'ordre $n$. On a aussi $\delta(\alpha)=v\alpha$ pour un $v\in \mathfrak{o}_{K}^{\times}$. On a $\delta(\alpha)^n=\delta(u)\varpi_{F}$, donc $v^n=\delta(u)u^{-1}$. Cette \'equation d\'efinissant une extension non ramifi\'ee de $E_{0}$, puisqu'elle a des solutions dans $K$, ces solutions sont dans $E_{0}$ donc $v\in \mathfrak{o}_{E_{0}}^{\times}$. On a $\delta\rho(\alpha)=\zeta^qv\alpha$, $\rho^{-1}\delta(\alpha)=\zeta^{-1}v\alpha$ donc $\zeta^{q+1}=1$, ce qui implique que $n$ divise $q+1$.  
  L'\'egalit\'e $\delta^2=1$ entra\^{\i}ne $v\delta(v)=1$ donc (Hilbert 90) il existe $x\in \mathfrak{o}_{E_{0}}^{\times}$ tel que $v=\delta(x)x^{-1}$. On a alors $\delta(ux^{-n})=ux^{-n}$, c'est-\`a-dire $ux^{-n}\in \mathfrak{o}_{F}^{\times}$. On peut remplacer $u$ par $ux^{-n}$ et supposer $u\in \mathfrak{o}_{F}^{\times}$. On a aussi $\mathfrak{o}_{F}^{\times}\subset \mathfrak{o}_{E_{0}}^{\times,n}$ car il en est ainsi en r\'eduction: la norme est surjective sur ${\mathbb F}_{q}^{\times}$ et, puisque $n$ divise $q+1$, toute norme est une puissance $n$-i\`eme. On peut donc remplacer $u$ par $1$ et alors $\alpha^n=\varpi_{F}$. Inversement, si $\delta_{n}(q+1)=1$, on voit facilement que l'extension $K$ indiqu\'ee v\'erifie 
  les conditions requises.

  Supposons maintenant que $E/F$ est ramifi\'ee. Alors
  
(4)  on a $n=2$ et $K=Q$.

Preuve. Notons $H$ la plus grande sous-extension non ramifi\'ee de $K$ et $p:\Gamma_{K/F}\to \Gamma_{H/F}$ la projection naturelle. Si $H=F$, $K/F$ est totalement ramifi\'ee donc cyclique. Puisque $[K:F]=2n$, $\Gamma_{K/F}$ contient un \'el\'ement d'ordre $2n$. Or un tel \'el\'ement ne peut pas \^etre dans $\Gamma_{K/E}$ et tout \'el\'ement de $\Gamma_{K/F}-\Gamma_{K/E}$ est d'ordre $2$. Donc $n=1$ ce qui est contraire \`a notre hypoth\`ese. Alors (2) entra\^{\i}ne que  $H/F$ est de degr\'e $2$, c'est-\`a-dire $H=E_{0}$.  On a $E_{0}\not=E$ puisque $E/F$ est ramifi\'ee donc $p(\rho)\not=1$, c'est-\`a-dire que $p(\rho)$ est l'unique \'el\'ement non trivial de $\Gamma_{E_{0}/F}$. Donc l'un des \'el\'ements $\delta$ ou $\delta\rho^{-1}$ appartient \`a $\Gamma_{K/E_{0}}$. Notons $\delta'$ cet \'el\'ement. L'\'el\'ement $\rho'=\rho^2$ appartient aussi \`a $\Gamma_{K/E_{0}}$. Ce dernier groupe est cyclique puisque $K/E_{0}$ est totalement ramifi\'ee. On a $\delta'\rho'={\rho'}^{-1}\delta'$ donc $\rho'={\rho'}^{-1}$ et ${\rho'}^2=1$. On a aussi ${\delta'}^2=1$ et $\delta'\not=1$ (car $\delta'$ agit non trivialement sur $E$). La cyclicit\'e entra\^{\i}ne que $\rho'$ est \'egal \`a $\delta'$ ou \`a $1$. Il n'est pas \'egal \`a $\delta'$ car $\rho'$ agit trivialement sur $E$. Donc $\rho'=1$, i.e. $\rho^2=1$. Donc $n=2$. Puisque $\rho^2=1$, le groupe $\Gamma_{K/F}$ est commutatif, isomorphe \`a $({\mathbb Z}/2{\mathbb Z})^2$. Inversement, bien s\^ur, l'extension biquadratique $Q$ v\'erifie nos hypoth\`eses.

\subsection{Extensions biquadratiques}\label{extensionsbiquadratiques}

On conserve les notations du paragraphe pr\'ec\'edent. Notons $Q_{0}$ l'extension biquadratique de $E_{0}$. On a

(1) l'extension $Q_{0}/F$ est galoisienne et $\Gamma_{Q_{0}/F}\simeq ({\mathbb Z}/4{\mathbb Z})\times ({\mathbb Z}/2{\mathbb Z})$.

En effet, $Q_{0}$ est la compos\'ee de l'extension non ramifi\'ee de degr\'e $4$ de $F$, dont le groupe de Galois est cyclique, et de l'extension $F(\sqrt{\varpi_{F}})/F$.  

Soit $E/F$ une extension quadratique ramifi\'ee. Notons $Q_{E}$ l'extension biquadratique de $E$. L'extension $Q_{E}/F$ est galoisienne. Elle contient l'extension quadratique $E_{0}$ non ramifi\'ee. Posons $\Gamma_{Q_{E}/E}=\{1,\rho',\rho_{0},\rho''\}$ o\`u $\rho_{0}$ est l'\'el\'ement tel que $\Gamma_{Q_{E}/E_{0}E}=\{1,\rho_{0}\}$. Remarquons que, puisque $E_{0}E/F$ est galoisienne, $\rho_{0}$ est central dans $\Gamma_{Q_{E}/F}$. Notons $K'$, resp. $K''$, la sous-extension de $Q_{E}$ telle que $\Gamma_{Q_{E}/K'}=\{1,\rho'\}$, resp. $\Gamma_{Q_{E}/K''}=\{1,\rho''\}$. Les extensions $K'/F$ et $K''/F$ sont totalement ramifi\'ees. 
Montrons que

(2) si $\delta_{4}(q-1)=1$,  $K'$ et $K''$ sont les deux extensions galoisiennes cycliques d'ordre $4$ de $F$ contenues dans  $Q_{E}$;   l'homomorphisme naturel 
$$\Gamma_{Q_{E}/F}\to \Gamma_{K'/F}\times \Gamma_{E_{0}/F}\simeq ({\mathbb Z}/4{\mathbb Z})\times {\mathbb Z}/2{\mathbb Z}$$
est un isomorphisme;

(3) si $\delta_{4}(q-1)=0$, on peut choisir un \'el\'ement  $\tau\in \Gamma_{Q_{E}/F}-\Gamma_{Q_{E}/E}$ v\'erifiant les conditions suivantes:  $\tau^2=1$; la conjugaison par $\tau$ fixe $\rho_{0}$ et \'echange $\rho'$ et $\rho''$; le sous-corps $E'$ de $Q_{E}$ tel que $\Gamma_{Q_{E}/E'}=\{1,\tau,\rho_{0},\tau\rho_{0}\}$ est l'unique extension quadratique ramifi\'ee de $F$ diff\'erente de $E$; de plus, il n'y a pas d'extensions galoisiennes cycliques d'ordre $4$ de $F$  contenant $E$.

Supposons $\delta_{4}(q-1)=1$. Ecrivons $E=F(\beta)$ avec $\beta^2=\varpi_{F}u$, o\`u $u\in \mathfrak{o}_{F}^{\times}$. Il y a  deux extensions galoisiennes cycliques d'ordre $4$ de $F$ qui sont totalement ramifi\'ees et qui contiennent $E$, \`a savoir $F(\alpha)$ o\`u $\alpha^4=\varpi_{F}u$ ou $\alpha^4=\varpi_{F}uv$, o\`u $v\in \mathfrak{o}_{F}^{\times,2}-\mathfrak{o}_{F}^{\times,4}$. Elles sont forc\'ement contenues dans $Q_{E}$ donc ce sont $K'$ et $K''$. Puisque $K'/F$ est totalement ramifi\'ee et $E_{0}/E$ est non ramifi\'ee, on a $Q_{E}=K'E_{0}$ et la derni\`ere assertion de (2) s'ensuit.  La structure de $\Gamma_{Q_{E}/F}$  entra\^{\i}ne que ce groupe n'admet que deux quotients isomorphes \`a ${\mathbb Z}/4{\mathbb Z}$, c'est-\`a-dire que $Q_{E}/F$ ne contient que   deux sous-extensions galoisiennes cycliques d'ordre $4$ qui sont $K'$ et $K''$. 

Supposons maintenant $\delta_{4}(q-1)=0$. Si $K/F$ est une extension galoisienne cyclique d'ordre $4$ et contenant $E$, elle ne peut pas contenir $E_{0}$, sinon l'homomorphisme naturel $\Gamma_{K/F}\to \Gamma_{E/F}\times \Gamma_{E_{0}/F}$ serait un isomorphisme et $\Gamma_{K/F}$ ne serait pas cyclique. Donc $K/F$ est totalement ramifi\'ee, ce qui est impossible puisque $\delta_{4}(q-1)=0$. Cela prouve la derni\`ere assertion de (3). Si $K'/F$ est galoisienne, elle est cyclique puisqu'elle est totalement ramifi\'ee. On vient de voir que c'est impossible donc $K'/F$ n'est pas galoisienne. 
Il y a donc un \'el\'ement $\tau\in \Gamma_{Q_{E}/F}-\Gamma_{Q_{E}/E}$ tel que $\tau\rho'\tau^{-1}\not=\rho'$. On a forc\'ement 
$\tau\rho'\tau^{-1}=\rho''$ ou $\rho_{0}$ et le cas $\rho_{0}$ est exclu car $\rho_{0}$ est central.   Puisque $\Gamma_{Q_{E}/E}$ est commutatif, l'\'egalit\'e $\tau\rho'\tau^{-1}=\rho''$ est en fait v\'erifi\'ee pour tout $\tau\in \Gamma_{Q_{E}/F}-\Gamma_{Q_{E}/E}$. Consid\'erons un tel $\tau$. L'extension $EE_{0}/F$ est l'extension biquadratique $Q$ de $F$ et son  groupe de Galois est  isomorphe \`a $({\mathbb Z}/2{\mathbb Z})^2$. L'image de $\tau^2$ dans ce groupe est triviale donc $\tau^2\in \Gamma_{Q_{E}/EE_{0}}=\{1,\rho_{0}\}$.  Si $\tau^2=\rho_{0}=\rho'\rho''$, il suffit de remplacer $\tau$ par $\rho'\tau$.  On a alors $\tau^2=1$. Le sous-corps $E'$ de $Q_{E}$ tel que $\Gamma_{Q_{E}/E'}=\{1,\tau,\rho_{0},\tau\rho_{0}\}$ est une extension quadratique de $F$ telle que $Q_{E}$ soit l'extension biquadratique de $E'$. Ce n'est pas $E_{0}$, sinon $\Gamma_{Q_{E}/F}$ serait ab\'elienne d'apr\`es (1). Ce n'est pas $E$ puisque $\tau\not\in \Gamma_{E}$. C'est donc l'autre extension quadratique ramifi\'ee de $F$.

\subsection{Sur certaines extensions de degr\'e $12$}\label{extensions12}

Montrons que

(1) il n'existe pas d'extensions galoisiennes $K/F$ non ab\'eliennes contenant une sous-extension galoisienne $E/F$ de sorte que   $\Gamma_{K/E}\simeq ({\mathbb Z}/2{\mathbb Z})^2$, $\Gamma_{E/F}\simeq {\mathbb Z}/3{\mathbb Z}$.

Preuve. Pour une telle extension, $K/E$ est l'extension biquadratique de $E$. Notons $L$ l'extension quadratique non ramifi\'ee de $E$.    Puisque $E/F$  est galoisienne, $L/F$ l'est aussi et $\Gamma_{K/L}$ est un sous-groupe distingu\'e   d'ordre $2$ de $\Gamma_{K/F}$. Un tel sous-groupe est central. Donc l'action  par conjugaison de $\Gamma_{E/F}$ sur $\Gamma_{K/E}$ fixe l'\'el\'ement non trivial de $\Gamma_{K/L}$. Ou bien elle fixe aussi les deux autres \'el\'ements non triviaux de $\Gamma_{K/E}$, ou bien elle les permute.  L'ordre de cette derni\`ere action ne divise pas $3$ donc l'action de $\Gamma_{E/F}$ sur $\Gamma_{K/E}$ est triviale. Mais alors le groupe $\Gamma_{K/F}$ est ab\'elien, contrairement \`a l'hypoth\`ese. $\square$

\section{Description des espaces $FC(\mathfrak{g}({\mathbb F}_{q}))$}

\subsection{G\'en\'eralit\'es}\label{generalites}
{\bf On suppose dans cette section que $G$ est un groupe semi-simple d\'efini sur ${\mathbb F}_{q}$ et absolument irr\'eductible.} Nous allons rappeler un certain nombre de faits bien connus concernant ces groupes. Une partie d'entre eux vaut pour les groupes de m\^eme type d\'efinis sur le corps $F$ et nous les utiliserons ult\'erieurement dans ce contexte. 

 On  suppose fix\'es un sous-groupe de Borel $B$ de $G$  et un sous-tore maximal $T\subset B$ tous deux conserv\'es par l'action galoisienne. On note $\Delta$ l'ensemble associ\'e de racines simples, ${\cal D}$ le diagramme de Dynkin et $Aut({\cal D})$ son groupe d'automorphismes. On utilise pour ces objets 
  les notations  de \cite{B}. On note $\check{\alpha}$ la coracine associ\'ee \`a une racine $\alpha$. 
  
  Nous rappellerons   la description d'une base de $FC(\mathfrak{g}({\mathbb F}_{q}))$ form\'ee d'\'el\'ements de $fc(\mathfrak{g}({\mathbb F}_{q}))$, c'est-\`a-dire de fonctions caract\'eristiques de faisceaux-caract\`eres cuspidaux \`a supports nilpotents et conserv\'es par l'action galoisienne. Cette description est due \`a Lusztig. On utilisera pour ces fonctions les notations introduites en \ref{groupessurFq}.

 On introduit les fonctions suivantes d\'efinies sur ${\mathbb N}_{>0}$:

la fonction d'Euler $\phi$ telle que, pour $n\in {\mathbb N}_{>0}$,  $\phi(n)$ est le nombre d'\'el\'ements de $\boldsymbol{\zeta}_{n,prim}({\mathbb C})$;

la fonction caract\'eristique $\delta_{\square}$ de l'ensemble des carr\'es dans ${\mathbb N}_{>0}$;

la fonction caract\'eristique $\delta_{\triangle}$ de l'ensemble des entiers de la forme $i(i+1)/2$ pour $i\in {\mathbb N}_{>0}$;

la fonction caract\'eristique $\delta_{2\triangle}$ de l'ensemble des entiers de la forme $i(i+1)$ pour $i\in {\mathbb N}_{>0}$. 

On note $sgn$ l'unique caract\`ere d'ordre $2$ de ${\mathbb F}_{q}^{\times}$.

\subsection{Type $A_{n-1}$ d\'eploy\'e}\label{An-1dep}
On suppose $G$ d\'eploy\'e de type $A_{n-1}$ avec $n\geq2$, c'est-\`a-dire $G_{SC}=SL(n)$. 

Supposons d'abord $G=G_{SC}$.   Pour $z\in \boldsymbol{\zeta}_{n}(\bar{{\mathbb F}}_{q})$, on pose $\underline{\zeta}(z)=\prod_{i=1,...,n-1}\check{\alpha}_{i}(z^{i})$. Alors $\underline{\zeta}$ est un isomorphisme de $\zeta_{n}(\bar{{\mathbb F}}_{q})$ sur $Z(G)$. 

Si $\delta_{n}(q-1)=0$, $FC(\mathfrak{g}({\mathbb F}_{q}))=\{0\}$. 

Supposons $\delta_{n}( q-1)=1$. Soit $N$ un \'el\'ement nilpotent r\'egulier de $\mathfrak{g}(F)$. On a $Z_{G}(N) /Z_{G}(N)^0\simeq Z(G)$. Alors $FC(\mathfrak{g}(F))$ a pour base la famille des fonctions $f_{N,\epsilon}$ quand $\epsilon$ d\'ecrit les caract\`eres d'ordre $n$ de $Z(G)$. En particulier,
 $dim(FC(\mathfrak{g}({\mathbb F}_{q})))=\phi(n)\delta_{n}( q-1)$.

 Consid\'erons maintenant un quotient $G=G_{SC}/{\cal Z}$, o\`u ${\cal Z}$ est un sous-groupe de $Z(G_{SC})$ non trivial. Alors $FC(\mathfrak{g}({\mathbb F}_{q}))=\{0\}$ car aucun des caract\`eres $\epsilon$ ne se factorise par ${\cal Z}$. 
 
 \subsection{Type $A_{n-1}$ unitaire}\label{An-1nondep}
 On suppose $G$ de type $A_{n-1}$ avec $n\geq 3$ et que  le Frobenius $Fr$ agit sur  ${\cal D}$ par l'unique automorphisme non trivial $\theta$ de ce diagramme. La seule diff\'erence avec le cas pr\'ec\'edent est que la condition d'invariance par l'action galoisienne est maintenant $\delta_{n}(q+1)=1$.
 
  Supposons $G=G_{SC}$. Si $\delta_{n}(q+1)=0$, $FC(\mathfrak{g}({\mathbb F}_{q}))=\{0\}$. Si $\delta_{n}( q+1)=1$, $FC(\mathfrak{g}({\mathbb F}_{q}))$ a pour base les $f_{N,\epsilon}$ pour un \'el\'ement nilpotent r\'egulier $N$ fix\'e, $\epsilon$ d\'ecrivant les caract\`eres d'ordre $n$ de $Z(G)$. D'o\`u $dim(FC(\mathfrak{g}({\mathbb F}_{q})))=\phi(n)\delta_{n}( q+1)$. 
  
  Si $G\not=G_{SC}$, $FC(\mathfrak{g}({\mathbb F}_{q}))=\{0\}$.
 
 \subsection{Type $B_{n}$ }\label{Bn}
 On suppose $G$ de type $B_{n}$ (donc d\'eploy\'e) avec $n\geq2$, c'est-\`a-dire $G_{SC}=Spin(2n+1)$.  
 
 Supposons d'abord  $G=G_{SC}$. On a $Z(G)=\{1,\check{\alpha}_{n}(-1)\}$. On a l'\'egalit\'e $dim(FC(\mathfrak{g}({\mathbb F}_{q})))=\delta_{\square}(2n+1)+\delta_{\triangle}(2n+1)$. 
 
 Supposons $\delta_{\square}(2n+1)=1$, c'est-\`a-dire $2n+1=h^2$ pour un $h\in {\mathbb N}$. Les classes de conjugaison par $G(\bar{{\mathbb F}}_{q})$ dans $\mathfrak{g}_{nil}(\bar{{\mathbb F}}_{q})$ sont param\'etr\'ees par les partitions orthogonales de $2n+1$. Consid\'erons un \'el\'ement $N_{\square}\in \mathfrak{g}_{nil}({\mathbb F}_{q})$ dont la partition associ\'ee est $(2h-1,2h-3,...,3,1)$. Le groupe $Z_{G_{AD}}(N_{\square})/Z_{G_{AD}}(N_{\square})^0$ est isomorphe au sous-groupe $({\mathbb Z}/2{\mathbb Z})_{0}^h$  des \'el\'ements de $({\mathbb Z}/2{\mathbb Z})^h$ de somme nulle. Plus pr\'ecis\'ement, \`a l'ensemble des entiers impairs $2l-1$ intervenant dans la partition est naturellement associ\'ee une base $(e_{2l-1})_{l=1,...,h}$ de  $({\mathbb Z}/2{\mathbb Z})^h$. On d\'efinit un caract\`ere $\epsilon_{\square}$ de ce groupe par $\epsilon_{\square}(\sum_{l=1,...,h}x_{2l-1}e_{2l-1})=\prod_{l=1,...,h}(-1)^{lx_{2l-1}}$, les $x_{2l-1}$ \'etant des \'el\'ements de ${\mathbb Z}/2{\mathbb Z}$. On restreint $\epsilon_{\square}$ au sous-groupe  
$Z_{G_{AD}}(N_{\square})/Z_{G_{AD}}(N_{\square})^0$ puis on le rel\`eve en un caract\`ere de $Z_{G}(N_{\square})/Z_{G}(N_{\square})^0$. Alors la fonction $f_{N_{\square},\epsilon_{\square}}$ appartient \`a $FC(\mathfrak{g}({\mathbb F}_{q}))$. 

Supposons $\delta_{\triangle}(2n+1)=1$, c'est-\`a-dire $2n+1=j(j+1)/2$ pour un entier $j\in {\mathbb N}$. On pose $m=[(j+1)/2]$.  On v\'erifie que $m$ est impair. Consid\'erons un \'el\'ement $N_{\triangle}\in \mathfrak{g}_{nil}({\mathbb F}_{q})$ dont la partition associ\'ee est $(2j+3,2j-1,2j-5,...)$. Cette partition a $m$ termes et son dernier terme est $1$ si $j$ est impair et $3$ si $j$ est pair. Le groupe $Z_{G}(N_{\triangle})/Z_{G}(N_{\triangle})^0$ est un groupe non commutatif  qui s'inscrit dans une suite exacte
$$1\to Z(G)\to Z_{G}(N_{\triangle})/Z_{G}(N_{\triangle})^0\to Z_{G_{AD}}(N_{\triangle})/Z_{G_{AD}}(N_{\triangle})^0\simeq ({\mathbb Z}/2{\mathbb Z})_{0}^{m}\to 0$$
Il a une unique repr\'esentation irr\'eductible $\epsilon_{\triangle}$ sur laquelle $Z(G)$ agit par le caract\`ere non trivial de ce groupe. Alors $f_{N_{\triangle},\epsilon_{\triangle}}\in FC(\mathfrak{g}({\mathbb F}_{q}))$. On peut pr\'eciser le support de cette fonction. R\'ealisons le groupe $G_{AD}$ comme le groupe sp\'ecial orthogonal d'un espace $V$ sur ${\mathbb F}_{q}$ muni d'une forme quadratique $Q$. Alors les classes de conjugaison par $G_{AD}({\mathbb F}_{q})$ dans l'orbite g\'eom\'etrique de $N_{\triangle}$ sont naturellement param\'etr\'ees par des familles de formes quadratiques de rang $1$, qui s'identifient \`a des familles $\underline{\eta}=(\eta_{l})_{l=1,...,m}$ d'\'el\'ements de ${\mathbb F}^{\times}_{q}/{\mathbb F}_{q}^{\times,2}$. Plus pr\'ecis\'ement, par les familles $\underline{\eta}=(\eta_{l})_{l=1,...,m}$ qui v\'erifient la condition 

(1) $(-1)^{(m-1)/2}\prod_{l=1,...,m}\eta_{l}=(-1)^{n}det(Q)$ 

\noindent dans ${\mathbb F}^{\times}_{q}/{\mathbb F}_{q}^{\times,2}$. Soit $N_{\underline{\eta}}$ un \'el\'ement dans la classe param\'etr\'ee par $\underline{\eta}$. La classe de conjugaison par $G_{AD}({\mathbb F}_{q})$ de $N_{\underline{\eta}}$ peut \^etre une seule classe de conjugaison par $G({\mathbb F}_{q})$ ou se couper en deux telles classes. Ce dernier cas se produit si et seulement si le fixateur de $N_{\underline{\eta}}$ dans $G_{AD}({\mathbb F}_{q})$ est contenu dans l'image naturelle de $G({\mathbb F}_{q})$ dans ce groupe. Ce fixateur \'etant facile \`a d\'ecrire, on voit que la condition est que tous les $\eta_{l}$ doivent \^etre \'egaux. Compte tenu  de l'\'egalit\'e (1) et de l'imparit\'e de $m$, la famille $\underline{\eta}$ est uniquement d\'etermin\'ee. Supposons que $N_{\triangle}$ appartient \`a la classe param\'etr\'ee par cette famille et notons $N'_{\triangle}$ un \'el\'ement conjugu\'e \`a $N_{\triangle}$ par $G_{AD}({\mathbb F}_{q})$ mais pas par $G({\mathbb F}_{q})$. Puisque $\epsilon$ se restreint en le caract\`ere non trivial de $Z(G)$, la fonction $f_{N_{\triangle},\epsilon_{\triangle}}$ se transforme par l'action de $G_{AD}({\mathbb F}_{q})$ selon le caract\`ere non trivial de ce groupe. Elle est forc\'ement nulle sur toute classe de conjugaison par $G_{AD}({\mathbb F}_{q})$ qui reste une unique classe de conjugaison par $G({\mathbb F}_{q})$. Par contre, ses valeurs aux points $N_{\triangle}$ et $N'_{\triangle}$ sont oppos\'ees. La dimension de $\epsilon_{\triangle}$ \'etant $2^{m-1}$, on peut normaliser $f_{N_{\triangle},\epsilon_{\triangle}}$ par
$f_{N_{\triangle},\epsilon_{\triangle}}(N_{\triangle})=2^{(m-1)/2}$, $f_{N_{\triangle},\epsilon_{\triangle}}(N'_{\triangle})=-2^{(m-1)/2}$ et $f_{N_{\triangle},\epsilon_{\triangle}}(N)=0$ pour tout $N\in \mathfrak{g}_{nil}({\mathbb F}_{q})$ qui n'est pas conjugu\'e \`a $N_{\triangle}$ ou $N'_{\triangle}$ par un \'el\'ement de $G({\mathbb F}_{q})$.

Supposons maintenant $G=G_{AD}$. Alors (quand elle existe) la fonction $f_{N_{\triangle},\epsilon_{\triangle}}$ dispara\^{\i}t puisqu'elle n'est pas invariante par $G_{AD}({\mathbb F}_{q})$ et il reste (quand elle existe) la fonction $f_{N_{\square},\epsilon_{\square}}$. D'o\`u  $dim(FC(\mathfrak{g}({\mathbb F}_{q})))=\delta_{\square}(2n+1)$.

\subsection{Type $C_{n}$}\label{Cn}
 On suppose $G$ de type $C_{n}$ (donc d\'eploy\'e) avec $n\geq2$, c'est-\`a-dire $G_{SC}=Sp(2n)$.  
 
 Supposons d'abord  $G=G_{SC}$. On a $Z(G)=\{1,z \}$, o\`u $z=\prod_{i=1,...,n}\check{\alpha}_{i}((-1)^{i})$. On a l'\'egalit\'e $dim(FC(\mathfrak{g}({\mathbb F}_{q})))=\delta_{\triangle}(n)$. 
 
 Supposons $\delta_{\triangle}(n)=1$, c'est-\`a-dire $2n=h(h+1)$ pour un $h\in {\mathbb N}$. Les classes de conjugaison par $G(\bar{{\mathbb F}}_{q})$ dans $\mathfrak{g}_{nil}(\bar{{\mathbb F}}_{q})$ sont param\'etr\'ees par les partitions symplectiques de $2n$. Consid\'erons un \'el\'ement $N_{\triangle}\in \mathfrak{g}_{nil}({\mathbb F}_{q})$ dont la partition associ\'ee est $(2h,2h-2,...,4,2)$. Le groupe $Z_{G}(N_{\triangle})/Z_{G}(N_{\triangle})^0$ est isomorphe \`a $({\mathbb Z}/2{\mathbb Z})^h$. Plus pr\'ecis\'ement, \`a l'ensemble des entiers pairs $2l$ intervenant dans la partition est naturellement associ\'ee une base $(e_{2l})_{l=1,...,h}$ de  $({\mathbb Z}/2{\mathbb Z})^h$. On d\'efinit un caract\`ere $\epsilon_{\triangle}$ de ce groupe par $\epsilon_{\triangle}(\sum_{l=1,...,h}x_{2l}e_{2l})=\prod_{l=1,...,h}(-1)^{lx_{2l}}$, les $x_{2l}$ \'etant des \'el\'ements de ${\mathbb Z}/2{\mathbb Z}$.   Alors la fonction $f_{N_{\triangle},\epsilon_{\triangle}}$ appartient \`a $FC(\mathfrak{g}({\mathbb F}_{q}))$. 
 
 Supposons maintenant $G=G_{AD}$. Alors (quand elle existe) la fonction $f_{N_{\triangle},\epsilon_{\triangle}}$ survit si et seulement si le caract\`ere $\epsilon_{\triangle}$ est trivial sur $Z(G_{SC})$, ce qui \'equivaut \`a $h\equiv 0,3\,\,mod\,\,4{\mathbb Z}$ ou encore \`a $n$ pair. D'o\`u $dim(FC(\mathfrak{g}({\mathbb F}_{q})))=\delta_{\triangle}(n)\delta_{2}(n)$.
 
 \subsection{Type $D_{n}$ d\'eploy\'e, $n$ pair}\label{Dndeppair}
 On suppose $G$ de type $D_{n}$ d\'eploy\'e, avec $n$ pair et $n\geq2$. On a $G_{SC}=Spin_{dep}(2n)$, l'indice $dep$ signifiant d\'eploy\'e. 
 
 Supposons d'abord $G=G_{SC}$.  On pose $z=\check{\alpha}_{n-1}(-1)\check{\alpha}_{n}(-1)$, $z'=\prod_{i=1,...,n}\check{\alpha}_{i}((-1)^{i})$, $z''=zz'$. On a $Z(G)=\{1,z,z',z''\}\simeq ({\mathbb Z}/2{\mathbb Z})^2$ et $G/\{1,z\}=SO_{dep}(2n)$. On note $\theta$  l'automorphisme d'ordre $2$ de $D$ qui \'echange $\alpha_{n-1}$ et $\alpha_{n}$ et on note de m\^eme l'automorphisme de $G$ qui s'en d\'eduit. On a $\theta(z)=z$, $\theta(z')=z''$. 
 On a l'\'egalit\'e $dim(FC(\mathfrak{g}({\mathbb F}_{q})))=\delta_{\square}(2n)+2\delta_{\triangle}(2n)$. 
 
 Supposons $\delta_{\square}(2n)=1$, c'est-\`a-dire $2n=h^2$ avec $h\in {\mathbb N}$. On construit une fonction $f_{N_{\square},\epsilon_{\square}}$ comme en \ref{Bn}. Cette fonction appartient \`a $FC(\mathfrak{g}({\mathbb F}_{q}))$.
 
 Supposons $\delta_{\triangle}(2n)=1$, c'est-\`a-dire $2n=j(j+1)/2$ pour un entier $j\in {\mathbb N}$. On pose $m=[(j+1)/2]$.  On v\'erifie que $m$ est divisible par $4$.   Notons $\mu'$ et $\mu''$ les deux caract\`eres de $Z(G)$  d\'efinis par $\mu'(z)=\mu''(z) =\mu'(z')=\mu''(z'')=-1$, $\mu'(z'')=\mu''(z')=1$. 
 Consid\'erons un \'el\'ement $N_{\triangle}\in \mathfrak{g}_{nil}({\mathbb F}_{q})$ dont la partition associ\'ee est $(2j+3,2j-1,2j-5,...)$. Cette partition a $m$ termes et son dernier terme est $1$ si $j$ est impair et $3$ si $j$ est pair. Le groupe $Z_{G}(N_{\triangle})/Z_{G}(N_{\triangle})^0$ est un groupe non commutatif  qui s'inscrit dans une suite exacte
$$1\to \{1,z\}\to Z_{G}(N_{\triangle})/Z_{G}(N_{\triangle})^0\to Z_{SO_{dep}(2n)}(N_{\triangle})/Z_{SO_{dep}(2n)}(N_{\triangle})^0\simeq ({\mathbb Z}/2{\mathbb Z})_{0}^{m}\to 0$$
Il a une unique repr\'esentation irr\'eductible $\epsilon'_{\triangle}$, resp. $\epsilon''_{\triangle}$, sur laquelle $Z(G)$ agit par le caract\`ere  $\mu'$, resp. $\mu''$. Alors $f_{N_{\triangle},\epsilon'_{\triangle}}$ et $f_{N_{\triangle},\epsilon''_{\triangle}}$ appartiennent \`a 
$FC(\mathfrak{g}({\mathbb F}_{q}))$. On peut pr\'eciser le support de cette fonction. R\'ealisons le groupe $SO_{dep}(2n)$ comme le groupe sp\'ecial orthogonal d'un espace $V$ sur ${\mathbb F}_{q}$ muni d'une forme quadratique $Q$. Alors les classes de conjugaison par $SO_{dep}(2n)({\mathbb F}_{q})$ dans l'orbite g\'eom\'etrique de $N_{\triangle}$ sont naturellement param\'etr\'ees par des familles de formes quadratiques de rang $1$, qui s'identifient \`a des familles $\underline{\eta}=(\eta_{l})_{l=1,...,m}$ d'\'el\'ements de ${\mathbb F}^{\times}_{q}/{\mathbb F}_{q}^{\times,2}$. Plus pr\'ecis\'ement, par les familles $\underline{\eta}=(\eta_{l})_{l=1,...,m}$ qui v\'erifient la condition 

(1) \quad $\prod_{l=1,...,m}\eta_{l}=1$ 

\noindent dans ${\mathbb F}^{\times}_{q}/{\mathbb F}_{q}^{\times,2}$. 
Soit $N_{\underline{\eta}}$ un \'el\'ement dans la classe param\'etr\'ee par $\underline{\eta}$. Comme  en \ref{Bn}, la
 classe de conjugaison par $SO_{dep}(2n;{\mathbb F}_{q})$ de $N_{\underline{\eta}}$ se coupe en deux classes de conjugaison par $G({\mathbb F}_{q})$ si et seulement si tous les $\eta_{l}$ sont \'egaux. Compte tenu de (1), il y a deux classes de conjugaison par $SO_{dep}(2n,{\mathbb F}_{q})$, d'o\`u $4$ classes de conjugaison par $G({\mathbb F}_{q})$. 
 Supposons que $N_{\triangle}$ appartient \`a l'une d'elles et notons $N^{-}_{\triangle}$, $N'_{\triangle}$, $N''_{\triangle}$ des repr\'esentants des trois autres, avec $N^{-}_{\triangle}$ conjugu\'e \`a $N_{\triangle}$ par un \'el\'ement de $SO_{dep}(2n)({\mathbb F}_{q})$. Alors, quitte \`a permuter $N'_{\triangle}$ et $N''_{\triangle}$, les fonctions $f_{N_{\triangle},\epsilon'_{\triangle}}$ et $f_{N_{\triangle},\epsilon''_{\triangle}}$  v\'erifient les \'egalit\'es
 $$f_{N_{\triangle},\epsilon'_{\triangle}}(N_{\triangle})=f_{N_{\triangle},\epsilon''_{\triangle}}(N_{\triangle})=f_{N_{\triangle},\epsilon'_{\triangle}}(N'_{\triangle})=f_{N_{\triangle},\epsilon''_{\triangle}}(N''_{\triangle})=2^{m/2-1}$$
 $$f_{N_{\triangle},\epsilon'_{\triangle}}(N^-_{\triangle})=f_{N_{\triangle},\epsilon''_{\triangle}}(N^-_{\triangle})=f_{N_{\triangle},\epsilon'_{\triangle}}(N''_{\triangle})=f_{N_{\triangle},\epsilon''_{\triangle}}(N'_{\triangle})=-2^{m/2-1},$$
 et $f_{N_{\triangle},\epsilon'_{\triangle}}(N)=f_{N_{\triangle},\epsilon''_{\triangle}}(N)=0$ pour tout $N\in \mathfrak{g}_{nil}({\mathbb F}_{q})$ qui n'est pas conjugu\'e \`a $N_{\triangle}$, $N^-_{\triangle}$, $N'_{\triangle}$ ou $N''_{\triangle}$ par un \'el\'ement de $G({\mathbb F}_{q})$. L'automorphisme $\theta$ conserve la classe de conjugaison de $N_{\triangle}$ par $SO_{dep}(2n)({\mathbb F}_{q})$ mais peut conserver ou non sa classe de conjugaison par $G({\mathbb F}_{q})$. Le premier cas se produit si $i$ est impair ou si $j$ est pair et $sgn(-1)=1$. Dans ce cas, $\theta(f_{N_{\triangle},\epsilon'_{\triangle}})=f_{N_{\triangle},\epsilon''_{\triangle}}$. Si $j$ est pair et $sgn(-1)=-1$, $\theta(f_{N_{\triangle},\epsilon'_{\triangle}})=-f_{N_{\triangle},\epsilon''_{\triangle}}$.

 Consid\'erons maintenant le cas d'un quotient $G=G_{SC}/ {\cal Z}$ o\`u ${\cal Z}$ est un sous-groupe de $Z(G_{SC})$. Les fonctions $f_{N,\epsilon}$ ci-dessus, quand elles existent, survivent si et seulement si $\epsilon$ est trivial sur ${\cal Z}$. Remarquons que $\epsilon_{\square}(z')=\epsilon_{\square}(z'')$ vaut $(-1)^{h/2}$ et que $h/2$ est pair si et seulement si $n$ est divisible par $4$. On obtient
 
  pour ${\cal Z}=\{1,z\}$, $dim(FC(\mathfrak{g}({\mathbb F}_{q})))=\delta_{\square}(2n)$; 
  
  pour ${\cal Z}=\{1,z'\}$ ou ${\cal Z}=\{1,z''\}$, $dim(FC(\mathfrak{g}({\mathbb F}_{q})))=\delta_{\square}(2n)\delta_{4}(n)+\delta_{\triangle}(2n)$;
  
  pour ${\cal Z}=Z(G_{SC})$, $dim(FC(\mathfrak{g}({\mathbb F}_{q})))=\delta_{\square}(2n)\delta_{4}(n)$.
 
 \subsection{Type $D_{n}$ non d\'eploy\'e,  $n$ pair}\label{Dnnondeppair}
  On suppose que $G$ de type $D_{n}$, avec $n$ pair et $n\geq2$, et que $Fr$ agit sur ${\cal D}$ par $\theta$. On a $G_{SC}=Spin_{{\mathbb F}_{q^2}/{\mathbb F}_{q}}(2n)$, l'indice ${\mathbb F}_{q^2}/{\mathbb F}_{q}$ indiquant que le groupe de Galois de cette extension agit non trivialement. Les faisceaux-caract\`eres sont les m\^emes que dans le cas pr\'ec\'edent mais il ne sont plus tous conserv\'es par l'action galoisienne, celle-ci \'echangeant les caract\`eres $\mu'$ et $\mu''$. 
  
  Pour $G=G_{SC}$, on a $dim(FC(\mathfrak{g}({\mathbb F}_{q})))=\delta_{\square}(2n)$ et l'espace $FC(\mathfrak{g}({\mathbb F}_{q}))$ est engendr\'e par une fonction $f_{N_{\square},\epsilon_{\square}}$ similaire \`a celle du cas d\'eploy\'e.
 
Pour un quotient  $G=G_{SC}/ {\cal Z}$, le sous-groupe ${\cal Z}$ doit \^etre conserv\'e par $\theta$. On obtient

pour ${\cal Z}=\{1,z\}$, $dim(FC(\mathfrak{g}({\mathbb F}_{q})))=\delta_{\square}(2n)$; 

pour ${\cal Z}=Z(G_{SC})$, $dim(FC(\mathfrak{g}({\mathbb F}_{q})))=\delta_{\square}(2n)\delta_{4}(n)$.

\subsection{Type $D_{n}$ d\'eploy\'e, $n$ impair}\label{Dndepimp}
On suppose $G$ de type $D_{n}$ d\'eploy\'e, avec $n$ impair et $n\geq3$. On  conserve la plupart des notations du cas $n$ pair. 

Supposons $G=G_{SC}$. La structure de $Z(G)$ n'est plus la m\^eme. Posons 
$$z'=\check{\alpha}_{n-1}(i)\check{\alpha}_{n}(-i)\prod_{j=1,...,n-2}\check{\alpha}_{j}((-1)^j),$$
 o\`u $i$ est ici une racine d'ordre $4$ de $1$ dans $\bar{{\mathbb F}}_{q}^{\times}$. Alors $Z(G)$ est le groupe cyclique d'ordre $4$ engendr\'e par $z'$. On a ${z'}^2=z$ o\`u $z$ est comme en \ref{Dndeppair}, $\theta(z')={z'}^{-1}$ et $Fr({z'})={z'}^{sgn(-1)}$. On a l'\'egalit\'e $dim(FC(\mathfrak{g}({\mathbb F}_{q})))= 2\delta_{\triangle}(2n)\delta_{4}(q-1)$.

  Supposons $\delta_{\triangle}(2n)=1$, c'est-\`a-dire $2n=j(j+1)/2$ pour un entier $j\in {\mathbb N}$. On pose $m=[(j+1)/2]$.  On v\'erifie que $m$ est pair et que $m/2$ est impair.   Notons $\mu'$ le caract\`ere de $Z(G)$ tel que $\mu'(z')=i$ o\`u cette fois $i$ est une racine d'ordre $4$ de $1$ dans ${\mathbb C}^{\times}$.  Notons $\mu''$ son inverse. Soit $N_{\triangle}\in \mathfrak{g}({\mathbb F}_{q})$ dont la classe de conjugaison par $SO_{dep}(2n)({\mathbb F}_{q})$ est param\'etr\'ee par la m\^eme partition qu'en \ref{Dndeppair}. De nouveau, le groupe $Z_{G}(N_{\triangle})/Z_{G}(N_{\triangle})^0$ a une unique repr\'esentation irr\'eductible $\epsilon'_{\triangle}$, resp. $\epsilon''_{\triangle}$, sur laquelle $Z(G)$ agit par le caract\`ere $\mu'$, resp. $\mu''$. Les classes de ces repr\'esentations sont conserv\'ees par l'action galoisienne si et seulement si $sgn(-1)=1$, c'est-\`a-dire $4$ divise $q-1$. Supposons que $4$ divise $q-1$. On construit alors les fonctions $f_{N_{\triangle},\epsilon'_{\triangle}}$ et $f_{N_{\triangle},\epsilon''_{\triangle}}$ qui appartiennent \`a $FC(\mathfrak{g}({\mathbb F}_{q}))$. On pr\'ecise le support de ces fonctions comme en \ref{Dndeppair}. La relation (1) de ce paragraphe devient
  
  (1) $\prod_{l=1,...,m}\eta_{l}=-1$.
  
  Les fonctions sont support\'ees par les  classes de conjugaison par $SO_{dep}(2n)({\mathbb F}_{q})$ param\'etr\'ees par les familles $\underline{\eta}=(\eta_{l})_{l=1,...,m}$ pour lesquelles $\eta_{l}$ est constant. Il y en a deux d'apr\`es (1), la parit\'e de $m$ et l'hypoth\`ese $sgn(-1)=1$. On choisit $N_{\triangle}$, $N^-_{\triangle}$, $N'_{\triangle}$, $N''_{\triangle}$ comme en \ref{Dndeppair}. Quitte \`a permuter $N'_{\triangle}$ et $N''_{\triangle}$, on a les \'egalit\'es
  $$f_{N_{\triangle},\epsilon'_{\triangle}}(N_{\triangle})=f_{N_{\triangle},\epsilon''_{\triangle}}(N_{\triangle})=2^{m/2-1},$$
   $$f_{N_{\triangle},\epsilon'_{\triangle}}(N^-_{\triangle})=f_{N_{\triangle},\epsilon''_{\triangle}}(N^-_{\triangle})=-2^{m/2-1},$$
    $$f_{N_{\triangle},\epsilon'_{\triangle}}(N'_{\triangle})=f_{N_{\triangle},\epsilon''_{\triangle}}(N''_{\triangle})=i2^{m/2-1},$$
 $$f_{N_{\triangle},\epsilon'_{\triangle}}(N''_{\triangle})=f_{N_{\triangle},\epsilon''_{\triangle}}(N'_{\triangle})=-i2^{m/2-1}.$$
On a $\theta(f_{N_{\triangle},\epsilon'_{\triangle}})=f_{N_{\triangle},\epsilon''_{\triangle}}$.

Pour un quotient $G=G_{SC}/{\cal Z}$, o\`u ${\cal Z}$ est un sous-groupe non trivial de $Z(G_{SC})$, on a $FC(\mathfrak{g}({\mathbb F}_{q}))=\{0\}$. 

\subsection{Type $D_{n}$  non d\'eploy\'e, $n$ impair}\label{Dnnondepimp}
On suppose que  $G$ est de type $D_{n}$, avec $n$ impair et $n\geq3$, et que $Fr$ agit sur ${\cal D}$ par $\theta$. 

Supposons $G=G_{SC}$.  La diff\'erence avec le cas pr\'ec\'edent est que l'on a maintenant $Fr(z')={z'}^{-sgn(-1)}$.   On a l'\'egalit\'e $dim(FC(\mathfrak{g}({\mathbb F}_{q})))= 2\delta_{\triangle}(2n)\delta_{4}(q+1)$. 

Supposons $2n=j(j+1)/2$ et $4$ divise $q+1$, c'est-\`a-dire $sgn(-1)=-1$. La relation (1) de \ref{Dndepimp} est remplac\'ee par

 (1) $\prod_{l=1,...,m}\eta_{l} =-\nu$,
 
 \noindent o\`u $\nu$ est l'\'el\'ement non trivial de ${\mathbb F}_{q}^{\times}/{\mathbb F}_{q}^{\times,2}$. Les fonctions $f_{N_{\triangle},\epsilon'_{\triangle}}$ et $f_{N_{\triangle},\epsilon''_{\triangle}}$ sont comme dans le cas pr\'ec\'edent. 
On a $\theta( f_{N_{\triangle},\epsilon'_{\triangle}})=(-1)^{j+1}f_{N_{\triangle},\epsilon''_{\triangle}}$. 
 
Pour un quotient $G=G_{SC}/{\cal Z}$, o\`u ${\cal Z}$ est un sous-groupe non trivial de $Z(G_{SC})$, on a $FC(\mathfrak{g}({\mathbb F}_{q}))=\{0\}$. 

\subsection{Type $D_{4}$ trialitaire}\label{D4trialitaire}

On suppose $G$ de type $D_{4}$. Le diagramme de Dynkin a des automorphismes suppl\'ementaires. On  a d\'ej\`a introduit  l'automorphisme   $\theta$ d'ordre $2$ et on note $\theta_{3}$ l'automorphisme d'ordre $3$ qui fixe $\alpha_{2}$ et envoie $\alpha_{1}$, $\alpha_{3}$ et $\alpha_{4}$ respectivement sur $\alpha_{3},\alpha_{4},\alpha_{1}$. On note de m\^eme l'automorphisme de $G_{SC}$ qui s'en d\'eduit. Il envoie $z',z,z''$ respectivement sur $z,z'',z'$. 
Le groupe $Aut({\cal D})$ est engendr\'e par $\theta$ et $\theta_{3}$ et est isomorphe \`a $\mathfrak{S}_{3}$.  On suppose que $Fr$ agit sur ${\cal D}$ par $\theta_{3}$ (on pourrait remplacer $\theta_{3}$ par $\theta_{3}^{-1}$, cela remplacerait $G$ par un groupe isomorphe).  Cela impose que $G=G_{SC}$ ou $G=G_{AD}$. On a $FC(\mathfrak{g}({\mathbb F}_{q}))=\{0\}$. 

\subsection{Type $E_{6}$ d\'eploy\'e}\label{E6dep}
On suppose $G$ d\'eploy\'e de type $E_{6}$. 

Supposons d'abord  $G=G_{SC}$. Pour $\zeta\in \boldsymbol{\zeta}_{3}(\bar{{\mathbb F}}_{q})$, posons $\underline{\zeta}(z)=\check{\alpha}_{1}(z)\check{\alpha}_{3}(z^2)\check{\alpha}_{5}(z)\check{\alpha}_{6}(z^2)$. Alors $\underline{\zeta}$ est un isomorphisme de $\zeta_{3}(\bar{{\mathbb F}}_{q})$ sur $Z(G_{SC})$. On note $\theta$ l'unique automorphisme non trivial  de ${\cal D}$. Il fixe $\alpha_{2}$ et $\alpha_{4}$ et \'echange $\alpha_{1}$ et $\alpha_{6}$ ainsi que $\alpha_{3}$ et $\alpha_{5}$. On a $\theta(\underline{\zeta}(z))=\underline{\zeta}(z^2)$. On a
$dim(FC(\mathfrak{g}({\mathbb F}_{q})))=2\delta_{3}(q-1)$.

Introduisons l'\'el\'ement $H$ de $\mathfrak{t}({\mathbb F}_{q})$ tel que $\alpha_{i}(H)=2$ pour $i=1,4,6$ et $\alpha_{i}(H)=0$ pour $i=2,3,5$. A l'aide de cet \'el\'ement, on d\'efinit comme en \ref{groupessurFq} des sous-espaces $\mathfrak{g}_{j}$ et $\mathfrak{g}_{\geq j}$ de $\mathfrak{g}$ pour tout $j\in {\mathbb Z}$, un Levi $M$  de $G$ d'alg\`ebre de Levi $\mathfrak{g}_{0}$ et    l'orbite ouverte   $\tilde{\mathfrak{g}}_{2}$ de l'action de $M$ dans $\mathfrak{g}_{2}$. On fixe $N\in \tilde{\mathfrak{g}}_{2}({\mathbb F}_{q})$. On a $Z_{G}(N)/Z_{G}(N)^0\simeq {\mathbb Z}/6{\mathbb Z}$ et l'image naturelle de $Z(G)$ dans ce groupe est $2{\mathbb Z}/6{\mathbb Z}$. Notons  $\epsilon'$ et $\epsilon''$ les deux caract\`eres   de $Z_{G}(N)/Z_{G}(N)^0$  d'ordre $6$. Ils donnent naissance \`a des faisceaux-caract\`eres cuspidaux. Ces faisceaux sont conserv\'es par l'action galoisienne si et seulement s'il en est de m\^eme des caract\`eres $\epsilon'$ et $\epsilon''$, c'est-\`a-dire si cette action galoisienne fixe tout \'el\'ement de $Z(G)$. La condition pour qu'il en soit ainsi est que $3$ divise $q-1$. Supposons cette condition v\'erifi\'ee. Alors les faisceaux-caract\`eres donnent naissance aux fonctions   $f_{N,\epsilon'}$ et $f_{N,\epsilon''}$ qui  engendrent  $FC(\mathfrak{g}({\mathbb F}_{q}))$. 

Supposons maintenant $G=G_{AD}$. Les restrictions des caract\`eres $\epsilon'$ et $\epsilon''$ \`a $Z(G_{SC})$ \'etant non triviales, on a $FC(\mathfrak{g}({\mathbb F}_{q}))=\{0\}$.

\subsection{Type $E_{6}$ non d\'eploy\'e}\label{E6nondep}
On suppose que $G$ est de type $E_{6}$ et que $Fr$ agit sur ${\cal D}$ par l'automorphisme $\theta$. Les faisceaux-caract\`eres sont les m\^emes que dans le cas pr\'ec\'edent. La diff\'erence est que l'action galoisienne sur le centre a chang\'e. Cette action fixe tout \'el\'ement central si et seulement si $3$ divise $q+1$. D'o\`u $dim(FC(\mathfrak{g}({\mathbb F}_{q})))=2\delta_{3}(q+1)$ si $G=G_{SC}$, $FC(\mathfrak{g}({\mathbb F}_{q}))=\{0\}$ si $G=G_{AD}$.

\subsection{Type $E_{7}$}\label{E7}
On suppose que $G$ est de type $E_{7}$, donc d\'eploy\'e. 

Supposons d'abord $G=G_{SC}$. Posons $z=\check{\alpha}_{2}(-1)\check{\alpha}_{5}(-1)\check{\alpha}_{7}(-1)$. Alors $Z(G)=\{1,z\}$. On a $dim(FC(\mathfrak{g}({\mathbb F}_{q})))=1$.

Introduisons l'\'el\'ement $H$ de $\mathfrak{t}({\mathbb F}_{q})$ tel que $\alpha_{i}(H)=2$ pour $i=4,7$ et $\alpha_{i}(H)=0$ pour $i=1,2,3,5,6$. A l'aide de cet \'el\'ement, on d\'efinit comme   en \ref{groupessurFq} le Levi $M$ et l'orbite ouverte $\tilde{\mathfrak{g}}_{2}$. Fixons un \'el\'ement $N\in \tilde{\mathfrak{g}}_{2}({\mathbb F}_{q})$. Il existe une unique repr\'esentation irr\'eductible $\epsilon$ de $Z_{G}(N)/Z_{G}(N)^0$ auquel est associ\'e un faisceau-caract\`ere cuspidal ($\epsilon$ est en fait un caract\`ere).  Le groupe $Z(G)$ agit sur cette repr\'esentation par le caract\`ere non trivial de ce groupe. D'o\`u une fonction $f_{N,\epsilon}\in 
FC(\mathfrak{g}({\mathbb F}_{q}))$.

Supposons maintenant $G=G_{AD}$.  Le groupe $Z(G)$ agit sur la repr\'esentation $\epsilon$ par le caract\`ere non trivial de ce groupe. D'o\`u $FC(\mathfrak{g}({\mathbb F}_{q}))=\{0\}$.

\subsection{Type $E_{8}$}\label{E8}
On suppose que $G$ est de type $E_{8}$, donc d\'eploy\'e.  On a forc\'ement $G=G_{SC}=G_{AD}$. On a  $dim(FC(\mathfrak{g}({\mathbb F}_{q})))=1$.

Introduisons l'\'el\'ement $H$ de $\mathfrak{t}({\mathbb F}_{q})$ tel que $\alpha_{i}(H)=2$ pour $i=5$ et $\alpha_{i}(H)=0$ pour $i\not=5$. A l'aide de cet \'el\'ement, on d\'efinit comme   en \ref{groupessurFq} le Levi $M$ et l'orbite ouverte $\tilde{\mathfrak{g}}_{2}$. Fixons un \'el\'ement $N\in \tilde{\mathfrak{g}}_{2}({\mathbb F}_{q})$. On a $Z_{G}(N)/Z_{G}(N)^0\simeq \mathfrak{S}_{5}$. Soit $\epsilon$ le caract\`ere signe usuel de ce groupe. Il lui est associ\'e un faisceau-caract\`ere cuspidal dont la fonction caract\'eristique $f_{N,\epsilon}$ appartient \`a $FC(\mathfrak{g}({\mathbb F}_{q}))$.

\subsection{Type $F_{4}$}\label{F4}
On suppose que $G$ est de type $F_{4}$, donc d\'eploy\'e.  On a forc\'ement $G=G_{SC}=G_{AD}$. On a  $dim(FC(\mathfrak{g}({\mathbb F}_{q})))=1$. 

Introduisons l'\'el\'ement $H$ de $\mathfrak{t}({\mathbb F}_{q})$ tel que $\alpha_{i}(H)=2$ pour $i=2$ et $\alpha_{i}(H)=0$ pour $i=1,3,4$. A l'aide de cet \'el\'ement, on d\'efinit comme   en  \ref{groupessurFq} le Levi $M$ et l'orbite ouverte $\tilde{\mathfrak{g}}_{2}$. Fixons un \'el\'ement $N\in \tilde{\mathfrak{g}}_{2}({\mathbb F}_{q})$. On a $Z_{G}(N)/Z_{G}(N)^0\simeq \mathfrak{S}_{4}$. Soit $\epsilon$ le caract\`ere signe usuel de ce groupe. Il lui est associ\'e un faisceau-caract\`ere cuspidal dont la fonction caract\'eristique $f_{N,\epsilon}$ appartient \`a $FC(\mathfrak{g}({\mathbb F}_{q}))$.

\subsection{Type $G_{2}$}\label{G2}
On suppose que $G$ est de type $F_{4}$, donc d\'eploy\'e.  On a forc\'ement $G=G_{SC}=G_{AD}$. On a  $dim(FC(\mathfrak{g}({\mathbb F}_{q})))=1$. 

Introduisons l'\'el\'ement $H$ de $\mathfrak{t}({\mathbb F}_{q})$ tel que $\alpha_{i}(H)=2$ pour $i=2$ et $\alpha_{i}(H)=0$ pour $i=1$. A l'aide de cet \'el\'ement, on d\'efinit comme   en  \ref{groupessurFq} le Levi $M$ et l'orbite ouverte $\tilde{\mathfrak{g}}_{2}$. Fixons un \'el\'ement $N\in \tilde{\mathfrak{g}}_{2}({\mathbb F}_{q})$. On a $Z_{G}(N)/Z_{G}(N)^0\simeq \mathfrak{S}_{3}$. 
Soit $\epsilon$ le caract\`ere signe usuel de ce groupe. Il lui est associ\'e un faisceau-caract\`ere cuspidal dont la fonction caract\'eristique $f_{N,\epsilon}$ appartient \`a $FC(\mathfrak{g}({\mathbb F}_{q}))$.

\section{Pr\'esentation des r\'esultats}

\subsection{Les r\'esultats}\label{resultats}
Pour toute la suite de l'article, $G$ est un groupe r\'eductif connexe d\'efini sur $F$, simplement connexe et absolument quasi-simple. On pr\'esente $G$ \`a l'aide d'une forme quasi-d\'eploy\'ee $G^*$ et d'un cocycle $\underline{n}_{G}$ comme en \ref{orbites}.  

On va \'etudier successivement tous les groupes $G$ possibles. Dans chaque cas, nous d\'efinirons les objets suivants:

(A) un ensemble fini ${\cal X}$ de nature combinatoire; pour chaque $x\in {\cal X}$, un sous-espace $FC_{x}\subset FC(\mathfrak{g}(F))$; en fait, ce sous-espace est muni d'une base dont l'inconv\'enient est que ses \'el\'ements ne sont canoniques qu'\`a un scalaire pr\`es; pour chaque $x\in {\cal X}$, on indiquera l'entier $d_{x}=dim(FC_{x})$; 

(B) un ensemble fini ${\cal Y}$ de nature combinatoire; pour chaque $y\in {\cal Y}$, un sous-espace 
$$FC^{{\cal E}}_{y}\subset FC^{{\cal E}}(\mathfrak{g}(F))=\oplus_{{\bf G}'\in Endo_{ell}(G)}FC^{st}(\mathfrak{g}'(F))^{Out({\bf G}')};$$
de nouveau $FC^{{\cal E}}_{y}$ est muni d'une base qui a le m\^eme inconv\'enient que ci-dessus;

(C) une bijection $\varphi:{\cal X}\to {\cal Y}$. 

Dans le cas o\`u $G$ est quasi-d\'eploy\'e, nous d\'efinirons de plus

(D) un sous-ensemble ${\cal X}^{st}\subset {\cal X}$.

Les propri\'et\'es que nous d\'emontrerons sont

(1) on a l'\'egalit\'e $FC(\mathfrak{g}(F))=\oplus_{x\in {\cal X}}FC_{x}$ et cette d\'ecomposition est orthogonale pour le produit scalaire elliptique;

(2) on a l'\'egalit\'e $ FC^{{\cal E}}(\mathfrak{g}(F))=\oplus_{y\in {\cal Y}}FC^{{\cal E}}_{y}$ et cette d\'ecomposition est orthogonale pour le produit scalaire elliptique d\'efini en \ref{donneesendoscopiques} ;

(3)  l'isomorphisme de transfert $transfert:FC(\mathfrak{g}(F))\to FC^{{\cal E}}(\mathfrak{g}(F))$ envoie $FC_{x}$ sur $FC^{{\cal E}}_{\varphi(x)}$ pour tout $x\in {\cal X}$. 

Dans le cas o\`u $G$ est quasi-d\'eploy\'e, nous prouverons de plus

(4) on a l'\'egalit\'e $FC^{st}(\mathfrak{g}(F))=\oplus_{x\in {\cal X}^{st}}FC_{x}$. 

\subsection{Quelques ingr\'edients des preuves}\label{ingredients}

La d\'emonstration se fait par r\'ecurrence sur le rang de $G$.  Pr\'ecis\'ement, le groupe $G$ \'etant fix\'e, on  suppose les r\'esultats connus pour tous les groupes $G'$ (simplement connexes et absolument quasi-simples) d\'efinis sur une extension finie $F'$ de $F$ tels que le rang de $G'$ sur $\bar{F}$ soit strictement inf\'erieur \`a celui de $G$.  De plus, si  $G$ n'est pas d\'eploy\'e, on suppose  les r\'esultats  connus pour sa forme quasi-d\'eploy\'ee $G^*$.

D\'ecrivons l'un des principaux arguments que nous utiliserons. Avec les notations du paragraphe pr\'ec\'edent, on  suppose donn\'es les objets (A), (B), (C)   v\'erifiant les conditions (1) et (2). On suppose donn\'es deux sous-ensembles ${\cal X}^{\star}\subset {\cal X}$ et ${\cal Y}^{\star}\subset {\cal Y}$ tels que $\varphi({\cal X}^{\star})={\cal Y}^{\star}$. On suppose d\'emontr\'e que

(1) $ transfert(\oplus_{x\in {\cal X}^{\star}}FC_{x})=\oplus_{y\in {\cal Y}^{\star}}FC^{{\cal E}}_{y}$.

On suppose donn\'ee une relation de pr\'eordre $\leq$ sur ${\cal X}^{\star}$, d'o\`u une relation d'\'equivalence $x\equiv x'$ si et seulement si $x\leq x'$ et $x'\leq x$.On note $(x)$ la classe d'\'equivalence de $x$ et $\underline{{\cal X}}^{\star}$ l'ensemble des classes d'\'equivalence. Le pr\'eordre induit un ordre sur $\underline{{\cal X}}^{\star}$. Par la bijection $\varphi$, on transporte le pr\'eordre et la relation d'\'equivalence sur ${\cal X}^{\star}$ en des  relations sur ${\cal Y}^{\star}$ et on adopte de m\^emes notations pour ces objets.  Evidemment, $\varphi$ se quotiente en une bijection encore not\'ee $\varphi$ de $\underline{{\cal X}}^{\star}$ sur $\underline{{\cal Y}}^{\star}$. Pour $(x)\in \underline{{\cal X}}^{\star}$, on pose $FC_{(x)}=\oplus_{x'\in (x)}FC_{x'}$. De m\^eme, pour $(y)\in \underline{{\cal Y}}^{\star}$, on pose $FC_{(y)}^{{\cal E}}=\oplus_{y'\in (y)}FC^{{\cal E}}_{y'}$.  On suppose

(2) $dim(FC_{(x)})=dim(FC^{{\cal E}}_{\varphi((x))})$ pour tout $(x)\in \underline{{\cal X}}^{\star}$.

On consid\`ere un sous-ensemble $\underline{{\cal Y}}^{\sharp}\subset \underline{{\cal Y}}^{\star}$ qui est soit $\underline{{\cal Y}}^{\star}$ tout entier, soit $\underline{{\cal Y}}^{\star}-\{(y_{min})\}$ o\`u $(y_{min})$ est le plus petit \'el\'ement de $\underline{{\cal Y}}^{\star}$ (l'existence de ce plus petit \'el\'ement est suppos\'ee dans ce cas).

Pour ${\bf G}'\in Endo_{ell}(G)$ et $Y\in \mathfrak{g}'_{reg}(F)$, la forme lin\'eaire $f'\mapsto S^{G'}(Y,f')$ est bien d\'efinie sur $SI(\mathfrak{g}'(F))$, a fortiori sur $FC^{st}(\mathfrak{g}'(F))^{Out({\bf G}')}$. On l'\'etend en une forme lin\'eaire sur $FC^{{\cal E}}(\mathfrak{g}(F))$ par $0$ sur  $FC^{st}(\mathfrak{g}''(F))^{Out({\bf G}'')}$ pour tout ${\bf G}''\in Endo_{ell}(G)$, ${\bf G}''\not={\bf G}'$.

Pour tout $(y)\in \underline{{\cal Y}}^{\sharp}$, posons $d_{(y)}=dim(FC^{{\cal E}}_{(y)})$ et supposons donn\'ee des familles $({\bf G}'_{i})_{i=1,...,d_{(y)}}$ et $(Y_{i})_{i=1,...,d_{(y)}}$ telles que: pour tout $i$, ${\bf G}'_{i}\in  Endo_{ell}(G)$, $Y_{i}\in \mathfrak{g}'_{i,ell}(F)$ et $Y_{i}$ est semi-simple et $G$-r\'egulier.
Remarquons que l'on ne suppose pas que les donn\'ees  ${\bf G}'_{i}$ sont distinctes. 
On suppose que, pour tout $(y)\in \underline{{\cal Y}}^{\sharp}$, les conditions suivantes sont satisfaites.

(3) Pour tout $i=1,...,d_{(y)}$ et tout $(y')\in \underline{{\cal Y}}^{\star}$ avec $(y')\not=(y)$, la forme lin\'eaire $f'\mapsto S^{G'_{i}}(Y_{i},f')$ est nulle sur $FC^{{\cal E}}_{(y')}$.

(4) Les restrictions \`a $FC^{{\cal E}}_{(y)}$ des   formes lin\'eaires $f'\mapsto S^{G'_{i}}(Y_{i},f')$, pour $i=1,...,d_{(y)}$, forment une base de l'espace dual $(FC^{{\cal E}}_{(y)})^*$.

(5) Soient $(x)\in \underline{{\cal X}}$, $f\in FC_{x}$, $i\in \{1,...,d_{(y)}\}$ et $X$ un \'el\'ement de $\mathfrak{g}_{reg}(F)$ dont la classe de conjugaison stable correspond \`a celle de $Y_{i}$. Supposons $I^G(X,f)\not=0$. Alors $(x)\geq \varphi^{-1}((y))$. 

\begin{lem}{Sous ces hypoth\`eses, on a l'\'egalit\'e $transfert(FC_{(x)})=FC^{{\cal E}}_{\varphi((x))}$ pour tout $(x)\in \underline{{\cal X}}^{\star}$.}\end{lem}

Preuve. On  suppose que $\underline{{\cal Y}}^{\sharp}=\underline{{\cal Y}}^{\star}-\{(y_{min})\}$ (quand $\underline{{\cal Y}}^{\sharp}=\underline{{\cal Y}}^{\star}$, il suffit de supprimer dans ce qui suit toute mention de l'\'el\'ement $(y_{min})$). On raisonne par r\'ecurrence en supposant l'\'egalit\'e  de l'\'enonc\'e d\'emontr\'ee pour tout $(x')<(x)$. Soit $f\in FC_{(x)}$. Gr\^ace \`a (1),  \'ecrivons
$$transfert(f)=\sum_{(y)\in \underline{{\cal Y}}^{\star}}f'_{(y)},$$
avec $f'_{y}\in FC^{{\cal E}}_{(y)}$ pour tout $(y)$. On va d'abord prouver

(6) soit $(y)\in \underline{{\cal Y}}^{\star}$, supposons $f'_{(y)}\not=0$, alors   $(x)\geq \varphi^{-1}((y))$. 

C'est tautologique si $(y)=(y_{min})$. Supposons $(y)\not=(y_{min})$. 
D'apr\`es (3) et (4), on peut fixer $i\in \{1,...,d_{(y)}\}$ tel que $S^{G'_{i}}(Y_{i},transfert(f))\not=0$, c'est-\`a-dire, d'apr\`es les d\'efinitions, $S^{G'_{i}}(Y_{i},transfert^{{\bf G}'_{i}}(f))\not=0$. Cette int\'egrale est par d\'efinition combinaison lin\'eaire d'int\'egrales $I^G(X,f)$ pour des \'el\'ements $X\in \mathfrak{g}_{reg}(F)$ dont la classe de conjugaison stable correspond \`a celle de $Y_{i}$. L'une de ces int\'egrales doit \^etre non nulle. Alors (5) entra\^{\i}ne que $(x)\geq \varphi^{-1}((y))$.  Cela prouve (6).

 Soit $(y)\in \underline{{\cal Y}}^{\star}$, supposons $f'_{(y)}\not=0$. 
Posons $(z)=\varphi^{-1}((y))$. On a $(x)\geq (z)$ d'apr\`es (6).   Supposons $(x)> (z)$. Par l'hypoth\`ese de r\'ecurrence, il existe $f_{(z)}\in FC_{(z)}$ de sorte que $transfert(f_{(z)})=f'_{(y)}$. Puisque l'application $transfert $ est une isom\'etrie, on a 
$$(f,f_{(z)})_{ell}=(transfert(f),transfert(f_{(z)}))_{ell}=(f'_{(y)},f'_{(y)})_{ell}.$$
Mais $(f,f_{(z)})_{ell}=0$ puisque $(x)\not=(z)$. Alors $f'_{(y)}=0$ contrairement \`a l'hypoth\`ese. Cela d\'emontre que $(y)=\varphi((x))$, donc $transfert(FC_{(x)})\subset FC^{{\cal E}}_{\varphi((x))}$.   Puisque $transfert$ est bijective, (2) entra\^{\i}ne que   cette inclusion est une \'egalit\'e. $\square$

A chaque fois que nous utiliserons ce lemme, il faudra construire les \'el\'ements $Y_{i}$ et d\'emontrer leurs propri\'et\'es. Nous le ferons compl\`etement dans le cas des groupes exceptionnels. Dans le cas des groupes (plus ou moins) classiques, notre construction utilise les descriptions des immeubles par l'alg\`ebre lin\'eaire. On a d\'ej\`a donn\'e en \cite{W2} de telles descriptions dans les cas sp\'eciaux orthogonaux ou symplectiques. Pour \'eviter d'allonger d\'emesur\'ement l'article, on ne donnera ici de description compl\`ete que dans le cas le plus nouveau qui est celui des groupes sp\'eciaux unitaires relatifs \`a une extension quadratique ramifi\'ee. Dans les autres cas, on se contentera de donner les d\'efinitions en laissant les d\'emonstrations au lecteur.  

En fait, on construira souvent des \'el\'ements $Y_{i}\in \mathfrak{g}'_{i,ell}(F)$ v\'erifiant (3) et (4) et pour lesquels on saurait d\'emontrer (5) si ces \'el\'ements \'etaient $G$-r\'eguliers. Mais ils ne seront pas toujours $G$-r\'eguliers. On utilisera alors l'argument suivant. On sait que les \'el\'ements $G$-r\'eguliers sont denses dans $\mathfrak{g}'_{i,ell}(F)$. On sait aussi que les int\'egrales orbitales sont localement constantes sur cet ensemble. On peut alors remplacer les $Y_{i}$ d'origine par des \'el\'ements tr\`es voisins de sorte qu'ils deviennent $G$-r\'eguliers et conservent les propri\'et\'es (3) et (4). On saura alors d\'emontrer (5) pour ces \'el\'ements.

\section{Type $A_{n-1}$}

\subsection{Type $A_{n-1}$ d\'eploy\'e}\label{An-1deppadique}
On suppose que $G$ est d\'eploy\'e de type $A_{n-1}$ avec $n\geq2$. On a $G=SL(n)$, $Z(G)\simeq \boldsymbol{\zeta}_{n}(F)$, $G_{AD}(F)/\pi(G(F))\simeq F^{\times}/F^{\times,n}$, $G_{AD}(F)_{0}/\pi(G(F))\simeq \mathfrak{o}_{F}^{\times}/\mathfrak{o}_{F}^{\times,2}$.   

Si $n$ ne divise pas $q-1$, on pose ${\cal X}=\emptyset$. Si $n$ divise $q-1$, notons $\Xi_{n}^{ram}$ l'ensemble des \'el\'ements de  $\Xi$, c'est-\`a-dire des caract\`eres de $F^{\times}/F^{\times,n}$, dont la restriction \`a $\mathfrak{o}_{F}^{\times}/\mathfrak{o}_{F}^{\times,n}$ est d'ordre $n$. Le nombre d'\'el\'ements de $\Xi_{n}^{ram}$ est $n\phi(n)$. On pose ${\cal X}=\Xi_{n}^{ram}$ et $d_{x}=1$ pour tout $x\in {\cal X}$.

L'ensemble $\underline{S}(G)$ compte $n$ sommets et forme une seule orbite pour l'action du groupe $G_{AD}(F)$. Fixons un sommet $s$.  Le groupe $G_{s}$ est le groupe $SL(n)$ d\'efini sur ${\mathbb F}_{q}$. Si $n$ ne divise pas $q-1$, $FC(\mathfrak{g}_{s}({\mathbb F}_{q}))=\{0\}$ d'apr\`es \ref{An-1dep}, donc $FC(\mathfrak{g}(F))=\{0\}$. Alors \ref{resultats} (1) est v\'erifi\'e avec ${\cal X}=\emptyset$.  

   Supposons que $n$ divise $q-1$.  Alors $dim(FC(\mathfrak{g}_{s}({\mathbb F}_{q})))=\phi(n)$ d'apr\`es \ref{An-1dep}. 
L'espace $FC(\mathfrak{g}_{s}({\mathbb F}_{q}))$ a pour base les fonctions $f_{N,\epsilon}$, o\`u $N$ est un nilpotent r\'egulier et $\epsilon$ parcourt les caract\`eres d'ordre $n$ de $Z(G_{s})$. On a $Z(G_{s})\simeq \boldsymbol{\zeta}_{n}({\mathbb F}_{q})\simeq \boldsymbol{\zeta}_{n}(F)$. On a $K_{AD,s}^{\dag}=K_{AD,s}^0$ et 
$K_{AD,s}^0/\pi(K_{s}^0)\simeq \mathfrak{o}_{F}^{\times}/\mathfrak{o}_{F}^{\times,n}$. Un caract\`ere $\epsilon$ de $Z(G)$ s'identifie \`a un caract\`ere de ce groupe, pr\'ecis\'ement au caract\`ere $x\mapsto \epsilon(Fr(y)y^{-1})$ pour $x\in \mathfrak{o}_{F}^{\times}$, o\`u $y$ est une racine $n$-i\`eme de $x$ dans $F^{nr,\times}$. D'apr\`es \ref{actionsurFC}, on a $\xi_{f_{N,\epsilon}}=\epsilon$ et  de $f_{N,\epsilon}$ se d\'eduit  une fonction $f_{\xi}\in FC(\mathfrak{g}(F))$ pour tout $\xi\in \Xi$ qui co\"{\i}ncide avec $\epsilon$ sur $\mathfrak{o}_{F}^{\times}/\mathfrak{o}_{F}^{\times,n}$. 
Un \'el\'ement $\xi\in \Xi$ co\"{\i}ncide avec un tel caract\`ere $\epsilon$ de $\mathfrak{o}_{F}^{\times}/\mathfrak{o}_{F}^{\times,n}$ si et seulement si $\xi\in \Xi_{n}^{ram}$.  Ainsi, pour tout $\xi\in {\cal X}=\Xi_{n}^{ram}$, on a d\'efini une fonction $f_{\xi}\in FC(\mathfrak{g}(F))$ et on note $FC_{\xi}$ la droite ${\mathbb C}f_{\xi}$. D'apr\`es \ref{actionsurFC}, la propri\'et\'e (1) de \ref{resultats} est v\'erifi\'ee. 

On pose ${\cal X}^{st}=\emptyset$. La propri\'et\'e \ref{resultats} (4) est imm\'ediate si $n$ ne divise pas $q-1$. Si $n$ divise $q-1$, une fonction $f_{\xi}$ appartient \`a $I_{cusp,\xi}(\mathfrak{g}(F))$ par d\'efinition. Puisque $\Xi_{n}^{ram}$ ne contient pas le caract\`ere trivial ${\bf 1}\in \Xi$, l'espace $FC(\mathfrak{g}(F))\cap I_{cusp,{\bf 1}}(\mathfrak{g}(F))$ est nul, a fortiori  l'espace $FC^{st}(\mathfrak{g}(F))$ est nul et \ref{resultats} (4) est v\'erifi\'ee. Explicitons-la:

(1) $FC^{st}{\mathfrak{g}}(F))=\{0\}$. 

 Si $n$ ne divise pas $q-1$, on pose ${\cal Y}=\emptyset$. Puisqu'on a vu que $FC(\mathfrak{g}(F))$ est nul, l'espace $FC^{{\cal E}}(\mathfrak{g}(F))$ l'est aussi et \ref{resultats} (2) et (3) sont v\'erifi\'es.
 
 Supposons que $n$ divise $q-1$. 
Soit $\xi\in \Xi_{n}^{ram}$. Via l'isomorphisme du corps de classes, $\xi$ s'identifie \`a un caract\`ere de $\Gamma_{F}$. Notons $K_{\xi}$ l'extension galoisienne ab\'elienne de degr\'e $n$ de $F$ telle que $\Gamma_{K_{\xi}}$ soit le noyau de $\xi$. C'est une extension totalement ramifi\'ee. La classification des donn\'ees endoscopiques de $G=SL(n)$ est bien connue. On peut la retrouver gr\^ace \`a \ref{description} et le calcul du caract\`ere associ\'e \`a une donn\'ee r\'esulte de \ref{donneesendoscopiques}. On voit qu'il y a une unique donn\'ee ${\bf G}'_{\xi}\in Endo_{ell}(G)$ telle que $\xi_{{\bf G}'}=\xi$.  
 Consid\'erons le groupe $Res_{K_{\xi}/F}GL(1)$.  Il y a un homomorphisme norme naturel  de ce groupe vers $GL(1)$ et $G'_{\xi}$ est son noyau.   
  En particulier, on a l'isomorphisme $G'_{\xi}(F)=\{x\in K_{\xi}^{\times}; norme_{K_{\xi}/F}(x)=1\}$. 
  Le groupe $G'_{\xi}$ est
un tore, l'immeuble de son groupe adjoint est r\'eduit \`a un point, notons-le $e'$. Le groupe $G'_{\xi,e'}$ est un tore et on a $X_{*}(G'_{\xi,e'})\simeq X_{*}(G'_{\xi})^{I_{F}}$. Puisque $K_{\xi}/F$ est totalement ramifi\'e, cette description entra\^{\i}ne que $G'_{\xi,e'}=\{1\}$. Alors $FC^{st}(\mathfrak{g}'_{\xi}(F))$ est la droite port\'ee par la fonction caract\'eristique de 
 $\mathfrak{k}'_{e'}=\{x\in \mathfrak{o}_{K_{\xi}}; trace_{K_{\xi}/F}(x)=0\}$. On note cette fonction $f'_{\xi}$. Le groupe d'automorphismes ext\'erieurs  de  ${\bf G}'_{\xi}$ est isomorphe \`a $\Gamma_{K_{\xi}/F}$. Il agit sur $G'_{\xi}(F)$ par l'action naturelle de ce groupe. Cette action fixe $f'_{\xi}$ donc  $FC^{st}(\mathfrak{g}'_{\xi}(F))^{Out({\bf G}'_{\xi})}$ est la droite port\'ee par $f'_{\xi}$. On note $FC_{\xi}^{{\cal E}} = FC^{st}(\mathfrak{g}'_{\xi}(F))^{Out({\bf G}'_{\xi})}$. 
 Cette construction fournit  un sous-espace 
 
  $\oplus_{\xi\in \Xi_{n}^{ram}}FC^{{\cal E}}_{\xi}\subset FC^{{\cal E}}(\mathfrak{g}(F))$. 
 
 \noindent Mais on sait que ce dernier espace a m\^eme dimension que $FC(\mathfrak{g}(F))$ et on a vu que cette dimension \'etait $\vert \Xi_{n}^{ram}\vert $. L'inclusion ci-dessus est donc une \'egalit\'e, ce qui d\'emontre \ref{resultats} (2). 
 
 On note $\varphi:{\cal X}=\Xi_{n}^{ram}\to {\cal Y}=\Xi_{n}^{ram}$ l'identit\'e. En utilisant \ref{donneesendoscopiques} (3), on voit que l'on a $transfert(FC_{\xi})\subset FC^{{\cal E}}_{\xi}$ pour tout $\xi\in \Xi_{n}^{ram}$ et cette inclusion est forc\'ement une \'egalit\'e puisque le transfert est un isomorphisme. Cela d\'emontre \ref{resultats} (3).

 \subsection{Forme int\'erieure du type $A_{n-1}$ d\'eploy\'e}
On suppose que $G$ n'est pas d\'eploy\'e mais que $G^*$ est d\'eploy\'e de type $A_{n-1}$.

Il y a une d\'ecomposition $n=md$ avec $m,d\in {\mathbb N}_{>0}$ et $d\geq2$ de sorte que $G$ soit de type $^dA_{n-1}$, cf. \cite{T} p.62. Alors, pour tout sommet $s\in \underline{S}(G)$,  $G_{s}$ est le groupe tel que $G_{s}({\mathbb F}_{q})=\{x\in GL(m,{\mathbb F}_{q^d}); norme_{{\mathbb F}_{q}^d/{\mathbb F}_{q}}(x)=1\}$. Parce que $d\geq2$, $Z(G_{s})^0\not=\{1\}$, donc $FC(\mathfrak{g}_{s}({\mathbb F}_{q}))=\{0\}$. Cela entra\^{\i}ne $FC(\mathfrak{g}(F))=\{0\}$. Alors les assertions de \ref{resultats} sont v\'erifi\'ees en posant ${\cal X}={\cal Y}=\emptyset$.

\subsection{Type $A_{n-1}$ quasi-d\'eploy\'e, $E/F$ non ramifi\'ee}\label{An-1quasidepnonram}
Soit $E$ une extension quadratique  de $F$. Notons $\tau$ l'\'el\'ement non trivial de $\Gamma_{E/F}$. On suppose que $G$ est quasi-d\'eploy\'e de type $A_{n-1}$, avec $n\geq3$, que $\Gamma_{E}$ agit trivialement sur ${\cal D}$ tandis que   $\tau$ y  agit par l'automorphisme non trivial $\theta$ du diagramme. Avec une notation \'evidente, $G$ est la forme quasi-d\'eploy\'ee de $SU_{E/F}(n)$.   Introduisons le groupe $U_{E/F}(1)$: on a $U_{E/F}(1,\bar{F})= \bar{F}^{\times}$ et il est muni de l'action galoisienne   $\sigma\mapsto \sigma_{E/F}$  telle que, pour $x\in \bar{F}^{\times}$ et $\sigma\in \Gamma_{F}$, on a $\sigma_{E/F}(x)=\sigma(x)$ si $\sigma\in \Gamma_{E}$ et $\sigma_{E/F}(x)=\sigma(x)^{-1}$ si $\sigma\not\in \Gamma_{E}$. On a la suite exacte
$$1\to Z(G)\to U_{E/F}(1)\stackrel{x\mapsto x^n}{\to}U_{E/F}(1)\to 1$$
 Il s'en d\'eduit une suite exacte de cohomologie. On calcule $H^1(\Gamma_{F},U_{E/F}(1))\simeq {\mathbb Z}/2{\mathbb Z}$. Compte tenu de l'isomorphisme $G_{AD}(F)/\pi(G(F)) \simeq H^1(\Gamma_{F},Z(G))$, on en d\'eduit une suite exacte
 $$1\to E^1/(E^1)^n\to G_{AD}(F)/\pi(G(F))\to Ker({\mathbb Z}/2{\mathbb Z}\stackrel{x\mapsto nx}{\to}{\mathbb Z}/2{\mathbb Z})\to 0$$
 o\`u $E^1$ est le noyau de l'homomorphisme $norme_{E/F}$ et $(E^1)^n$ le sous-groupe des puissances $n$-i\`emes. Le groupe $E^1/(E^1)^n$ est l'image naturelle dans $G_{AD}(F)/\pi(G(F))$ du groupe unitaire $U_{E/F}(n,F)$. Donc ce groupe s'envoie surjectivement sur $G_{AD}(F)$ si $n$ est impair et son image est d'indice $2$ si $n$ est pair.

  Supposons dans ce paragraphe que $E/F$ est non ramifi\'ee.  Si $n$ ne divise pas $q+1$, posons ${\cal X}=\emptyset$. Supposons que $n$ divise $q+1$. Notons $\Xi_{n}^{ram}$
 le sous-ensemble des \'el\'ements de $\Xi$ dont la restriction \`a 
$E^1/(E^1)^n$ est d'ordre $n$. Le nombre d'\'el\'ements de $\Xi^{ram}_{n}$ est $(1+\delta_{2}(n))\phi(n)$. On pose ${\cal X}=\Xi_{n}^{ram}$ et $d_{x}=1$ pour tout $x\in {\cal X}$. 
 
 En se r\'ef\'erant aux tables de Tits (le groupe est de type $^2A'_{n-1}$) ou en d\'ecrivant l'immeuble par l'alg\`ebre lin\'eaire, on voit que $\underline{S}(G)$ est en bijection avec l'ensemble des  couples $(a,b)\in {\mathbb N}\times {\mathbb N}$ tels que $b$ est pair et $a+b=n$.  Pour un sommet $s$ param\'etr\'e par $(a,b) $, le groupe $G_{s}$ est celui pour lequel $G_{s}({\mathbb F}_{q})=\{(x,y)\in U_{{\mathbb F}_{q^2}/{\mathbb F}_{q}}(a,{\mathbb F}_{q})\times U_{{\mathbb F}_{q^2}/{\mathbb F}_{q}}(b,{\mathbb F}_{q}); det(x)det(y)=1\}$, avec des notations \'evidentes. 
Si $n$ est impair, l'action de $G_{AD}(F)$ d\'efinie en \ref{actionsurFC} conserve chaque sommet. Si $n$ est pair, cette action \'echange les sommets param\'etr\'es par $(n,0)$ et $(0,n)$ et conserve les autres. Puisque $a$ est forc\'ement non nul quand $n$ est impair, on voit  l'ensemble des orbites  $\underline{S}(G_{AD})$ est en tout cas en bijection avec le sous-ensemble des couples  $(a,b)$ tels que $a\not=0$. 

Consid\'erons un sommet $s$ param\'etr\'e par $(a,b)$, avec $a\not=0$. Si $b\not=0$, on voit que $Z(G_{s})^0\not=\{1\}$, donc  $FC(\mathfrak{g}_{s}({\mathbb F}_{q}))=\{0\}$. Supposons $b=0$. Si $n$ ne divise pas $q+1$, on a $FC(\mathfrak{g}_{s}({\mathbb F}_{q}))=\{0\}$  d'apr\`es \ref{An-1nondep}. Cela entra\^{\i}ne $FC(\mathfrak{g}(F))=\{0\}$ et l'assertion (1)  de \ref{resultats} est v\'erifi\'ee avec  notre d\'efinition ${\cal X}=\emptyset$. Supposons  que $n$ divise $q+1$. Alors  une base de $FC(\mathfrak{g}_{s}({\mathbb F}_{q}))$ est form\'ee des fonctions $f_{N,\epsilon}$ o\`u $N$ est un nilpotent r\'egulier et $\epsilon$ est un caract\`ere d'ordre $n$ de $Z(G_{s})\simeq Z(G)$. Un tel caract\`ere s'identifie \`a un caract\`ere du groupe $E^1/(E^1)^n$ et on a $\xi_{f_{N,\epsilon}}=\epsilon$. D'apr\`es \ref{orbites}, pour tout $\xi\in \Xi$ de restriction $\epsilon$ \`a $E^1/(E^1)^n$, la fonction $f_{N,\epsilon}$ donne naissance \`a un \'el\'ement de  $FC(\mathfrak{g}(F))$. On note $FC_{\xi}$ la droite  engendr\'ee par cet \'el\'ement. Les \'el\'ements de $\Xi$ dont la restriction \`a $E^1/(E^1)^n$ est d'ordre $n$ sont pr\'ecis\'ement les \'el\'ements de $\Xi^{ram}_{n}$. On a ainsi associ\'e \`a tout \'el\'ement  $\xi\in \bar{\cal X}=\Xi_{n}^{ram}$ une droite $FC_{\xi}\subset FC(\mathfrak{g}(F))$. D'apr\`es \ref{orbites}, l'assertion  \ref{resultats}(1) est v\'erifi\'ee.  

On pose ${\cal X}^{st}=\emptyset$. L'assertion \ref{resultats} (4) est triviale si $n$ ne divise pas $q+1$. Dans le cas o\`u $n$ divise $q+1$, elle r\'esulte comme en \ref{An-1deppadique} du fait que $\Xi_{n}^{ram}$ ne contient pas l'\'el\'ement neutre de $\Xi$. Explicitons-la:

(1) $FC^{st}{\mathfrak{g}}(F))=\{0\}$.

Si $n$ ne divise pas $q+1$, on pose ${\cal Y}=\emptyset$. Puisqu'on a d\'ej\`a vu que $FC(\mathfrak{g}(F))=\{0\}$, les assertions \ref{resultats} (2) et (3) sont triviales. 

Supposons que $n$ divise $q+1$. On pose ${\cal Y}=\Xi_{n}^{ram}$. 
On identifie $\hat{\Delta}_{a}$ \`a ${\mathbb Z}/n{\mathbb Z}$, la racine $\hat{\alpha}_{0}$ s'identifiant \`a $0$ et l'automorphisme non trivial  $\theta$ du diagramme $\hat{{\cal D}}$ \`a $j\mapsto -j$. Le groupe $\hat{\Omega}$ est celui des translations par ${\mathbb Z}/n{\mathbb Z}$. L'action galoisienne $\sigma\mapsto \sigma_{G}$ attach\'ee \`a $G$ est $\sigma_{G}=1$ si $\sigma\in \Gamma_{E}$, $\sigma_{G}=\theta$ si $\sigma\in \Gamma_{F}-\Gamma_{E}$. Introduisons l'extension $K/E$ de degr\'e $n$ de l'assertion \ref{extensionsdiedrales}(3). 
  Fixons un g\'en\'erateur $\rho$ de $\Gamma_{K/E}$ et un \'el\'ement $\delta\in \Gamma_{K/F}-\Gamma_{K/E}$. On a $\delta^2=1$ et $\delta\rho=\rho^{-1}\delta$. Fixons un entier $u\in \{1,...,n-1\}$ premier \`a $n$. Si $n$ est impair, posons $v=0$. Si $n$ est pair, soit $v\in \{0,1\}$. D\'efinissons une action $\sigma\mapsto \sigma_{G'}$ de $\Gamma_{F}$ sur ${\mathbb Z}/n{\mathbb Z}$, triviale sur $\Gamma_{K}$, par $\rho_{G'}(j) =j+u$, $\delta_{G'}(j)=-j+u(1-v)$. Posons ${\cal O}={\mathbb Z}/n{\mathbb Z}$. Le couple $(\sigma\mapsto \sigma_{G'},{\cal O})$ appartient \`a ${\cal E}_{ell}(G)$ et il lui est associ\'ee une donn\'ee endoscopique ${\bf G}'\in Endo_{ell}(G)$. Remarquons que $(\rho^v\delta)_{G'}(0)=u$ et $(\rho^v\delta)_{G'}(u)=0$. Notons $L$ le sous-corps des points fixes de $\rho^v\delta$ dans $K$. L'extension $K/L$ est quadratique et $L$ ne contient pas $E$ puisque $\rho^v\delta$ n'appartient pas \`a $\Gamma_{K/E}$. Donc $K$ est la compos\'ee des deux extensions $L$ et $E$ de $F$. On v\'erifie que le groupe $G'$ est un tore et que $G'(F)=\{x\in U_{K/L}(1,L); norme_{K/E}(x)=1\}$.   L'immeuble $Imm(G'_{AD})$ est r\'eduit \`a un unique sommet $e'$. Le sous-espace des  invariants par l'action de $\rho_{G'}$ dans $X^*(G')$ est nul. Puisque 
 $K/E$ est totalement ramifi\'ee,  on a donc $X^*(G)^{I_{F}}=\{0\}$, do\`u  $G'_{e'}=\{1\}$. L'espace $FC(\mathfrak{g}'(F))$ est une droite port\'ee par la fonction caract\'eristique de $\mathfrak{k}'_{e'}$. Elle est \'evidemment stable et invariante par automorphisme. Donc $FC^{st}(\mathfrak{g}'(F))^{Out({\bf G}')}$ est la droite port\'ee par la fonction pr\'ec\'edente. On calcule le caract\`ere $\xi_{{\bf G}'}$ en explicitant la construction de \ref{donneesendoscopiques}. On voit que $\rho_{G'}(s_{sc})s_{sc}^{-1}$ est \'egal \`a $\zeta_{\star}^{u}$, o\`u $\zeta_{\star}$ est un g\'en\'erateur de $Z(\hat{G}_{SC})$ ind\'ependant de $u$ et $v$. Il en r\'esulte  
  qu'il existe un caract\`ere $\xi_{\star}$ d'ordre $n$ de $E^1/(E^1)^n$, ind\'ependant de $u$ et $v$, de sorte que la restriction de $\xi_{{\bf G}'}$ \`a ce groupe soit $\xi_{\star}^{u}$. Il en r\'esulte que $\xi_{{\bf G}'}$ appartient \`a $\Xi_{n}^{ram}$. Faisons maintenant varier $u$ et $v$: consid\'erons des donn\'ees ${\bf G}'_{0}$ et ${\bf G}'_{1}$ associ\'ees \`a des couples distincts $(u_{0},v_{0})$ et $(u_{1},v_{1})$. Si $u_{0}\not=u_{1}$, on a $\xi_{\star}^{u_{0}}\not=\xi_{\star}^{u_{1}}$, a fortiori $\xi_{{\bf G}'_{0}}\not=\xi_{{\bf G}'_{1}}$. Si $u_{0}=u_{1}$, on a $v_{0}\not=v_{1}$ donc $n$ pair et on peut supposer que $v_{0}=0$ et $v_{1}=1$. On a alors $\delta_{G'_{0}}=\rho_{G'_{1}}\delta_{G'_{1}}$ d'o\`u $\delta_{G'_{0}}(s_{sc})s_{sc}^{-1}=\rho_{G'_{1}}(\delta_{G'_{1}}(s_{sc})s_{sc}^{-1})\rho_{G'_{1}}(s_{sc})s_{sc}^{-1}=\delta_{G'_{1}}(s_{sc})s_{sc}^{-1}\zeta_{\star}^{u_{1}}$.  Puisque $n$ est pair et que $u_{1}$ est premier \`a $n$, $\zeta_{\star}^{u_{1}}$ n'est pas un carr\'e dans $Z(\hat{G}_{SC})$. Puisque $\delta$ agit sur $Z(\hat{G}_{SC})$ par inversion, les deux cocycles associ\'es \`a ${\bf G}'_{0}$ et ${\bf G}'_{1}$ ne peuvent \^etre cohomologues que si $\delta_{G'_{0}}(s_{sc})s_{sc}^{-1}$ et $\delta_{G'_{1}}(s_{sc})s_{sc}^{-1}$ diff\`erent par un carr\'e de $Z(\hat{G}_{SC})$. Ce n'est pas le cas donc les cocycles ne sont pas cohomologues et $\xi_{{\bf G}'_{0}}\not=\xi_{{\bf G}'_{1}}$. En r\'esum\'e, l'application ${\bf G}'\mapsto \xi_{{\bf G}'}$ est une bijection de l'ensemble des  donn\'ees ${\bf G}'$ associ\'ees \`a nos couples $(u,v)$ sur l'ensemble $\Xi_{n}^{ram}$. Pour $\xi\in {\cal Y}=\Xi_{n}^{ram}$, on note ${\bf G}'_{\xi}$ celle de nos donn\'ees pour laquelle $\xi_{{\bf G}'}=\xi$ et on pose $FC^{{\cal E}}_{\xi}=FC^{st}(\mathfrak{g}'_{\xi}(F))^{Out({\bf G}')}$.  
  Les m\^emes arguments qu'en \ref{An-1deppadique} montrent que l'assertion (2)  de \ref{resultats} est v\'erifi\'ee.
  
  On note $\varphi:{\cal X}=\Xi_{n}^{ram}\to {\cal Y}=\Xi_{n}^{ram}$ l'identit\'e. De nouveau, les m\^emes arguments qu'en \ref{An-1deppadique} montrent que l'assertion (3)  de \ref{resultats} est v\'erifi\'ee.

\subsection{Forme int\'erieure du type $A_{n-1}$ quasi-d\'eploy\'e, $E/F$ non ramifi\'ee}\label{formeinterieurenonramifiee}
On suppose que $G^*$ est du type pr\'ec\'edent mais que $G$ n'est pas quasi-d\'eploy\'e. On a $H^1(\Gamma_{F},G^*_{AD})\simeq Z(\hat{G}_{SC})^{\Gamma_{F}}$ et ce groupe est trivial si $n$ est impair, isomorphe \`a ${\mathbb Z}/2{\mathbb Z}$ si $n$ est pair. Notre hypoth\`ese sur $G$ implique donc que $n$ est pair et $G$ est l'unique  forme non quasi-d\'eploy\'ee de $SU_{E/F}(n)$. 

On pose ${\cal X}=\emptyset$.  
 La situation  est tr\`es voisine du pr\'ec\'edent mais, cette fois, les \'el\'ements de $\underline{S}(G)$ sont en bijection avec les couples $(a,b)\in {\mathbb N}^2$ tels que $a+b=n$ et $a,b$ impairs (cela r\'esulte de la description de l'immeuble par l'alg\`ebre lin\'eaire ou des tables de Tits, le groupe \'etant du type $^2A''_{n-1}$). Le m\^eme calcul que dans le paragraphe pr\'ec\'edent montre que $FC(\mathfrak{g}_{s}({\mathbb F}_{q}))=\{0\}$ pour tout sommet $s$. D'o\`u $FC(\mathfrak{g}(F))=\{0\}$ et l'assertion \ref{resultats} (1). 

On pose ${\cal Y}=\emptyset$. La nullit\'e de $FC(\mathfrak{g}(F))$ entra\^{\i}ne trivialement les assertions (2) et (3) de \ref{resultats}.

\subsection{Type $A_{n-1}$ quasi-d\'eploy\'e, $E/F$  ramifi\'ee}\label{An-1quasidepram}
Soit $E/F$ une extension quadratique ramifi\'ee. On consid\`ere un groupe $G$ quasi-d\'eploy\'e de type $A_{n-1}$, avec $n\geq3$,  tel que l'action galoisienne sur le diagramme de Dynkin ${\cal D}$ soit triviale sur $\Gamma_{E}$ mais que tout \'el\'ement  de $\Gamma_{F}-\Gamma_{E}$ agisse par l'automorphisme $\theta$ de ${\cal D}$. Le groupe $E^1$ est une extension d'un pro-$p$-groupe par ${\mathbb Z}/2{\mathbb Z}$, donc $E^1/(E^1)^n=\{1\}$ si $n$ est impair et $E^1/(E^1)^{n}\simeq {\mathbb Z}/2{\mathbb Z}$ si $n$ est pair. D'apr\`es le calcul fait en \ref{An-1quasidepnonram}, $G_{AD}(F)/\pi(G(F))$ est trivial si $n$ est impair et a $4$ \'el\'ements si $n$ est pair. Supposons que $n$ est pair.   Le groupe $G_{AD}(F)_{0}/\pi(G(F))$ est isomorphe \`a $(Z(G)^{I_{F}})_{\Gamma_{{\mathbb F}_{q}}}$. Or $I_{F}$ agit sur $Z(G)$ par $\sigma(z)=z$ pour $\sigma\in I_{E}$ et $\sigma(z)=z^{-1}$ pour $\sigma\in I_{F}-I_{E}$. Le groupe pr\'ec\'edent a donc $2$ \'el\'ements (remarquons que les deux sous-groupes $G_{AD}(F)_{0}/\pi(G(F))$ et $E^1/(E^1)^n$ de $G_{AD}(F)/\pi(G(F)$ ont tous deux $2$ \'el\'ements mais ne sont pas \'egaux).  On note $\Xi_{n}$ le sous-ensemble des \'el\'ements de $\Xi$, c'est-\`a-dire des caract\`eres de $G_{AD}(F)_{0}/\pi(G(F))$ dont la restriction \`a $G_{AD}(F)_{0}/\pi(G(F))$ est triviale si $n$ est divisible par $4$, non triviale sinon. On revient maintenant au cas g\'en\'eral o\`u $n$ n'est pas suppos\'e pair.

Notons ${\mathbb X}$ l'ensemble des couples $(k,h)\in {\mathbb N}^2$ tels que $h^2+k(k+1)=n$. Si  $\delta_{2\triangle}(n)=0$, on pose ${\cal X}={\mathbb X}$. Si $\delta_{2\triangle}(n)=1$ (ce qui implique que $n$ est pair), l'ensemble ${\mathbb X}$ poss\`ede un \'el\'ement $(k,h)$ tel que $h=0$, on le note $(k_{0},0)$.   On note ${\cal X}$ la r\'eunion de ${\mathbb X}-\{(k_{0},0)\}$ et de $\{(k_{0},0,\xi); \xi\in \Xi_{n}\}$. Quel que soit $\delta_{2\triangle}(n)$,  on pose $d_{x}=1$ pour tout $x\in {\cal X}$.  

On d\'etermine l'ensemble $\underline{S}(G)$ par l'alg\`ebre lin\'eaire ou les tables de Tits (le groupe est de type $C-BC_{(n-1)/2}$ si $n$ est impair, $B-C_{n/2}$ si $n$ est pair). Si $n$ est impair, $\underline{S}(G)$ est en bijection avec l'ensemble des couples $(a,b)\in {\mathbb N}^2$ tels que $a+b=n$, $a\not=2$  et $b$ est pair. Si $n$ est pair, il y a une surjection de $\underline{S}(G)$ sur cet ensemble dont les fibres ont un seul \'el\'ement, sauf la fibre au-dessus de $(0,n)$ qui en a deux. 
L'action de $G_{AD}(F)$ est triviale si $n$ est impair. Si $n$ est pair, ses orbites sont exactement les fibres de l'application pr\'ec\'edente. Pour un sommet $s$ param\'etr\'e par $(a,b)$, on a $G_{s}=SO_{dep}(a)\times Sp(b)$.  
On a vu en \ref{Bn} que $FC(\mathfrak{so}_{dep}(a,{\mathbb F}_{q}))$ n'est non nul que si $a$ est de la forme $a=h^2$.   On a vu en \ref{Cn} que $FC(\mathfrak{sp}(b,{\mathbb F}_{q}))$ n'est non nul que si  $b$ est de la forme $b=k(k+1)$. Supposons ces conditions  v\'erifi\'ees. Alors  ces espaces sont des droites, donc $FC(\mathfrak{g}_{s}({\mathbb F}_{q}))$ est lui-aussi une droite.   Si $a\not=0$, de cette droite se d\'eduit d'apr\`es \ref{orbites} une droite dans $FC(\mathfrak{g}(F))$ que l'on note $FC_{k,h}$. Supposons $a=0$ et $b=k(k+1)$. Ces conditions se produisent si et seulement si $\delta_{2\triangle}(n)=1$ et alors $k=k_{0}$.  L'espace  $ FC(\mathfrak{sp}(k_{0}(k_{0}+1)),{\mathbb F}_{q}))$ est la droite port\'ee  par une fonction $f_{N,\epsilon}$. Le caract\`ere $\epsilon$ en question est trivial sur $Z(Sp(k_{0}(k_{0}+1))$ si $k_{0}(k_{0}+1)/2$ est pair et non trivial sinon. Ces conditions sont \'equivalentes \`a $n$ divisible ou non par $4$. La fonction se transforme donc selon le caract\`ere trivial de $G_{AD}(F)_{0}$ si $4$ divise $n$ et selon le caract\`ere non trivial de ce groupe si $4$ ne divise pas $n$. D'apr\`es \ref{orbites}, pour tout caract\`ere $\xi\in \Xi_{n}$, la fonction $f_{N,\epsilon}$ donne donc naissance \`a un \'el\'ement de $FC(\mathfrak{g}(F))$ se transformant selon le caract\`ere $\xi$ de $G_{AD}(F)/\pi(G(F))$. On note $FC_{0,k_{0},\xi}$ la droite port\'ee par cette fonction.  On a ainsi associ\'e une droite $FC_{x}$ \`a tout $x\in {\cal X}$  de sorte que \ref{resultats}(1) soit  v\'erifi\'ee.

Notons ${\mathbb Y}$ l'ensemble des couples $(i,j)$ tels que $i,j\in {\mathbb N}$, $i\leq j$, $ i(i+1)/2+j(j+1)/2=n$. Si $\delta_{2\triangle}(n)=0$, on pose ${\cal Y}={\mathbb Y}$. Supposons $\delta_{2\triangle}(n)=1$.    Le couple   $(k_{0},k_{0})$ appartient \`a ${\mathbb Y}$. 
 On note ${\cal Y}$ la r\'eunion de ${\mathbb Y}-\{(k_{0},k_{0})\}$ et de $\{(k_{0},k_{0},\xi); \xi\in \Xi_{n}\}$.

   D\'eterminons les donn\'ees endoscopiques elliptiques ${\bf G}$ de $G$ telles que $FC^{st}(\mathfrak{g}'(F))^{Out({\bf G}')}$ $\not=\{0\}$. On identifie $\hat{\Delta}_{a}$ \`a ${\mathbb Z}/n{\mathbb Z}$ comme en \ref{An-1quasidepnonram}, l'action $\sigma\mapsto \sigma_{G}$ et le groupe $\hat{\Omega}$ sont comme dans ce paragraphe, l'extension $E/F$ \'etant maintenant ramifi\'ee. On consid\`ere une donn\'ee $(\sigma\mapsto \sigma_{G'},{\cal O})\in {\cal E}_{ell}(G)$. A \'equivalence pr\`es, c'est-\`a-dire modulo conjugaison par l'action de $\hat{\Omega}$, on peut supposer et on suppose que $0\in {\cal O}$.  Rappelons que l'on a not\'e   $E_{G'}/F$ l'extension telle que  $\Gamma_{E_{G'}}$ soit le noyau de l'action   $\sigma\mapsto \sigma_{G'}$. On sait que $E\subset E_{G'}$ et que $\Gamma_{E_{G'}}/\Gamma_{E}$ s'envoie injectivement dans $\hat{\Omega}$. Les propri\'et\'es suivantes en r\'esultent:
   $E_{G'}/E$ est une  extension cyclique de degr\'e divisant $n$;  pour $\rho\in \Gamma_{E_{G'}/E}$ et $\delta\in \Gamma_{E_{G'}/F}-\Gamma_{E_{G'}/E}$, on a les \'egalit\'es $\delta^2=1$ et $\delta\rho=\rho^{-1}\delta$. D'apr\`es \ref{extensionsdiedrales}(4), cela entra\^{\i}ne $E_{G'}=E$ ou $K=Q$ (l'extension biquadratique de $F$) et, dans ce dernier cas, $n$ est pair.

Supposons d'abord $E_{G'}=E$. L'orbite ${\cal O}$ est r\'eduite au point $0$ ou \`a deux points $0,c$, avec $0<c<n$. Le premier cas donne la donn\'ee endoscopique principale ${\bf G}$ dont on ne conna\^{\i}t pas encore l'espace $FC^{st}(\mathfrak{g}(F))$ associ\'e. Laissons-la provisoirement en suspens. Dans le second cas, on remarque d'abord que l'on peut supposer $0<c\leq n/2$. En effet,  si $c>n/2$, on conjugue nos donn\'ees par l'action de  l'application $j\mapsto j-c$ qui appartient \`a $\hat{\Omega}$. Cela  remplace $c$ par $-c=n-c<n/2$. 
On a $G'(F)=\{(u,v)\in U_{E/F}(c,F)\times U_{E/F}(n-c,F); det(u)det(v)=1\}$.  D'apr\`es le lemme \ref{centre}, $FC^{st}(\mathfrak{g}'(F))$ est le m\^eme espace que celui associ\'e au groupe $G'_{SC}$. On peut utiliser pour ce groupe l'assertion  (5) ci-dessous par r\'ecurrence puisque $c,n-c<n$. L'espace $FC^{st}(\mathfrak{g}'(F))$ n'est non nul que si $c$ et $n-c$ sont de la forme $ i(i+1)/2$ et $j(j+1)/2$. Supposons ces conditions v\'erifi\'ees.  L'espace $FC^{st}(\mathfrak{g}'(F))$ est alors  une droite.
 La donn\'ee ${\bf G}'$ n'a pas d'automorphisme non trivial sauf si $n$ est pair et $c=n/2$ ou encore $i=j$. Dans ce cas, il y a un automorphisme non trivial qui \'echange les deux facteurs $U_{E/F}(n/2)$. Mais cet automorphisme agit trivialement sur $FC^{st}(\mathfrak{g}'(F))$ donc $FC^{st}(\mathfrak{g}'(F))^{Out({\bf G}')}$ est une droite. 
 Si $i\not=j$, on note  ${\bf G}'_{i,j}$ notre donn\'ee endoscopique et on pose $FC^{\cal E}_{i,j}=FC^{st}(\mathfrak{g}_{i,j}'(F))^{Out({\bf G}_{i,j}')}$.  Supposons $i=j$. Alors $\delta_{2\triangle}(n)=1$ et que $i=k_{0}$. On pose alors $\xi=\xi_{{\bf G}'}$, on note ${\bf G}'_{k_{0},k_{0},\xi}$ notre donn\'ee endoscopique et on pose $FC^{\cal E}_{k_{0},k_{0},\xi}=FC^{st}(\mathfrak{g}_{k_{0},k_{0},\xi}'(F))^{Out({\bf G}_{k_{0},k_{0},\xi}')}$. Remarquons qu'\`a ce point, on ne sait pas encore que $\xi\in \Xi_{n}$.

Supposons maintenant que $E_{G'}=Q$, donc $n$ est pair. On note $\Gamma_{Q/E}=\{1,\rho\}$   et  $\delta$ l'\'el\'ement non trivial de $\Gamma_{Q/F}$ qui fixe tout \'el\'ement de l'extension quadratique non ramifi\'ee $E_{0}$ de $F$. On a forc\'ement $\rho_{G'}(j)=j+n/2$ pour tout $j\in {\mathbb Z}/n{\mathbb Z}$ et il existe $u\in \{0,...,n-1\}$ tel que $\delta_{G'}(j)=-j+u$. L'orbite ${\cal O}$, qui est celle de $0$ d'apr\`es notre hypoth\`ese $0\in {\cal O}$, peut avoir $2$ ou $4$ \'el\'ements. Si elle en a $4$, le fixateur d'une racine dans ${\cal O}$ est $\Gamma_{Q}$. Puisque $Q/F$ n'est pas totalement ramifi\'ee, la donn\'ee endoscopique correspondante v\'erifie $FC(\mathfrak{g}'(F))=\{0\}$ d'apr\`es le lemme \ref{centre}. On peut exclure ce cas. Supposons donc que ${\cal O}$ n'a que deux \'el\'ements. Alors forc\'ement ${\cal O}=\{0,n/2\}$ donc $u=0$ ou $u=n/2$. Posons $e=0$ si $u=0$, $e=1$ si $u=n/2$. Alors $(\rho^{e}\delta)_{G'}(j)=-j$ pour tout $j$. Le fixateur dans $\Gamma_{Q/F}$ de la racine $\hat{\alpha}_{0}\in {\cal O}$ est $\{1,\rho^{e}\delta\}$. D'apr\`es le lemme \ref{centre}, on peut supposer que le sous-corps des points fixes de ce groupe est une extension ramifi\'ee de $F$. D'apr\`es la d\'efinition de $\delta$, cela impose $e=1$ et $u=n/2$.  
   On voit que  $G'(F)=\{x\in U_{Q/E_{0}}(n/2,E_{0}); norme_{Q/E}(det(x))=1\}$.  L'extension $Q/E_{0}$ est ramifi\'ee. On peut comme dans le cas $E_{G'}=E$ utiliser (5) ci-dessous par r\'ecurrence, en rempla\c{c}ant le corps de base $F$ par $E_{0}$: l'espace $FC^{st}(\mathfrak{g}'(F))$ n'est non nul que si $n/2$ est de la forme $i(i+1)/2$. On a alors n\'ecessairement $i=k_{0}$.
 Supposons qu'il en soit ainsi. L'espace $FC^{st}(\mathfrak{g}'(F))$  est   une droite qui poss\`ede un g\'en\'erateur naturel.  La donn\'ee ${\bf G}'$ a un unique automorphisme non trivial, qui est l'action de l'\'el\'ement  $j\mapsto j+n/2$ de $\hat{\Omega}$, ou encore celle de $\rho$.  C'est-\`a-dire que cet automorphisme agit sur $G'(F)\subset U_{Q/E_{0}}(n/2,E_{0})$ par l'action naturelle de $\rho$ sur ce dernier groupe. On doit prouver que l'action d\'eduite sur  $FC^{st}(\mathfrak{g}'(F))$ est triviale. On peut remplacer $G'$ par $G'_{SC}$. Posons $h=[(k_{0}+1)/2]$, $k=[k_{0}/2]$.
  On a $n/2=h^2+k(k+1)$. On utilise les r\'esultats du pr\'esent paragraphe par r\'ecurrence, pr\'ecis\'ement la description ci-dessous de ${\cal X}^{st}$.  Le
  g\'en\'erateur de $FC^{st}(\mathfrak{g}'(F))$ est la fonction  issue d'un sommet $s'$ de l'immeuble de $G'_{AD}$ pour lequel $G'_{SC,s'}({\mathbb F}_{q})= SO_{dep}(h^2,{\mathbb F}_{q^2})\times  Sp(k(k+1),{\mathbb F}_{q^2})$. Ce groupe est muni d'un automorphisme naturel d\'eduit du Frobenius $Fr$. 
On voit que l'action de $\rho$ conserve le sommet $s'$ de l'immeuble de $G'_{AD}$ et agit sur $G'_{SC,s'} $ comme ce Frobenius. En se reportant \`a la description des g\'en\'erateurs $f_{N,\epsilon}$ pour les deux composantes de l'espace $\mathfrak{g}'_{SC,s'}({\mathbb F}_{q})$, on voit que les supports de ces fonctions contiennent des points fixes par l'action de Frobenius donc elles ne peuvent qu'\^etre fix\'ees par cette action. On obtient l'assertion souhait\'ee: $Out({\bf G}')$ agit trivialement sur $FC^{st}(\mathfrak{g}'(F))$ donc $FC^{st}(\mathfrak{g}'(F))^{Out({\bf G}')}$ est une droite. Posons $\xi=\xi_{{\bf G}'}$. 
On note ${\bf G}'_{k_{0},k_{0},\xi}$ notre donn\'ee endoscopique et on pose  $FC^{\cal E}_{k_{0},k_{0},\xi}=FC^{st}(\mathfrak{g}_{k_{0},k_{0},\xi}'(F))^{Out({\bf G}_{k_{0},k_{0},\xi}')}$.  Comme plus haut, on ne sait pas encore que $\xi\in \Xi_{n}$. 

A ce point, on a obtenu une description de $FC^{{\cal E}}(\mathfrak{g}(F))$ qui ressemble \`a celle de \ref{resultats} (2). Il y a trois diff\'erences. On n'a pas trait\'e   la  donn\'ee principale ${\bf G}$. Dans le cas o\`u il existe un \'el\'ement $y=(i,j)\in {\mathbb Y}$ tel que $i=0$ (cet \'el\'ement est alors unique), on n'a pas d\'efini d'espace $FC^{{\cal E}}_{y}$ (dans la construction du cas $E_{G'}=E$, on a suppos\'e $0<c$ ce qui impose $i>0$). Dans le cas o\`u $\delta_{2\triangle}(n)=1$, on a bien construit deux droites $FC^{\cal E}_{k_{0},k_{0},\xi}$ mais on n'a prouv\'e ni que les caract\`eres $\xi$ intervenant \'etaient distincts, ni qu'ils appartenaient \`a  $\Xi_{n}$.  

Traitons ce dernier point. 
   On suppose donc $n$ pair et $\delta_{2\triangle}(n)=1$. On utilise les g\'en\'erateurs $\rho$ et $\delta$ de $\Gamma_{Q/F}$ introduits ci-dessus. Notons ${\bf G}'_{ 1}$ et ${\bf G}'_{ 2}$ les deux donn\'ees endoscopiques  not\'ees ${\bf G}'_{k_{0},k_{0},\xi}$ ci-dessus.  Elles sont associ\'ees 
  \`a la m\^eme orbite ${\cal O}=\{0,n/2\}$.  La diff\'erence est dans les actions galoisiennes.    On a $\delta_{G'_{1}}(j)=\delta_{G'_{2}}(j)=-j+n/2$  mais $\rho_{G'_{1}}(j)=j$ et $\rho_{G'_{2}}(j)=j+n/2$ pour tout $j\in {\mathbb Z}/n{\mathbb Z}$. On calcule les \'el\'ements de $H^1(W_{F},Z(\hat{G}_{SC}))$ associ\'es \`a ces deux actions. En identifiant $Z(\hat{G}_{SC})$ \`a $\boldsymbol{\zeta}_{n}({\mathbb C})$, ce sont les deux cocycles  triviaux sur $W_{Q}$ et tels que $\delta\mapsto (-1)^{n/2}$, $\rho\mapsto  1$ pour ${\bf G}'_{1}$, $\rho\mapsto -1$ pour ${\bf G}'_{2}$. En se rappelant la d\'efinition de $\delta$, on obtient les deux cocycles dont la restriction \`a $I_{F}$ est triviale si $n/2$ est pair, non triviale si $n/2$ est impair. En utilisant les suites exactes que l'on a \'ecrites dans \ref{actionsurFC}, on voit que cet ensemble de cocycles s'identifie \`a $\Xi_{n}$, comme on le voulait.

 On d\'efinit une application 
$$\begin{array}{cccc}\phi:&{\mathbb X}&\to&{\mathbb Y}\\ &(k,h)&\mapsto&\varphi(k,h)=(i,j)\\ \end{array}$$
par les formules suivantes:

si $k\geq h$, $i=k-h$, $j=k+h$;

si $k<h$, $i=h-k-1$, $j=h+k$.

On v\'erifie que $\phi$ est une bijection. Si $\delta_{2\triangle}(n)=0$, on pose $\varphi=\phi$. Supposons $\delta_{2\triangle}(n)=1$. On voit que $\phi(k_{0},0)=(k_{0},k_{0})$. Alors $\phi$ se rel\`eve en une bijection $\varphi:{\cal X}\to {\cal Y}$ qui envoie $(k_{0},0,\xi)$ sur $(k_{0},k_{0},\xi)$ pour $\xi\in \Xi_{n}$. 
On note ${\cal X}^{st}$ le sous-ensemble des $(k,h)\in {\mathbb X}$ tels que $h=k$ ou $h=k+1$. On voit que ${\cal X}^{st}$ a au plus un \'el\'ement.  Il correspond par $\varphi$ au sous-ensemble ${\cal Y}^{st}$  des $(i,j)\in {\mathbb Y}$ tels que $i=0$, qui a lui aussi au plus un \'el\'ement (c'est l'\'el\'ement que l'on n'a pas obtenu dans la description ci-dessus de $FC^{{\cal E}}(\mathfrak{g}(F))$).  Si ${\cal X}^{st}$ est non vide, on note son unique \'el\'ement $(k^{st},h^{st})$ (remarquons qu'il est forc\'ement diff\'erent de l'\'eventuel \'el\'ement $(k_{0},0)$) et on note $(0,j^{st})$ celui de ${\cal Y}^{st}$.    Utilisons  l'\'egalit\'e
$$dim(FC(\mathfrak{g}(F)))= dim(FC^{{\cal E}}(\mathfrak{g}(F))$$
$$=dim(FC^{st}(\mathfrak{g}(F)))+\sum_{{\bf G}'\in Endo_{ell}(G), {\bf G}'\not={\bf G}}dim(FC^{st}(\mathfrak{g}'(F))^{Out({\bf G}')}).$$
Puisque tous les sous-espaces $FC_{x}$ et $FC^{{\cal E}}_{y}$ que nous avons construits sont des droites, les r\'esultats d\'ej\`a obtenus entra\^{\i}nent que $dim(FC(\mathfrak{g}(F)))=\vert {\cal X}\vert=\vert {\cal Y}\vert  $ tandis que la derni\`ere somme de l'expression ci-dessus vaut $\vert {\cal Y}\vert -\vert {\cal Y}^{st}\vert  $. Par diff\'erence, on obtient $dim(FC^{st}(\mathfrak{g}(F)))=\vert {\cal Y}^{st}\vert $. Si   ${\cal Y}^{st}=\emptyset$, nos constructions pr\'ec\'edentes ont donc d\'ecrit tout $FC^{{\cal E}}(\mathfrak{g}(F))$. Si ${\cal Y}^{st}$ a un \'el\'ement   $(0,j^{st})$, il suffit de poser $FC^{{\cal E}}_{0,j^{st}}=FC^{st}(\mathfrak{g}(F))$ pour compl\'eter ces constructions et obtenir \ref{resultats} (2).

  Munissons ${\mathbb X}$ de la relation $(k,h)\leq (k',h')$ si et seulement si $h+k\leq h'+k'$. Montrons que

(1) c'est une relation d'ordre total.

Soit $x=(k,h)\in {\mathbb X}$. Posons $a(x)=2h+2k+1$, $b(x)=2k+1-2h$. On v\'erifie que $a(x)^2+b(x)^2=8h^2+8k(k+1)+2=8n+2$. Soit $x'=(k',h')\in {\mathbb X}$, supposons $x\leq x'$ et $x'\leq x$. Alors $a(x)=a(x')$. L'\'egalit\'e pr\'ec\'edente entra\^{\i}ne $b(x)=\pm b(x')$. Si $b(x)=-b(x')$, on a $ 4k+2=a(x)+b(x)=a(x')-b(x')=4h'$, ce qui est impossible, les deux termes extr\^emes n'\'etant pas congrus modulo $4$. Donc $b(x)=b(x')$ puis $x'=x$.   Cela prouve (1).

Par la surjection naturelle ${\cal X}\to {\mathbb X}$, on rel\`eve notre relation d'ordre sur ${\mathbb X}$ en une relation de pr\'eordre sur ${\cal X}$ que l'on note encore $\leq$. 
  Dans les constructions ci-dessus, on a associ\'e \`a tout $y\in {\cal Y}$ une donn\'ee endoscopique ${\bf G}'_{y}$ et l'application $y\mapsto {\bf G}'_{y}$ est injective.   Soit $y\in {\cal Y}$.  On introduira dans le paragraphe \ref{preuve}  un \'el\'ement  $G$-r\'egulier $Y_{ y}\in \mathfrak{g}'_{ y, ell}(F)$ qui a les propri\'et\'es suivantes:

(2) pour tout \'el\'ement non nul $f'\in FC^{{\cal E}}_{y}=FC^{st}(\mathfrak{g}'_{y}(F))^{Out({\bf G}'_{y})}$, on a $S^{G'_{ y}}(Y_{y},f_{{\bf G}'_{y}})\not=0$;

(3) soient $x\in {\cal X}$ et $f\in FC_{x}$;  soit $X$ un \'el\'ement de $\mathfrak{g}_{ell}(F)$ correspondant \`a $Y_{ y}$; supposons $I^G(X,f)\not=0$; alors $\varphi^{-1}(y)\leq x$. 

Alors les hypoth\`eses (1) \`a (5) de \ref{ingredients} sont satisfaites pour $\underline{{\cal Y}}^{\sharp}=\underline{{\cal Y}}$ avec les notations de ce paragraphe. Cela entra\^{\i}ne

$$(4) \qquad transfert(FC_{(x)})=FC^{{\cal E}}_{\varphi((x))}$$
pour tout $((x))\in \underline{{\cal X}}$. Nos classes d\'equivalence $((x))$ sont presque toutes r\'eduites \`a un \'el\'ement $x$. Dans ce cas, la relation ci-dessus \'equivaut \`a $transfert(FC_{x})=FC^{{\cal E}}_{\varphi(x)}$. Si $\delta_{2\triangle}(n)=1$, il y a la classe$ \{(k_{0},0,\xi); \xi\in \Xi_{n}\}$ qui a deux \'el\'ements. Mais les deux droites $FC_{k_{0},0,\xi}$ pour $\xi\in \Xi_{n}$ se distinguent par l'action de $G_{AD}(F)/\pi(G(F))$ qui agit par $\xi$ sur $FC_{0,k_{0},\xi}$. Du c\^ot\'e endoscopique, les deux droites $FC^{\cal E}_{k_{0},k_{0},\xi}$ de l'espace correspondant se distinguent aussi par le caract\`ere $\xi$ associ\'e \`a la donn\'ee ${\bf G}'_{k_{0},k_{0},\xi}$. Puisque le transfert est compatible avec l'action de $G_{AD}(F)/\pi(G(F))$, l'\'egalit\'e (4) se raffine en les \'egalit\'es 
$$transfert(FC_{k_{0},0,\xi})=FC^{{\cal E}}_{k_{0},k_{0},\xi}$$
pour $\xi\in \Xi_{n}$. Cela d\'emontre \ref{resultats}(3).

La relation \ref{resultats} (4) r\'esulte de \ref{resultats} (3) et des constructions: $FC^{st}(\mathfrak{g}(F))$ est l'image r\'eciproque par l'application $transfert$ du sous-espace $FC^{st}(\mathfrak{g}(F))\subset FC^{{\cal E}}(\mathfrak{g}(F))$. Ce dernier espace est par d\'efinition $FC^{{\cal E}}_{0,j^{st}}$, donc $FC^{st}(\mathfrak{g}(F))=FC_{\varphi^{-1}(0,j^{st})}=FC_{k^{st},h^{st}}$. Remarquons que ${\cal Y}^{st}\not=\emptyset$ si et seulement si $\delta_{\triangle}(n)=1$.  L'assertion suivante est alors cons\'equence de \ref{resultats} (4):

(5) $dim(FC^{st}(\mathfrak{g}(F)))=\delta_{\triangle}(n)$.

\subsection{Forme int\'erieure du type $A_{n-1}$ quasi-d\'eploy\'e, $E/F$ ramifi\'ee}
On suppose que $G^*$ est du type pr\'ec\'edent mais que $G$ n'est pas quasi-d\'eploy\'e. Comme en \ref{formeinterieurenonramifiee}, cela impose que $n$ est pair et que $G$ est l'unique forme non quasi-d\'eploy\'ee du groupe $SU_{E/F}(n)$. 

On conserve les objets ${\mathbb X}$, ${\mathbb Y}$ et $\phi$ du cas pr\'ec\'edent.  On note ${\cal X}$ le sous-ensemble des $(k,h)\in {\mathbb X}$ tels que $h\not=0$ et on pose $d_{x}=1$ pour tout $x\in {\cal X}$.  

L'ensemble $\underline{S}(G)$ est maintenant param\'etr\'e par les couples $(a,b)\in {\mathbb N}^2$ tels que $a+b=n$, $b$ est pair et $a\not=0$. Pour un sommet $s$ param\'etr\'e par $(a,b)$, on a $G_{s}=SO_{{\mathbb F}_{q^2}/{\mathbb F}_{q}}(a)\times Sp(b)$, o\`u $SO_{{\mathbb F}_{q^2}/{\mathbb F}_{q}}(a)$ est la forme non d\'eploy\'ee du groupe sp\'ecial orthogonal. La construction de la droite $FC_{x}$ pour $x\in {\cal X}$ est alors identique \`a celle du paragraphe pr\'ec\'edent. Les petites difficult\'es caus\'ees par les sommets qui \'etaient conjugu\'es par le groupe $G_{AD}(F)$ mais pas par $G(F)$ disparaissent car ces sommets disparaissent.

 On note  ${\cal Y}$ le sous-ensemble des $(i,j)\in {\mathbb Y}$ tels que $i<j$. La construction de la droite $FC^{{\cal E}}_{y}$ pour $y\in {\cal Y}$ est la m\^eme que dans le paragraphe pr\'ec\'edent. Les donn\'ees endoscopiques un peu exceptionnelles param\'etr\'ees par l'\'eventuel couple $(k_{0},k_{0})$ ne contribuent pas: ces donn\'ees sont les seules pour lesquelles le groupe d'automorphismes ext\'erieurs  $Out({\bf G}')$ est non trivial. Maintenant, ce groupe agit par son unique caract\`ere non trivial sur $FC^{st}(\mathfrak{g}'(F))$ donc $FC^{st}(\mathfrak{g}'(F))^{Out({\bf G}'}=\{0\}$.
 
 On voit que $\phi({\cal X})={\cal Y}$ et on note $\varphi$ l'application $\phi$ restreinte \`a ${\cal X}$.  On prouve l'assertion (3) de \ref{resultats} comme dans le paragraphe pr\'ec\'edent.

  \section{Calcul d'int\'egrales orbitales, type $A_{n-1}$ quasi-d\'eploy\'e, $E/F$ ramifi\'ee}

\subsection{Description explicite des \'el\'ements de $FC(\mathfrak{g}(F))$}\label{descriptionexplicite}
Soit $E/F$ une extension quadratique ramifi\'ee. On fixe une uniformisante $\varpi_{E}\in E^{\times}$ telle que $\varpi_{E}^2\in F^{\times}$. Soit $V$ un espace vectoriel sur $E$ de dimension $n\geq 2$, muni d'une forme hermitienne $q$ (relativement \`a $E/F$). 
 On note ${\bf  G}$ le groupe unitaire de $(V,q)$ et $G$ le sous-groupe sp\'ecial unitaire. On note $\boldsymbol{\mathfrak{g}}$ et $\mathfrak{g}$  leurs alg\`ebres de Lie. On a  ${\bf G}(E)=GL_{E}(V)$ et $\boldsymbol{\mathfrak{g}}(E)=End_{E}(V)$.   Pour $X\in \boldsymbol{\mathfrak{g}}(E)$, on note $X^{\star}$ l'\'el\'ement de $\boldsymbol{\mathfrak{g}}(E)$ tel que $q(Xv,v')=q(v,X^{\star}v')$ pour tous $v,v'\in V$. Alors $\boldsymbol{\mathfrak{g}}(F)$ est l'ensemble des $X\in \boldsymbol{\mathfrak{g}}(E)$ tels que $X=-X^{\star}$ et $\mathfrak{g}(F)$ est le sous-ensemble des $X\in \boldsymbol{\mathfrak{g}}(F)$ tels que $trace(X)=0$. Si $V'$ et $V''$ sont deux sous-$E$-espaces de $V$  en dualit\'e pour la forme $q$ et si $R$ est un sous-$\mathfrak{o}_{E}$-r\'eseau de $V'$, on note $R^*=\{v''\in V''; \forall v'\in R', q(v'',v')\in \mathfrak{o}_{E}\}$.

 L'ensemble  $\boldsymbol{\mathfrak{g}}(E)$    est non seulement une alg\`ebre de Lie mais une alg\`ebre tout court pour le produit matriciel habituel. Quand on parlera de sous-alg\`ebre de $\boldsymbol{\mathfrak{g}}(E)$, il s'agira d'une sous-alg\`ebre pour ce produit matriciel. Nous d\'efinirons divers sous-$\mathfrak{o}_{E}$-modules de $\boldsymbol{\mathfrak{g}}(E)$ invariants par l'application $X\mapsto X^{\star}$ et nous utiliserons leurs intersections avec $\boldsymbol{\mathfrak{g}}(F)$ et $\mathfrak{g}(F)$. Dans ce cas, pour simplifier les notations, nous noterons le sous-module de $\boldsymbol{\mathfrak{g}}(E)$  par une lettre gothique grasse affect\'ee d'un exposant $E$, son intersection avec $\boldsymbol{\mathfrak{g}}(F)$ par la  m\^eme lettre grasse sans exposant et son intersection avec $\mathfrak{g}(F)$ par la lettre gothique fine sans exposant. Par exemple, $\boldsymbol{\mathfrak{h}}^{E}\subset \boldsymbol{\mathfrak{g}}(E)$, $\boldsymbol{\mathfrak{h}}=\boldsymbol{\mathfrak{h}}^{E}\cap \boldsymbol{\mathfrak{g}}(F)$ et $\mathfrak{h}=\boldsymbol{\mathfrak{h}}\cap \mathfrak{g}(F)$.

   On note $sgn$ l'unique caract\`ere d'ordre $2$ de ${\mathbb F}_{q}^{\times}$, ou son rel\`evement en un caract\`ere de $\mathfrak{o}_{F}^{\times}$. Pour tout entier $m\in {\mathbb N}$, on note $\chi$ le caract\`ere de $({\mathbb F}_{q}^{\times})^m$  d\'efini par $\chi(\underline{\gamma})=\prod_{i=1,...,m}sgn(\gamma_{i})^{i}$ pour tout $\underline{\gamma}=(\gamma_{i})_{i=1,...,m}\in ({\mathbb F}_{q}^{\times})^m$. On note encore $\chi$ le caract\`ere similaire de $(\mathfrak{o}_{F}^{\times})^m$.
  
  Soient deux entiers $h^{o},h^s\in {\mathbb N}$, supposons $n={h^{o}}^2+h^s(h^s+1)$. Les exposants $o$ et $s$ \'evoquent "orthogonal" et "symplectique" pour une raison qui va appara\^{\i}tre et seront trait\'es comme des symboles math\'ematiques.  Supposons d'abord  $h^{o}\not=0$. On peut choisir
   une d\'ecomposition en sous-espaces
 $$V= \oplus_{i=-2h^s,...2h^{o}-2}V_{i}$$

\noindent et, pour tout   indice $i\in \{-2h^{s},...,2h^{o}-2\}$, un $\mathfrak{o}_{E}$-sous-r\'eseau $R_{i}\subset V_{i}$
 de sorte que
 
 pour $i,j\in \{-2h^{s},...,2h^{o}-2\}$, $V_{i}$ et $V_{j}$ sont  orthogonaux sauf si  $i,j\geq0$ et $i+j=2h^{o}-2$ ou si  $i,j\leq-1$ et $i+j=-1-2h^s$, auquel cas ils sont en dualit\'e; en particulier,  si $h^{o}\geq1$, l'espace $V_{h^{o}-1}$ est non d\'eg\'en\'er\'e;

 pour $i=0,...,h^{o}-1$, $dim_{E}(V_{i})= i+1 $ et, pour $i=-h^s,...,-1$, $dim_{E}(V_{i})=\vert i\vert$;
 
 pour $i\geq0$,    $R_{i}=R^*_{2h^{o}-2-i}$ pour la dualit\'e entre $V_{i}$ et $V_{2h^{o}-2-i}$;
 
 pour $i\leq-1$,  $R_{i}=\mathfrak{p}_{E}^{-1}R^*_{-2h^s-1-i}$ pour  la dualit\'e entre $V_{i}$ et $V_{-2h^s-1-i}$ .

 Posons $d=2h^{o}+2h^s-1$. On a d\'efini les $\mathfrak{o}_{E}$-modules $R_{i}$ pour $i\in \{-2h^s,...,2h^{o}-2\}$. On prolonge cette d\'efinition en posant $R_{i}=\mathfrak{p}_{E}^mR_{j}$ pour $i=j+md\in {\mathbb Z}$ avec $j\in \{-2h^s,...,2h^{o}-2\}$ et $m\in {\mathbb Z}$.  
 On pose $R_{\geq i}=\sum_{j\geq i}R_{j}$. 
 On a $R_{\geq i}\supset R_{\geq i+1}$, $R_{\geq -2h^s}=R_{\geq 0}^*\supset R_{\geq 0}\supset R_{\geq 2h^0-1}=\mathfrak{p}_{E}R_{\geq 0}^*$. La forme $q$ se r\'eduit en une forme orthogonale   non d\'eg\'en\'er\'ee $q^{o}$ sur le ${\mathbb F}_{q}$-espace  $V^{o}=R_{\geq 0}/\mathfrak{p}_{E}R_{\geq 0}^*$  et la forme $\varpi_{E}q$ se r\'eduit en une forme symplectique $q^s$ sur  le ${\mathbb F}_{q}$-espace $V^s=R_{\geq 0}^*/R_{\geq 0}$. On note $G^{o}$, resp. $G^s$,   le groupe sp\'ecial orthogonal de $(V^{o},q^{o})$, resp. le groupe symplectique de $(V^s,q^s)$ (ces groupes sont d\'efinis sur ${\mathbb F}_{q}$). 
  Notons $\boldsymbol{\mathfrak{k}}^{E}$ la sous-$\mathfrak{o}_{E}$-alg\`ebre de $\boldsymbol{\mathfrak{g}}(E)$ form\'ee des \'el\'ements qui conservent les r\'eseaux $R_{\geq 0}$ et $R_{\geq 0}^*$. Notons $(\boldsymbol{\mathfrak{k}}^{E})^{\perp}$   l'id\'eal de $\boldsymbol{\mathfrak{k}}^{E}$ form\'e des \'el\'ements $X$ tels que $X(R_{\geq 0}^*)\subset R_{\geq 0}$ et $X(R_{\geq 0})\subset \mathfrak{p}_{E}R_{\geq 0}^*$.   On a les isomorphismes $\boldsymbol{\mathfrak{k}}/\boldsymbol{\mathfrak{k}}^{\perp} =\mathfrak{k}/\mathfrak{k}^{\perp}\simeq \mathfrak{g}^{o}({\mathbb F}_{q})\oplus \mathfrak{g}^s({\mathbb F}_{q})$. 
 
 Pour $i\in {\mathbb Z}$, on note $\boldsymbol{\mathfrak{q}}^{E}_{i/d}$ le sous-espace de $\boldsymbol{\mathfrak{g}}(E)$ form\'e des \'el\'ements $X$ tels que $X(R_{j})\subset R_{i+j}$ pour tout $j\in {\mathbb Z}$. On pose $\boldsymbol{\mathfrak{p}}^{E}_{ i/d}=\sum_{j\geq i}\boldsymbol{\mathfrak{q}}^{E}_{j/d}$.  Pour $i,j\in {\mathbb Z}$ le produit matriciel envoie $\boldsymbol{\mathfrak{p}}^{E}_{ i/d}\times \boldsymbol{\mathfrak{p}}^{E}_{ j/d}$ dans $\boldsymbol{\mathfrak{p}}^{E}_{ (i+j)/d}$. En particulier $\boldsymbol{\mathfrak{p}}^{E}_{ 0}$ est une alg\`ebre et $\boldsymbol{\mathfrak{p}}^{E}_{ i/d}$ en est un id\'eal pour tout $i\geq0$. De plus, $\boldsymbol{\mathfrak{p}}^{E}_{0}\subset \boldsymbol{\mathfrak{k}}^{E}$.   Remarquons que $\boldsymbol{\mathfrak{q}}^{E}_{i/d}$ et $\boldsymbol{\mathfrak{p}}^{E}_{ i/d}$ sont invariants par l'application $X\mapsto X^{\star}$. Notons que la "p\'eriodicit\'e" des applications $i\mapsto \boldsymbol{\mathfrak{p}}^{E}_{i/d}$ et $i\mapsto \boldsymbol{\mathfrak{p}}_{i/d}$ n'est pas tout-\`a-fait la m\^eme: on a $\boldsymbol{\mathfrak{p}}^{E}_{1+i/d}=\mathfrak{p}_{E}\boldsymbol{\mathfrak{p}}^{E}_{i/d}$ et $\boldsymbol{\mathfrak{p}}_{2+i/d}=\mathfrak{p}_{F}\boldsymbol{\mathfrak{p}}_{i/d}$.

   Pour $i\geq0$, notons $\mathfrak{g}^{o}_{2i}$ et $\mathfrak{g}^{o}_{ \geq 2i}$, resp. $\mathfrak{g}^{s}_{2i}$ et $\mathfrak{g}^{s}_{ \geq 2i}$, les projections dans $\mathfrak{g}^{o}({\mathbb F}_{q})$, resp. $\mathfrak{g}^s({\mathbb F}_{q})$, de $(\boldsymbol{\mathfrak{q}}_{i/d}+\boldsymbol{\mathfrak{k}}^\perp)/\boldsymbol{\mathfrak{k}}^\perp$ et $(\boldsymbol{\mathfrak{p}}_{ i/d}+\boldsymbol{\mathfrak{k}}^\perp)/\boldsymbol{\mathfrak{k}}^\perp$. Pour $a\in \{o,s\}$, l'ensemble $\mathfrak{g}^{a}_{\geq  0}$ est une sous-alg\`ebre parabolique de $\mathfrak{g}^{a}({\mathbb F}_{q})$ et $\mathfrak{g}^{a}_{\geq2}$ en est son radical nilpotent. L'espace $\mathfrak{g}^{a}_{0}$ en est une sous-alg\`ebre de Levi et on note $G^{a}_{0}$ le sous-groupe associ\'e de $G^{a}$. L'action de $G^{a}_{0}({\mathbb F}_{q})$ dans l'espace $\mathfrak{g}^{a}_{2}$ poss\`ede un nombre fini d'orbites ouvertes.   Pour $a=s$, chaque orbite est form\'ee d'\'el\'ements nilpotents de $\mathfrak{g}^{s}({\mathbb F}_{q})$ param\'etr\'es par la partition $(2,4,...,2h^s)$ et par toutes les familles $(q_{2},q_{4},...,q_{2h^s})$ de classes d'isomorphie de formes quadratiques sur ${\mathbb F}_{q}$, non d\'eg\'en\'er\'ees de rang $1$. L'ensemble de ces familles s'identifie \'evidemment \`a $({\mathbb F}_{q}^{\times}/{\mathbb F}_{q}^{\times,2})^{h^s}$. Pour $a=o$, chaque orbite est form\'ee d'\'el\'ements nilpotents de $\mathfrak{g}^{o}({\mathbb F}_{q})$ param\'etr\'es par la partition $(1,3,...,2h^{o}-1)$ et par  les familles $(q_{1},q_{3},...,q_{2h^{o}-1})$ de classes d'isomorphie de formes quadratiques sur ${\mathbb F}_{q}$, non d\'eg\'en\'er\'ees de rang $1$, telles que $q^{o}$ soit isomorphe \`a la somme $\oplus_{i=1,...,h^o}q_{2i-1}$ \`a laquelle on ajoute un certain nombre de plans isotropes. Il existe un \'el\'ement $\gamma(V)\in {\mathbb F}_{q}^{\times}/{\mathbb F}_{q}^{\times,2}$ tel que l'ensemble de ces familles soit naturellement isomorphe \`a l'ensemble des $\underline{\gamma}=(\gamma_{i})_{i=1,...,h^{o}}\in ({\mathbb F}_{q}^{\times}/{\mathbb F}_{q}^{\times,2})^{h^o}$ tels que $(-1)^{[h^{o}/2]}\prod_{i}\gamma_{i}= \gamma(V) $.  On note cet ensemble $({\mathbb F}_{q}^{\times}/{\mathbb F}_{q}^{\times,2})^{h^{o}}(V)$.   Pour $a=o,s$, notons $\tilde{\mathfrak{g}}^{a}_{2}$ la r\'eunion de ces orbites ouvertes. 
 On d\'efinit une fonction $\tilde{f}^{a}$ sur $\mathfrak{g}^{a}({\mathbb F}_{q})$ de la fa\c{c}on suivante. Elle est \`a support dans $\tilde{\mathfrak{g}}^{a}_{2}+\mathfrak{g}^{a}_{\geq 4}$ et est invariante par translations par $\mathfrak{g}^{a}_{\geq  4}$. Pour un \'el\'ement $X\in \tilde{\mathfrak{g}}^{a}_{2}$ param\'etr\'e par un \'el\'ement $\underline{\gamma}=(\gamma_{i})_{i=1,...,h^{a}}$ de $({\mathbb F}_{q}^{\times}/{\mathbb F}_{q}^{\times,2})^{h^{o}}(V)$ si $a=o$, de $({\mathbb F}_{q}^{\times}/{\mathbb F}_{q}^{\times,2})^{h^s}$ si $a=s$, $\tilde{f}^{a}(X)= \chi(\underline{\gamma})$. Notons $P^{a}$ le sous-groupe parabolique de $G^{a}$ d'alg\`ebre de Lie $\mathfrak{g}^{a}_{\geq 0}$. La fonction $\tilde{f}^{a}$ est invariante par conjugaison par $P^{a}({\mathbb F}_{q})$. On d\'efinit la fonction $f^{a}$ sur $\mathfrak{g}^{a}({\mathbb F}_{q})$ par 
 $$f^{a}(X)=\sum_{g\in G^{a}({\mathbb F}_{q})/P^{a}({\mathbb F}_{q})}\tilde{f}(g^{-1}Xg)$$
 pour tout $X\in \mathfrak{g}^{a}({\mathbb F}_{q})$. En fait, on sait que, pour tout $X$, il existe au plus un $g\in G^{a}({\mathbb F}_{q})/P^{a}({\mathbb F}_{q})$ tel que $\tilde{f}(g^{-1}Xg)\not=0$, la somme est donc r\'eduite \`a au plus un \'el\'ement. 
 On pose $\tilde{f}=\tilde{f}^{o}\oplus \tilde{f}^s$ et $f=f^{o}\oplus f^{s}$. On note $f_{G}$ la fonction sur $\mathfrak{g}(F)$, \`a support dans $\mathfrak{k}$ et invariante par $\mathfrak{k}^{\perp}$, telle que $f_{G}(X)=f(\bar{X})$ pour tout $X\in \mathfrak{k}$, o\`u $\bar{X}$ est la r\'eduction de $X$ dans $\mathfrak{k}/\mathfrak{k}^{\perp}=\mathfrak{g}^{a}({\mathbb F}_{q})\oplus \mathfrak{g}^{s}({\mathbb F}_{q})$. On d\'efinit de m\^eme la fonction  $\tilde{f}_{G}$.  
 
 {\bf Remarque.} Les espaces $\mathfrak{g}^{a}_{2i}$ sont les m\^emes qu'en \ref{groupessurFq}.
  \bigskip
 
 Il existe un sommet $s$ de l'immeuble $Imm(G_{AD})$ (avec un autre emploi du symbole $s$) tel que 
$\mathfrak{k}$ s'identifie \`a $\mathfrak{k}_{s}$. La fonction $f$ engendre la droite $ FC(\mathfrak{g}_{s}({\mathbb F}_{q}))$.
La fonction $f_{G}$ est celle not\'ee $f_{x}$ en \ref{An-1quasidepram}, o\`u $x=(h^{s},h^{o})$. On voit que les suites de r\'eseaux $(R_{\geq i})_{i\in {\mathbb Z}}$ que l'on a introduites sont uniquement associ\'ees \`a $(h^{s},h^{o})$,  \`a conjugaison pr\`es par $G(F)$. 

Consid\'erons maintenant le cas $h^{o}=0$. On suppose alors que $G$ est quasi-d\'eploy\'e. On peut  alors choisir
   une d\'ecomposition en sous-espaces
 $$V= \oplus_{i=-2h^s,...-1}V_{i}$$
 et, pour tout   indice $i\in \{-2h^{s},...,-1\}$, un $\mathfrak{o}_{E}$-sous-r\'eseau $R_{i}\subset V_{i}$ de sorte que les propri\'et\'es pr\'ec\'edentes soient v\'erifi\'ees. Tout ce que l'on a dit ci-dessus reste vrai, aux changements suivants pr\`es:
 
  la p\'eriodicit\'e $d=2h^{o}+2h^{s}-1$ de nos suites est remplac\'ee par $d=2h^{s}$; 
  
  l'espace $V^{o}$ est nul et dispara\^{\i}t (on peut dire que notre hypoth\`ese que $G$ est quasi-d\'eploy\'e \'equivaut \`a l'\'egalit\'e $\gamma(V)=1$); 
  
  les suites de r\'eseaux $(R_{\geq i})_{i\in {\mathbb Z}}$
 sont uniques \`a conjugaison pr\`es par ${\bf G}(F)$ mais il y a deux suites possibles \`a conjugaison pr\`es par $G(F)$: \`a partir de la suite que l'on a fix\'ee, on peut remplacer les modules $R_{-1}
$ et $R_{-2h^{s}}$ respectivement par $\mathfrak{p}_{E}R_{-2h^{s}}$ et $\mathfrak{p}_{E}^{-1}R_{-1}$. 

On obtient deux fonctions $f_{G}$. Pour $x=(h^{s},0,\xi)$ comme en  \ref{An-1quasidepram}, la fonction $f_{x}$ de ce paragraphe  est combinaison lin\'eaire de ces deux fonctions. 

Toujours dans le cas $h^{o}=0$, on introduit la variante suivante de nos suites de r\'eseaux. Pour $i\in \{-2h^{s}+1,...,-1\}$, on pose $R^{\natural}_{i}=R_{i}$ si $i\not=-1$, $R^{\natural}_{-1}=R_{-1}\oplus \mathfrak{p}_{E}R_{-2h^{s}}$. Cette suite se prolonge comme pr\'ec\'edemment en une suite $(R^{\natural}_{i})_{i\in {\mathbb Z}}$ dont la p\'eriodicit\'e est $d^{\natural}=2h^{s}-1$. On d\'efinit ensuite les modules $\boldsymbol{\mathfrak{p}}^{E}_{i/d^{\natural}}$ etc.. en utilisant les nouvelles suites de r\'eseaux.

On revient au cas g\'en\'eral o\`u $h^{o}$ est quelconque. 

\begin{lem}{Soit $X$ un \'el\'ement du support de $f_{G}$. Soit $\alpha\in \bar{F}^{\times}$ une valeur propre de $X$, o\`u $X$ est consid\'er\'e comme un \'el\'ement de $End_{E}(V)$. Alors on a $val_{E}(\alpha)\geq \frac{1}{2h^{o}+2h^{s}-1}$.}\end{lem}

Preuve. A conjugaison pr\`es par $G(F)$, on peut supposer que $X$ appartient au support de $\tilde{f}_{G}$.   Alors $X$ appartient \`a $\mathfrak{p}_{1/d}\subset \boldsymbol{\mathfrak{p}}^{E}_{1/d}$. Donc $X^{dm}\in \boldsymbol{\mathfrak{p}}^{E}_{m}=\mathfrak{p}_{E}^m\boldsymbol{\mathfrak{p}}^{E}_{0}$ pour tout entier $m\geq 1$. Cela entra\^{\i}ne que $val_{E}(\alpha)\geq \frac{1}{d}$. Si $h^{o}\not=0$, on a $d=2h^{o}+2h^{s}-1$, d'o\`u  la conclusion de l'\'enonc\'e.

Supposons maintenant $h^{o}=0$. Montrons que $\mathfrak{p}_{1/d}=\mathfrak{p}^{\natural}_{1/d^{\natural}}$.  L'inclusion du membre de droite dans celui de gauche est \'evidente. Soit $Y\in \mathfrak{p}_{1/d}$. On veut prouver que $Y\in \mathfrak{p}^{\natural}_{1/d^{\natural}}$, c'est-\`a-dire que $Y(R^{\natural}_{\geq i})\subset R^{\natural}_{\geq( i+1)}$ pour tout $i=-2h^s+1,...,-1$. Pour tout tel $i$, on a $R^{\natural}_{\geq i}=R_{\geq i}$ et, si $i\not=-1$, $R^{\natural}_{\geq (i+1)}=R_{\geq(i+1)}$. Si $i\not=-1$, l'inclusion \`a d\'emontrer r\'esulte donc de l'hypoth\`ese $Y\in \mathfrak{p}_{1/d}$. Pour $i=-1$, on a $R^{\natural}_{\geq( i+1)}=R_{\geq 1}$. On doit prouver l'inclusion $Y(R_{\geq -1})\subset R_{\geq 1}$ alors que l'hypoth\`ese  
$Y\in \mathfrak{p}_{1/d}$ nous dit seulement que $Y(R_{\geq -1})\subset R_{\geq 0}$ et $Y(R_{\geq 0})\subset R_{\geq 1}$. 
 Or l'espace $R_{\geq -1}/R_{\geq 1}$ de dimension $2$ sur ${\mathbb F}_{q}$ est naturellement muni d'une forme quadratique d\'eploy\'ee. L'\'el\'ement $Y$ se r\'eduit en un \'el\'ement $\bar{Y}$ de l'alg\`ebre de Lie de son groupe sp\'ecial orthogonal. De plus, les relations impos\'ees \`a $Y$ entra\^{\i}nent que $\bar{Y}$ est nilpotent. Or le groupe sp\'ecial orthogonal est un tore, le seul nilpotent dans son alg\`ebre de Lie est nul. Donc $\bar{Y}=0$, c'est-\`a-dire $Y(R_{\geq -1})\subset R_{\geq 1}$, ce qui d\'emontre l'assertion. En particulier $X\in \mathfrak{p}^{\natural}_{1/d^{\natural}}$. Le m\^eme raisonnement que dans le cas $h^{o}\not=0$ montre que $val_{E}(\alpha)\geq \frac{1}{d^{\natural}}$. Or $d^{\natural}=2h^{o}+2h^{s}-1$, d'o\`u  encore la conclusion de l'\'enonc\'e. $\square$
 
 \subsection{Le cas stable}\label{lecasstable}
On conserve les hypoth\`eses du cas pr\'ec\'edent et on suppose maintenant que $G$ est quasi-d\'eploy\'e et que  $h^{s}=h^{o}$ ou $h^{s}=h^{o}-1$. On pose simplement $h=h^{o}$ et $\eta=h-h^{s}$. Notre hypoth\`ese $n\geq2$ entra\^{\i}ne $h\geq h^{s}\geq 1$.   On fixe des suites de r\'eseaux $(R_{i})_{i=-2h+2\eta,...2h-2}$ comme dans le paragraphe pr\'ec\'edent et tous les objets qui s'en d\'eduisent. On identifie $\mathfrak{o}_{F}^{\times}/\mathfrak{o}_{F}^{\times,2}$ \`a un ensemble de repr\'esentants ${\cal C}\subset \mathfrak{o}_{F}^{\times}$.

Soit $m\in \{1,...,h\}$. On pose $d_{m}=4m-1-2\eta$. On fixe un \'el\'ement $\alpha_{m}\in \bar{F}^{\times}$ tel que $\alpha_{m}^{d_{m}}=\varpi_{E}$. On v\'erifie que l'extension $E(\alpha_{m})/F(\alpha_{m}^2)$ est quadratique, plus pr\'ecis\'ement $E(\alpha_{m})$ est la compos\'ee des deux extensions $E$ et $F(\alpha_{m}^2)$ de $F$.  On consid\`ere  $E(\alpha_{m})$ comme un espace sur $E$ et, pour $\gamma\in {\cal C}$, on le munit de la forme hermitienne $q_{m,\gamma}$ (relative \`a l'extension $E/F$) d\'efinie par $q_{m,\gamma}(v,v')=(-1)^{m-1}\gamma d_{m}^{-1}trace_{E(\alpha_{m})/E}(\bar{v}v'\alpha_{m}^{2-2m})$ pour $v,v'\in E(\alpha_{m})$, o\`u $v\mapsto \bar{v}$ est la conjugaison galoisienne associ\'ee \`a l'extension $E(\alpha_{m})/F(\alpha_{m}^2)$. 

Soit $\underline{\gamma}=(\gamma_{m})_{m=1,...,h}\in {\cal C}^h$. On  pose $V_{\underline{\gamma}}=\oplus_{m=1,...,h}E(\alpha_{m})$ et on munit cet espace  de la forme hermitienne $q_{\underline{\gamma}}=\oplus_{m=1,...,h}q_{m,\gamma_{m}}$. On v\'erifie que cet espace hermitien est isomorphe \`a $(V,q)$ si et seulement si $\underline{\gamma}$ appartient \`a ${\cal C}^h(V)$ (c'est l'ensemble $({\mathbb F}_{q}^{\times}/{\mathbb F}_{q}^{\times,2})^{h}(V)$ du paragraphe pr\'ec\'edent). Supposons cette condition v\'erifi\'ee. On voit que l'on peut fixer, et on fixe,  un isomorphisme d'espaces hermitiens $\iota_{\underline{\gamma}}:V_{\underline{\gamma}}\to V$ v\'erifiant les conditions suivantes:

pour $i\in \{1-h,...,h-1\}$, $R_{i}$ est l'image par $\iota_{\underline{\gamma}}$ de la somme sur $m\in \{\vert i\vert +1,..;,h\}$ des modules $\mathfrak{o}_{E}\alpha_{m}^{m-1+i}$;

pour $i\in \{1,... h-\eta\}$, $R_{-h-1+\eta+i}$ est l'image par $\iota_{\underline{\gamma}}$ de la somme sur $m\in \{\eta+i,..;,h\}$ des modules $\mathfrak{o}_{E}\alpha_{m}^{-m-1+\eta+i}$;

pour $i\in \{1,... h-\eta\}$, $R_{-h+\eta-i}$ est l'image par $\iota_{\underline{\gamma}}$ de la somme sur $m\in \{\eta+i,..;,h\}$ des modules $\mathfrak{o}_{E}\alpha_{m}^{-m+\eta-i}$. 

Notons $X'_{\underline{\gamma}}$ l'endomorphisme de $V_{\underline{\gamma}}$ qui agit sur chaque $E(\alpha_{m})$ par multiplication par $\alpha_{m}$, \`a l'exception suivante: si $\eta=1$ et $m=1$, auquel cas $d_{1}=1$, $X'_{\underline{\gamma}}$ agit sur $E_{\alpha_{1}}=E$ par $0$. On note $X_{\underline{\gamma}}$ l'endomorphisme de $V$ qui s'en d\'eduit via $\iota_{\underline{\gamma}}$. On v\'erifie que $X_{\underline{\gamma}}\in \mathfrak{g}_{ell}(F)\cap \mathfrak{k}$. Notons $\bar{X}_{\underline{\gamma}}=\bar{X}_{\underline{\gamma}}^{o}\oplus \bar{X}_{\underline{\gamma}}^{s}$ la r\'eduction de $X_{\underline{\gamma}}$ dans $\mathfrak{k}/\mathfrak{k}^{\perp}=\mathfrak{g}^{o}({\mathbb F}_{q})\oplus \mathfrak{g}^{s}({\mathbb F}_{q})$. Pour $a\in \{o,s\}$, on v\'erifie que $\bar{X}_{\underline{\gamma}}^{a}$ appartient \`a $\tilde{\mathfrak{g}}^{a}_{2}+\mathfrak{g}^{a}_{\geq4}$ et que cet \'el\'ement nilpotent est param\'etr\'e 

dans le cas $a=o$, par la partition $(1,3,...,h)$ et l'\'el\'ement $\underline{\gamma}\in {\cal C}^{h}(V)$;

dans le cas $a=s$, par la partition $(2,4,...,2h-2\eta)$ et l'\'el\'ement

\noindent $((-1)^{1+\eta}\gamma_{1+\eta},(-1)^{1+\eta}\gamma_{2+\eta},...,(-1)^{1+\eta}\gamma_{h})$ de ${\cal C}^{h-\eta}$. On en d\'eduit l'\'egalit\'e

$$(1) \qquad \tilde{f}_{G}(X_{\underline{\gamma}})=\left\lbrace\begin{array}{cc}sgn(-1)^{h(h+1)/2},& \,\,si\,\,\eta=0;\\ sgn((-1)^{[h/2]}\gamma(V)),&\,\,si\,\,\eta=1.\\ \end{array}\right.$$

Cette valeur ne d\'epend pas de $\underline{\gamma}$. 

Il est bien connu que les \'el\'ements $X_{\underline{\gamma}}$ sont tous stablement conjugu\'es, plus pr\'ecis\'ement que la famille $(X_{\underline{\gamma}})_{\underline{\gamma}\in {\cal C}^h(V)}$ est un ensemble de repr\'esentants des classes de conjugaison par ${\bf G}(F)$ dans leur classe de conjugaison stable commune. 
 C'est aussi un ensemble de repr\'esentants des classes de conjugaison par $G(F)$, c'est-\`a-dire que, pour tout $\underline{\gamma}\in {\cal C}^{h}(V)$, la classe de conjugaison de $X_{\underline{\gamma}}$ par $G(F)$ est \'egale \`a sa classe de conjugaison par ${\bf G}(F)$. C'est \'evident si $n$ est impair puisqu'alors $\pi(G(F))=G_{AD}(F)$. Si $n$ est pair, le commutant dans ${\bf G}(F)$ de $X_{\underline{\gamma}}$ contient l'image via $\iota_{\underline{\gamma}}$ du sous-groupe isomorphe \`a $\{\pm 1\}^h$ form\'e des \'el\'ements qui agissent par $\pm 1$ sur chaque $E(\alpha_{m})$. Ce sous-groupe s'envoie surjectivement sur l'image $E^1/(E^1)^n$ de ${\bf G}(F)$ dans $\pi(G(F))\backslash G_{AD}(F)$ et l'assertion en r\'esulte.

\bigskip
  On a donc pour tout $\gamma\in {\cal C}^h(V)$ une \'egalit\'e
$$ (2)\qquad S^G(X_{\underline{\gamma}},\tilde{f}_{G})=\sum_{\underline{\gamma}'\in {\cal C}^h(V)}I^{G}(X_{\underline{\gamma}'},\tilde{f}_{G}).$$
 
  Notons $K_{R}$ le sous-groupe des \'el\'ements $g\in G(F)$ tels que $g(R_{\geq i})\subset R_{\geq i}$ pour tout $i\in {\mathbb Z}$.  L'action de ce groupe   fixe la fonction $\tilde{f}_{G}$.  Montrons que
  
  (3) l'ensemble des $g\in G(F)$ tels que $g^{-1}X_{\underline{\gamma}}g$ appartienne au support de $ \tilde{f}_{G}$ est \'egal \`a $ K_{R}$.
  
  Le temps de la preuve de (3), on identifie $V$ \`a $V_{\underline{\gamma}}$ par l'isomorphisme $\iota_{\underline{\gamma}}$. L'\'el\'ement $\underline{\gamma}$ \'etant fix\'e, on le supprime de la notation pour all\'eger celle-ci. On pose $V^{<h}=\oplus_{m=1,...,h-1}E(\alpha_{m})$. 
  Fixons $g\in G(F)$ tel que $g^{-1}Xg$  appartienne au support de $\tilde{f}_{G}$. Pour $i\in {\mathbb Z}$, posons $S_{i}=gR_{i}$, $S_{\geq i}=gS_{\geq i}$. Ces $\mathfrak{o}_{E}$-modules v\'erifient les m\^emes conditions que les $R_{i}$ ou $R_{\geq i}$. Par construction, le support de $\tilde{f}_{G}$ est contenu dans $\mathfrak{p}_{1/d}$ (o\`u $d=d_{h}$), c'est-\`a-dire que tout \'el\'ement $Y$ de ce support satisfait la relation $Y(R_{\geq i})\subset R_{\geq i+1}$ pour tout $i$. Puisque $g^{-1}Xg$ appartient \`a ce support, on a donc 
  
  (4) $XS_{\geq i}\subset S_{\geq i+1}$ pour tout $i$. 
  
  Nous allons prouver par r\'ecurrence sur $n$ que cette relation entra\^{\i}ne $S_{\geq i}=R_{\geq i}$ pour tout $i$. La relation (4) entra\^{\i}ne $X^dS_{\geq i}\subset S_{\geq i+d}=\mathfrak{p}_{E}S_{\geq i}$. Pour $m=1,...,h$, l'\'el\'ement $\varpi_{E}^{-1}X^d\in \boldsymbol{\mathfrak{g}}(E)$ agit sur $E(\alpha_{m})$ par multiplication par $\alpha_{m}^d\varpi_{E}^{-1}$ (\`a l'exception du cas $\eta=1$, $m=1$ o\`u l'\'el\'ement agit par $0$ sur $E(\alpha_{1})$). On a $val_{E}(\alpha_{m}^d\varpi_{E}^{-1})>0$ si $m<h$ tandis que $val_{E}(\alpha_{h}^d\varpi_{E}^{-1})=0$. Un sous-$\mathfrak{o}_{E}$-r\'eseau de $V$ qui est conserv\'e par un tel \'el\'ement est forc\'ement somme de ses intersections avec les deux sous-espaces $E(\alpha_{h})$ et $V^{<h}$. Ecrivons conform\'ement $S_{\geq i}=S_{\geq i}^h\oplus S_{\geq i}^{<h}$. La relation (4) entra\^{\i}ne d'abord que $S_{\geq i}^h$ est un sous-$\mathfrak{o}_{E(\alpha_{h})}$-module de $E(\alpha_{h})$. Cela implique que $S_{\geq i}^h=\mathfrak{p}_{E(\alpha_{h})}^{j(i)}$ pour un entier $j(i)\in {\mathbb Z}$. L'application $i\mapsto j(i)$ est forc\'ement croissante. L'\'egalit\'e $S_{\geq i+d}^h=\mathfrak{p}_{E}S_{\geq i}^h$ implique que $j(i+d)=j(i)+d$. 
  La relation (4) entra\^{\i}ne aussi $\alpha_{h}S_{\geq i}^h\subset S_{\geq i+1}^h$, c'est-\`a-dire $\mathfrak{p}_{E(\alpha_{h})}^{j(i)+1}\subset \mathfrak{p}_{E(\alpha_{h})}^{j(i+1)}$, d'o\`u $j(i+1)\leq j(i)+1$.  Avec la relation $j(i+d)=j(i)+d$, cela entra\^{\i}ne $j(i+1)=j(i)+1$ pour tout $i$. On se rappelle l'\'egalit\'e $\mathfrak{p}_{E}S_{\geq 0}^{\star}=S_{\geq 2h-1}$, d'o\`u   $\mathfrak{p}_{E}S_{\geq 0}^{h,\star}=S_{\geq 2h-1}^h$. La d\'efinition de la forme $q_{h,\gamma_{h}}$ entra\^{\i}ne que, pour tout $j\in {\mathbb Z}$, $(\mathfrak{p}_{E(\alpha_{h})}^j)^{\star}=\mathfrak{p}_{E(\alpha_{h})}^{-2h+2\eta-j}$. Donc $\mathfrak{p}_{E}S_{\geq 0}^{h,\star}=\mathfrak{p}_{E(\alpha_{h})}^{-2h+2\eta+d-j(0)}$ et l'\'egalit\'e pr\'ec\'edente entra\^{\i}ne $j(2h-1)=-2h+2\eta+d-j(0)$, c'est-\`a-dire $2j(0)=-4h+1+2\eta+d=0$. Or donc $S_{\geq i}^h=\mathfrak{p}_{E(\alpha_{h})}^{i}$ pour tout $i\in {\mathbb Z}$. La suite de r\'eseaux $(S^h_{\geq i})_{i\in {\mathbb Z}}$ est donc uniquement d\'etermin\'ee et on voit que l'on a $S^h_{\geq i}=R_{\geq i}\cap E(\alpha_{h})$ (cela r\'esulte soit de la d\'efinition de la suite de r\'eseaux $(R_{\geq i})_{i\in {\mathbb Z}}$, soit du fait que cette suite v\'erifie (4)). 
  
  On doit maintenant d\'eterminer les r\'eseaux $S_{\geq i}^{<h}$. Supposons d'abord $h\geq 2+\eta$.    Pour $i\in \{-1,0\}\cup \{2h-3,2h-2\}$, il r\'esulte des d\'efinitions que
  $$dim_{{\mathbb F}_{q}}(S_{\geq i}/S_{\geq i+1})=1=dim_{{\mathbb F}_{q}}(S^h_{\geq i}/S^h_{\geq i+1}).$$
  Il en r\'esulte que $S_{\geq i}^{<h}=S_{\geq i+1}^{<h}$. On peut r\'eindicer la suite $(S_{\geq i}^{<h})_{i\in {\mathbb Z}}$ en d\'efinissant
  
  $\bar{S}_{\geq i}^{<h}=S_{\geq i+1}^{<h}$ pour $i\in \{0,...,2h-4\}$, $\bar{S}_{\geq i}^{<h}=S_{\geq i-1}^{h}$ pour $i\in \{-2h+2+2\eta,-1\}$ et $\bar{S}^{<h}_{\geq i+d_{h-1}r}=\mathfrak{p}_{E}^r\bar{S}_{\geq i}^{<h}$ pour $i\in \{-2h+2+2\eta,2h-4\}$ et $r\in {\mathbb Z}$. 
  
  On a encore la relation $X\bar{S}_{\geq i}^{<h}\subset \bar{S}_{\geq i+1}$ pour tout $i$, ce qui est l'analogue de (4).
  Les donn\'ees $V^{<h}$, $\bar{S}_{\geq i}^{<h}$ et $X_{\vert V^{<h}}$ v\'erifient les m\^emes hypoth\`eses que $V$, $S_{\geq i}$ et $X$. On applique l'hypoth\`ese de r\'ecurrence: les r\'eseaux $\bar{S}_{\geq i}^{<h}$ sont uniquement d\'etermin\'es et \'egaux aux analogues pour nos donn\'ees des r\'eseaux $R_{\geq i}$. Ces analogues sont les r\'eseaux $\bar{R}_{\geq i}^{<h}$ d\'eduits de $R_{\geq i}\cap V^{<h}$ de la m\^eme fa\c{c}on que $\bar{S}_{\geq i}^{<h}$ a \'et\'e d\'eduit de $S_{\geq i}^{<h}$. Cela entra\^{\i}ne la relation voulue $S_{\geq i}=R_{\geq i}$ pour tout $i$. Il reste \`a lever l'hypoth\`ese $h\geq 2+\eta$. 
Si $\eta=0$ et $h=1$, on voit que $V^{<1}=\{0\}$ donc $S_{\geq i}^1=\{0\}=R_{\geq i}\cap V^{<1}$. Si $\eta=1$, on a d\'ej\`a exclu le cas $h=1$ (qui impose $n=1$). Si $\eta=1$ et $h=2$, $V^{<2}$ se r\'eduit \`a $E(\alpha_{2})=E$. On a alors $S_{\geq i}^{<2}=\mathfrak{p}_{E}^{j'(i)}$ pour un $j'(i)\in {\mathbb Z}$. Les m\^emes arguments de dimension et de dualit\'e utilis\'es plus haut entra\^{\i}nent que $j'(i+5r)=r$ pour $i\in \{-3,...1\}$ et $r\in {\mathbb Z}$. On obtient encore $S_{\geq i}^{<2}=R_{\geq i}\cap V^{<2}$ et la conclusion. 

Puisque $gR_{\geq i}=S_{\geq i}=R_{\geq i}$ pour tout $i$, on a $g \in K_{R}$. Cela prouve (3).

\begin{lem}{On a $S^G(X_{\underline{\gamma}},f_{G})\not=0$ pour tout $\underline{\gamma}\in {\cal C}^h(V)$.} \end{lem}

Preuve. Puisque $K_{R}$ fixe $\tilde{f}_{G}$, il r\'esulte de (3) que $I^{G}(X_{\underline{\gamma}},\tilde{f}_{G})=c_{\underline{\gamma}}\tilde{f}_{G}(X_{\underline{\gamma}})$ o\`u $c_{\underline{\gamma}}>0$ (on v\'erifie qu'en fait, cette constante ne d\'epend pas de $\underline{\gamma}$ mais peu nous importe). On d\'eduit de (1) et (2)  que $S^G(X_{\underline{\gamma}},\tilde{f}_{G})\not=0$ pour tout $\underline{\gamma}\in {\cal C}^h(V)$. Mais, par construction, $S^G(X_{\underline{\gamma}},f_{G})=cS^G(X_{\underline{\gamma}},\tilde{f}_{G})$ o\`u $c$ est un entier strictement positif. $\square$

La propri\'et\'e suivante r\'esulte de la construction:

(5) pour tout $\underline{\gamma}\in {\cal C}^h(V)$, $X_{\underline{\gamma}}$ poss\`ede une valeur propre $\alpha\in \bar{F}^{\times}$ telle que $val_{E}(\alpha)=\frac{1}{2h^{o}+2h^s-1}$.

\subsection{Preuve des assertions (2) et (3) de \ref{An-1quasidepram}}\label{preuve}
On reprend les notations de ce paragraphe. Soit $y\in {\cal Y}$, notons $(i,j)$ son image dans ${\mathbb Y}$ et posons $(k,h)=\phi^{-1}(i,j)$. Rappelons que, par d\'efinition de l'ensemble ${\mathbb Y}$, on a $i\leq j$. On a associ\'e \`a $y$ une donn\'ee endoscopique ${\bf G}'_{y}$. Celle-ci v\'erifie $G'_{y,SC}(F)= SU_{E/F}(i(i+1)/2,F)\times  SU_{E/F}(j(j+1)/2,F)$, sauf dans le cas o\`u $i=j$, auquel cas   $G'_{y,SC}(F)$ peut \^etre soit le groupe pr\'ec\'edent, soit $SU_{Q/E_{0}}(n/2,E_{0})$ (o\`u $E_{0}/F$ est l'extension quadratique non ramifi\'ee).    Notons $G''$ l'un des facteurs de ces groupes. On applique \`a ce groupe $G''$ les constructions du paragraphe pr\'ec\'edent, quitte pour le dernier groupe \`a remplacer le corps de base $F$ par $E_{0}$. Elles nous fournissent un \'el\'ement $f_{G''}\in FC(\mathfrak{g}''(F))$ et divers \'el\'ements not\'es $X_{\underline{\gamma}}$ dans ce paragraphe.  On fixe un tel \'el\'ement $X''\in \mathfrak{g}_{ell}''(F)$. Le lemme \ref{lecasstable} affirme que  $S^{G''}(X'',f_{G''})\not=0$. Quand le groupe $G'_{y,SC}(F)$ n'a qu'un seul facteur, on note $f'_{y}$ la fonction $f_{G''}$ et $Y_{y}=X''$. 
  Quand il  a deux facteurs, on note $f'_{y}$ le produit tensoriel des deux fonctions $f_{G''}$. On a envie de d\'efinir  $Y_{y}$ comme la somme des deux \'el\'ements $X''$ mais l'\'el\'ement   ainsi d\'efini n'est pas forc\'ement $G$-r\'egulier car une m\^eme valeur propre peut appara\^{\i}tre dans chacun des deux facteurs. On utilise l'argument de \ref{ingredients}: on remplace les \'el\'ements $X''$ par des \'el\'ements suffisamment voisins de sorte que ces derniers v\'erifient les m\^emes propri\'et\'es que les \'el\'ements initiaux (c'est-\`a-dire le lemme et la relation (5) du paragraphe pr\'ec\'edent) mais que leur somme soit $G$-r\'eguli\`ere. On note alors $Y_{y}$ leur somme.

  On peut pr\'eciser une propri\'et\'e de cet \'el\'ement $Y_{y}$. Supposons $G'_{y,SC}(F)= SU_{E/F}(i(i+1)/2,F)\times  SU_{E/F}(j(j+1)/2,F)$ et consid\'erons le groupe $G''=SU_{E/F}(j(j+1)/2)$. D'apr\`es \ref{lecasstable}(5), l'\'el\'ement $X''$ poss\`ede une  valeur propre $\alpha\in \bar{F}^{\times}$ telle que $val_{E}(\alpha)=\frac{1}{2j-1}$. 
 On v\'erifie \`a l'aide de la d\'efinition de $\phi$ que $2j-1=2h+2k-1$.  On en d\'eduit

(1) pour tout \'el\'ement $X\in \mathfrak{g}(F)$ correspondant \`a $Y_{y}$, $X$ poss\`ede une valeur propre $\alpha$ telle que $val_{E}(\alpha)=\frac{1}{2h+2k-1}$.

Le m\^eme r\'esultat vaut dans le cas o\`u $G'_{y,SC}(F)=SU_{Q/E_{0}}(n/2,E_{0})$. Le fait que l'on doit changer de corps de base ne change pas la valuation puisque $E_{0}/F$ est non ramifi\'ee.

D\'emontrons \ref{An-1quasidepram} (2), c'est-\`a-dire

(2) pour tout \'el\'ement non  nul $f'\in FC^{st}(\mathfrak{g}'_{y}(F))^{Out({\bf G}_{y}')}$, on a $S^{G'_{y}}(Y_{y},f')\not=0$.

Les groupes $G''$ consid\'er\'es ci-dessus   sont
  du type $A_{n''-1}$ quasi-d\'eploy\'e, relatif \`a une extension ramifi\'ee. Utilisons pour ces groupes les notations de \ref{An-1quasidepram} en y ajoutant des $''$. On a vu en \ref{descriptionexplicite} que la fonction $f_{G''}$ appartenait \`a la droite $FC''_{x''}$ pour l'unique \'el\'ement $x''\in {\cal X}^{'',st}$. Excluons d'abord le cas $i=0$. Alors $n''<n$ et on peut appliquer  par r\'ecurrence l'assertion (4) de \ref{resultats}: $FC_{x''}=FC^{st}(\mathfrak{g}''(F))$. Il en r\'esulte que $f'_{y}\in FC^{st}(\mathfrak{g}'_{y}(F))$ et on a vu que cette droite \'etait fixe par $Out({\bf G}'_{y})$. L'in\'egalit\'e $S^{G'_{y}}(Y_{y},f'_{y})\not=0$ entra\^{\i}ne alors (2). 
 Consid\'erons maintenant le cas $i=0$. La donn\'ee ${\bf G}'_{y}$ est la donn\'ee principale ${\bf G}$ et $(k,h)$ est l'unique \'el\'ement de ${\cal X}^{st}$. La fonction $f'_{y}$  est alors la fonction $f_{k,h}$ mais on n'a pas encore d\'emontr\'e qu'elle appartenait \`a $FC^{st}(\mathfrak{g}(F))$. On a toutefois d\'ej\`a prouv\'e en  \ref{An-1quasidepram} que cet espace \'etait de dimension $1$. Notons $f_{{\bf G}}$ un g\'en\'erateur de cet espace. On peut \'ecrire $f'_{y}=cf_{{\bf G}}+f'$, o\`u $c\in {\mathbb C}$ et $f'\in \sum_{{\bf G}'\in Endo_{ell}(G), {\bf G}'\not={\bf G}}FC(\mathfrak{g}(F) ,{\bf G}')$.  L'\'el\'ement $Y_{y}$ est elliptique. Par d\'efinition des espaces $ I_{cusp}(\mathfrak{g}(F),{\bf G}')$, on a  donc $S^G(Y_{y},f')=0$. Donc $S^G(Y_{y},f'_{y})=cS^G(Y_{y},f_{{\bf G}})$ et on conclut $S^G(Y_{y},f_{{\bf G}})\not=0$.  Cela d\'emontre (2). 

  D\'emontrons maintenant \ref{An-1quasidepram} (3). Soit $x'\in {\cal X}$, notons $(k',h')$ son  image dans ${\mathbb X}$.  La droite $FC_{x'}$ est port\'ee par la fonction not\'ee $f_{G}$ en \ref{descriptionexplicite} associ\'ee au couple $(h^s,h^{o})=(k',h')$, notons-la $f_{x'}$. Soit $X\in \mathfrak{g}_{reg}(F)$ un \'el\'ement correspondant \`a $Y_{y}$. Supposons $I^G(X,f_{x'})\not=0$. D'apr\`es (1), $X$ poss\`ede une valeur propre $\alpha$ telle que $val_{E}(\alpha)=\frac{1}{2h+2k-1}$. Le lemme \ref{descriptionexplicite} nous dit que $val_{E}(\alpha)\geq \frac{1}{2h'+2k'-1}$. Il en r\'esulte que $k'+h'\geq k+h$. Par d\'efinition, cela signifie que $x'\geq \varphi^{-1}(y)$. Cela d\'emontre  \ref{An-1quasidepram} (3).
  
  \subsection{Action d'un automorphisme}\label{actiondunautomorphisme}
On reprend les hypoth\`eses de \ref{lecasstable}. Pour exhiber un \'epinglage du  groupe $G$, posons $\alpha=1$ si $n$ est impair et $\alpha=\varpi_{E}^{-1}$ si $n$ est pair. Quitte \`a multiplier $q$ par un \'el\'ement de $\mathfrak{o}_{F}^{\times}$, on peut supposer qu'il existe une base $(e_{i})_{i=1,...,n}$ de $V$ v\'erifiant les \'egalit\'es suivantes

$q(e_{i},e_{n+1-i})=(-1)^{i}\alpha$ pour tout $i=1,...,n$,

$q(e_{i},e_{j})=0$ pour tous $i,j=1,...,n$ tels que $i+j\not=n+1$. 

On identifie $GL(V\otimes_{E}\bar{F})$ \`a $GL(n,\bar{F})$ gr\^ace \`a cette base. Notons $B$ le sous-groupe de Borel triangulaire sup\'erieur de $GL(n)$ et $T$ le sous-tore diagonal. Introduisons l'\'epinglage habituel form\'e des matrices $(E_{i})_{i=1,...,n-1}$, o\`u $E_{i}$ est l'\'el\'ement de $\mathfrak{gl}(n)$ dont tous les coefficients sont nuls sauf le $(i,i+1)$-i\`eme qui vaut $1$. Notons $\theta$ l'automorphisme de $GL(n)$ qui pr\'eserve $B$ et $T$ et envoie $E_{i}$ sur $E_{n-i}$ pour tout $i=1,...,n-1$. On note $\sigma\mapsto \sigma_{GL(n)}$ l'action galoisienne habituelle (action sur les coefficients des matrices) et on introduit l'action galoisienne $\sigma\mapsto \sigma_{G}$ d\'efinie par $\sigma_{G}=\sigma_{GL(n)}$ si $\sigma\in \Gamma_{E}$ et $\sigma_{G}=\theta\circ \sigma_{GL(n)}$ si $\sigma\in \Gamma_{F}-\Gamma_{E}$. On v\'erifie qu'alors $G$ est \'egal \`a $GL(n)$ muni de l'action galoisienne $\sigma\mapsto \sigma_{G}$. L'automorphisme $\theta$ se restreint en un automorphisme de $G$ d\'efini sur $F$. Remarquons que, modulo le plongement $G(F)\subset G(E)=GL(n,E)$, l'action de $\theta$ sur $G(F)$ s'identifie \`a l'action $g\mapsto \tau_{GL(n)}(g)$, o\`u  $\tau$ est l'element non trivial de $\Gamma_{E/F}$. 

  Posons $n^{o}=h^2$, $n^{s}=(h-\eta)(h+1-\eta)$. Supposons d'abord $n$ impair donc $\alpha=1$.On peut choisir notre suite de r\'eseaux $(R_{i})_{i\in {\mathbb Z}}$ de sorte que $R_{\geq 0}$ soit engendr\'e sur $\mathfrak{o}_{E}$ par $e_{1},... e_{n^{s}/2+n^{o}}, \varpi_{E}e_{n^{s}/2+n^{o}+1},...,\varpi_{E}e_{n}$. Alors $R_{\geq 0}^*$ est engendr\'e sur $\mathfrak{o}_{E}$ par $\varpi_{E}^{-1}e_{1},...,\varpi_{E}^{-1}e_{n^{s}/2},$

\noindent$e_{n^{s}/2+1},...,e_{n}$. Posons $I^{o}=\{n^{s}/2+1,...,n^{s}/2+n^{o}\}$, $I^{s}_{+}=\{1,...,n^{s}/2\}$, $I^{s}_{-}=\{n^{s}/2+n^{o}+1,...,n\}$, $I^{s}=I^{s}_{+}\cup I^{s}_{-}$. Pour $i\in I^{o}$, resp. $i\in I^{s}_{+}$, $i\in I^{s}_{-}$, notons $\underline{e_{i}}$  la r\'eduction dans $V^{o}=R_{\geq 0}/\mathfrak{p}_{E}R_{\geq 0}^*$, resp $V^{s}=R_{\geq 0}^*/R_{\geq 0}$, de $e_{i}$, resp. $\varpi_{E}^{-1}e_{i}$, resp. $e_{i}$.  
L'espace $V^{o}$ est muni de la base $(\underline{e}_{i})_{i\in I^{o}}$ et l'espace $V^{s}$ est muni de la base $(\underline{e}_{i})_{i\in I^{s}}$. Pour $i,j\in I^{s}$, posons $a_{i,j}=0$ si $i,j\in I^{s}_{+}$ ou $i,j\in I^{s}_{-}$, $a_{i,j}=1$ si $i\in I^{s}_{+}$ et $j\in I^{s}_{-}$ et $a_{i,j}=-1$ si $i\in I^{s}_{-}$ et $j\in I^{s}_{+}$. 
Pour $X=(x_{i,j})_{i,j=1,...,n}\in \boldsymbol{\mathfrak{g}}(E)$, notons $X^{o}$ la matrice extraite $(x_{i,j})_{i,j\in I^{o}}$ et $X^{s}$ la matrice   $(\varpi_{E}^{a_{i,j}}x_{i,j})_{i,j\in I^{s}}$. Si $X\in \boldsymbol{\mathfrak{k}}^{E}$, on v\'erifie que ces matrices $X^{o}$ et $X^{s}$ sont \`a coefficients dans $\mathfrak{o}_{E}$ et que les images de $X$ dans $\mathfrak{g}^{o}({\mathbb F}_{q})$, resp. $\mathfrak{g}^{s}({\mathbb F}_{q})$, sont les r\'eductions naturelles de $X^{o}$, resp. $X^{s}$, c'est-\`a-dire que l'on envoie chaque coefficient sur sa r\'eduction dans ${\mathbb F}_{q}$.   A cause des termes $\varpi_{E}^{a_{i,j}}$ dans les formules ci-dessus et parce que $\tau(\varpi_{E})=-\varpi_{E}$, cette description montre que l'action de $\theta$ pr\'eserve $\boldsymbol{\mathfrak{k}}^{E}$ et se r\'eduit en l'identit\'e de $\mathfrak{g}^{o}({\mathbb F}_{q})$ et en l'action sur $\mathfrak{g}^{s}({\mathbb F}_{q})$ de la similitude symplectique qui agit par $-1$ sur les $\underline{e}_{i}$ pour $i\in I^{s}_{+}$ et par $1$ sur les $\underline{e}_{i}$ pour $i\in I^{s}_{-}$. On v\'erifie que cette similitude multiplie la fonction $f^{s}$ par $sgn(-1)^{n^{s}/2}$. Posons $j=2h-\eta$. On a $n= j(j+1)/2$.  On a suppos\'e $n$ impair. Puisque $n=h^2+n^{s}$ et que $n^s$ est pair,   $h$ est impair. On voit alors que $n^{s}/2$ est de la m\^eme parit\'e que $[(j+2)/4]$. D'o\`u

(1) $\theta(f_{G})=sgn(-1)^{[(j+2)/4]}f_{G}$.

Supposons maintenant $n$ pair donc $\alpha=\varpi_{E}^{-1}$. On peut choisir notre suite de r\'eseaux $(R_{i})_{i\in {\mathbb Z}}$ de sorte que $R_{\geq 0}$ soit engendr\'e sur $\mathfrak{o}_{E}$ par $\varpi_{E}e_{1},...,\varpi_{E}e_{n/2},e_{n/2+1},... e_{n^{s}/2+n^{o}}, $

\noindent$\varpi_{E}e_{n^{s}/2+n^{o}+1},...,\varpi_{E}e_{n}$. Alors $R_{\geq 0}^*$ est engendr\'e sur $\mathfrak{o}_{E}$ par $e_{1},...,e_{n^{s}/2},\varpi_{E}e_{n^{s}/2+1},...,\varpi_{E}e_{n/2},$

\noindent$e_{n/2+1},...,e_{n}$. Le calcul se poursuit comme ci-dessus, les r\^oles de $\mathfrak{g}^{o}$ et $\mathfrak{g}^{s}$ \'etant invers\'es. L'action de $\theta$ pr\'eserve $\boldsymbol{\mathfrak{k}}^{E}$ et se r\'eduit en l'identit\'e de $\mathfrak{g}^{s}({\mathbb F}_{q})$ et en l'action sur $\mathfrak{g}^{o}({\mathbb F}_{q})$ de la similitude de rapport $-1$. On v\'erifie que cette similitude multiplie la fonction $f^{o}$ par $sgn(-1)^{[h(h+1)/2]}$. Cette fois, $h$ est pair et $h(h+1)/2$ est  encore de m\^eme parit\'e que $[(j+2)/4]$. D'o\`u encore (1).

  \section{Les groupes (quasi)-classiques}

\subsection{Type $B_{n}$ d\'eploy\'e}\label{Bndeppadique}
On suppose $G$ d\'eploy\'e de type $B_{n}$ avec $n\geq 2$. C'est-\`a-dire que $G=Spin_{dep}(2n+1)$ est la forme d\'eploy\'ee du groupe spinoriel et on a $G_{AD}=SO_{dep}(2n+1)$. La  norme spinorielle se quotiente en un isomorphisme $\pi(G(F))\backslash G_{AD}(F)\to F^{\times}/F^{\times,2}$.  L'image de $G_{AD}(F)_{0}$ dans ce groupe est le sous-groupe $\mathfrak{o}_{F}^{\times}/\mathfrak{o}_{F}^{\times,2}$. On note $\Xi^{nr}$, resp. $\Xi^{ram}$, l'ensemble des caract\`eres de $F^{\times}/F^{\times,2}$ qui sont triviaux, resp. non triviaux,  sur $\mathfrak{o}_{F}^{\times}/\mathfrak{o}_{F}^{\times,2}$.

Notons ${\mathbb X}^{nr}$ l'ensemble des couples $(k,h)\in {\mathbb N}^2$ tels que $k^2+h^2=2n+1$, $k$ est pair et $h$ est impair. S'il n'existe pas de tel couple avec $k=0$, on pose ${\cal X}^{nr}={\mathbb X}^{nr}$. Si un tel couple $(0,h_{0})$ existe, on note ${\cal X}^{nr}$  la r\'eunion de ${\mathbb X}^{nr}-\{(0,h_{0})\}$ et de $\{(0,h_{0},\xi); \xi\in \Xi^{nr}\}$. Remarquons que le couple $(0,h_{0})$ n'existe que si $n$ est pair. Notons ${\mathbb X}^{ram}$ l'ensemble des couples $(k,h)\in {\mathbb N}^2$ tels que $k(k+1)/2+h(h+1)/2=2n+1$, $k\geq h$ et, si $\delta_{4}(q-1)=0$, l'un des termes $k(k+1)/2$ ou $h(h+1)/2$ est divisible par $4$ (remarquons que l'\'egalit\'e v\'erifi\'ee par $k$ et $h$ entra\^{\i}ne en tout cas qu'un et un seul des deux termes pr\'ec\'edents est pair). On note ${\cal X}^{ram}$ l'ensemble des triplets $(k,h,\xi)$ avec $(k,h)\in {\mathbb X}^{ram}$ et $\xi\in \Xi^{ram}$. On pose ${\cal X}={\cal X}^{nr}\sqcup {\cal X}^{ram}$ et $d_{x}=1$ pour tout $x\in {\cal X}$. 

Pour tout groupe $Spin(N)$, avec $N\geq3$,  on note $z$ l'el\'ement non trivial du noyau de la projection naturelle $Spin(N)\to SO(N)$. En d\'ecrivant l'immeuble \`a l'aide de l'alg\`ebre lin\'eaire ou en utilisant les tables de Tits, on obtient la description suivante. 
L'ensemble $\underline{S}(G)$ s'envoie surjectivement sur l'ensemble des couples $(a,b)\in {\mathbb N}^2$ tels que $a+b=n$ et $a\not=1$. Les fibres de cette surjection ont un seul \'el\'ement sauf celle au-dessus de $(0,n)$ qui en a deux. L'action de $G_{AD}(F)$ pr\'eserve les fibres et permute les deux \'el\'ements de la fibre au-dessus de $(0,n)$. Pour $s\in \underline{S}(G)$ param\'etr\'e par le couple $(a,b)$, on a $G_{s}=Spin(2n+1)$ si $(a,b)=(0,n)$, $G_{s}=Spin_{dep}(2n)$ si $(a,b)=(n,0)$, $G_{s}=(Spin_{dep}(2a)\times Spin(2b+1))/\{1,(z,z)\}$ si $ab\not=0$. Consid\'erons un sommet $s$ param\'etr\'e par un couple $(a,b)$ avec $ab\not=0$.  Consid\'erons deux fonctions  $f_{N^{a},\epsilon^{a}}\in fc(\mathfrak{Spin}_{dep}(2a,{\mathbb F}_{q}))$ et  $f_{N^{b},\epsilon^{b}} \in fc(\mathfrak{Spin}(2b+1,{\mathbb F}_{q}))$. Parce que $G_{s}$ est le quotient par $\{1, (z,z)\}$ de $G_{s,SC}$, le produit tensoriel des deux  fonctions est un \'el\'ement de $FC(\mathfrak{g}_{s}({\mathbb F}_{q}))$ si et seulement si $\epsilon^{a}(z)=\epsilon^{b}(z)$. On suppose cette condition v\'erifi\'ee. 
L'image dans  $G_{s}$ de l'\'el\'ement $z$ de $G$ est l'image dans $G_{s}$ de $(1,z)$ o\`u, ici, $z$ est l'\'el\'ement de $Spin(2b+1)$. L'action de $G_{AD}(F)_{0}$ sur $f_{N^{a},\epsilon^{a}}\otimes  f_{N^{b},\epsilon^{b}}$ est  donc triviale si $\epsilon^{b}(z)=1$, non triviale sinon.   L'action du groupe $G_{AD}(F)$ tout entier est d\'etermin\'ee par l'action suppl\'ementaire de l'automorphisme $\theta$ de $Spin_{dep}(2a)$. On utilise alors les descriptions de \ref{Bn}, \ref{Dndeppair}  et \ref{Dndepimp}.  Il y a une fonction $f_{N^{a},\epsilon^{a}}\otimes  f_{N^{b},\epsilon^{b}}$ avec $\epsilon^{a}(z)=\epsilon^{b}(z)=1$ si et seulement si $(2a,2b+1)$ est de la forme $(k^2,h^2)$. Si cette condition est v\'erifi\'ee,  le couple $(k,h)$ appartient \`a ${\cal X}^{nr}$. La fonction $f_{N^{a},\epsilon^{a}}$ est  invariante par $\theta$. D'apr\`es \ref{actionsurFC}, 
$f_{N^{a},\epsilon^{a}}\otimes  f_{N^{b},\epsilon^{b}}$   donne naissance \`a un \'el\'ement de $FC(\mathfrak{g}(F))$  et on note $FC_{k,h}$ la droite port\'ee par cette fonction.    Il y a des fonctions $f_{N^{a},\epsilon^{a}}\otimes  f_{N^{b},\epsilon^{b}}$ avec $\epsilon^{a}(z)=\epsilon^{b}(z)=-1$ si et seulement si $(2a,2b+1)$ est de la forme $(k'(k'+1)/2,h'(h'+1)/2)$ et, de plus, $a$ est pair ou $\delta_{4}(q-1)=0$. Supposons ces conditions v\'erifi\'ees. Il y a alors deux fonctions $f_{N^{a},\epsilon^{a}}$ qui, en les normalisant correctement, sont permut\'ees par $\theta$.  On obtient deux fonctions $f_{N^{a},\epsilon^{a}}\otimes  f_{N^{b},\epsilon^{b}}$ sur lesquelles $G_{AD}(F)_{0}$ agit par son caract\`ere non trivial et qui sont permut\'ees par $\theta$. Pour chaque caract\`ere $\xi\in \Xi^{ram}$, la somme ou la diff\'erence de ces deux fonctions se transforme selon le caract\`ere $\xi$ de $G_{AD}(F)/\pi(G(F))$. D'apr\`es \ref{actionsurFC}, il se d\'eduit  pour tout $\xi\in \Xi^{ram}$ un \'el\'ement de $FC(\mathfrak{g}(F))$ et on note  $FC_{k,h,\xi}$ la droite qu'elle porte,
 o\`u $(k,h)$ est le seul des couples $(k',h')$ ou $(h',k')$ tel que $k\geq h$.   Remarquons que la condition $a$ est pair ou $\delta_{4}(q-1)=0$ que l'on a impos\'ee \'equivant \`a: si $\delta_{4}(q-1)=0$, $k(k+1)/2$ ou $h(h+1)/2$ est divisible par $4$. Alors, le triplet $(k,h,\xi)$ appartient \`a ${\cal X}^{ram}$. 
Consid\'erons maintenant le sommet $s$ param\'etr\'e par le couple $(n,0)$. On a $G_{s}=Spin_{dep}(2n)$. La description est la m\^eme que dans le cas pr\'ec\'edent, en oubliant le facteur $Spin(2b+1)$, apr\`es avoir remarqu\'e que l'image dans $G_{s}$ de l'\'el\'ement $z$ de $G$ est l'\'el\'ement $z$ de $Spin_{dep}(2n)$. Cela se v\'erifie en explicitant ces \'el\'ements gr\^ace \`a l'ensemble de racines  affines $\Delta_{a}$ de $G$: le $z$ de $G$ est $\check{\alpha}_{n}(-1)$, celui de $Spin_{dep}(2n)$ est l'image naturelle de $\check{\alpha}_{0}(-1)\check{\alpha}_{1}(-1)$ mais ces deux \'el\'ements sont en fait \'egaux car $\check{\alpha}_{0}\check{\alpha}_{1}\check{\alpha}_{n}\prod_{i=2,...,n-1}\check{\alpha}_{i}^2=1$. 
Consid\'erons enfin un sommet $s$ param\'etr\'e par le couple $(0,n)$. Le stabilisateur de ce sommet dans $G_{AD}(F)$ est le sous-groupe $G_{AD}(F)_{0}$. Tout \'el\'ement de $FC(\mathfrak{g}_{s}({\mathbb F}_{q}))$ qui se transforme selon
un caract\`ere $\xi_{0}$ de $G_{AD}(F)_{0}$ donne naissance \`a deux \'el\'ements de $FC(\mathfrak{g}(F))$ qui se transforment selon les deux \'el\'ements de $\Xi$ dont la restriction \`a $G_{AD}(F)_{0}$ est $\xi_{0}$. Il  y a une fonction $f\in 
 FC(\mathfrak{g}_{s}({\mathbb F}_{q}))$ se transformant selon le caract\`ere $\xi_{0}=1$ si et seulement si $2n+1$ est de la forme $h_{0}^2$. Dans ce cas, on note naturellement $FC_{0,h_{0},\xi}$ les droites port\'ees par ces fonctions,  avec $\xi\in \Xi^{nr}$. On a $(0,h_{0},\xi)\in {\cal X}^{nr}$. Il y a une fonction $ f\in 
 FC(\mathfrak{g}_{s}({\mathbb F}_{q}))$ se transformant selon l'unique caract\`ere $\xi_{0}\not=1$ si et seulement si $2n+1$ est de la forme $k(k+1)/2$.  Dans ce cas, les deux droites d\'eduites se notent naturellement $FC_{k,0,\xi}$ pour $\xi\in \Xi^{ram}$. On a $(k,0,\xi)\in {\cal X}^{ram}$. 
 On voit que cette description d\'emontre \ref{resultats} (1).

Notons ${\mathbb Y}^{nr}$ l'ensemble des couples $(i,j)\in {\mathbb N}^2$ tels que $i(i+1)+j(j+1)=n$ et $i\leq j$. S'il n'existe pas de tel couple avec $i=j$, on pose ${\cal Y}^{nr}={\mathbb Y}^{nr}$. Si un tel couple $(i_{0},i_{0})$ existe, on note ${\cal Y}^{nr}$ la r\'eunion de ${\mathbb Y}^{nr}-\{(i_{0},i_{0})\}$ et de $\{(i_{0},i_{0},\xi); \xi\in \Xi^{nr}\}$. Remarquons que ce couple $(i_{0},i_{0})$ n'existe que si $n$ est pair. Notons ${\mathbb Y}^{ram}$ l'ensemble des couples $(i,j)$ avec $i,j\in {\mathbb N}$, v\'erifiant les conditions suivantes:

$2i(i+1)+j(j+1)/2=n$;

 si $\delta_{4}(q-1)=0$, $[(i+1)/2]+[(j+2)/4]$ est pair.
 
 On note ${\cal Y}^{ram}$ l'ensemble des triplets $(i,j,\xi)$ pour $(i,j)\in {\mathbb Y}^{ram}$ et $\xi\in \Xi^{ram}$. 
  On pose ${\cal Y}={\cal Y}^{nr}\sqcup {\cal Y}^{ram}$.
  
  D\'eterminons les donn\'ees endoscopiques elliptiques ${\bf G}'$ de $G$ telles que $FC^{st}(\mathfrak{g}'(F))^{Out({\bf G}')}\not=\{0\}$. On consid\`ere un couple $(\sigma\mapsto \sigma_{G'},{\cal O})\in {\cal E}_{ell}(G)$. Rappelons   que $\hat{D}_{a}$ est le diagramme de Dynkin compl\'et\'e  du groupe dual de $G$, c'est-\`a-dire d'un groupe de type $C_{n}$. Le groupe $\hat{\Omega}$ co\"{\i}ncide avec le groupe d'automorphismes de ce diagramme, lequel est $\{1,\omega\}$, o\`u $\omega$ permute $\hat{\alpha}_{i}$ et $\hat{\alpha}_{n-i}$ pour $i=0,...,n$.  
  
Supposons d'abord que l'action $\sigma\mapsto \sigma_{G'}$ soit triviale. L'orbite ${\cal O}$ est r\'eduite \`a une seule racine $\hat{\alpha}_{m}$.  Deux donn\'ees \'etant \'equivalentes si et seulement si elles se d\'eduisent l'une de l'autre par l'action de $\omega$, on peut supposer $m\leq n/2$.  Le diagramme $\hat{D}_{a}-{\cal O}$ est le produit de deux diagrammes de type $C_{m}$ et $C_{n-m}$.  La donn\'ee n'a pas d'automorphisme non trivial  si $m\not=n/2$ et elle en a un, l'action de $\omega$, si $m=n/2$. Dualement, on obtient que $G'$ est semi-simple et $G'_{SC}=Spin(2m+1)\times Spin(2(n-m)+1)$. L'action galoisienne \'etant triviale, il est clair que $\xi_{{\bf G}'}=1$.   Si $m\not=0$, on peut appliquer par r\'ecurrence l'assertion  (4) ci-dessous: on a $FC^{st}(\mathfrak{g}'(F))\not=\{0\}$ si et seulement si   $(m,n-m)$ est de la forme $(i(i+1),j(j+1))$. Dans  ce cas, $(i,j)$ appartient \`a ${\cal Y}^{nr}$, l'espace $FC^{st}(\mathfrak{g}'(F))$ est  une droite.  Si de plus $m<n/2$, il n'y a pas d'automorphisme non trivial. On pose ${\bf G}'={\bf G}'_{i,j}$ et  $FC^{{\cal E}}_{i,j}=FC^{st}(\mathfrak{g}_{i,j}'(F))^{Out({\bf G}'_{i,j})}$. 
  Dans le cas o\`u $m=n/2$, on a $i=j$ et l'action de l'automorphisme \'echange les deux facteurs $Spin(n+1)$. On v\'erifie qu'elle agit trivialement sur $FC^{st}(\mathfrak{g}'(F))$.  En notant ${\bf 1}$  l'\'el\'ement neutre de $\Xi$, le triplet $(i,i,{\bf 1})$ appartient \`a ${\cal Y}^{nr}$. On pose ${\bf G}'={\bf G}'_{i,i,{\bf 1}}$ et  $FC^{{\cal E}}_{i,i,{\bf 1}}=FC^{st}(\mathfrak{g}_{i,i,{\bf 1}}'(F))^{Out({\bf G}_{i,i,{\bf 1}}')}$. 
   Enfin, si $m=0$, on a ${\bf G}'={\bf G}$ et  on ne peut encore rien dire  de l'espace $FC^{st}(\mathfrak{g}(F))$.

Supposons maintenant que l'action $\sigma\mapsto \sigma_{G'}$ soit non triviale. L'extension $E_{G'}/F$ est quadratique et on a  $\sigma_{G'}=1$ pour $\sigma\in \Gamma_{E_{G'}}$ et $\sigma_{G'}=\omega$ pour $\sigma\in \Gamma_{F}-\Gamma_{E_{G'}}$. L'orbite ${\cal O}$ est de la forme $\{\hat{\alpha}_{m},\hat{\alpha}_{n-m}\}$ pour un $m<n/2$ ou $\{\hat{\alpha}_{n/2}\}$ dans le cas o\`u $n$ est pair. Traitons d'abord le premier cas. Le  stabilisateur de $\hat{\alpha}_{m}$ dans $\Gamma_{F}$ est $\Gamma_{E_{G'}}$. Le lemme \ref{centre} exclut le cas o\`u $E_{G'}/F$ est non ramifi\'ee. Supposons $E_{G'}/F$ ramifi\'ee.  Le diagramme $\hat{D}_{a}-{\cal O}$ est r\'eunion de deux diagrammes de type $C_{m}$ et d'un diagramme de type $A_{n-2m-1}$ (ce dernier disparaissant si $n=2m+1$).  Pour $\sigma\in \Gamma_{F}-\Gamma_{E_{G'}}$, $\sigma_{G'}$ \'echange les deux premiers diagrammes et agit sur celui de type $A_{n-2m-1}$ par l'automorphisme non trivial. Dualement, on a $G'_{SC}=Res_{E_{G'}/F} (Spin(2m+1)_{\vert E_{G'}})\times SU_{E_{G'}/F}(n-2m)$, o\`u $Spin(2m+1)_{\vert E_{G'}} $ est le groupe $Spin(2m+1)$ d\'efini sur $E_{G'}$, avec la convention $SU_{E_{G'}/F}(1)=\{1\}$. On utilise par r\'ecurrence  les r\'esultats du pr\'esent paragraphe et de \ref{An-1quasidepram} pour les deux groupes $Spin(2m+1)$ et $SU_{E_{G'}/F}(n-2m)$: on a $FC^{st}(\mathfrak{g}'(F))\not=\{0\}$ si et seulement si  $m$ est de la forme $i(i+1)$ et $n-2m$ est de la forme $j(j+1)/2$ avec $i,j\in {\mathbb N}$   (le cas $j=0$ est exclu puisque l'on a suppos\'e $m<n/2$). Si ces conditions sont v\'erifi\'ees, l'espace $FC^{st}(\mathfrak{g}'(F))$ est  une droite. La donn\'ee a un automorphisme non trivial. Cet automorphisme agit sur chacun des facteurs de $G'_{SC}$  par l'automorphisme galoisien associ\'e \`a un \'el\'ement $\sigma\in \Gamma_{F}-\Gamma_{E}$. Nous montrerons en \ref{variance}  que l'automorphisme du premier facteur agit par multiplication par $sgn(-1)^{[(i+1)/2]}$ sur $FC^{st}(\mathfrak{g}'(F))$ et  nous avons montr\'e en \ref{actiondunautomorphisme} que celui du second facteur agissait par multiplication par $sgn(-1)^{[( j+2)/4]}$. L'action de $Out({\bf G}')$ est donc triviale si $sgn(-1)=1$, c'est-\`a-dire si $\delta_{4}(q-1)=1$, ou si $[(i+1)/2]+[( j+2)/4]$ est pair. Elle est non triviale si $\delta_{4}(q-1)=0$ et $[(i+1)/2]+[( j+2)/4]$ est impair.  Dans ce dernier cas, ${\bf G}'$ ne nous int\'eresse pas. Dans le premier cas, on calcule le caract\`ere $\xi_{{\bf G}'}$ en utilisant \ref{donneesendoscopiques} et on obtient que $\xi_{{\bf G}'}$ est le caract\`ere dont le noyau est l'image de $E_{G'}^{\times}$ dans $F^{\times}/F^{\times,2}$ par l'application norme.  C'est un \'el\'ement de $\Xi^{ram}$.   Remarquons que dans cette construction,  $E_{G'}$ peut \^etre l'une ou l'autre des deux extensions quadratiques ramifi\'ees de $F$ et  $\xi_{{\bf G}'}$ d\'ecrit alors les deux \'el\'ements de $\Xi^{ram}$. Pour $\xi\in \Xi^{ram}$, le triplet $(i,j,\xi)$ appartient \`a ${\cal Y}^{ram}$. On note ${\bf G}'_{i,j,\xi}$ la donn\'ee ${\bf G}'$ associ\'ee \`a l'extension $E_{G'}$ telle que $\xi_{{\bf G}'}=\xi$ et on pose  $FC^{{\cal E}}_{i,j,\xi}=FC^{st}(\mathfrak{g}'_{i,j,\xi}(F))^{Out({\bf G}'_{i,j,\xi})}$. 
Traitons maintenant le cas o\`u $n$ est pair et $m=n/2$. Le groupe unitaire dispara\^{\i}t. Si $E_{G'}/F$ est ramifi\'ee, on voit que tout se passe comme pr\'ec\'edemment \`a condition de prendre  $j=0$. Mais maintenant, le cas $E_{G'}=E_{0}$ est autoris\'e (o\`u $E_{0}$ est l'extension quadratique non ramifi\'ee de $F$). On a alors $G'_{SC}=Res_{E_{0}/F} (Spin(n+1)_{\vert E_{0}})$. De nouveau, $FC^{st}(\mathfrak{g}'(F))\not=\{0\}$ si et seulement si  $n/2$ est de la forme $i_{0}(i_{0}+1)$ et, si cette condition est v\'erifi\'ee, cet espace est une droite. De nouveau, $Out({\bf G}')$ a deux \'el\'ements. Nous montrerons en \ref{variance}  que ce groupe agit trivialement sur $FC^{st}(\mathfrak{g}'(F))$.  En utilisant \ref{donneesendoscopiques}, on voit que $\xi_{{\bf G}'}$ est l'\'el\'ement non trivial $\xi_{0}$ de $\Xi^{nr}$. On a $(i_{0},i_{0},\xi_{0})\in {\cal Y}^{nr}$, on pose ${\bf G}'={\bf G}'_{i_{0},i_{0},\xi_{0}}$ et   $FC^{{\cal E}}_{i_{0},i_{0},\xi_{0}}=FC^{st}(\mathfrak{g}_{i_{0},i_{0},\xi_{0}}'(F))^{Out({\bf G}_{i_{0},i_{0},\xi_{0}}')}$.

A ce point, on a obtenu une description de $FC^{{\cal E}}(\mathfrak{g}(F))$ qui ressemble \`a celle de \ref{resultats} (2), aux deux diff\'erences suivantes pr\`es. On n'apas trait\'e la donn\'ee principale ${\bf G}$. Dans le cas o\`u il existe un \'el\'ement de ${\cal Y}^{nr}$ de la forme $y=(0,j)$, on ne lui a  pas associ\'e d'espace $FC^{{\cal E}}_{y}$ (parce que, dans le cas d'une action galoisienne triviale, on n'a trait\'e que le cas $m>0$). 

On d\'efinit deux applications $\phi^{nr}:{\mathbb X}^{nr}\to {\mathbb Y}^{nr}$ et $\phi^{ram}={\mathbb X}^{ram}\to {\mathbb Y}^{ram}$ par les formules suivantes:

 pour $x=(k,h)\in {\mathbb X}^{nr}$, $\phi^{nr}(x)=((\vert k-h\vert -1)/2,(k+h-1)/2)$;

pour $x=(k,h)\in {\mathbb X}^{ram}$,
$$\phi^{ram}(x)=\left\lbrace\begin{array}{cc}((k+h-1)/4,(k-h-1)/2),&\,\,si\,\, k \not\equiv h\,\,mod\,\,2{\mathbb Z},\\   ((k-h-2)/4,(k+h)/2),&\,\,si\,\, k \equiv h\,\,mod\,\,2{\mathbb Z}.\\ \end{array} \right.$$

On v\'erifie que ce sont des bijections. Elles se rel\`event de fa\c{c}on \'evidente en des bijections $\varphi^{nr}:{\cal X}^{nr}\to {\cal Y}^{nr}$ et $\varphi^{ram}:{\cal X}^{ram}\to {\cal Y}^{ram}$. Par exemple, on voit qu'il existe un  couple $(0,h_{0})\in {\mathbb X}^{nr}$ si et seulement s'il existe un couple $(i_{0},i_{0})\in {\mathbb Y}^{nr}$. Si ces couples existent, on voit que $\phi^{nr}(0,h_{0})=(i_{0},i_{0})$.  On pose $\varphi^{nr}(0,h_{0},\xi)=(i_{0},i_{0},\xi)$ pour tout $\xi\in \Xi^{nr}$. On note $\varphi:{\cal X}\to {\cal Y}$ la bijection dont les restrictions \`a ${\cal X}^{nr}$ et ${\cal X}^{ram}$ sont $\varphi^{nr}$ et $\varphi^{ram}$. On voit que ${\cal X}^{nr}$ poss\`ede un \'el\'ement $(k,h)$ tel que $\vert k-h\vert =1$ si et seulement si $\delta_{2\triangle}(n)=1$. Dans ce cas, ${\cal X}^{nr}$ poss\`ede un unique tel \'el\'ement que l'on note $(k^{st},h^{st})$ et on pose ${\cal X}^{st}=\{(k^{st},h^{st})\}$. De m\^eme, ${\cal Y}^{nr}$ poss\`ede un \'el\'ement de la forme $(0,j)$
si et seulement si $\delta_{2\triangle}(n)=1$. Dans ce cas, ${\cal Y}^{nr}$ poss\`ede un unique tel \'el\'ement que l'on note $(0,j^{st})$ et on pose ${\cal Y}^{st}=\{(0,j^{st})\}$. On v\'erifie que $\varphi({\cal X}^{st})={\cal Y}^{st}$. 

Maintenant, le m\^eme argument de comparaison des dimensions  qu'en \ref{An-1quasidepram} prouve que $FC^{st}(\mathfrak{g}(F))=\{0\}$ si $\delta_{2\triangle}(n)=0$ tandis que $FC^{st}(\mathfrak{g}(F))$ est une droite si $\delta_{2\triangle}(n)=1$. Dans ce dernier cas, on compl\`ete notre description de l'espace $FC^{{\cal E}}(\mathfrak{g}(F))$ en posant ${\bf G}'_{0,j^{st}}={\bf G}$ et $FC^{{\cal E}}_{0,j^{st}}=FC^{st}(\mathfrak{g}(F))$. On a obtenu \ref{resultats}(2).

On a vu que l'action de $G_{AD}(F)_{0}/\pi(G(F))$ sur un espace  $FC_{x}$ \'etait triviale si $x\in {\cal X}^{nr}$ et non triviale si $x\in {\cal X}^{ram}$. De m\^eme, la restriction de $\xi_{{\bf G}'_{y}}$ \`a $G_{AD}(F)_{0}$ est triviale si $y\in {\cal Y}^{nr}$, non triviale si $y\in {\cal Y}^{ram}$. Pour $\star=nr$ ou $ram$, cela entra\^{\i}ne l'\'egalit\'e

$transfert(\oplus_{x\in {\cal X}^{\star}}FC_{x})=\oplus_{y\in {\cal Y}^{\star}}FC^{{\cal E}}_{y}$.
 
On munit l'ensemble ${\mathbb X}^{\star}$ de la relation $(k,h)\leq (k',h')$ si et seulement si $k+h\leq k'+h'$. On prouve comme en \ref{An-1quasidepram} que  c'est une relation  d'ordre total.
On rel\`eve notre relation d'ordre sur ${\mathbb X}^{\star}$ en une relation de pr\'eordre sur ${\cal X}^{\star}$ qu'on note encore $\leq $. Dans les constructions ci-dessus, on a associ\'e \`a tout $y\in {\cal Y}$ une donn\'ee endoscopique ${\bf G}'_{y}$ et l'application $y\mapsto {\bf G}'_{y}$ est injective.   Soit $y\in {\cal Y}^{\star}$.  On introduira  en \ref{elementsYy}   un \'el\'ement $Y_{ y}\in \mathfrak{g}'_{ y, ell}(F)$ qui a les propri\'et\'es suivantes:

(1) pour tout \'el\'ement non nul $f'\in FC^{{\cal E}}_{y}=FC^{st}(\mathfrak{g}'_{y}(F))^{Out({\bf G}'_{y})}$, on a $S^{G'_{ y}}(Y_{y},f_{{\bf G}'_{y}})\not=0$;

(2) soient $x\in {\cal X}^{\star}$ et $f\in FC_{x}$;  soit $X$ un \'el\'ement de $\mathfrak{g}_{ell}(F)$ correspondant \`a $Y_{ y}$; supposons $I^G(X,f)\not=0$; alors $\varphi^{-1}(y)\leq x$. 

Alors les hypoth\`eses (1) \`a (5) de \ref{ingredients} sont satisfaites pour $\underline{{\cal Y}}^{\sharp}=\underline{{\cal Y}}^{\star}$ avec les notations de ce paragraphe. Cela entra\^{\i}ne

$$(3) \qquad transfert(FC_{(x)})=FC^{{\cal E}}_{\varphi((x))}$$
pour tout $((x))\in \underline{{\cal X}}^{\star}$. Comme en \ref{An-1quasidepram}, on raffine cette \'egalit\'e en tenant compte de l'action de $G_{AD}(F)/\pi(G(F))$ et on obtient \ref{resultats}(3). La relation \ref{resultats} (4) s'en d\'eduit comme en \ref{An-1quasidepram}. Explicitons la cons\'equence de \ref{resultats}(4):

(4) on a $dim(FC^{st}(\mathfrak{g}(F)))=\delta_{2\triangle}(n)$.

   \subsection{Forme int\'erieure du type $B_{n}$ d\'eploy\'e}
On suppose que $G^*$ d\'eploy\'e de type $B_{n}$ avec $n\geq2$ et que $G$ n'est pas d\'eploy\'e. C'est-\`a-dire que $G=Spin_{ndep}(2n+1)$ est la forme non d\'eploy\'ee du groupe spinoriel. Ce cas est presque le m\^eme que le pr\'ec\'edent. 

Notons  ${\mathbb X}^{nr}={\cal X}^{nr}$ l'ensemble des couples $(k,h)\in {\mathbb N}^2$ tels que $k^2+h^2=2n+1$, $k$ est pair, $h$ est impair et $k\not=0$. Si $\delta_{4}(q-1)=1$, on pose ${\cal X}^{ram}=\emptyset$. Si $\delta_{4}(q-1)=0$, on note ${\mathbb X}^{ram}$ l'ensemble des couples $(k,h)\in {\mathbb N}^2$ tels que $k(k+1)/2+h(h+1)/2=2n+1$, $k\geq h$ et l'un des termes $k(k+1)/2$, $h(h+1)/2$ est $\equiv 2\,\,mod\,\, 4{\mathbb Z}$. On  note ${\cal X}^{ram}$ l'ensemble des triplets $(k,h,\xi)$ avec $(k,h)\in {\mathbb X}^{ram}$ et $\xi\in \Xi^{ram}$, o\`u $\Xi^{ram}$ est le m\^eme ensemble que dans le paragraphe pr\'ec\'edent. On pose ${\cal X}={\cal X}^{nr}\cup {\cal X}^{ram}$ et $d_{x}=1$ pour tout $x\in {\cal X}$.

L'ensemble $\underline{S}(G)$ est en bijection avec l'ensemble des couples $(a,b)\in {\mathbb N}^2$ tels que $a+b=n$ et $a\not=0$. Pour $s\in \underline{S}(G)$ param\'etr\'e par $(a,b)$, on a $G_{s}=(Spin_{ndep}(2a)\times Spin(2b+1))/ \{1,(z,z)\}$ si $b\not=0$, o\`u $Spin_{ndep}(2a)$ est la forme non d\'eploy\'ee de ce groupe et $G_{s}=Spin_{ndep}(2a)$ si $b=0$. 
 La preuve du cas pr\'ec\'edent s'applique et conduit \`a la relation \ref{resultats}(1). Il y a deux diff\'erences. Les deux sommets qui \'etaient conjugu\'es par $G_{AD}(F)$ mais pas par $G(F)$ disparaissent. Parce que les groupes $Spin(2a)$ apparaissant dans les groupes $G_{s}$ sont maintenant non d\'eploy\'es, les conditions d'invariance par $\Gamma_{{\mathbb F}_{q}}$ des faisceaux-caract\`eres changent, ce qui conduit \`a la modification de l'ensemble ${\mathbb X}^{ram}$. 

On note ${\mathbb Y}^{nr}={\cal Y}^{nr}$ l'ensemble des couples $(i,j)\in {\mathbb N}^2$ tels que $i(i+1)+j(j+1)=n$ et $i<j$. Si $\delta_{4}(q-1)=1$, on pose ${\cal Y}^{ram}=\emptyset$. Si $\delta_{4}(q-1)=0$, on note ${\mathbb Y}^{ram}$ l'ensemble des couples $(i,j)$ avec $i,j\in {\mathbb N}$, v\'erifiant les conditions

$2i(i+1)+j(j+1)/2=n$;

$[(i+1)/2]+[( j+2)/4]$ est impair.

On note ${\cal Y}^{ram}$ l'ensemble des triplets $(i,j,\xi)$ pour $(i,j)\in {\mathbb Y}^{ram}$ et $\xi\in \Xi^{ram}$. On pose ${\cal Y}={\cal Y}^{nr}\cup {\cal Y}^{ram}$. 

 La construction du paragraphe pr\'ec\'edent s'applique et conduit \`a la relation \ref{resultats}(2). 
Pour une donn\'ee endoscopique  ${\bf G}'$ telle que $Out({\bf G}')\not=\{1\}$,   l'action de ce groupe d'automorphismes ext\'erieurs  est tordue par le caract\`ere non trivial de ce groupe parce que $G$ n'est plus d\'eploy\'e.   Les telles donn\'ees qui apparaissaient dans le cas pr\'ec\'edent disparaissent, c'est celles qui \'etaient \'elimin\'ees qui interviennent maintenant.

On d\'efinit la bijection $\varphi$ comme dans le paragraphe pr\'ec\'edent. La preuve de \ref{resultats}(3) est similaire \`a celle de ce paragraphe.

\subsection{Type $C_{n}$ d\'eploy\'e}\label{Cndeppadique}
On suppose que $G$ est d\'eploy\'e de type $C_{n}$ avec $n\geq2$, c'est-\`a-dire $G=Sp(2n)$. On d\'efinit un homomorphisme $T_{ad}(F)\to F^{\times}/F^{\times,2}$ par $\prod_{l=1,...,n}\check{\varpi}_{l}(x_{l})\mapsto \prod_{l=1,...,n}x_{l}^l$, o\`u  $\check{\varpi}_{l}$ est le copoids associ\'e \`a la racine $\alpha_{l}$. On a $G_{AD}(F)/\pi(G(F))=T_{ad}(F)/\pi(T(F))$. De l'homomorphisme pr\'ec\'edent se d\'eduit un isomorphisme $G_{AD}(F)/\pi(G(F))\to F^{\times}/F^{\times,2}$. L'image de $G_{AD}(F)_{0}$ est $\mathfrak{o}_{F}^{\times}/\mathfrak{o}_{F}^{\times,2}$. Si $n$ est pair, resp. impair, on note $\Xi_{n}$ l'ensemble des \'el\'ements de $\Xi$ triviaux, resp. non triviaux, sur $G_{AD}(F)_{0}$.

On note ${\mathbb X}$ l'ensemble des couples $(k,h)\in {\mathbb N}^2$ tels que $k\geq h$ et $2n=k(k+1)+h(h+1)$. Supposons $\delta_{2\triangle}(n)=0$. On note ${\cal X}$ l'ensemble des triplets $(k,h,\xi)$ avec $(k,h)\in {\mathbb X}$ et $\xi\in \Xi_{n}$. Supposons maintenant   $\delta_{2\triangle}(n)=1$. Alors il existe un unique couple $(k,h)\in {\mathbb X}$ tel que $k=h$. On le note $(k^{st},k^{st})$.  On note ${\cal X}$  la r\'eunion de $\{(k^{st},k^{st})\}$ et de l'ensemble des triplets $(k,h,\xi)$ avec $(k,h)\in {\mathbb X}$, $k>h$ et $\xi\in \Xi_{n}$.  On pose $d_{x}=1$ pour tout $x\in {\cal X}$. 

L'ensemble $\underline{S}(G)$ s'envoie surjectivement sur l'ensemble des couples $(a,b)\in {\mathbb N}^2$ tels que $a\geq b$ et $a+ b=n$. Les fibres de cette application ont deux \'el\'ements sauf, dans le cas o\`u $n$ est pair, au-dessus du couple $(n/2,n/2)$ o\`u la fibre n'a qu'un \'el\'ement. L'action de $G_{AD}(F)/\pi(G(F))$ conserve chaque fibre et y agit transitivement. Pour un sommet $s$ param\'etr\'e par $(a,b)$, on a $G_{s}\simeq Sp(2a)\times Sp(2b)$. D'apr\`es \ref{Cn}, l'espace $FC(\mathfrak{g}_{s}({\mathbb F}_{q}))$ est non nul si et seulement si $(2a,2b)$ est de la forme $(k(k+1),h(h+1))$. Si cette condition est v\'erifi\'ee, l'espace $FC(\mathfrak{g}_{s}({\mathbb F}_{q}))$ est une droite port\'ee par une fonction $f_{N^{a},\epsilon^{a}}\times f_{N^{b},\epsilon^{b}}$. On note $z$ l'\'el\'ement central non trivial de tout groupe symplectique. L'\'el\'ement $z$ de $G$ s'envoie sur le couple $(z,z)$ de $G_{s}$ (ou l'unique $z\in G_{s}$ si $b=0$). On sait que $\epsilon^{a}(z)=(-1)^{a}$ et $\epsilon^{b}(z)=(-1)^b$. Donc l'\'el\'ement $z$ de $G$ agit sur notre fonction par multiplication par $(-1)^{a+b}=(-1)^n$. Puisque l'action de $G_{AD}(F)_{0}$ sur cette fonction est d\'etermin\'ee par cette action de $z$, on voit que    $G_{AD}(F)_{0}$ agit sur notre fonction par le caract\`ere trivial si $n$ est pair et non trivial si $n$ est impair. Si $a>b$, de la fonction $f_{N^{a},\epsilon^{a}}\times f_{N^{b},\epsilon^{b}}$  sont issus deux \'el\'ements de $FC(\mathfrak{g}(F))$ correspondant aux deux \'el\'ements de $\Xi$ dont les restrictions \`a $G_{AD}(F)_{0}$ sont le caract\`ere pr\'ec\'edent, c'est-\`a-dire aux deux \'el\'ements de $\Xi_{n}$. On note $FC_{k,h,\xi}$ la droite port\'ee par la fonction associ\'e \`a $\xi\in \Xi_{n}$. Si $a=b$, ce qui se produit si et seulement si  $\delta_{2\triangle}(n)=1$, on a $(k,h)=(k^{st},k^{st})$. L'action de $G_{AD}(F)$ tout entier sur $G_{s}$ est r\'ecup\'er\'ee par la permutation des deux facteurs et notre fonction $f_{N^{a},\epsilon^{a}}\times f_{N^{b},\epsilon^{b}}$ est fix\'ee par cette action. De cette fonction est issue un unique \'el\'ement de $FC(\mathfrak{g}(F))$  et on note $FC_{k^{st},k^{st}}$ la droite qu'il engendre. Cette description prouve \ref{resultats}(1).

On note ${\mathbb Y}$ l'ensemble des couples $(i,j)\in {\mathbb N}^2$ tels que $n=i(i+1)+j^2$. Remarquons que $j$ est de la m\^eme parit\'e que $n$. Supposons $\delta_{2\triangle}(n)=0$. On note ${\cal  Y}$ l'ensemble des triplets $(i,j,\xi)$ avec $(i,j)\in {\mathbb Y}$ et $\xi\in \Xi_{n}$. Supposons maintenant $\delta_{2\triangle}(n)=1$. 
 Alors il existe un unique couple $(i,j)\in {\mathbb Y}$ avec $j=0$, \`a savoir le couple $(k^{st},0)$.  On note  ${\cal Y}$ la r\'eunion de $ \{(k^{st},0)\}$ et de l'ensemble des triplets $(i,j,\xi)$ avec $(i,j)\in {\mathbb Y}$, $j\not=0$ et $\xi\in \Xi_{n}$.  
 
  D\'eterminons les donn\'ees endoscopiques elliptiques ${\bf G}'$ de $G$ telles que $FC^{st}(\mathfrak{g}'(F))^{Out({\bf G}')}\not=\{0\}$. On consid\`ere un couple $(\sigma\mapsto \sigma_{G'},{\cal O})\in {\cal E}_{ell}(G)$.  Le diagramme $\hat{D}_{a}$ est le diagramme de Dynkin compl\'et\'e d'un groupe de type $B_{n}$. Le groupe $\hat{\Omega}$ co\"{\i}ncide avec le groupe d'automorphismes de ce diagramme, lequel est $\{1,\omega\}$, o\`u $\omega$ permute $\hat{\alpha}_{0}$ et $\hat{\alpha}_{1}$ et fixe les autres racines.

Supposons que l'action galoisienne soit triviale. Si ${\cal O}=\{\hat{\alpha}_{0}\}$ ou $\{\hat{\alpha}_{1}\}$, ces deux cas \'etant conjugu\'es par $\omega$, on a ${\bf G}'={\bf G}$ et, \`a ce point, on ne peut rien dire de $FC^{st}(\mathfrak{g}(F))$. Si ${\cal O}=\{\hat{\alpha}_{m}\}$ avec $m\geq 2$, on voit que $G'_{SC}\simeq Spin_{dep}(2m)\times Sp(2(n-m))$. En appliquant par r\'ecurrence (1) ci-dessous et \ref{Dndeppairpadique} (5), l'espace $FC^{st}(\mathfrak{g}'(F))$ est non nul si et seulement si $m=j^2$  et $n-m=i(i+1)$  avec $i,j\in {\mathbb N}$ et $j$ pair. Si ces conditions sont v\'erifi\'ees, l'espace $FC^{st}(\mathfrak{g}'(F))$ est une droite munie d'un g\'en\'erateur. Le groupe  $Out({\bf G}')$ est $\{1,\omega\}$ et $\omega$ agit par l'automorphisme ext\'erieur non trivial de $Spin_{dep}(2m)$ dont on voit qu'il fixe le g\'en\'erateur. Donc $FC^{st}(\mathfrak{g}'(F))^{Out({\bf G}')}$ est encore une droite. Les actions galoisiennes \'etant triviale, le caract\`ere $\xi_{{\bf G}'}$ est l'\'el\'el\'ement neutre ${\bf 1}$ de $\Xi$. Remarquons que ${\bf 1}$ appartient \`a $\Xi_{n}$ car la parit\'e de $j$ entra\^{\i}ne celle de $n$. Donc $(i,j,1)\in {\cal Y}$.  On note ${\bf G}'_{i,j,{\bf 1}}$ la donn\'ee endoscopique et on pose  $FC^{{\cal E}}_{i,j,{\bf 1}}=FC^{st}(\mathfrak{g}_{i,j,{\bf 1}}'(F))^{Out({\bf G}_{i,j,{\bf 1}}')}$.

Supposons que l'action galoisienne se factorise en une bijection $\Gamma_{E_{0}/F}\to \hat{\Omega}$ (o\`u $E_{0}/F$ est l'extension quadratique non ramifi\'ee). Le lemme \ref{centre} exclut le cas o\`u ${\cal O}=\{\hat{\alpha}_{0},\hat{\alpha}_{1}\}$. Supposons que ${\cal O}=\{\hat{\alpha}_{m}\}$ avec $m\geq 2$. On voit que $G'_{SC}\simeq Spin_{E_{0}/F}(2m)\times Sp(2(n-m))$. Le r\'esultat est le m\^eme que ci-dessus: l'espace $FC^{st}(\mathfrak{g}'(F))^{Out({\bf G}')}$  est non nul   si et seulement si $ m=j^2$  et $n-m=i(i+1)$  avec $i,j\in {\mathbb N}$ et $j$ pair. Si ces conditions sont v\'erifi\'ees, c'est une droite. La diff\'erence est 
 que $\xi_{{\bf G}'}$ n'est plus \'egal \`a ${\bf 1}$. Un calcul facile montre que $\xi_{{\bf G}'}$ est l'\'el\'ement non trivial  de $\Xi_{n}$. Notons-le ici $\xi$. On a $(i,j,\xi)\in {\cal Y}$. On note ${\bf G}'_{i,j,\xi}$ la donn\'ee endoscopique ${\bf G}'$ et  on pose  $FC^{{\cal E}}_{i,j,\xi}=FC^{st}(\mathfrak{g}_{i,j,\xi}'(F))^{Out({\bf G}_{i,j,\xi}')}$. 
 
Soit $E$ une extension quadratique ramifi\'ee de $F$. Supposons que l'action galoisienne se factorise en une bijection $\Gamma_{E/F}\to \hat{\Omega}$. Supposons que ${\cal O}=\{\hat{\alpha}_{m}\}$ avec $m\geq 2$. On voit que $G'_{SC}\simeq Spin_{E/F}(2m)\times Sp(2(n-m))$. Alors l'espace $FC^{st}(\mathfrak{g}'(F))$  est  non nul si et seulement si $ m=j^2$  et $n-m=i(i+1)$  avec $i,j\in {\mathbb N}$ et $j$ impair, cf. (1) ci-dessous, \ref{Dnpairpadiqueram}(1) et \ref{Dnimppadiqueram} (1). Si ces conditions sont v\'erifi\'ees, l'espace $FC^{st}(\mathfrak{g}'(F))$ est une droite. Comme ci-dessus, le groupe $Out({\bf G}')$ agit trivialement sur cet espace. On calcule facilement $\xi_{{\bf G}'}$. On voit que, quand $E$ d\'ecrit les deux extensions quadratiques ramifi\'ees de $F$, le caract\`ere $\xi_{{\bf G}'}$ associ\'e d\'ecrit $\Xi_{n}$ (remarquons que l'hypoth\`ese $j$ impair implique que $n$ l'est aussi). 
Pour $\xi\in \Xi_{n}$ on note ${\bf G}'_{i,j,\xi}$ la donn\'ee d\'etermin\'ee par l'extension $E$ telle que $\xi_{{\bf G}'_{i,j,\xi}}=\xi$. On note  $FC^{{\cal E}}_{i,j,\xi}=FC^{st}(\mathfrak{g}_{i,j,\xi}'(F))^{Out({\bf G}'_{i,j,\xi})}$. Revenons \`a une extension $E$ ramifi\'ee fix\'ee et supposons maintenant que ${\cal O}=\{\hat{\alpha}_{0},\hat{\alpha}_{1}\}$ (ce cas n'est plus exclu par le lemme \ref{centre} puisque $E/F$ est ramifi\'ee).  On   peut remplacer $G'$ par $G'_{SC}\simeq Sp(2n-2)$. Les r\'esultats sont les m\^emes que pr\'ec\'edemment en consid\'erant que l'entier $j$ vaut $1$.

 A ce point, on a associ\'e une donn\'ee ${\bf G}'_{y}$ et une droite $FC^{{\cal E}}_{y}$ \`a tout \'el\'ement $y\in {\cal Y}$, sauf dans le cas $\delta_{2\triangle}(n)=1$, auquel cas on n'a rien associ\'e \`a l'\'el\'ement $(k^{st},0)$. On n'a pas trait\'e non plus la donn\'ee endoscopique principale ${\bf G}$.

 On d\'efinit une application $\phi:{\mathbb X}\to {\mathbb Y}$ par les formules suivantes:

pour $(k,h)\in {\mathbb X}$ avec $k\equiv h\,\,mod\,\,2{\mathbb Z}$, $i=(k+h)/2$, $j=(k-h)/2$;

pour $(k,h)\in {\mathbb X} $ avec $k\equiv h+1\,\,mod\,\,2{\mathbb Z}$, $i=(k-h-1)/2$, $j=(k+h+1)/2$.

C'est une bijection qui se rel\`eve naturellement en une bijection $\varphi={\cal X}\to {\cal Y}$.  Si $\delta_{2\triangle}(n)=0$, on pose ${\cal X}^{st}={\cal Y}^{st}=\emptyset$. Si $\delta_{2\triangle}(n)=1$, on pose ${\cal X}^{st}=\{(k^{st},k^{st})\}$ et ${\cal Y}^{st}=\{(k^{st},0)\}$. On voit que $\varphi({\cal X}^{st})={\cal Y}^{st}$. 

Maintenant, le m\^eme argument de comparaison des dimensions  qu'en \ref{An-1quasidepram} prouve que $FC^{st}(\mathfrak{g}(F))=\{0\}$ si $\delta_{2\triangle}(n)=0$ tandis que $FC^{st}(\mathfrak{g}(F))$ est une droite si $\delta_{2\triangle}(n)=1$. Dans ce dernier cas, on compl\`ete notre description de l'espace $FC^{{\cal E}}(\mathfrak{g}(F))$ en posant ${\bf G}'_{k^{st},0}={\bf G}$ et $FC^{{\cal E}}_{k^{st},0}=FC^{st}(\mathfrak{g}(F))$. On a obtenu \ref{resultats}(2).

On munit l'ensemble ${\mathbb X}$ de la relation $(k,h)\leq (k',h')$ si et seulement si $k+h\leq k'+h'$. On prouve comme en \ref{An-1quasidepram} que   cette relation  est un ordre total. On le rel\`eve en un pr\'eordre sur ${\cal X}$ par la surjection \'evidente ${\cal X}\mapsto {\mathbb X}$.  Soit $y\in {\cal Y}$.  On introduira dans le paragraphe \ref{elementsYy}  un \'el\'ement $Y_{ y}\in \mathfrak{g}'_{ y, ell}(F)$ qui a les propri\'et\'es (2) et (3) de \ref{An-1quasidepram}. 
Alors les hypoth\`eses (1) \`a (5) de \ref{ingredients} sont satisfaites pour $\underline{{\cal Y}}^{\sharp}=\underline{{\cal Y}}$ avec les notations de ce paragraphe. On en d\'eduit \ref{resultats} (3) et (4) comme en \ref{An-1quasidepram}. Explicitons la cons\'equence de \ref{resultats}(4):

(1) on a $dim(FC^{st}(\mathfrak{g}(F)))=\delta_{2\triangle}(n)$.

  \subsection{Forme int\'erieure du type $C_{n}$}
  On suppose que $G^*$ est du type pr\'ec\'edent et que $G$ en est la forme int\'erieure non d\'eploy\'ee. Dans les tables de Tits, le groupe est de type $^2C_{n}$. 
  
  Si   $\delta_{2\triangle}(n)=0$, on pose ${\cal X}= \emptyset$. Supposons  $\delta_{2\triangle}(n)=1$. On note $k^{st}$ l'entier tel que $n=k^{st}(k^{st}+1)$. On pose ${\cal X}=\{(k^{st},k^{st})\}$ et $d_{x}=1$ pour (tout) $x\in {\cal X}$.
  
   Le "local index" des tables de Tits est le diagramme affine ${\cal D}_{a}$ de type $C_{n}$ muni de l'action de $\Gamma_{{\mathbb F}_{q}}$ qui est triviale sur $\Gamma_{{\mathbb F}_{q^2}}$ et telle qu'un \'el\'ement $\sigma\in \Gamma_{{\mathbb F}_{q}}-\Gamma_{{\mathbb F}_{q^2}}$ agisse par l'unique automorphisme $\omega$ de ce diagramme: $\omega$ envoie $\alpha_{m}$ sur $\alpha_{n-m}$. Les \'el\'ements de $\underline{S}(G)$ correspondent aux orbites de cette action dans ${\cal D}_{a}$. Pour un sommet $s$ correspondant \`a une orbite \`a deux \'el\'ements, le lemme \ref{orbites}  montre que  $FC(\mathfrak{g}_{s}({\mathbb F}_{q}))=\{0\}$. Si $n$ est impair, toutes les orbites ont deux \'el\'ements, donc $FC(\mathfrak{g}(F))=\{0\}$. Supposons $n$ pair.  Il n'y a qu'une orbite galoisienne poss\'edant un seul \'el\'ement, \`a savoir $\{\alpha_{n/2}\}$. Notons $s$ le sommet associ\'e. Son unicit\'e entra\^{\i}ne qu'il est conserv\'e par l'action de $G_{AD}(F)$.  Le groupe $G_{s}$ est alors $Res_{{\mathbb F}_{q^2}/{\mathbb F}_{q}}(Sp(n)_{\vert {\mathbb F}_{q^2}})$.  D'apr\`es \ref{Cn}, on sait que $FC(\mathfrak{g}_{s}({\mathbb F}_{q}))$ est non nul si et seulement si $\delta_{2\triangle}(n)=1$. Si cette condition n'est pas v\'erifi\'ee, on a donc $FC(\mathfrak{g}(F))=\{0\}$. Supposons $\delta_{2\triangle}(n)=1$. Alors $FC(\mathfrak{g}_{s}({\mathbb F}_{q}))$ est une droite. On voit qu'elle est fix\'ee par $G_{AD}(F)$. Il s'en d\'eduit une droite dans $FC(\mathfrak{g}(F))$ que l'on note $FC_{k^{st},k^{st}}$. Cela d\'emontre \ref{resultats}(1). 
   
   Si   $\delta_{2\triangle}(n)=0$, on pose $ {\cal Y}=\emptyset$. Supposons  $\delta_{2\triangle}(n)=1$.  On pose  ${\cal Y}=\{(k^{st},0)\}$.

    Puisqu'on a d\'ej\`a d\'etermin\'e ${\cal X}$, on voit que $FC^{{\cal E}}(\mathfrak{g}(F))$ est nul si $\delta_{2\triangle}(n)=0$ et est une droite si $\delta_{2\triangle}(n)=1$. Pour d\'emontrer \ref{resultats} (2), il suffit dans ce dernier cas de d\'eterminer cette droite.  Or,  d'apr\`es \ref{Cndeppadique}, l'espace $FC^{st}(\mathfrak{g}^*(F))$ est une droite. On la note $FC^{{\cal E}}_{k^{st},0}$ et c'est forc\'ement $FC^{{\cal E}}(\mathfrak{g}(F))$ tout entier.   
    
    On note $\varphi$ l'unique bijection de ${\cal X}$ sur ${\cal Y}$. L'assertion \ref{resultats}(3) est triviale.

    \subsection{Type $D_{n}$ d\'eploy\'e, $n$ pair}\label{Dndeppairpadique}
Introduisons d'abord quelques notations g\'en\'erales pour les groupes de type $D_{n}$, avec $n\geq4$. On suppose que $\hat{G}$ est de ce type. 
 On note 

$\theta$ l'automorphisme de $\hat{{\cal D}}_{a}$ qui \'echange $\hat{\alpha}_{n-1}$ et $\hat{\alpha}_{n}$ et qui fixe les autres racines;

 $\theta'$ l'automorphisme de $\hat{{\cal D}}_{a}$ qui \'echange $\hat{\alpha}_{0}$ et $\hat{\alpha}_{1}$ et qui fixe les autres racines;
 
 $\delta$ l'automorphisme de $\hat{{\cal D}}_{a}$ qui envoie $\hat{\alpha}_{i}$ sur $\hat{\alpha}_{n-i}$ pour tout $i=0,...,n$.

 Si $n>4$, le groupe $Aut(\hat{{\cal D}}_{a})$ est engendr\'e par $\theta$, $\theta'$ et $\delta$ et le groupe $Aut(\hat{{\cal D}})$ est r\'eduit \`a $\{1,\theta\}$. Si $n=4$, il y a l'automorphisme suppl\'ementaire $\theta_{3}\in Aut(\hat{{\cal D}})$ introduit en \ref{D4trialitaire} et $Aut(\hat{{\cal D}})\simeq \mathfrak{S}_{3}$. 
 
 Si $n$ est pair,  $\hat{\Omega} =\{1,\delta,\theta\theta',\delta\theta\theta'\}\simeq ({\mathbb Z}/2{\mathbb Z})^2$. Si $n$ est impair, $\hat{\Omega}=\{1,\delta\theta,\theta\theta',\delta\theta'\}\simeq {\mathbb Z}/4{\mathbb Z}$.

On suppose que $G$ est d\'eploy\'e de type $D_{n}$ avec $n\geq4$, c'est-\`a-dire $G=Spin_{dep}(2n)$, avec une notation conforme \`a celles d\'ej\`a utilis\'ees. On suppose $n$ pair. Comme en \ref{Dndeppair}, on a $Z(G)=\{1,z,z',z''\}\simeq ({\mathbb Z}/2{\mathbb Z})^2$.  D\'efinissons un homomorphisme $T_{ad}(F)\to (F^{\times}/F^{\times,2})^2$ par $\prod_{i=1,...,n}\check{\varpi}(x_{i})\mapsto (\prod_{i=1,...,n}x_{i}^{i},x_{n-1}x_{n})$.  Il s'en  d\'eduit un isomorphisme 
$G_{AD}(F)/\pi(G(F))\simeq T_{ad}(F)/\pi(T(F))\simeq (F^{\times}/F^{\times,2})^2$.  Notons $\iota_{1}$ et $\iota_{2}$ les deux plongements \'evidents de $F^{\times}/F^{\times,2}$ dans $(F^{\times}/F^{\times,2})^2$: $\iota_{1}(x)=(x,1)$, $\iota_{2}(x)=(1,x)$. 
L'image de $SO_{dep}(2n,F)$ est   $\iota_{1}(F^{\times}/F^{\times,2})$ et l'image de $G_{AD}(F)_{0}$ est $(\mathfrak{o}_{F}^{\times}/\mathfrak{o}_{F}^{\times,2})^2$. 

Notons $\tilde{\Xi}$ le sous-ensemble des $\xi\in \Xi$ de la forme $(\chi_{E/F},1) $ ou $(\chi_{E/F},\chi_{E/F})$ o\`u $E/F$ est une extension quadratique ramifi\'ee de $F$ et $\chi_{E/F}$ est le caract\`ere associ\'e  de $F^{\times}$.  De l'inclusion  $(F^{\times}/F^{\times,2}\times \mathfrak{o}_{F}^{\times}/\mathfrak{o}_{F}^{\times,2})\to (F^{\times}/F^{\times,2})^2$ se d\'eduit un homomorphisme de restriction de $\Xi$ dans le groupe $(F^{\times}/F^{\times,2}\times \mathfrak{o}_{F}^{\times}/\mathfrak{o}_{F}^{\times,2})^{\vee}$. On voit que

(1) cet homomorphisme de restriction est injectif sur $\tilde{\Xi}$ et  l'image de $\tilde{\Xi}$ est l'ensemble des caract\`eres de $F^{\times}/F^{\times,2}\times \mathfrak{o}_{F}^{\times}/\mathfrak{o}_{F}^{\times,2}$ qui sont non triviaux sur   $\iota_{1}(\mathfrak{o}_{F}^{\times}/\mathfrak{o}_{F}^{\times,2})$.

Notons ${\mathbb X}^{+}$ l'ensemble des couples $(k,h)\in {\mathbb N}^2$ tels que $k^2+h^2=2n$ et $k\geq h$. Remarquons que $k$ et $h$ sont forc\'ement pairs puisque $n$ l'est.  
  On note $\tilde{{\cal X}}^{+}$ l'ensemble des triplets $(k,h,\xi)$ o\`u $(k,h)\in {\mathbb X}^{+}$ et $\xi\in \Xi$ v\'erifient

 $k>h$;

si $k>h>0$, $\xi$ vaut $1$ sur $\iota_{1}(F^{\times}/F^{\times,2})$ et vaut $sgn^{(k+h)/2}$ sur  $\iota_{2}(\mathfrak{o}_{F}^{\times}/\mathfrak{o}_{F}^{\times,2})$;

si $h=0$, $\xi=(1,sgn^{n/2})$ sur $(\mathfrak{o}_{F}^{\times}/\mathfrak{o}_{F}^{\times,2})^2$.

Si $\delta_{\square}(n)=0$, on a $k>h$ pour tout $(k,h)\in {\mathbb X}$. On pose ${\cal X}^+=\tilde{{\cal X}}^+
$. Si $\delta_{\square}(n)=1$, il existe un couple $(k,h)\in {\mathbb X}^+$ tel que $k=h$. On le note $(k^{st},k^{st})$ et on pose ${\cal X}^+=\{(k^{st},k^{st})\}\cup \tilde{{\cal X}}^+$.

Notons ${\mathbb X}^{-}$ l'ensemble des couples $(k,h)\in {\mathbb N}^2$ tels que $k(k+1)/2+h(h+1)/2=2n$, $k\geq h$ et

$k(k+1)/2$ et $h(h+1)/2$ sont pairs si $\delta_{4}(q-1)=1$, $k(k+1)/2$ et $h(h+1)/2$ sont divisibles par $4$ si $\delta_{4}(q-1)=0$.

On note ${\cal X}^{-}$ l'ensemble des triplets $(k,h,\xi)$ o\`u $(k,h)\in {\mathbb X}^{-}$ et $\xi\in \Xi$ v\'erifie:

si $k>h$, la restriction de $\xi$ \`a   $\iota_{1}(\mathfrak{o}_{F}^{\times}/\mathfrak{o}_{F}^{\times,2})$ est non triviale;

si $k=h$, $\xi\in \tilde{\Xi}$.

On pose ${\cal X}={\cal X}^{+}\sqcup {\cal X}^{-}$ et $d_{x}=1$ pour tout $x\in {\cal X}$.

L'ensemble $\underline{S}(G)$ s'envoie surjectivement sur l'ensemble des couples $(a,b)\in {\mathbb N}^2$ tels que $a\geq b$,  $a+b=n$ et $b\not=1$. Les fibres de cette surjection ont un \'el\'ement au-dessus de $(n/2,n/2)$, deux \'el\'ements au-dessus de $(a,b)$ pour $a\not=b$ et $b\not=0$, $4$ \'el\'ements au-dessus de $(n,0)$. L'action de $G_{AD}(F)$ pr\'eserve les fibres et est transitive sur chacune d'elles. Pour $s\in \underline{S}(G)$ param\'etr\'e par $(a,b)$, on a $G_{s}=(Spin_{dep}(2a)\times Spin_{dep}(2b))/  \{1,(z,z)\}$ si $b\not=0$, $G_{s}=Spin_{dep}(2n)$ si $(a,b)=(n,0)$. 
 Consid\'erons un sommet $s\in \underline{S}(G)$ param\'etr\'e par un couple $(a,b)$ avec   $b\geq 2$.   Comme  en \ref{Bndeppadique}, on doit d\'eterminer les couples de fonctions  $f_{N^{a},\epsilon^{a}}\times f_{N^{b},\epsilon^{b}}$ tels que $\epsilon^{a}(z)=\epsilon^{b}(z)$. On utilise les r\'esultats de \ref{Dndeppair} et \ref{Dndepimp}. Consid\'erons ceux pour lesquels $\epsilon^{a}(z)=\epsilon^{b}(z)=1$. Il y a un tel couple si et seulement si 
 $(2a,2b)$ est de la forme $(k^2,h^2)$ et alors, le couple est unique. Le stabilisateur du sommet $s$ dans $G_{AD}(F)$ est l'image r\'eciproque de $F^{\times}/F^{\times,2}\times \mathfrak{o}_{F}^{\times}/\mathfrak{o}_{F}^{\times,2}$. L'action du premier facteur est r\'ecup\'er\'ee par celle du groupe $\{(x,y)\in O_{dep}(2a,{\mathbb F}_{q})\times O_{dep}(2b,{\mathbb F}_{q}); det(x)=det(y)\}$. On voit que cette action fixe  $f_{N^{a},\epsilon^{a}}\times f_{N^{b},\epsilon^{b}}$. L'action du second facteur est triviale ou non selon que l'action sur les faisceaux-caract\`eres  de l'image de $Z(G)$ dans $G_{s}$ est triviale ou non. Remarquons que l'hypoth\`ese $(2a,2b) =(k^2,h^2)$ implique que $a$ et $b$ sont pairs. On voit que l'application de $Z(G)$ dans $Z(G_{s})$ est $z'\mapsto (z',z')$ et $z\mapsto (1,z)=(z,1)$ (puisque $(z,z)=1$). L'action de  $(1,z)$ est triviale. Celle de chaque composante $z'$ est la multiplication par $(-1)^{k/2}$, resp. $(-1)^{h/2}$, donc l'action de $\iota_{2}(\mathfrak{o}_{F}^{\times}/\mathfrak{o}_{F}^{\times,2})$ se fait par   le caract\`ere $sgn^{(k+h)/2}$. De cette fonction $f_{N^{a},\epsilon^{a}}\times f_{N^{b},\epsilon^{b}}$ sont donc issus deux \'el\'ements de $FC(\mathfrak{g}(F))$ correspondant aux deux caract\`eres de $(F^{\times}/F^{\times,2})^2$ prolongeant le caract\`ere que l'on vient de d\'eterminer du sous-groupe $F^{\times}/F^{\times,2}\times \mathfrak{o}_{F}^{\times}/\mathfrak{o}_{F}^{\times,2}$. Pour chacun de ces caract\`eres $\xi$, le triplet  $(k,h,\xi)$ appartient \`a ${\cal X}^{+}$, on note  $FC_{k,h,\xi}$ la droite port\'ee par la fonction correspondante.   Consid\'erons maintenant les fonctions $f_{N^{a},\epsilon^{a}}\times f_{N^{b},\epsilon^{b}}$ telles que $\epsilon^{a}(z)=\epsilon^{b}(z)=-1$. Il y en a si et seulement si $(2a,2b)$ est de la forme $(k(k+1)/2,h(h+1)/2)$ et, de plus, $a$ et $b$ sont pairs dans le cas o\`u $\delta_{4}(q-1)=0$. Il y  en a alors $4$ car chacun des groupes $Spin_{dep}(2a)$ et $ Spin_{dep}(2b)$ porte deux fonctions v\'erifiant $\epsilon^{a}(z)=-1$, resp. $\epsilon^{b}(z)=-1$. Si $a$ et $b$ sont pairs, notons ces fonctions $f_{N^{a},\epsilon^{a},\pm 1}$, resp. $f_{N^{a},\epsilon^{a},\pm 1}$, o\`u le signe est celui par lequel agit $z'$ sur le faisceau-caract\`ere. Si $a$ et $b$ sont impairs, notons ces fonctions $f_{N^{a},\epsilon^{a},\pm i}$, resp. $f_{N^{b},\epsilon^{b},\pm i}$ selon le scalaire $\pm i$ par lequel agit $z'$.  On voit que  les caract\`eres de $G_{AD}(F)_{0}$ selon lesquels se transforment nos fonctions sont les deux caract\`eres de $(\mathfrak{o}_{F}^{\times}/\mathfrak{o}_{F}^{\times,2})^2$ qui sont non triviaux sur la premi\`ere composante. Il y a deux fonctions pour chaque caract\`ere et elles sont permut\'ees par l'action du groupe $\{(x,y)\in O_{dep}(2a,{\mathbb F}_{q})\times O_{dep}(2b,{\mathbb F}_{q}); det(x)=det(y)\}$ (par exemple, si $a$ est pair,  un \'el\'ement de $O_{dep}(2a,{\mathbb F}_{q})$ de d\'eterminant $-1$ envoie $f_{N^{a},\epsilon^{a},1}$ sur un multiple de $f_{N^{a},\epsilon^{a},-1}$). Pour tout caract\`ere $\xi_{0}$ de $F^{\times}/F^{\times,2}\times \mathfrak{o}_{F}^{\times}/\mathfrak{o}_{F}^{\times,2}$ dont la restriction  \`a $\iota_{1}(\mathfrak{o}_{F}^{\times}/\mathfrak{o}_{F}^{\times,2})$ est non triviale, il y a donc une combinaison lin\'eaire de nos fonctions, unique \`a un scalaire pr\`es, qui se transforme selon le caract\`ere $\xi_{0}$. 
 Ensuite, chaque telle fonction donne naissance \`a deux \'el\'ements de $FC(\mathfrak{g}(F))$ correspondant aux deux caract\`eres $\xi$ de $(F^{\times}/F^{\times,2})^2$ prolongeant $\xi_{0}$. Pour tout tel caract\`ere $\xi$, le triplet  $(k,h,\xi)$ appartient \`a ${\cal X}^{-}$. On note $FC_{k,h,\xi}$ la droite port\'ee par la fonction se transformant  selon $\xi$.
 
 Consid\'erons maintenant un sommet $s$ param\'etr\'e par le couple $(n/2,n/2)$. La diff\'erence avec le cas pr\'ec\'edent est que $s$ est fix\'e par $G_{AD}(F)$ tout entier. L'action suppl\'ementaire  sur $G_{s}$ est la permutation des facteurs. Supposons que $n$ est de la forme  $n=2a=2b=(k^{st})^2$. Il y a une fonction   $f_{N^{a},\epsilon^{a}}\times f_{N^{a},\epsilon^{a}}$ avec $\epsilon^{a}(z)=1$. Le m\^eme calcul que ci-dessus montre qu'elle est fix\'ee par l'image r\'eciproque de $F^{\times}/F^{\times,2}\times \mathfrak{o}_{F}^{\times}/\mathfrak{o}_{F}^{\times,2}$ dans $G_{AD}(F)$ et, puisqu'elle est fix\'ee par permutation des facteurs, elle est fix\'ee par tout $G_{AD}(F)$. On note $FC_{k^{st},k^{st}}$  la droite de $FC(\mathfrak{g}(F))$ issue de cette fonction (on a  $(k^{st},k^{st})\in {\cal X}^+$).  Supposons que $n$ est de la forme $n=2a=2b=k(k+1)/2$ et que $n/2$ est pair si $\delta_{4}(q-1)=0$. Il y a alors quatre fonctions $f_{N^{a},\epsilon^{a}}\times f_{N^{b},\epsilon^{b}}$ pour lesquelles $\epsilon^{a}(z)=\epsilon^{b}(z)=-1$. On a d\'ej\`a associ\'e une combinaison lin\'eaire de telles fonctions \`a tout caract\`ere $\xi_{0}$ de $F^{\times}/F^{\times,2}\times \mathfrak{o}_{F}^{\times}/\mathfrak{o}_{F}^{\times,2}$ dont la restriction  \`a $\iota_{1}(\mathfrak{o}_{F}^{\times}/\mathfrak{o}_{F}^{\times,2})$ est non triviale. Son  unicit\'e entra\^{\i}ne qu'elle   se transforme selon un caract\`ere $\xi$ du groupe $G_{AD}(F)/\pi(G(F))$ prolongeant le caract\`ere $\xi_{0}$. Notons $\tilde{\Xi}^X$ l'ensemble des caract\`eres $\xi$ obtenus ainsi. Pour $\xi\in \tilde{\Xi}^X$, on note $FC_{k,k,\xi}$ la droite   de $FC(\mathfrak{g}(F))$ issue de cette fonction. Si $\tilde{\Xi}^X=\tilde{\Xi}$, le triplet $(k,k,\xi)$ appartient \`a ${\cal X}^-$. A ce point, on ne d\'emontre pas cette \'egalit\'e $\tilde{\Xi}^X=\tilde{\Xi}$. Il r\'esulte toutefois de la construction de $\tilde{\Xi}^X$ que cet ensemble v\'erifie la propri\'et\'e (1) donc qu'il a m\^eme nombre d'\'el\'ements que $\tilde{\Xi}$. 
 
   Consid\'erons enfin le cas d'un sommet $s$ param\'etr\'e par $(n,0)$. On a alors $G_{s}=Spin_{dep}(2n)$ et le stabilisateur de $s$ dans $G_{AD}(F)$ est $G_{AD}(F)_{0}$. Supposons $2n=k^2$. Il y a alors une fonction caract\'eristique de faisceau-caract\`ere $f_{N,\epsilon}$ telle que $\epsilon(z)=1$. En remarquant que $k/2$ est de m\^eme parit\'e que $n/2$, le m\^eme calcul que ci-dessus montre que cette fonction se transforme selon le caract\`ere $\xi_{0}=(1,sgn^{n/2})$  de $(\mathfrak{o}_{F}^{\times}/\mathfrak{o}_{F}^{\times,2})^2$. Ensuite, cette fonction donne naissance \`a quatre \'el\'ements de $FC(\mathfrak{g}(F))$  associ\'es aux caract\`eres $\xi$ de $G_{AD}(F)$ qui prolongent $\xi_{0}$. On note $FC_{k,0,\xi}$  la droite port\'ee par la fonction  correspondante.  On a 
     $(k,0,\xi)\in {\cal X}^{+}$. Supposons maintenant  $2n=k(k+1)/2$. Il y a deux fonctions caract\'eristiques de faisceaux-caract\`eres $f_{N,\epsilon}$ telles que $\epsilon(z)=-1$, qui se transforment par les deux caract\`eres de $(\mathfrak{o}_{F}^{\times}/\mathfrak{o}_{F}^{\times,2})^2$ qui sont non triviaux sur le premier facteur. Il s'en d\'eduit huit \'el\'ements de $FC(\mathfrak{g}(F))$ param\'etr\'es par les caract\`eres $\xi$ de $(F^{\times}/F^{\times,2})^2$ dont la restriction  \`a $\iota_{1}(\mathfrak{o}_{F}^{\times}/\mathfrak{o}_{F}^{\times,2})$ est non triviale. On note  $FC_{k,0,\xi}$ les droites port\'ees par ces fonctions,  o\`u  $(k,0,\xi)\in {\cal X}^-$.   Cela d\'emontre  l'assertion \ref{resultats}(1), \`a ceci pr\`es que l'on a remplac\'e l'ensemble $\tilde{\Xi}$ par $\tilde{\Xi}^X$ dans la description de ${\cal X}$.

Notons ${\mathbb Y}^{+}$ l'ensemble des couples $(i,j)\in {\mathbb N}^2$ tels que $i^2+j^2=n$ et $i\geq j$.  
 On note $\tilde{{\cal Y}}^{+}$ l'ensemble des triplets $(i,j,\xi)$ o\`u $(i,j)\in {\mathbb Y}^+$ et $\xi\in \Xi$ v\'erifient
 
 $j>0$;
 
 si $i>j>0$, $\xi$ vaut $1$ sur  $\iota_{1}(F^{\times}/F^{\times,2})$ et vaut $sgn^{i}$ sur   $\iota_{2}(\mathfrak{o}_{F}/\mathfrak{o}_{F}^{\times})$;
 
 si $i=j$, $\xi=(1,sgn^{n/2})$ sur $(\mathfrak{o}_{F}/\mathfrak{o}_{F}^{\times})^2$.
 
 Si $\delta_{\square}(n)=0$, on a $j>0$ pour tout $(i,j)\in {\mathbb Y}^+$ et on pose ${\cal Y}^+=\tilde{{\cal Y}}^+$. Si $\delta_{\square}(n)=1$, il y a un couple $(i,0)\in {\mathbb Y}^+$, \`a savoir le couple $(k^{st},0)$. On pose ${\cal Y}^+=\{(k^{st},0)\}\cup \tilde{{\cal Y}}^+$. 
 
  Notons ${\mathbb Y}^{-}$ l'ensemble des couples $(i,j)$ avec $i,j\in {\mathbb N}$,   tels que $2i^2+j(j+1)/2=n$,
$i$ est pair, 
  $i/2+[(j+2)/4]$ est pair ou $\delta_{4}(q-1)=0$.

 On note ${\cal Y}^{-}$ l'ensemble des triplets $(i,j,\xi)$ o\`u $(i,j)\in {\mathbb Y}^{-}$ et $\xi\in \Xi$ v\'erifie:
 
 si $i\not=0$, la restriction de $\xi$ \`a   $\iota_{1}(\mathfrak{o}_{F}^{\times}/\mathfrak{o}_{F}^{\times,2})$ est non triviale;
 
 si $i=0$, $\xi\in \tilde{\Xi}$.  

On pose ${\cal Y}={\cal Y}^{+}\sqcup {\cal Y}^{-}$.

 Avant de d\'eterminer les donn\'ees endoscopiques de $G$, faisons une remarque sur l'identification du caract\`ere que d\'etermine une telle donn\'ee ${\bf G}'$. Identifions $Z(\hat{G}_{SC})$ \`a $\{\pm 1\}^2$ par $z'\mapsto (-1,1)$ et $z\mapsto (1,-1)$ (les \'el\'ements $z$ et $z'$ sont ici ceux du groupe $\hat{G}_{SC}$). Une donn\'ee ${\bf G}'$ d\'etermine un \'el\'ement de $H^1(W_{F},Z(\hat{G}_{SC}))$ qui, par l'isomorphisme pr\'ec\'edent et celui du corps de classes, d\'etermine un  couple de caract\`eres quadratiques de $F^{\times}$. Les diff\'erentes identifications ont \'et\'e choisies de sorte que ce couple soit pr\'ecis\'ement le caract\`ere $\xi_{{\bf G}'}$ de $G_{AD}(F)/\pi(G(F))\simeq (F^{\times}/F^{\times,2})^2$. 
 
 Calculons les donn\'ees endoscopiques elliptiques ${\bf G}'$ de $G$ telles que $FC^{st}(\mathfrak{g}'(F))^{Out({\bf G}')}\not=\{0\}$.  On consid\`ere un couple $(\sigma\mapsto \sigma_{G'},{\cal O})\in {\cal E}_{ell}(G)$. Puisque $G$ est d\'eploy\'e, l'action $\sigma\mapsto \sigma_{G'}$ se fait par des \'el\'ements de $\hat{\Omega}$. 
 
 Supposons d'abord que l'action $\sigma\mapsto \sigma_{G'}$ soit triviale. L'orbite ${\cal O}$ est r\'eduite \`a une racine. Si celle-ci est une extr\'emit\'e du diagramme $\hat{{\cal D}}_{a}$,  on a ${\bf G}'={\bf G}$. A ce point, on ne peut rien dire de  l'espace $FC^{st}(\mathfrak{g}(F))$. Consid\'erons le cas o\`u ${\cal O}=\{\hat{\alpha}_{m}\}$ pour $m\in \{2,...,n-2\}$. Le diagramme $\hat{{\cal D}}_{a}-{\cal O}$ est le produit de deux diagrammes de type $D_{m}$ et $D_{n-m}$ et $G'_{SC}\simeq Spin_{dep}(2m)\times Spin_{dep}(2n-2m)$. L'action des \'el\'ements de $ \hat{\Omega}$ d\'efinit des \'equivalences entre donn\'ees endoscopiques. En faisant agir $\delta$, on peut supposer $m\geq n/2$. En raisonnant par r\'ecurrence, on peut appliquer les assertions    (5) ci-dessous et \ref{Dndepimppadique} (1): on a $FC^{st}(\mathfrak{g}'(F)) \not=\{0\}$ si et seulement si $m=i^2$ et $n-m=j^2$  avec $i,j$ pairs.  Supposons ces conditions v\'erifi\'ees. Alors l'espace $FC^{st}(\mathfrak{g}'(F))$ est une droite. La donn\'ee ${\bf G}'$ a un automorphisme non trivial, qui est l'action de $\theta\theta'$. En utilisant la description  de $FC^{st}(\mathfrak{g}'(F))$  fournie par   le pr\'esent paragraphe, on voit que cet automorphisme agit trivialement sur  $FC^{st}(\mathfrak{g}'(F))$. Si $m\not=n/2$, l'automorphisme pr\'ec\'edent est le seul non trivial. Si $m=n/2$, le groupe d'automorphismes ext\'erieurs  $Out({\bf G}')$ est $\hat{\Omega}$ tout entier. Mais l'action de $\delta$ consiste \`a \'echanger les deux facteurs de $G'$ et il est clair que cette action est l'identit\'e sur $FC^{st}(\mathfrak{g}'(F))$. Donc $FC^{st}(\mathfrak{g}'(F))^{Out({\bf G}')}$ est une droite. L'action galoisienne \'etant triviale, on a $\xi_{{\bf G}'}=1$.   On a $(i,j,{\bf 1})\in {\cal Y}^+$, on pose $ {\bf G}'_{i,j,{\bf 1}}={\bf G}'$ et $FC^{{\cal E}}_{i,j,{\bf 1}}=FC^{st}(\mathfrak{g}'(F))^{Out({\bf G}_{i,j,{\bf 1}}')}$.

 Supposons que l'action $\sigma\mapsto \sigma_{G'}$ ait pour image dans $\hat{\Omega}$ un sous-groupe d'ordre $2$. Il y a alors une extension quadratique $E/F$ et un \'el\'ement $\omega\in \hat{\Omega}-\{1\}$ tels que $\sigma_{G'}=1$ si $\sigma\in \Gamma_{E}$ et $\sigma_{G'}=\omega$ si $\sigma\in \Gamma_{F}-\Gamma_{E}$. Supposons d'abord que $\omega=\theta\theta'$. L'orbite ${\cal O}$ peut avoir un seul \'el\'ement $\hat{\alpha}_{m}$ pour $m\in \{2,...,n-2\}$. De nouveau, on peut supposer $m\geq n/2$. Ce cas est voisin du pr\'ec\'edent, le groupe   $G_{SC}'$ \'etant cette fois $Spin_{E/F}(2m)\times Spin_{E/F}(2n-2m)$. Si $E/F$ est non ramifi\'ee, on a $FC^{st}(\mathfrak{g}'(F)) \not=\{0\}$ si et seulement si $m=i^2$ et $n-m=j^2$  avec $i,j$ pairs. Si $E/F$ est ramifi\'ee,  on a $FC^{st}(\mathfrak{g}'(F)) \not=\{0\}$ si et seulement si $m=i^2$ et $n-m=j^2$  avec $i,j$ impairs (remarquons que l'on a $j\geq2$ puisque $m\leq n-2$). Supposons ces conditions v\'erifi\'ees. L'espace $FC^{st}(\mathfrak{g}'(F))$ est une droite et on voit comme ci-dessus que le groupe $Out({\bf G}')$ agit trivialement sur cette droite. Le caract\`ere $\xi_{{\bf G}'}$ est trivial sur $\iota_{1}(F^{\times}/F^{\times,2})$ et est \'egal au caract\`ere  $\chi_{E/F}$   sur $\iota_{2}(F^{\times}/F^{\times,2})$. Notons $\xi$ ce caract\`ere. On a $(i,j,\xi)\in {\cal Y}^+$, on pose ${\bf G}'_{i,j,\xi}={\bf G}'$ et  $FC^{{\cal E}}_{i,j,\xi}=FC^{st}(\mathfrak{g}_{i,j,\xi}'(F))^{Out({\bf G}_{i,j,\xi}')}$.
  Consid\'erons maintenant une orbite ${\cal O}$ \`a deux \'el\'ements. Cette orbite est alors \'egale 
  \`a $\{\hat{\alpha}_{0},\hat{\alpha}_{1}\}$ ou $\{\hat{\alpha}_{n-1},\hat{\alpha}_{n}\}$. Ces deux ensembles se d\'eduisent l'un de l'autre par l'action de $\delta$, on peut consid\'erer seulement le cas ${\cal O}= \{\hat{\alpha}_{n-1},\hat{\alpha}_{n}\}$.  Le lemme \ref{centre} exclut le cas o\`u $E/F$ est non ramifi\'ee.  Supposons $E/F$ ramifi\'ee. On a $G'_{SC}\simeq Spin_{E/F}(2n-2)$.  On a $FC^{st}(\mathfrak{g}'(F)) \not=\{0\}$ si et seulement si $n-1=i^2$ avec $i$ forc\'ement impair. La discussion se poursuit comme pr\'ec\'edemment,  on  pose ${\bf G}'_{i,1,\xi}={\bf G}'$, $FC^{{\cal E}}_{i,1,\xi}=FC^{st}(\mathfrak{g}'_{i,1,\xi}(F))^{Out({\bf G}_{i,1,\xi}')}$, o\`u $\xi$ est comme ci-dessus.
  
  Supposons maintenant que $\sigma_{G'}=1$ pour $\sigma\in \Gamma_{E}$ mais que $\sigma_{G'}=\delta$ pour $\sigma\in \Gamma_{F}-\Gamma_{E}$. Consid\'erons d'abord le cas o\`u ${\cal O}=\{\hat{\alpha}_{m},\hat{\alpha}_{n-m}\}$, avec $2\leq m<n/2$. De nouveau, le lemme \ref{centre} exclut  le cas o\`u $E/F$ est non ramifi\'ee.  On suppose donc $E/F$ ramifi\'ee. On a $G'_{SC}\simeq Res_{E/F}(Spin_{dep}(2m))\times SU_{E/F}(n-2m)$. En utilisant  par r\'ecurrence (5) ci-dessous et  \ref{An-1quasidepram} (5), on a $FC^{st}(\mathfrak{g}'(F))\not=\{0\}$ si et seulement si $m=i^2$ avec $i$ pair et $n-2m=j(j+1)/2$. Puisque $n$ est pair, $[(j+1)/2]$ est forc\'ement  pair. Supposons ces conditions sont v\'erifi\'ees. L'espace $FC^{st}(\mathfrak{g}'(F))$ est une droite.
   Le groupe $Out({\bf G}')$ est $\hat{\Omega}$ tout entier. L'action de $\theta\theta'$ est l'automorphisme usuel du facteur $Spin_{dep}(2m)$ qui agit  trivialement sur la droite  $FC^{st}(\mathfrak{g}'(F))$. L'action de $\delta$ n'est autre que l'action galoisienne naturelle de $\Gamma_{E/F}$ sur chacun des facteurs. On montrera en \ref{variance} que cette action est la multiplication par $sgn(-1)^{i/2}$ sur le premier facteur et on a montr\'e en \ref{actiondunautomorphisme} que c'\'etait la  multiplication par $sgn(-1)^{[(j+2)/4]}$ sur le second.  Donc $FC^{st}(\mathfrak{g}'(F))^{Out({\bf G}')}$ est une droite si $sgn(-1)=1$, c'est-\`a-dire $\delta_{4}(q-1)=1$, ou si $i/2+[(j+2)/4]$ est pair, et est nul sinon. Supposons ces conditions v\'erifi\'ees.  On calcule le caract\`ere $\xi_{{\bf G}'}$ comme en \ref{donneesendoscopiques}. Pour cela, donnons une formule g\'en\'erale pour l'\'el\'ement $s_{sc}$ associ\'e \`a notre donn\'ee, qui est valable sans hypoth\`ese de parit\'e sur  $m$ ou $n$. On voit que l'on peut choisir
$$(2) \qquad s_{sc}=(\prod_{l=1,...,m}\check{\hat{\alpha}}_{l}((-1)^l))(\prod_{l=m+1,...,n-m}\check{\hat{\alpha}}_{l}(i^{l+m})(\prod_{l=n-m+1,...,n-2}\check{\hat{\alpha}}_{l}(i^n))\check{\hat{\alpha}}_{n-1}(\zeta_{8}^n)\check{\hat{\alpha}}_{n}(\zeta_{8}^n),$$
 o\`u $\zeta_{8}$ est une racine carr\'ee de $i$ dans ${\mathbb C}^{\times}$ et o\`u $\check{\hat{\alpha}}_{l}$ est la coracine associ\'ee \`a la racine $\hat{\alpha}_{l}$. 
 On calcule $\delta(s_{sc})s_{sc}^{-1}=z'z^{n/2}$ et $\theta\theta'(s_{sc})s_{sc}^{-1}=z$. 
  Pour $\tau\in \Gamma_{E}-\Gamma_{F}$,  on a $\tau_{G'}=\delta$. On en d\'eduit que $\xi_{{\bf G}'}$
    est \'egal  \`a $(\chi_{E/F},\chi_{E/F}^{n/2})$. Notons $\xi$ ce caract\`ere.  On a $(i,j,\xi)\in {\cal Y}^{-}$. On pose ${\bf G}'_{i,j,\xi}={\bf G}'$ et $FC^{{\cal E}}_{i,j,\xi}=FC^{st}(\mathfrak{g}_{i,j,\xi}'(F))^{Out({\bf G}_{i,j,\xi}')}$.
   Consid\'erons maintenant le cas o\`u ${\cal O}=\{\hat{\alpha}_{n/2}\}$. Le groupe unitaire dispara\^{\i}t. Si $E/F$ est ramifi\'ee, ce cas ne diff\`ere pas du pr\'ec\'edent, on a simplement $j=0$. Mais, cette fois, le cas o\`u $E/F$ est non ramifi\'ee n'est plus exclu. Dans ce cas, on a $FC^{st}(\mathfrak{g}'(F))\not=\{0\}$ si et seulement si $n/2=i^2$ avec $i$ pair et on montrera en \ref{variance} qu'alors, $Out({\bf G}')$ agit trivialement sur cet espace. Le caract\`ere $\xi_{{\bf G}'}$ est encore \'egal \`a $(\chi_{E/F},\chi_{E/F}^{n/2})$. Notons $\xi$ ce caract\`ere.  On a $(i,i,\xi)\in {\cal Y}^+$. On pose ${\bf G}'_{i,i,\xi}={\bf G}'$ et $FC^{{\cal E}}_{i,i,\xi}=FC^{st}(\mathfrak{g}_{i,i,\xi}'(F))^{Out({\bf G}_{i,i,\xi}')}$.
  Consid\'erons enfin le cas o\`u ${\cal O}$ est \'egale \`a $\{\hat{\alpha}_{0},\hat{\alpha}_{n}\}$ ou $\{\hat{\alpha}_{1},\hat{\alpha}_{n-1}\}$. Ces deux orbites se d\'eduisent l'une de l'autre par $\theta\theta'\in\hat{ \Omega}$, on peut ne consid\'erer que la premi\`ere. De nouveau, l'extension $E/F$ non ramifi\'ee est exclue. Pour $E/F$ ramifi\'ee, la situation est la m\^eme que dans le premier cas trait\'e (${\cal O}=\{\hat{\alpha}_{m},\hat{\alpha}_{n-m}\}$, avec $2\leq m<n/2$), l'entier $i$ devenant nul.
  
  Supposons que $\sigma_{G'}=1$  pour $\sigma\in \Gamma_{E}$ mais que $\sigma_{G'}=\delta\theta\theta'$ pour $\sigma\in \Gamma_{F}-\Gamma_{E}$. Ce cas est analogue au pr\'ec\'edent, la seule diff\'erence est le calcul du caract\`ere $\xi_{{\bf G}'}$: on a maintenant $\tau_{G'} =\delta\theta\theta'$. Le caract\`ere $\xi_{{\bf G}'}$ vaut alors $(\chi_{E/F},\chi_{E/F}^{1+n/2})$. On pose les m\^emes d\'efinitions que pr\'ec\'edemment, au changement pr\`es de ce caract\`ere. 
  
  Supposons enfin que l'action $\sigma\mapsto \sigma_{G'}$ ait pour image $\hat{\Omega}$ tout entier. Puisque ce groupe est isomorphe \`a $ ({\mathbb Z}/2{\mathbb Z})^2$, le noyau de cette action est forc\'ement le groupe $\Gamma_{Q}$ o\`u $Q/F$ est l'extension biquadratique. Le cas o\`u  ${\cal O}=\{\hat{\alpha}_{0},\hat{\alpha}_{1},\hat{\alpha}_{n-1},\hat{\alpha}_{n}\}$ est exclu par le lemme \ref{centre} puisque $Q/F$ n'est pas totalement ramifi\'ee.   Consid\'erons maintenant le cas o\`u ${\cal O}=\{\hat{\alpha}_{m},\hat{\alpha}_{n-m}\}$, avec $2\leq m<n/2$. Notons $E/F$ l'extension quadratique telle que l'image de $\Gamma_{E}$ par $\sigma\mapsto \sigma_{G'}$ soit $\{1,\theta\theta'\}$. Le stabilisateur d'un \'el\'ement de l'orbite est $\Gamma_{E}$. Le lemme \ref{centre} exclut le cas o\`u $E/F$ est non ramifi\'ee. Supposons donc $E/F$ ramifi\'ee. 
  Fixons $\tau\in \Gamma_{F}-\Gamma_{E}$. On a $\tau_{G'}=\delta$ ou $\delta\theta\theta'$. On a $G'_{SC}\simeq Res_{E/F}(Spin_{Q/E}(2m))\times SU_{E/F}(n-2m)$. Puisque $E/F$ est ramifi\'ee, $Q/E$ ne l'est pas et la situation est \`a peu pr\`es la m\^eme que dans le cas o\`u l'action galoisienne se factorise par $\Gamma_{E}$. On a $FC^{st}(\mathfrak{g}'(F))\not=\{0\}$ si et seulement si $m=i^2$ avec $i$ pair et $n-2m=j(j+1)/2$. Supposons ces conditions v\'erifi\'ees. On voit encore que $FC^{st}(\mathfrak{g}'(F))^{Out({\bf G}')}$ est une droite si $\delta_{4}(q-1)=1$ ou si $i/2+[(j+2)/4]$ est pair, et est nul sinon. L'\'el\'ement $s_{sc}$ est comme ci-dessus et on obtient que le caract\`ere $\xi_{{\bf G}'}$ vaut $(\chi_{E/F},\chi_{E'/F})$, o\`u $E'/F$ est l'extension quadratique telle que l'image de $\Gamma_{E'}$ par l'application $\sigma\mapsto \sigma_{G'}$ soit $\{1,\delta\}$ si $n/2$ est pair, $\{1,\delta\theta\theta'\}$ si $n/2$ est impair. Les deux actions possibles de $\tau_{G'}$ donnent en tout cas les deux extensions $E'/F$ distinctes de $E$, donc deux caract\`eres $\xi\in \Xi$. Pour chacun d'eux, on a $(i,j,\xi)\in {\cal Y}^-$. On pose ${\bf G}'_{i,j,\xi}={\bf G}'$ et $FC^{{\cal E}}_{i,j,\xi}=FC^{st}(\mathfrak{g}_{i,j,\xi}'(F))^{Out({\bf G}_{i,j,\xi}')}$.
 Consid\'erons enfin le cas o\`u ${\cal O}=\{\hat{\alpha}_{n/2}\}$. On d\'efinit $E/F$ comme ci-dessus. Si $E/F$ est ramifi\'ee, la situation est la m\^eme que dans le cas que l'on vient de traiter, $j$ devenant nul. Maintenant, le cas $E/F$ non ramifi\'ee n'est plus exclu. Supposons que $E/F$ soit non ramifi\'ee.    On a $G'_{SC}\simeq Res_{E/F}(Spin_{Q/E}(n))$. L'extension $Q/E$ est maintenant ramifi\'ee. Alors $FC^{st}(\mathfrak{g}'(F))$ est non nul si et seulement si $n/2=i^2$ avec $i$ impair. Supposons qu'il en soit ainsi.  L'espace $FC^{st}(\mathfrak{g}'(F))$ est une droite et on voit que $Out({\bf G}')$ y agit trivialement. L'\'el\'ement $s_{sc}$ est comme pr\'ec\'edemment et  $n/2$ est  impair. Alors $\xi_{{\bf G}'}$ est \'egal \`a $(\chi_{E/F},\chi_{E'/F})$, o\`u $E'/F$ est l'extension quadratique telle que l'image de $\Gamma_{E'}$ par l'application $\sigma\mapsto \sigma_{G'}$ soit  $\{1,\delta\theta\theta'\}$.
  Selon la valeur de $\tau_{G'}$ qui peut \^etre \'egal \`a $\delta$ ou $\delta\theta\theta'$, on obtient  les deux extensions $E'/F$ ramifi\'ees donc deux caract\`eres $\xi$.   Pour chacun d'eux, on  a $(i,i,\xi)\in {\cal Y}^+$. On pose ${\bf G}'_{i,i,\xi}={\bf G}'$ et $FC^{{\cal E}}_{i,i,\xi}=FC^{st}(\mathfrak{g}_{i,i,\xi}'(F))^{Out({\bf G}_{i,i,\xi}')}$.
 
   En rassemblant ces calculs, on a obtenu la description \ref{resultats}(2) sauf sur deux points. On n'a pas trait\'e la donn\'ee principale ${\bf G}$. Si $\delta_{\square}(n)=1$,   aucune des donn\'ees d\'ecrites n'est param\'etr\'ee par le triplet $(k^{st},0)\in {\cal Y}^{+}$.     
      
      On d\'efinit des applications $\phi^+:{\mathbb X}^+\to {\mathbb Y}^+$ et  $\phi^-:{\mathbb X}^-\to {\mathbb Y}^-$   par les formules suivantes:

pour $(k,h)\in {\mathbb X}^{+}$, $\phi^+(k,h)=((k+h)/2,(k-h)/2)$;

pour $(k,h)\in {\mathbb X}^{-}$,
$$\phi^-(k,h)=\left\lbrace\begin{array}{cc}((k-h)/4,(k+h)/2),& \,\,si\,\,k \equiv h\,\,mod\,\,2{\mathbb Z},\\ ((k+h+1)/4,(k-h-1)/2),&\,\,si \,\,k \not\equiv h\,\,mod\,\,2{\mathbb Z}.\\   \end{array}\right.$$
On v\'erifie que  ce sont des bijections.     Si $\delta_{\square}(n)=0$, on pose ${\cal X}^{st}={\cal Y}^{st}=\emptyset$. Si $\delta_{\square}(n)=1$, on pose 
${\cal X}^{st}=\{(k^{st},k^{st})\}$ et ${\cal Y}^{st}=\{(k^{st},0)\}$.

Posons $\star=\pm$. Si $\star=+$ ou si $\star=-$ et qu'il n'existe pas d'\'el\'ement de ${\mathbb X}^-$ de la forme $(k,k)$, la bijection $\phi^{\star}$  se rel\`eve  naturellement en une bijection $\varphi^{\star}:{\cal X}^{\star}\to {\cal Y}^{\star}$. Par contre, si $\star=-$ et qu'il existe un tel \'el\'ement $(k,k)\in {\mathbb X}^-$, on ne peut pas encore relever naturellement la bijection $\phi^{\star}$  en une bijection $\varphi^{\star}:{\cal X}^{\star}\to {\cal Y}^{\star}$ parce que l'on 
 n'a pas encore d\'emontr\'e l'\'egalit\'e $\tilde{\Xi}^X=\tilde{\Xi}$.  Mais on a vu que ces deux ensembles avaient m\^eme nombre d'\'el\'ements. Fixons de fa\c{c}on provisoire  une bijection arbitraire entre ces deux ensembles. Sous les hypoth\`eses pr\'ec\'edentes, on  peut alors relever $\phi^{\star}$ en une bijection $\varphi^{\star}:{\cal X}^{\star}\to {\cal Y}^{\star}$. On r\'eunit $\varphi^+$ et $\varphi^-$ en une bijection $\varphi:{\cal X}\to {\cal Y}$. Remarquons que $\varphi({\cal X}^{st})={\cal Y}^{st}$. 
 
  Le m\^eme argument de dimensions qu'en \ref{An-1quasidepram}   montre alors que $FC^{st}(\mathfrak{g}(F))$ est  une droite si  $\delta_{\square}(n)=1$ et est nul sinon. Si $\delta_{\square}(n)=1$, on pose ${\bf G}'_{k^{st},0}={\bf G}$, $FC^{{\cal E}}_{k^{st},0}=FC^{st}(\mathfrak{g}(F))$. Cela ach\`eve la preuve de \ref{resultats}(2).  
  
    Soit de nouveau $\star=\pm$. Montrons que
 
 (3) $transfert(\oplus_{x\in {\cal X}^{\star}}FC_{x})=\oplus_{y\in {\cal Y}^{\star}}FC^{{\cal E}}_{y}$.

 Pour $x\in {\cal X}$, notons $\xi_{x}\in \Xi$ le caract\`ere tel que $FC_{x}\subset I_{cusp,\xi_{x}}(\mathfrak{g}(F))$. Il r\'esulte de nos descriptions que $\xi_{x}$ est trivial, resp. non trivial, sur $\iota_{1}(\mathfrak{o}_{F}^{\times}/\mathfrak{o}_{F}^{\times,2})$ si et seulement si $x\in {\cal X}^+$, resp. $x\in {\cal X}^-$. De m\^eme, pour $y\in {\cal Y}$, le caract\`ere $\xi_{{\bf G}'_{y}}$ est trivial, resp. non trivial, sur $\iota_{1}(\mathfrak{o}_{F}^{\times}/\mathfrak{o}_{F}^{\times,2})$ si et seulement si $y\in {\cal Y}^+$, resp. $x\in {\cal Y}^-$. Alors (3) r\'esulte de la compatibilit\'e du transfert avec les actions de $G_{AD}(F)/\pi(G(F))$.

    On munit l'ensemble ${\mathbb X}^{\star}$ de la relation $(k,h)\leq (k',h')$ si et seulement si $k+h\leq k'+h'$. On prouve comme en \ref{An-1quasidepram} que c'est un ordre total. On applique alors les constructions de \ref{ingredients}. Remarquons que la bijection  encore not\'ee $\varphi:\underline{{\cal X}}^{\star}\to \underline{{\cal Y}}^{\star}$ qui se d\'eduit de $\varphi$ s'identifie \`a $\phi^{\star}:{\mathbb X}^{\star}\to {\mathbb Y}^{\star}$ et est donc  ind\'ependante du choix arbitraire que l'on a fait ci-dessus.  
   Soit $y\in {\cal Y}^{\star}$.  On introduira dans le paragraphe \ref{elementsYy}  un \'el\'ement $Y_{ y}\in \mathfrak{g}'_{ y, ell}(F)$ qui a les propri\'et\'es  (1) et (2) de \ref{Bndeppadique}.
 Alors les hypoth\`eses (1) \`a (5) de \ref{ingredients} sont satisfaites pour $\underline{{\cal Y}}^{\sharp}=\underline{{\cal Y}}^{\star}$. Cela entra\^{\i}ne

$$(4) \qquad transfert(FC_{(x)})=FC^{{\cal E}}_{\varphi^{\star}((x))}$$
pour tout $((x))\in \underline{{\cal X}}^{\star}$. Supposons  qu'il existe un \'el\'ement de ${\mathbb X}^-$ de la forme $(k,k)$. Appliquons l'\'egalit\'e (4) pour $\star=-$ en prenant pour $(x)$ la classe qui se projette sur cet \'el\'ement. On a $FC_{(x)}=\oplus_{\xi\in \tilde{\Xi}^X}FC_{k,k,\xi}$ et  $FC_{k,k,\xi}\subset I_{cusp,\xi}(\mathfrak{g}(F))$. L'ensemble des caract\`eres $\xi_{{\bf G}'_{y}}$ pour $y\in \varphi^{-}((x))$ est $\tilde{\Xi}$. La compatibilit\'e du transfert \`a l'action de $G_{AD}(F)/\pi(G(F))$ entra\^{\i}ne l'\'egalit\'e $\tilde{\Xi}^X=\tilde{\Xi}$, ce qui ach\`eve la preuve de \ref{resultats}(1). Il y avait plus haut un arbitraire dans la d\'efinition de $\varphi$. Puisque $\tilde{\Xi}^X=\tilde{\Xi}$, il y a un choix canonique de telle bijection  et c'est celui que l'on fait. 
On prouve les relations \ref{resultats}(3) et (4) comme en \ref{An-1quasidepram}. Explicitons la cons\'equence de \ref{resultats}(4):

(5) on a $dim(FC^{st}(\mathfrak{g}(F))=\delta_{\square}(n)$. 

 \subsection{Forme int\'erieure classique du type $D_{n}$ d\'eploy\'e, $n$ pair}
 On suppose que $G^*$ est du type pr\'ec\'edent, en particulier  $n$ est pair. Le diagramme de Dynkin compl\'et\'e ${\cal D}_{a}$ de $G^*$ est identique au diagramme $\hat{{\cal D}}_{a}$ de $\hat{G}$. Identifions ces deux diagrammes. Le groupe $N$ de \ref{orbites}   s'identifie au groupe $\hat{\Omega}$.  On suppose ici que $G$ est la forme int\'erieure de $G^*$ associ\'ee \`a l'\'el\'ement $\theta\theta'\in \hat{\Omega}\simeq N$, ou encore au caract\`ere de $Z(\hat{G}_{SC})$ dont le noyau est $\{1,z\}$. Alors $G$ est la forme non d\'eploy\'ee "classique" du groupe $Spin(2n)$, autrement dit le groupe $Spin$ associ\'e \`a un espace de dimension $2n $ sur $F$ muni d'une forme quadratique de d\'eterminant $1$ telle que les sous-espaces totalement isotropes maximaux sont de dimension $n-2$. Dans les tables de Tits, le groupe est de type $^2D'_{n}$. On a encore $G_{AD}(F)/\pi(G(F))\simeq (F^{\times}/F^{\times,2})^2$, le quotient $G_{AD}(F)/\pi(G(F))$ ne changeant pas par passage \`a une forme int\'erieure.

 Notons ${\mathbb X}^{+}$ l'ensemble des couples $(k,h)\in {\mathbb N}^2$ tels que $k^2+h^2=2n$ et $k\geq h\geq2$. Remarquons que $k$ et $h$ sont forc\'ement pairs puisque $n$ l'est.  
  On note $\tilde{{\cal X}}^{+}$ l'ensemble des triplets $(k,h,\xi)$ o\`u $(k,h)\in {\mathbb X}^{+}$ et $\xi\in \Xi$ v\'erifient

 $k>h$;

 $\xi$ vaut $1$ sur $\iota_{1}(F^{\times}/F^{\times,2})$ et vaut $sgn^{(k+h)/2}$ sur  $\iota_{2}(\mathfrak{o}_{F}^{\times}/\mathfrak{o}_{F}^{\times,2})$.

Si $\delta_{\square}(n)=0$, on pose ${\cal X}^+=\tilde{{\cal X}}^+$. Si $\delta_{\square}(n)=1$, il y a un \'el\'ement $(k,h)\in {\mathbb X}^+$ tel que $k=h$, on le note $(k^{st},k^{st})$. On pose ${\cal X}^+=\{(k^{st},k^{st})\}\cup \tilde{{\cal X}}^+$. 

 Si $\delta_{4}(q-1)=1$, posons ${\mathbb X}^-=\emptyset$. Si $\delta_{4}(q-1)=0$, ce qui \'equivaut \`a $\delta_{4}(q+1)=1$, notons 
 ${\mathbb X}^{-}$ l'ensemble des couples $(k,h)\in {\mathbb N}^2$ tels que $k(k+1)/2+h(h+1)/2=2n$, $k\geq h$ et
 $k(k+1)/2$ et $h(h+1)/2$ sont  congrus \`a $2$ modulo $4$.

On note ${\cal X}^{-}$ l'ensemble des triplets $(k,h,\xi)$ o\`u $(k,h)\in {\mathbb X}^{-}$ et $\xi\in \Xi$ v\'erifient:

si $k>h$, la restriction de $\xi$ \`a   $\iota_{1}(\mathfrak{o}_{F}^{\times}/\mathfrak{o}_{F}^{\times,2})$ est non triviale;

si $k=h$, $\xi\in \tilde{\Xi}$, avec le m\^eme ensemble $\tilde{\Xi}$ qu'en \ref{Dndeppairpadique}.

On pose ${\cal X}={\cal X}^{+}\sqcup {\cal X}^{-}$ et $d_{x}=1$ pour tout $x\in {\cal X}$.  

 L'ensemble $\underline{S}(G)$ s'envoie surjectivement sur l'ensemble des couples $(a,b)\in {\mathbb N}^2$ tels que $a\geq b$,  $a+b=n$ et $b\geq1$. Les fibres de cette surjection ont un \'el\'ement au-dessus de $(n/2,n/2)$, deux \'el\'ements au-dessus de $(a,b)$ pour $a\not=b$. L'action de $G_{AD}(F)$ pr\'eserve les fibres et est transitive sur chacune d'elles. Pour $s\in \underline{S}(G)$ param\'etr\'e par $(a,b)$, on a $G_{s}=(Spin_{ndep}(2a)\times Spin_{ndep}(2b))/ \{1,(z,z)\}$ (il s'agit des formes non d\'eploy\'ees des groupes en question). La preuve de l'assertion \ref{resultats}(1) est alors similaire \`a celle du paragraphe pr\'ec\'edent. Les diff\'erences sont d'une part que les sommets param\'etr\'es par la couple $(n,0)$ disparaissent. D'autre part, parce que les groupes $Spin$ apparaissant dans les groupes $G_{s}$ sont maintenant non d\'eploy\'es,  les fonctions  $f_{N^{a},\epsilon^{a}}\times f_{N^{b},\epsilon^{b}}$ telles que $\epsilon^{a}(z)=\epsilon^{b}(z)=-1$ existent  si et seulement si $\delta_{4}(q+1)=1$,  $(2a,2b)$ est de la forme $(k(k+1)/2,h(h+1)/2)$ et $a$ et $b$ sont impairs.

Notons ${\mathbb Y}^{+}$ l'ensemble des couples $(i,j)\in {\mathbb N}^2$ tels que $i^2+j^2=n$ et $i> j$.  
  On note $\tilde{{\cal Y}}^{+}$ l'ensemble des triplets $(i,j,\xi)$ o\`u $(i,j)\in {\mathbb Y}^{+}$ et $\xi\in \Xi$ v\'erifie
 
 $j>0$;
 
 $\xi$ vaut $1$ sur  $\iota_{1}(F^{\times}/F^{\times,2})$ et vaut $sgn^{i}$ sur   $\iota_{2}(\mathfrak{o}_{F}/\mathfrak{o}_{F}^{\times})$.
 
 Si $\delta_{\square}(n)=0$, on pose ${\cal Y}^+=\tilde{{\cal Y}}^+$. Si $\delta_{\square}(n)=1$, il y a un \'el\'ement $(i,0)\in {\mathbb Y}^+$, \`a savoir l'\'el\'ement  $(k^{st},0)$. On pose ${\cal Y}^+=\{(k^{st},0)\}\cup \tilde{{\cal Y}}^+$. 
 
  Si $\delta_{4}(q-1)=1$ posons ${\mathbb Y}^-=\emptyset$. Si $\delta_{4}(q+1)=1$, notons ${\mathbb Y}^{-}$ l'ensemble des couples $(i,j)$ avec $i,j\in {\mathbb N}$,  tels que $2i^2+j(j+1)/2=n$,
$i$ est pair   et $i/2+[(j+2)/4] $ est impair.

 On note ${\cal Y}^{-}$ l'ensemble des triplets $(i,j,\xi)$ o\`u $(i,j)\in {\mathbb Y}^{-}$ et $\xi\in \Xi$ v\'erifie:
 
 si $i\not=0$, la restriction de $\xi$ \`a   $\iota_{1}(\mathfrak{o}_{F}^{\times}/\mathfrak{o}_{F}^{\times,2})$ est non triviale;
 
 si $i=0$, $\xi\in \tilde{\Xi}$.  

On pose ${\cal Y}={\cal Y}^{+}\sqcup {\cal Y}^{-}$.

   La description des donn\'ees endoscopiques elliptiques ${\bf G}'$ telles que $FC^{st}(\mathfrak{g}'(F))$ est bien s\^ur la m\^eme que dans le paragraphe pr\'ec\'edent. Ce qui change est l'action du groupe d'automorphismes ext\'erieurs. En reprenant tous les cas un \`a un, on s'aper\c{c}oit que cette action est tordue par la restriction \`a $Out({\bf G}')\subset \hat{\Omega}$ du caract\`ere de ce groupe qui vaut $1$ sur $\theta\theta'$ mais $-1$ sur $\delta$. Ainsi, les donn\'ees telles que $Out({\bf G}')\subset \{1,\theta\theta'\}$ se conservent sans changement. Pour les donn\'ees telles que $Out({\bf G}')$ contient $\delta$ ou $\delta\theta\theta'$, les donn\'ees qui \'etaient admissibles dans le paragraphe pr\'ec\'edent ne le sont plus tandis que celles que l'on avait exclues \`a cause de l'action de ces \'el\'ements $\delta$ ou $\delta\theta\theta'$ deviennent admissibles. On obtient ainsi l'assertion \ref{resultats}(2).
   
   Pour $\star=\pm$, on d\'efinit une application $\phi^{\star}: {\mathbb X}^{\star}\to {\mathbb Y}^{\star}$  par les m\^emes formules qu'en \ref{Dndeppairpadique}.  La fin de la d\'emonstration est identique \`a celle de ce paragraphe. 
 
 \subsection{ Forme int\'erieure non classique du type $D_{n}$ d\'eploy\'e, $n$ pair}\label{Dndeppairnonclassique}
 On suppose que $G^*$ est comme en \ref{Dndeppairpadique}, en particulier  $n$ est pair. On suppose ici que $G$ est la forme int\'erieure de $G^*$ associ\'ee  \`a l'\'el\'ement $\delta\in \hat{\Omega}\simeq N$ ou encore au caract\`ere de $Z(\hat{G}_{SC})$ dont le noyau est $\{1,z'\}$. Alors $G$ est une forme non d\'eploy\'ee "non classique" de $Spin(2n)$, c'est-\`a-dire que ce n'est pas le groupe associ\'e \`a un espace muni d'une forme quadratique. Dans les tables de Tits, le groupe est de type $2D''_{n}$. 
 
 Notons $\Xi'$ le sous-ensemble de $\Xi$ form\'e 
 des deux  caract\`eres $(\chi_{E/F},\chi_{E/F}^{\delta_{4}(q+1)(\delta_{4}(n)+1)})$ o\`u $E/F$ est une extension quadratique ramifi\'ee de $F$.  Si $\delta_{\square}(n)=0$, posons ${\cal X}^+ =\emptyset$. Si $\delta_{\square}(n)=1$, c'est-\`a-dire $n=k^2$, pour un entier $k\in {\mathbb N}$, posons ${\cal X}^+=\{(k,k)\}$. Si $\delta_{\triangle}(n)=0$, posons ${\cal X}^-=\emptyset$. Si $\delta_{\triangle}(n)=1$, c'est-\`a-dire si $n=k(k+1)/2$ pour un entier $k\in {\mathbb N}$, posons ${\cal X}^-=\{(k,k,\xi); \xi\in \Xi'\}$. On pose ${\cal X}={\cal X}^+\sqcup {\cal X}^-$ et $d_{x}=1$ pour tout $x\in {\cal X}$.
  
 Le diagramme ${\cal D}_{a}^{nr}$ des tables de Tits est le diagramme ${\cal D}_{a}$ muni de l'action de $\Gamma_{{\mathbb F}_{q}}$ qui est triviale sur $\Gamma_{{\mathbb F}_{q^2}}$ et telle qu'un \'el\'ement $\sigma\in \Gamma_{{\mathbb F}_{q}}-\Gamma_{{\mathbb F}_{q^2}}$ agisse par $\delta$. Les \'el\'ements de $\underline{S}(G)$ correspondent aux orbites de cette action dans ${\cal D}_{a}$. Pour un sommet $s$ correspondant \`a une orbite \`a deux \'el\'ements,  le lemme \ref{orbites} montre que   $FC(\mathfrak{g}_{s}({\mathbb F}_{q}))=\{0\}$. Il n'y a qu'une orbite galoisienne poss\'edant un seul \'el\'ement, \`a savoir $\{\alpha_{n/2}\}$. Notons $s$ le sommet associ\'e. Son unicit\'e entra\^{\i}ne qu'il est conserv\'e par l'action de $G_{AD}(F)$.  Le groupe $G_{s}$ est alg\'ebriquement le m\^eme qu'en \ref{Dndeppairpadique}, c'est-\`a-dire $G_{s}=(Spin(n)\times Spin(n))/\{1,(z,z)\}$. Notons $\sigma\mapsto \sigma_{dep}$ l'action galoisienne d\'eploy\'ee sur $Spin(n)$. L'action galoisienne de $\sigma\in \Gamma_{{\mathbb F}_{q}}$ sur $G_{s}$ est $(g_{1},g_{2})\mapsto (\sigma_{dep}(g_{1}),\sigma_{dep}(g_{2}))$ pour $\sigma\in \Gamma_{{\mathbb F}_{q^2}}$ et $(g_{1},g_{2})\mapsto (\sigma_{dep}(g_{2}),\sigma_{dep}(g_{1}))$ pour $\sigma\not\in \Gamma_{{\mathbb F}_{q^2}}$.   Les faisceaux-caract\`eres cuspidaux sur $\mathfrak{g}_{s}$ sont les m\^emes qu'en \ref{Dndeppairpadique} puisque ces objets sont alg\'ebriques. Ce qui change est l'action galoisienne. Si $n=k^2$ pour un entier $k\in {\mathbb N}$, il y  a sur chaque facteur $\mathfrak{spin}(n)$ un faisceau-caract\`ere associ\'e \`a un unipotent $N$ et un caract\`ere $\epsilon$ tel que $\epsilon(z)=1$. Leur produit est clairement fixe par l'action galoisienne et donne naissance \`a un \'el\'ement de $FC(\mathfrak{g}(F))$. On note $FC_{k,k}$ la droite port\'ee par cet \'el\'ement (on a $(k,k)\in {\cal X}^+$). L'action de $G_{AD}(F)/\pi(G(F))$ n'est pas claire pour l'auteur. Mais l'action du sous-groupe $\iota_{1}(\mathfrak{o}_{F}^{\times}/\mathfrak{o}_{F}^{\times,2})$ est facile \`a d\'ecrire. C'est  l'action par conjugaison de $SO_{dep}(n,{\mathbb F}_{q^2})$. Cette action est triviale puisque les caract\`eres $\epsilon$ valent $1$ sur $z$. Si $n=k(k+1)/2$, il y  a sur chaque facteur $\mathfrak{spin}(n)$ deux faisceaux-caract\`eres associ\'es \`a un m\^eme unipotent $N$ et aux deux caract\`eres $\epsilon$ tel que $\epsilon(z)=-1$. Notons ${\cal F}_{1}'$ et ${\cal F}_{1}''$ les deux faisceaux-caract\`eres de la premi\`ere composante et ${\cal F}_{2}'$ et ${\cal F}_{2}''$ ceux de la seconde composante. Ils sont chacun fix\'es par $Fr^2$ car $4$ divise $q^2-1$. L'action du Frobenius envoie par exemple ${\cal F}_{1}'$ sur ${\cal F}_{2}'$ ou ${\cal F}_{2}''$ selon que $4$ divise $q-1$ ou $q+1$. En tout cas, parmi les quatre produits ${\cal F}'_{1}\otimes {\cal F}'_{2}$, etc... il y a deux faisceaux-caract\`eres invariants par l'action galoisienne tandis que les deux autres sont permut\'es. Donc $dim(FC(\mathfrak{g}_{s}({\mathbb F}_{q})))=2$. Cette fois, le groupe $\iota_{1}(\mathfrak{o}_{F}^{\times}/\mathfrak{o}_{F}^{\times,2})$ agit  sur chaque fonction par le caract\`ere non trivial puisque $\epsilon(z)=-1$. D'apr\`es \ref{actionsurFC}, ce caract\`ere se prolonge en deux caract\`eres $\xi$  de $G_{AD}(F)/\pi(G(F))$ de sorte que, pour chacun d'eux, une combinaison lin\'eaire des \'el\'ements de $ FC(\mathfrak{g}_{s}({\mathbb F}_{q}))$ donne naissance \`a un \'el\'ement de $FC(\mathfrak{g}(F))$ se transformant par $\xi$. On note $FC_{k,k,\xi}$ la droite port\'ee par cet \'el\'ement. On obtient une assertion similaire  \`a \ref{resultats}(1):   puisqu'on n'a pas montr\'e que les deux caract\`eres $\xi$ pr\'ec\'edents \'etaient les \'el\'ements de $\Xi'$, on obtient un param\'etrage par un ensemble $\tilde{{\cal X}}={\cal X}^+\cup\tilde{{\cal X}}^-$ d\'efini en rempla\c{c}ant $\Xi'$ par l'ensemble de ces deux caract\`eres dans la d\'efinition de ${\cal X}$. 
 
 Si $\delta_{\square}(n)=0$, posons ${\cal Y}^+ =\emptyset$. Si $\delta_{\square}(n)=1$, c'est-\`a-dire $n=k^2$, pour un entier $k\in {\mathbb N}$, posons ${\cal Y}^+=\{(k,0)\}$. Si $\delta_{\triangle}(n)=0$, posons ${\cal Y}^-=\emptyset$. Si $\delta_{\triangle}(n)=1$, c'est-\`a-dire si $n=j(j+1)/2$ pour un entier $j\in {\mathbb N}$, posons ${\cal Y}^-=\{(0,j,\xi); \xi\in \Xi'\}$. On pose ${\cal Y}={\cal Y}^+\sqcup {\cal Y}^-$.

 Les donn\'ees endoscopiques elliptiques ${\bf G}'$ de $G$ telles que $FC^{st}(\mathfrak{g}'(F))\not=\{0\}$ sont  les m\^emes qu'en \ref{Dndeppairpadique}. Ce qui change est l'action du groupe d'automorphismes ext\'erieurs. En reprenant tous les cas un \`a un, on s'aper\c{c}oit que pour toutes les donn\'ees telles que $Out({\bf G}')$ contient $\theta\theta'$, ce groupe agit maintenant non trivialement sur la droite $FC^{st}(\mathfrak{g}'(F))$. Il  y a au plus  $5$ donn\'ees telles que $\theta\theta'\not\in Out({\bf G}')$. D'abord la donn\'ee principale ${\bf G}$. On a vu que son espace $FC^{st}(\mathfrak{g}^*(F))$ \'etait non nul si et seulement si $n=k^2$ pour un entier $k\in {\mathbb N}$. Si cette condition est v\'erifi\'ee, l'espace pr\'ec\'edent est une droite que l'on note $FC^{{\cal E}}_{k,0}$ (on a $(k,0)\in {\cal Y}^+$). Les quatre autres donn\'ees sont d\'etermin\'ees par une extension quadratique $E/F$ ramifi\'ee, par l'action $\sigma\mapsto \sigma_{G'}$ triviale sur $\Gamma_{E}$ et telle que $\sigma_{G'}=\delta$, resp. $\delta\theta\theta'$, pour $\sigma\in \Gamma_{F}-\Gamma_{E}$ et par l'orbite ${\cal O}=\{\hat{\alpha}_{0},\hat{\alpha}_{n}\}$, resp. ${\cal O}=\{\hat{\alpha}_{0},\hat{\alpha}_{n-1}\}$. Elles n'interviennent que si $n= j(j+1)/2$ pour un entier $j\in {\mathbb N}$. Supposons cette condition v\'erifi\'ee. 
  
 Supposons $\sigma_{G'}=\delta$ pour $\sigma\in \Gamma_{F}-\Gamma_{E}$ et ${\cal O}=\{\hat{\alpha}_{0},\hat{\alpha}_{n}\}$. Alors $Out({\bf G}')=\{1,\delta\}$. On a vu en \ref{Dndeppairpadique} que l'action naturelle de $\delta$ sur $FC^{st}(\mathfrak{g}'(F))$ \'etait la multiplication par $sgn(-1)^{[(j+2)/4]}$. L'\'egalit\'e $n= j(j+1)/2$ et la parit\'e de $n$ entra\^{\i}nent  que  $[(j+2)/4] $ est de m\^eme parit\'e que $n/2$, donc $sgn(-1)^{[(j+2)/4]}=sgn(-1)^{n/2}$. On a calcul\'e $\delta(s_{sc})s_{sc}^{-1}=z'z^{n/2}$. Puisque le caract\`ere de $Z(\hat{G}_{SC})$ associ\'e \`a $G$ est trivial sur $z'$ mais pas sur $z$, l'action de $\delta$ sur le facteur de transfert est la multiplication par $(-1)^{n/2}$. Donc l'action de $\delta$ sur $FC^{st}(\mathfrak{g}'(F))$ (qui est l'action naturelle multipli\'ee par la constante pr\'ec\'edente) est la multiplication par $(-sgn(-1))^{n/2}$. Elle est triviale si $sgn(-1)=-1$, c'est-\`a-dire $\delta_{4}(q+1)=1$, ou si $4$ divise $n$. Supposons ces conditions v\'erifi\'ees. Le caract\`ere $\xi_{{\bf G}'}$ est le m\^eme qu'en \ref{Dndeppairpadique}, c'est-\`a-dire le caract\`ere $\xi=(sgn_{E/F},sgn_{E/F}^{n/2})$. On a   $(0,j,\xi)\in {\cal Y}^-$. On pose ${\bf G}'_{0,j,\xi}={\bf G}'$ et $FC^{{\cal E}}=FC^{st}(\mathfrak{g}'_{0,j,\xi}(F))^{Out({\bf G}'_{0,j,\xi})}$.

 Supposons maintenant que $\sigma_{G'}=\delta\theta\theta'$ pour $\sigma\in \Gamma_{F}-\Gamma_{E}$ et ${\cal O}=\{\hat{\alpha}_{0},\hat{\alpha}_{n-1}\}$. Alors $Out({\bf G}')=\{1,\delta\theta\theta'\}$. L'action naturelle de $\delta\theta\theta'$ sur $FC^{st}(\mathfrak{g}'(F))$ est encore la multiplication par $sgn(-1)^{n/2}$. On a calcul\'e $\delta\theta\theta'(s_{sc})s_{sc}^{-1}=z'z^{1+n/2}$. Donc  l'action de $\delta\theta\theta'$ sur $FC^{st}(\mathfrak{g}'(F))$ est la multiplication par $(-1)^{n/2+1}sgn(-1)^{n/2}$. Elle est triviale si $n/2$ est impair et $sgn(-1)=1$, c'est-\`a-dire $\delta_{4}(q-1)=1$. Ce sont les cas oppos\'es a ceux ci-dessus. Supposons ces conditions v\'erifi\'ees. Le caract\`ere $\xi_{{\bf G}'}$ est maintenant $\xi=(sgn_{E/F},sgn_{E/F}^{n/2+1})=(sgn_{E/F},1)$.  On a   $(0,j,\xi)\in {\cal Y}^-$. On pose ${\bf G}'_{0,j,\xi}={\bf G}'$ et $FC^{{\cal E}}=FC^{st}(\mathfrak{g}'_{0,j,\xi}(F))^{Out({\bf G}'_{0,j,\xi})}$. 
 
 Cette description d\'emontre l'assertion \ref{resultats}(2).

 Les droites $FC_{x}$ pour $x\in \tilde{{\cal X}}$ 
  se distinguent par le caract\`ere de $G_{AD}(F)/\pi(G(F))$ par lequel agit ce groupe: le groupe $\iota_{1}(\mathfrak{o}_{F}^{\times}/\mathfrak{o}_{F}^{\times,2})$ agit trivialement sur  $FC_{x}$ pour $x\in {\cal X}^+$ et non trivialement pour $x\in {\cal X}^-$; le groupe $G_{AD}(F)/\pi(G(F))$ agit par deux caract\`eres distincts sur  les deux droites  $FC_{x}$ pour $x\in \tilde{{\cal X}}^-$ (quand elles existent). De m\^eme, les caract\`eres $\xi_{{\bf G}'_{y}}$ pour $y\in {\cal Y}$ sont tous distincts. Puisque le transfert est compatible \`a l'action de $G_{AD}(F)/\pi(G(F))$, cela entra\^{\i}ne d'abord que les deux caract\`eres intervenant dans la  d\'efinition de $\tilde{{\cal X}}^-$ sont bien les deux \'el\'ements de $\Xi'$, ce qui ach\`eve la preuve de \ref{resultats}(1). Il y a alors une bijection \'evidente $\varphi:{\cal X}\to {\cal Y}$. 
 Pour $x\in {\cal X}$, on a forc\'ement $transfert(FC_{x})=FC^{{\cal E}}_{\varphi(x)}$ car $\varphi(x)$ est l'unique \'el\'ement  $y\in{\cal Y}$ tel que $\xi_{{\bf  G}'_{y}}$ puisse co\"{\i}ncider avec le caract\`ere par lequel $G_{AD}(F)/\pi(G(F))$ agit sur $FC_{x}$. Cela d\'emontre \ref{resultats}(3).

 {\bf Remarque.} On peut aussi consid\'erer la forme int\'erieure $G$ de $G^*$ associ\'ee  \`a l'\'el\'ement $\delta\theta\theta'\in \hat{\Omega}\simeq N$. Ce groupe est isomorphe au pr\'ec\'edent, on a simplement compos\'e le torseur int\'erieur par l'automorphisme $\theta$ de $G^*$. On obtient \'evidemment des r\'esultats similaires. De fa\c{c}on naturelle, l'ensemble $\Xi'$ est remplac\'e par l'ensemble $\Xi''$ des caract\`eres $(sgn_{E/F},sgn_{E/F}^{1+\delta_{4}(q+1)(\delta_{4}(n)+1)})$.

\subsection{Type $D_{n}$ d\'eploy\'e, $n$ impair}\label{Dndepimppadique}
 On suppose que $G$ est d\'eploy\'e de type $D_{n}$ avec $n\geq4$, c'est-\`a-dire $G=Spin_{dep}(2n)$, et on suppose $n$ impair. Comme en \ref{Dndepimp}, on a $Z(G)=\{1,z',z,z''\}\simeq {\mathbb Z}/4{\mathbb Z}$, avec $(z')^2=z$. D\'efinissons un homomorphisme $T_{ad}(F)\to F^{\times}/F^{\times,4}$ par $\prod_{i=1,...,n}\check{\varpi}_{i}(x_{i})\mapsto (\prod_{i=1,...,n-2}x_{i}^2) x_{n-1}x_{n}^{-1}$.  Il s'en d\'eduit un isomorphisme $G_{AD}(F)/\pi(G(F))\simeq T_{ad}(F)/\pi(T(F))\simeq  F^{\times}/F^{\times,4}$. L'image de $SO_{2n}(F)$ dans ce groupe est $F^{\times,2}/F^{\times,4}$. L'image de $G_{AD}(F)_{0}$ est $\mathfrak{o}_{F}^{\times}/\mathfrak{o}_{F}^{\times,4}$. Remarquons que, dans le cas o\`u $\delta_{4}(q-1)=1$ ce dernier groupe est cyclique d'ordre $4$. 
 
 Si $\delta_{4}(q-1)=0$, on pose ${\cal X}=\emptyset$. 
 Supposons $\delta_{4}(q-1)=1$. On note ${\mathbb X}$ l'ensemble des couples $(k,h)\in {\mathbb N}^2$ tels que $k(k+1)/2+h(h+1)/2=2n$, $k>h$ et $k(k+1)/2$ et $h(h+1)/2$ sont pairs. On note ${\cal X}$ l'ensemble des triplets $(k,h,\xi)$ o\`u $(k,h)\in {\mathbb X}$ et $\xi$ est un \'el\'ement de $\Xi$ dont la restriction \`a  $\mathfrak{o}_{F}^{\times}/\mathfrak{o}_{F}^{\times,4}$ est d'ordre $4$.  On pose $d_{x}=1$ pour tout $x\in {\cal X}$.

 L'ensemble $\underline{S}(G)$ s'envoie surjectivement sur l'ensemble des couples $(a,b)\in {\mathbb N}^2$ tels que $a>b$, $a+b=n$ et $b\not=1$. Remarquons que les deux nombres $a$ et $b$ sont de parit\'e distincte. Les fibres de cette surjection ont deux \'el\'ements sauf au-dessus de $(n,0)$ o\`u la fibre a $4$ \'el\'ements. L'action de $G_{AD}(F)$ pr\'eserve les fibres et est transitive sur chacune d'elles. Pour un sommet $s$ param\'etr\'e par $(a,b)$, on a $G_{s}=(Spin_{dep}(2a)\times Spin_{dep}(2b))/\{1,(z,z)\}$ si $b\not=0$, $G_{s}=Spin_{dep}(2n)$ si $(a,b)=(n,0)$. Consid\'erons un sommet $s$ param\'etr\'e par $(a,b)$ avec $b\geq2$. On doit d\'eterminer les couples de fonctions  $f_{N^{a},\epsilon^{a}}\times f_{N^b,\epsilon^b}$ tels que $\epsilon^{a}(z)=\epsilon^b(z)$. Il n'y en a pas pour lesquels $\epsilon^{a}(z)=\epsilon^b(z)=1$. En effet, d'apr\`es \ref{Dndeppair} et \ref{Dndepimp}, un tel couple n'existe que si $(2a,2b)$ est de la forme $(k^2,h^2)$. Or cette condition implique que $a$ et $b$ sont tous deux pairs et on a d\'ej\`a dit que ce n'\'etait pas possible. Consid\'erons les couples pour lesquels $\epsilon^{a}(z)=\epsilon^b(z)=-1$. Pour fixer les id\'ees, supposons $a$ impair et $b$ pair. Il y a  une telle fonction $f_{N^b,\epsilon^b}$ si et seulement si $2b$ est de la forme $h(h+1)/2$ et, dans ce cas, il y en a deux. Il y a une telle fonction $f_{N^{a},\epsilon^{a}}$ si et seulement si $\delta_{4}(q-1)=1$ et $2a$ est de la forme $k(k+1)/2$. Dans ce cas, il y en a deux. Supposons ces conditions v\'erifi\'ees. On obtient alors $4$ g\'en\'erateurs de $FC(\mathfrak{g}_{s}({\mathbb F}_{q}))$.   Le stabilisateur de $s$ dans $G_{AD}(F)/\pi(G(F))$ est l'image r\'eciproque de $2{\mathbb Z}/4{\mathbb Z}$ par l'homomorphisme $F^{\times}/F^{\times,4}\to {\mathbb Z}/4{\mathbb Z}$ issu de la valuation. L'action du sous-groupe
  $ \mathfrak{o}_{F}^{\times,2}/\mathfrak{o}_{F}^{\times,4}$  de $G_{AD}(F)_{0}/\pi(G(F))$ est facile \`a d\'eterminer: c'est l'action naturelle de $\{(x,y)\in O_{dep}(2a,{\mathbb F}_{q})\times O_{dep}(2b,{\mathbb F}_{q}); det(x)=det(y)\}$. L'action du sous-groupe $SO_{dep}(2a,{\mathbb F}_{q})\times SO_{dep}(2b,{\mathbb F}_{q})$  est non triviale puisque  $\epsilon^{a}(z)=\epsilon^b(z)=-1$. Donc les caract\`eres de $G_{AD}(F)_{0}$ par lesquels se transforment nos fonctions sont d'ordre $4$. Il est clair que la situation est invariante par conjugaison complexe. On obtient donc que les deux caract\`eres d'ordre $4$ de $G_{AD}(F)_{0}/\pi(G(F))$ interviennent et que, pour chacun de ces caract\`eres, il y a deux g\'en\'erateurs de $FC(\mathfrak{g}_{s}({\mathbb F}_{q}))$ qui se transforment selon ce caract\`ere.  Comme en \ref{Dndeppairpadique}, l'action  d'un \'el\'ement  $(x,y)\in O_{dep}(2a,{\mathbb F}_{q})\times O_{dep}(2b,{\mathbb F}_{q})$ tel que $det(x)=det(y)=-1$ permute les deux g\'en\'erateurs associ\'es \`a chaque caract\`ere de $G_{AD}(F)_{0}$. Pour tout caract\`ere $\xi_{0}$ du stabilisateur de $s$ dont la restriction \`a $\mathfrak{o}_{F}^{\times}/\mathfrak{o}_{F}^{\times,4}$ est d'ordre $4$, il y a donc une combinaison lin\'eaire de nos fonctions qui se transforme selon ce caract\`ere.   Ensuite, chacune de ces fonctions donne naissance \`a deux \'el\'ements de $FC(\mathfrak{g}(F))$  se transformant selon deux caract\`eres $\xi$ de $F^{\times}/F^{\times,4}$ prolongeant $\xi_{0}$. On note $FC_{k,h,\xi}$ la droite port\'ee par la fonction correspondant \`a $\xi$. On a $(k,h,\xi)\in {\cal X}$.  Consid\'erons maintenant un sommet $s$ param\'etr\'e par $(n,0)$. On a $G_{s}=Spin_{dep}(2n)$. De nouveau, parce que $n$ est impair, il n'y a pas de fonction $f_{N,\epsilon}$ telle que $\epsilon(z)=1$. Il y a une telle fonction avec $\epsilon(z)=-1$ si et seulement si $\delta_{4}(q-1)=1$ et $2n=k(k+1)/2$ pour un entier $k\in {\mathbb N}$. Supposons ces conditions v\'erifi\'ees. Il y a alors deux telles fonctions. De nouveau, elles se transforment selon les deux caract\`eres d'ordre $4$ de $\mathfrak{o}_{F}^{\times}/\mathfrak{o}_{F}^{\times,4}$. Maintenant, ce groupe est le stabilisateur de $s$ dans $G_{AD}(F)/\pi(G(F))$ et, des deux fonctions pr\'ec\'edentes se d\'eduisent $8$ \'el\'ements de $FC(\mathfrak{g}(F))$ qui se transforment selon les caract\`eres de $G_{AD}(F)/\pi(G(F))$  qui prolongent les deux caract\`eres pr\'ec\'edents de $\mathfrak{o}_{F}^{\times}/\mathfrak{o}_{F}^{\times,4}$. On note $FC_{k,0,\xi}$ la droite port\'ee par la fonction correspondant \`a $\xi$. On a $(k,0,\xi)\in {\cal X}$.

   Si $\delta_{4}(q-1)=0$, on pose ${\cal Y}=\emptyset$. Les assertions \ref{resultats}(2) et (3) sont triviales puisqu'on a d\'ej\`a prouv\'e que $FC(\mathfrak{g}(F))=\{0\}$. On suppose d\'esormais que $\delta_{4}(q-1)=1$. 
On note ${\mathbb Y}$ l'ensemble des couples $(i,j)$ o\`u $i,j\in {\mathbb N}$,  tels que $2i^2+j(j+1)/2=n$ et $i$  est  impair. On note ${\cal Y}$ l'ensemble des triplets $(i,j,\xi)$ o\`u $(i,j)\in {\cal Y}$ et $\xi$ est un \'el\'ement de $\Xi$ dont la restriction \`a  $\mathfrak{o}_{F}^{\times}/\mathfrak{o}_{F}^{\times,4}$ est d'ordre $4$.

 Commen\c{c}ons par une remarque concernant l'identification des caract\`eres $\xi_{{\bf G}'}$.  Identifions $Z(\hat{G}_{SC})$ \`a $\boldsymbol{\zeta}_{4}({\mathbb C})$ par $z'\mapsto i$ (il s'agit ici de l'\'el\'ement $z'\in Z(\hat{G}_{SC})$). Un \'el\'ement de $H^1(W_{F},Z(\hat{G}_{SC}))$ devient un caract\`ere d'ordre $4$ de $W_{F}$ qui, par l'isomorphisme du corps de classes, s'identifie \`a un caract\`ere de $F^{\times}/F^{\times,4}$. Une donn\'ee endoscopique elliptique ${\bf G}'$ de $G$ d\'etermine un \'el\'ement de $H^1(W_{F},Z(\hat{G}_{SC}))$ donc un caract\`ere de $F^{\times}/F^{\times,4}\simeq G_{AD}(F)/\pi(G(F))$. Celui-ci n'est autre que $\xi_{{\bf G}'}$. 
 
 Consid\'erons un couple $(\sigma\mapsto \sigma_{G'},{\cal O})\in {\cal E}_{ell}(G)$. Rappelons que, $n$ \'etant impair,  $\hat{\Omega}=\{1,\delta\theta,\theta\theta',\delta\theta'\}$ est cyclique d'ordre $4$, engendr\'e par $\delta\theta$.

   Supposons d'abord  que l'action $\sigma\mapsto \sigma_{G'}$ soit triviale. L'orbite ${\cal O}$ est r\'eduite \`a une seule racine. Si celle-ci est une extr\'emit\'e du diagramme, on a ${\bf G}'={\bf G}$ et on ne peut \`a pr\'esent rien dire de l'espace $FC^{st}(\mathfrak{g}(F))$. Consid\'erons le cas o\`u ${\cal O}=\{\hat{\alpha}_{m}\}$ avec $m\in \{2,...,n-2\}$. On a alors $G'_{SC}\simeq Spin_{dep}(2m)\times Spin_{dep}(2n-2m)$. Les deux nombres $m$ et $n-m$ sont de parit\'e distincte. Supposons par exemple que $m$ est impair. On applique (1) ci-dessous   par r\'ecurrence: l'espace $FC^{st}(\mathfrak{spin}_{dep}(2m,F))$ est nul. Donc $FC^{st}(\mathfrak{g}'(F))=\{0\}$.
 
   Supposons qu'il existe une extension quadratique $E/F$ telle que $\sigma_{G'}=1$ pour $\sigma\in \Gamma_{E}$ et $\sigma_{G'}=\theta\theta'$ pour $\sigma\in \Gamma_{F}-\Gamma_{E}$. L'orbite ${\cal O}$ peut avoir un seul \'el\'ement $\hat{\alpha}_{m}$ pour $m\in \{2,...,n-2\}$. On a alors $G'_{SC}\simeq Spin_{E/F}(2m)\times Spin_{E/F}(2n-2m)$. On applique \ref{resultats}(4) par r\'ecurrence.  L'espace $FC^{st}(\mathfrak{g}'(F))$ n'est non nul que si $(m,n-m)$ est de la forme $(i^2,j^2)$ avec $i,j$ pairs si $E/F$ est non ramifi\'ee, $i,j$ impairs si $E/F$ est ramifi\'ee. Ces conditions ne peuvent pas \^etre v\'erifi\'ees puisque $m$ et $n-m$ sont de parit\'e distincte. L'orbite ${\cal O}$ peut aussi avoir deux \'el\'ements: $\{\hat{\alpha}_{0},\hat{\alpha}_{1}\}$ ou $\{\hat{\alpha}_{n-1},\hat{\alpha}_{n}\}$. Ces deux cas sont d'ailleurs \'equivalents (par conjugaison par $\delta\theta\in \hat{\Omega}$). Le stabilisateur d'un sommet de l'orbite est $\Gamma_{E}$. Le lemme \ref{centre} exclut le cas o\`u $E/F$ est non ramifi\'ee. Supposons $E/F$ ramifi\'ee. Alors $G'_{SC}\simeq Spin_{E/F}(2n-2)$.  L'espace $FC^{st}(\mathfrak{g}'(F))$ n'est non nul que si $n-1=j^2$ avec $j$ impair, ce qui est impossible puisque $n-1$ est pair. 
 
 Supposons que l'image de $\Gamma_{F}$ par l'homomorphisme $\sigma\mapsto \omega_{G'}$ soit $\hat{\Omega}$ tout entier. L'action se factorise alors par une extension galoisienne $E_{G'}/F$ cyclique d'ordre $4$. On fixe un g\'en\'erateur $\rho$ de $\Gamma_{E_{G'}/F}$. On peut avoir $\rho_{G'}=\delta\theta$ ou $\rho_{G'}=\delta\theta'$. Supposons $\rho_{G'}=\delta\theta$. On note $E$ l'extension quadratique de $F$ contenue dans $E_{G'}$, c'est-\`a-dire telle que $\Gamma_{E_{G'}/E}=\{1,\rho^2\}$. On peut avoir ${\cal O}=\{\hat{\alpha}_{m},\hat{\alpha}_{n-m}\}$ avec $m\in \{2,...,(n-1)/2\}$. Le fixateur de $\hat{\alpha}_{m}$ dans $\Gamma_{F}$ est $\Gamma_{E}$. D'apr\`es le lemme \ref{centre}, on peut supposer que $E/F$ est ramifi\'ee. Il en r\'esulte que $E_{G'}/F$ est totalement ramifi\'ee (si $E_{G'}$ contient l'extension quadratique $E_{0}/F$ non ramifi\'ee, l'homomorphisme naturel $\Gamma_{E_{G'}/F}\to \Gamma_{E/F}\times \Gamma_{E_{0}/F}$ est un  isomorphisme et $\Gamma_{E_{G'}/F}$ n'est pas cyclique). Cela force $\delta_{4}(q-1)=1$ et il y a $4$ extensions $E_{G'}/F$ possibles. On voit que $G'_{SC}\simeq Res_{E/F}(Spin_{E_{G'}/E}(2m))\times SU_{E/F}(n-2m)$. L'extension $E_{G'}/E$ est ramifi\'ee. En utilisant \ref{Dnpairpadiqueram} (1), \ref{Dnimppadiqueram} (1) et \ref{An-1quasidepram} (5) par r\'ecurrence, on obtient que  $FC^{st}(\mathfrak{g}'(F))$ est non nul si et seulement si $m=i^2$ avec $i$ impair et $n-2m=j(j+1)/2$.   Le groupe $Out({\bf G}')$ est $\hat{\Omega}$ tout entier. L'action de $\delta\theta$ est l'action naturelle de $\rho$ sur chacun des facteurs. On montrera en \ref{variance} que c'est l'identit\'e sur le premier facteur et on a vu en \ref{actiondunautomorphisme} que c'\'etait la multiplication par $sgn(-1)^{[(j+2)/4]}$ sur le second. Mais $\delta_{4}(q-1)=1$ donc $sgn(-1)=1$ et l'action est triviale.  Donc $FC^{st}(\mathfrak{g}'(F))^{Out({\bf G}')}$ est une droite. Calculons le caract\`ere $\xi_{{\bf G}'}$. On peut choisir $s_{sc}$ comme en \ref{Dndeppairpadique}(2). 
  On calcule $\delta\theta(s_{sc})s_{sc}^{-1}=(z')^n$.  Le cocycle de $W_{F}$ dans $Z(\hat{G}_{SC})$ qui  est trivial sur $W_{E_{G'}}$ et envoie $\rho$ sur $z'$ correspond \`a un caract\`ere $\chi_{E_{G'}/F}$ de $F^{\times}$ dont le noyau est le groupe des normes de l'extension $E_{G'}/F$. Alors $\xi_{{\bf G}'}=\chi_{E_{G'}/F}^n$. Notons $\xi$ ce caract\`ere. On a $(i,j,\xi)\in {\cal Y}$. On pose ${\bf G}'_{i,j,\xi}={\bf G}'$ et $FC^{{\cal E}}_{i,j,\xi}= FC^{st}(\mathfrak{g}_{i,j,\xi}'(F))^{Out({\bf G}_{i,j,\xi}')}$.  On peut aussi avoir ${\cal O}=\{\hat{\alpha}_{0},\hat{\alpha}_{1},\hat{\alpha}_{n-1},\hat{\alpha}_{n}\}$. Le fixateur de $\hat{\alpha}_{0}$ dans $\Gamma_{F}$ est $\Gamma_{E_{G'}}$. Le lemme \ref{centre} permet de supposer que $E_{G'}/F$ est totalement ramifi\'ee, donc encore $\delta_{4}(q-1)=1$. On a $G'_{SC}\simeq SU_{E/F}(n-2)$, o\`u $E$ est comme ci-dessus. La suite de la discussion est la m\^eme, en consid\'erant que $i=1$ (cas qui \'etait exclu pr\'ec\'edemment puisqu'on avait $i^2=m\geq2$).
 
 Le cas o\`u $\rho_{G'}=\delta\theta'$ est similaire. Le seul changement est le calcul de $\xi_{{\bf G}'}$. L'\'el\'ement $s_{sc}$ est le m\^eme que ci-dessus mais on doit calculer $\delta\theta'(s_{sc})s_{sc}^{-1}$. Ce terme vaut $(z')^{-n}$ et $\xi_{{\bf G}'}=\chi_{E_{G'}/F}^{-n}$.  Modulo ce changement, on pose les m\^emes d\'efinitions que ci-dessus. Remarquons que la construction associe \`a chaque extension $E_{G'}$, que l'on a munie d'un g\'en\'erateur $\rho$ de $\Gamma_{E_{G'}/F}$, un unique caract\`ere $\chi_{E_{G'}/F}$. Quand $E_{G'}/F$ d\'ecrit les $4$ extensions cycliques ramifi\'ees de  degr\'e $4$ de $F$, les caract\`eres $\chi_{E_{G'}/F}^{\pm n}$ d\'ecrivent tous les caract\`eres de $F^{\times}/F^{\times,4}$ dont les restrictions \`a $\mathfrak{o}_{F}^{\times}/\mathfrak{o}_{F}^{\times,4}$ sont d'ordre $4$. 
 
 A ce point, on a associ\'e \`a tout \'el\'ement $y\in {\cal Y}$ une donn\'ee endoscopique elliptique ${\bf G}'_{y}\not={\bf G}$ et une droite $FC^{{\cal E}}_{y}$.  On n'a pas trait\'e le cas de la donn\'ee principale ${\bf }$.  On d\'efinit une bijection $\phi:{\mathbb X}\to {\mathbb Y}$ par la m\^eme formule qu'en \ref{Dndeppairpadique}. Elle se rel\`eve de fa\c{c}on \'evidente en une bijection $\varphi:{\cal X}\to {\cal Y}$. On applique alors l'argument de comparaison des dimensions. Il implique que l'espace $FC^{st}(\mathfrak{g}(F))$ que l'on n'a pas encore d\'etermin\'e est nul. Cela ach\`eve la preuve de \ref{resultats}(2). En posant ${\cal X}^{st}=\emptyset$, cela d\'emontre en m\^eme temps \ref{resultats}(4).

 On munit ${\mathbb X}$ de la relation $(k',h')\leq (k,h)$ si et seulement si $k'+h'\leq k+h$. On montre que c'est une relation d'ordre et on la rel\`eve en une relation de pr\'eordre sur ${\cal X}$. Pour tout $y\in {\cal Y}$, on introduira en \ref{elementsYy} un \'el\'ement $Y_{y}\in \mathfrak{g}'_{y,ell}(F)$ de sorte que les propri\'et\'es \ref{An-1quasidepram}(2) et (3) soient v\'erifi\'ees. L'assertion \ref{resultats}(3) s'en d\'eduit. Explicitons la cons\'equence de \ref{resultats}(4):
 
 (1) $FC^{st}(\mathfrak{g}(F))=\{0\}$. 
 
 \subsection{Forme int\'erieure classique du type $D_{n}$ d\'eploy\'e, $n$ impair}
 On suppose que $G^*$ est du type pr\'ec\'edent, en particulier $n$ est impair. On suppose que $G$ est la forme int\'erieure de $G^*$ associ\'ee \`a l'\'el\'ement $\theta\theta'\in \hat{\Omega}\simeq N$, ou encore au caract\`ere de $Z(\hat{G}_{SC})$ dont le noyau est $\{1,z\}$. Alors $G$ est le groupe $Spin(2n)$ associ\'e \`a un espace de dimension $2n$ sur $F$ muni d'une forme quadratique de d\'eterminant $-1$ (c'est-\`a-dire de d\'eterminant normalis\'e $1$ puisque $n$ est impair) telle que les sous-espaces totalement isotropes maximaux sont de dimension $n-2$. Dans les tables de Tits, le groupe est de type $^2D'_{n}$.
 
 On pose ${\cal X}=\emptyset$.   
 L'ensemble $\underline{S}(G)$ s'envoie surjectivement sur l'ensemble des couples $(a,b)\in {\mathbb N}^2$ tels que $a>b\geq1$ et $a+b=n$. Les fibres de cette surjection ont deux \'el\'ements.  L'action de $G_{AD}(F)$ pr\'eserve les fibres et permute leurs deux \'el\'ements. Pour un sommet $s$ param\'etr\'e par $(a,b)$, avec $b\geq2$,   on a $G_{s}=(Spin_{{\mathbb F}_{q^2}/{\mathbb F}_{q}}(2a)\times  Spin_{{\mathbb F}_{q^2}/{\mathbb F}_{q}}(2b))/\{1,(z,z)\}$. Supposons par exemple $a$ pair donc $b$ impair. D'apr\`es \ref{Dnnondeppair}, il n'y a pas de fonction $f_{N^{a},\epsilon^{a}}$ sur $\mathfrak{spin}(2a,{\mathbb F}_{q})$ avec $\epsilon^{a}(z)=-1$. D'apr\`es \ref{Dnnondepimp}, il n'y a pas de fonction $f_{N^b,\epsilon^b}$ sur $\mathfrak{spin}(2b,{\mathbb F}_{q})$ avec $\epsilon^{b}(z)=1$. Donc il n'y a pas de 
 fonctions de la forme $f_{N^{a},\epsilon^{a}}\times f_{N^b,\epsilon^b}$  sur $\mathfrak{g}_{s}({\mathbb F}_{q})$ telles que $\epsilon^{a}(z)=\epsilon^b(z)$. Donc $FC(\mathfrak{g}_{s}({\mathbb F}_{q}))=\{0\}$. Pour un sommet $s$ param\'etr\'e par $(n-1,1)$, le groupe  $Spin_{{\mathbb F}_{q^2}/{\mathbb F}_{q}}(2b)$ ci-dessus devient $GL(1)$ muni de l'action non triviale de $\Gamma_{{\mathbb F}_{q^2}/{\mathbb F}_{q}}$. Donc $Z(G_{s})^0\not=\{1\}$ et $FC(\mathfrak{g}_{s}({\mathbb F}_{q}))=\{0\}$.  Cela prouve \ref{resultats}(1). 
 
 En posant ${\cal Y}=\emptyset$, les assertions \ref{resultats}(2) et (3) s'ensuivent.

 \subsection{Forme int\'erieure non classique du type $D_{n}$ d\'eploy\'e, $n$ impair}
 On conserve le m\^eme groupe $G^*$ et on suppose que $G$ est la forme int\'erieure de $G^*$ associ\'ee \`a l'\'el\'ement $\delta\theta\in \hat{\Omega}\simeq N$, ou encore \`a un certain caract\`ere d'ordre $4$ de $Z(\hat{G}_{SC})$. Dans les tables de Tits, le groupe est de type $^4D_{n}$. 
 
 On pose ${\cal X} =\emptyset$.   
 
 Comme en \ref{Dndeppairnonclassique}, on consid\`ere le diagramme ${\cal D}_{a}$ muni  cette fois de l'action de $\Gamma_{{\mathbb F}_{q}}$ qui est triviale sur $\Gamma_{{\mathbb F}_{q^4}}$ et telle que le Frobenius agisse par $\delta\theta$. 
 Les \'el\'ements $\underline{S}(G)$ sont param\'etr\'es par les orbites de cette action. Toute orbite a au moins $2$ \'el\'ements.  Le lemme \ref{orbites} montre que, pour tout $s\in \underline{S}(G)$, on a  $FC(\mathfrak{g}_{s}({\mathbb F}_{q}))=\{0\}$. Donc $FC(\mathfrak{g}(F))=\{0\}$ et \ref{resultats}(1) est v\'erifi\'ee.
 
En posant ${\cal Y}=\emptyset$, les assertions \ref{resultats}(2) et (3) s'ensuivent.

 \subsection{Type $D_{n}$ quasi-d\'eploy\'e, $n$ pair, $E_{0}/F$ non ramifi\'ee}\label{Dnpairpadiquenonram}
 On suppose que $G$ est quasi-d\'eploy\'e de type $D_{n}$ avec $n\geq4$, $n$ pair, et que $\Gamma_{F}$ agit sur le diagramme ${\cal D}$ de $G$ par l'action $\sigma\mapsto \sigma_{G}$ triviale sur $\Gamma_{E_{0}}$ et telle que $\sigma_{G}=\theta$ pour $\sigma\in \Gamma_{F}-\Gamma_{E_{0}}$ (on rappelle que $E_{0}/F$ est l'extension quadratique non ramifi\'ee). On note $\tau$ un \'el\'ement de Frobenius de $\Gamma_{F}$, qui agit donc non trivialement sur $E_{0}$.  L'action de $\tau$ sur $Z(G)$ fixe $z$ et \'echange $z'$ et $z''$. 
 Un \'el\'ement de $T_{AD}(F)$ s'\'ecrit $\prod_{l=1,...,n}\check{\varpi}_{l}(x_{l})$ avec $x_{l}\in F^{\times}$ pour $l=1,...,n-2$, $x_{n-1},x_{n}\in E_{0}^{\times}$ et $x_{n}=\tau(x_{n-1})$. D\'efinissons un homomorphisme $T_{ad}(F)\to E_{0}^{\times}/E_{0}^{\times,2}$ par $\prod_{l=1,...,n}\check{\varpi}_{l}(x_{l})\mapsto \prod_{l=1,...,n}x_{l}^{l}$.  Il s'en d\'eduit un isomorphisme $G_{AD}(F)/\pi(G(F))\simeq E_{0}^{\times}/E_{0}^{\times,2}$. L'image de $SO_{E_{0}/F}(2n,F)$ est celle de $F^{\times}/F^{\times,2}$ (cette image a deux \'el\'ements car $\mathfrak{o}_{F}^{\times}\subset \mathfrak{o}_{E_{0}}^{\times,2}$). L'image de $G_{AD}(F)_{0}$ est $\mathfrak{o}_{E_{0}}^{\times}/\mathfrak{o}_{E_{0}}^{\times,2}$. On note $\Xi_{0}$ le sous-ensemble des \'el\'ements de $\Xi$ dont la restriction \`a    $\mathfrak{o}_{E_{0}}^{\times}/\mathfrak{o}_{E_{0}}^{\times,2}$ est triviale si $\delta_{4}(n)=1$, non triviale si $\delta_{4}(n)=-1$.
 
 Notons ${\mathbb X}$ l'ensemble des couples $(k,h)\in {\mathbb N}^2$ tels que $k^2+h^2=2n$ et $k\geq h$. Puisque $n$ est pair, $k$ et $h$ sont forc\'ement pairs. Dans le cas o\`u $\delta_{\square}(2n)=0$, on pose ${\cal X}={\mathbb X}$ et $d_{x}=2$ pour tout $x=(k,h)\in {\cal X}$ avec $k>h$, $d_{x}=1$ si $x=(k,h)\in {\cal X}$ avec $k=h$.  Supposons $\delta_{\square}(2n)=1$. Il y a alors un unique couple $(k,h)\in {\mathbb X}$ tel que $h=0$. Notons-le $(k_{0},0)$. On pose ${\cal X}=({\mathbb X}-\{(k_{0},0)\})\sqcup \{(k_{0},0,\xi); \xi\in \Xi_{0}\}$ et $d_{x}=2$ si $x\in {\mathbb X}-\{(k_{0},0)\}$, $d_{x}=1$ si $x=(k_{0},0,\xi)$.

 L'ensemble $\underline{S}(G)$ s'envoie surjectivement sur l'ensemble des couples $(a,b)\in {\mathbb N}^2$ tels que $a+b=n$, $b\not=n$ et $b\not=1$. Les fibres ont un seul \'el\'ement sauf au-dessus de $(n,0)$ o\`u la fibre a deux \'el\'ements. L'action de $G_{AD}(F)$ est triviale sauf sur cette fibre \`a deux \'el\'ements, lesquels  sont permut\'es par cette action. Pour un sommet $s$ param\'etr\'e par $(a,b)$, avec $b\geq 2$,  on a $G_{s}\simeq (Spin_{{\mathbb F}_{q^2}/{\mathbb F}_{q}}(2a)\times Spin_{dep}(2b))/\{1,(z,z)\}$, avec la convention suivante: quand $a=1$, $Spin_{{\mathbb F}_{q^2}/{\mathbb F}_{q}}(2a)$ devient $GL(1)$ muni de l'action non triviale de $\Gamma_{{\mathbb F}_{q^2}/{\mathbb F}_{q}}$ et on  note $z$   l'\'el\'ement $-1\in GL(1)$. Pour un sommet $s$ param\'etr\'e par $(n,0)$, on a $G_{s}\simeq Spin_{{\mathbb F}_{q^2}/{\mathbb F}_{q}}(2n)$. Consid\'erons un sommet param\'etr\'e par $(a,b)$ avec $b\geq2$. Comme dans les paragraphes pr\'ec\'edents, on cherche les couples de fonctions $f_{N^{a},\epsilon^{a}}\times f_{N^b,\epsilon^b}$ telles que $\epsilon^{a}(z)=\epsilon^{b}(z)$. Supposons d'abord $\epsilon^{a}(z)=\epsilon^{b}(z)=-1$. D'apr\`es \ref{Dndeppair} et \ref{Dndepimp}, $b$ est pair ou $\delta_{4}(q-1)=1$. D'apr\`es \ref{Dnnondeppair} et \ref{Dnnondepimp}, on a $a$ impair et $\delta_{4}(q+1)=1$. Ces conditions sont contradictoires puisque l'\'egalit\'e $a+b=n$ impose que $a$ et $b$ sont de m\^eme parit\'e. Supposons maintenant $\epsilon^{a}(z)=\epsilon^{b}(z)=1$. Il y a de telles fonctions si et seulement si $(2a,2b)$ est de la forme $(k^2,h^2)$, avec $k,h$ forc\'ement pairs. Supposons ces conditions v\'erifi\'ees.  Il y a alors une unique fonction  $f_{N^{a},\epsilon^{a}}\times f_{N^b,\epsilon^b}$ qui donne naissance \`a un \'el\'ement de $FC(\mathfrak{g}(F))$. On note $\tilde{FC}_{k,h}$ la droite port\'ee par cette fonction. On n'a pas forc\'ement $(k,h)\in {\cal X}$ car on n'a pas forc\'ement $k\geq h$, condition que l'on a impos\'ee aux \'el\'ements de ${\cal X}$. Mais consid\'erons un \'el\'ement de ${\cal X}$ de la forme $(k,h)$. On a forc\'ement $k\not=0$ car $k\geq h$ et $h\not=0$ car un couple $(k,0)$ peut appartenir \`a ${\mathbb X}$ mais pas \`a ${\cal X}$. Nos constructions pour les couples $(a,b)=(k^2,h^2)$ et $(h^2,k^2)$ d\'efinissent des droites $\tilde{FC}_{k,h}$ et $\tilde{FC}_{h,k}$. Si $h=k$, elles sont \'egales et on pose $FC_{k,h}=\tilde{FC}_{k,h}$. Si $k>h$, on pose $FC_{k,h}=\tilde{FC}_{k,h}\oplus \tilde{FC}_{h,k}$. Cet espace est de dimension $2$.

   Consid\'erons maintenant un sommet $s$ param\'etr\'e par $(n,0)$. D'apr\`es \ref{Dnnondeppair}, $FC(\mathfrak{g}_{s}({\mathbb F}_{q}))$ est non nul si et seulement si $2n=k^2$ pour un $k\in {\mathbb N}$. Supposons cette condition  v\'erifi\'ee. Alors  $FC(\mathfrak{g}_{s}({\mathbb F}_{q}))$ est une droite  engendr\'ee par une fonction $f_{N,\epsilon}$  telle que $\epsilon(z)=1$. Le stabilisateur de $s$ dans $G_{AD}(F)$ est $G_{AD}(F)_{0}$. L'action de ce groupe sur $f_{N,\epsilon}$ est celle du groupe $G_{s,AD}({\mathbb F}_{q})$. Celle-ci  est triviale si $\delta_{4}(n)=1$, non triviale si $\delta_{4}(n)=-1$, cf \ref{Dnnondeppair}. 
 La fonction $f_{N,\epsilon}$ donne naissance \`a deux \'el\'ements de $FC(\mathfrak{g}(F))$ qui se transforment selon les deux caract\`eres $\xi$ de $G_{AD}(F)/\pi(G(F))$ qui prolongent le caract\`ere ainsi d\'etermin\'e  de $G_{AD}(F)_{0}$. On note $FC_{k,0,\xi}$ la droite port\'ee par l'\'el\'ement  qui se transforme selon  $\xi$.  On a $(k,0,\xi)\in {\cal X}$. Cela d\'emontre \ref{resultats}(1).

 Notons ${\mathbb Y}$ l'ensemble des couples $(i,j)\in {\mathbb N}^2$ tels que $i^2+j^2=n$ et  $i\geq j$. Dans le cas o\`u $\delta_{\square}(2n)=0$, on pose ${\cal Y}={\mathbb Y}$. Dans le cas o\`u $\delta_{\square}(2n)=1$, il y a un unique couple $(i,j)\in {\mathbb Y}$ tel que $i=j$, on le note $(i_{0},i_{0})$. On pose ${\cal Y}=({\mathbb Y}-\{(i_{0},i_{0})\})\sqcup \{(i_{0},i_{0},\xi);\xi\in \Xi_{0}\}$.

 Consid\'erons un \'el\'ement $(\sigma\mapsto \sigma_{G'},{\cal O})\in {\cal E}_{ell}(G)$.  Rappelons que l'application $\sigma\mapsto \omega_{G'}(\sigma)$ est un homomorphisme injectif de $\Gamma_{E_{G'}/E_{0}}$ dans $\hat{\Omega}$. Cela entra\^{\i}ne que $E_{G'}$ est inclus dans l'extension biquadratique $Q_{0}$ de $E_{0}$. On sait que $Q_{0}/F$ est galoisienne ab\'elienne et que $\Gamma_{Q_{0}/F}\simeq ({\mathbb Z}/4{\mathbb Z})\times ({\mathbb Z}/2{\mathbb Z})$ (cf. \ref{extensionsbiquadratiques}(1)).  Soit $\sigma\in \Gamma_{E_{G'}/E_{0}}$. Puisque $Q_{0}/F$ est ab\'elienne, on a $(\sigma\tau)_{G'}=(\tau
\sigma)_{G'}$, c'est-\`a-dire $\omega_{G'}(\sigma)\omega_{G'}(\tau)\theta=\omega_{G'}(\tau)\theta\omega_{G'}(\sigma)$. Donc $\omega_{G'}(\sigma)$ commute \`a $\theta$. Les seuls \'el\'ements de $\hat{\Omega}$ qui commutent \`a $\theta$ sont $1$ et $\theta\theta'$. Il en r\'esulte que $E_{G'}=E_{0}$ ou $E_{G'}$ est une extension quadratique de $E_{0}$ et que, dans ce cas,  l'image de $\Gamma_{E_{G'}/E_{0}}$ par l'action galoisienne est $\{1,\theta\theta'\}$. Rappelons que $Q_{0}$ est la compos\'ee de l'extension non ramifi\'ee de degr\'e  $4$ de $F$ et de l'extension $F(\sqrt{\varpi_{F}})$ de $F$. Pour fixer la notation, on suppose que $\tau$ agit trivialement sur $F(\sqrt{\varpi_{F}})$ et on note $\rho$ l'\'el\'ement de $\Gamma_{Q_{0}/F}$ qui agit trivialement sur l'extension non ramifi\'ee de degr\'e $4$ et non trivialement sur $F(\sqrt{\varpi_{F}})$. 
Il y a deux extensions quadratiques $K$ de $E_{0}$
telles que $\Gamma_{K}/\Gamma_{F}\simeq {\mathbb Z}/4{\mathbb Z}$, l'une non ramifi\'ee (son fixateur est $\{1,\rho\}$), l'autre ramifi\'ee (son fixateur est $\{1,\tau^2\rho\}$). Il y a une extension quadratique $K$ de $E_{0}$ telles que $\Gamma_{K}/\Gamma_{F}\simeq ({\mathbb Z}/2{\mathbb Z})^2$, \`a savoir $K=Q$ l'extension biquadratique de $F$ (son fixateur est $\{1,\tau^2\}$). 

Supposons d'abord $E_{G'}=E_{0}$. On doit avoir $\tau_{G'}^2=1$ c'est-\`a-dire $\omega_{G'}(\tau)\theta\omega_{G'}(\tau)\theta=1$. De nouveau, cela implique $\omega_{G'}(\tau)=1$ ou $\theta\theta'$. On voit que les deux cas sont conjugu\'es par $\delta\in \Omega$, donc que les donn\'ees endoscopiques associ\'ees sont \'equivalentes. Supposons donc $\tau_{G'}=\theta$. Si l'orbite ${\cal O}$ a deux \'el\'ements, on a ${\cal O}=\{\hat{\alpha}_{n-1},\hat{\alpha}_{n}\}$. Le fixateur de chacune de ces deux racines est $\Gamma_{E_{0}}$. Puisque $E_{0}/F$ est non ramifi\'ee, ce cas est exclu par le lemme \ref{centre}. Si ${\cal O} $ est r\'eduite \`a $\hat{\alpha}_{0}$ ou $\hat{\alpha}_{1}$, ces deux cas \'etant d'ailleurs conjugu\'es, on a ${\bf G}'={\bf G}$ et, \`a ce point, on ne peut rien dire de l'espace $FC^{st}(\mathfrak{g}(F))$. Supposons ${\cal O}=\{\hat{\alpha}_{m}\}$ pour un $m\in \{2,...,n-2\}$. On voit alors que $G'_{SC}\simeq Spin_{dep}(2m)\times Spin_{E_{0}/F}(2n-2m)$. En utilisant \ref{resultats}(4) par r\'ecurrence, on a $FC^{st}(\mathfrak{g}'(F))\not=\{0\}$ si et seulement si $n-m=i^2$ et $m=j^2$ avec $i,j\in {\mathbb N}$ et $i,j$ pairs. Supposons ces conditions v\'erifi\'ees. L'espace $FC^{st}(\mathfrak{g}'(F))$ est une droite. Le groupe d'automorphismes ext\'erieurs  de ${\bf G}'$ est $\{1,\theta\theta'\}$. Il agit trivialement sur $FC^{st}(\mathfrak{g}'(F))$.  Si $i\not=j$,  on pose $ {\bf G}'_{i,j}={\bf G}'$ et $\tilde{FC}^{{\cal E}}_{i,j}=FC^{st}(\mathfrak{g}_{i,j}'(F))^{Out({\bf G}'_{i,j})}$. 
 Supposons $i=j$. On a $n/2=i^2$ qui est pair. Calculons $\xi_{{\bf G}'}$. On peut choisir
$$(1) \qquad s_{sc}=(\prod_{l=1,...,n/2}\check{\hat{\alpha}}_{l}((-1)^l))(\prod_{l=n/2+1,...,n-2}\check{\hat{\alpha}}_{l}((-1)^{n/2}))\check{\hat{\alpha}}_{n-1}(i^{n/2})\check{\hat{\alpha}}_{n}(i^{n/2}).$$
  On calcule $\theta(s_{sc})=s_{sc}$, $\delta(s_{sc})s_{sc}^{-1}=z'z^{n/2}$. 
On a \'ecrit ces formules sous une forme g\'en\'erale car elle nous serviront plus loin. Ici, $n/2$ est pair et les formules se simplifient. Dans notre cas, $\sigma_{G'}(s_{sc})s_{sc}^{-1}=1$ pour tout $\sigma\in \Gamma_{F}$. Puisque $n/2$ est pair, le caract\`ere ${\bf 1}\in \Xi$ appartient \`a $\Xi_{0}$ et $(i,i,{\bf 1})$ appartient \`a ${\cal Y}$. On pose $ {\bf G}'_{i,i,{\bf 1}}={\bf G}'$ et $FC^{{\cal E}}_{i,i,{\bf 1}}=FC^{st}(\mathfrak{g}_{i,i,{\bf 1}}'(F))^{Out({\bf G}'_{i,i,{\bf 1}})}$. 

Supposons que $E_{G'}=Q$. L'extension $E_{G'}/E_{0}$ est ramifi\'ee et $\Gamma_{E_{G'}/E_{0}}$ est engendr\'e par l'image de $\rho$. On a $\rho_{G'}=\theta\theta'$. De nouveau, on a $\tau_{G'}^2=1$ donc $\omega_{G'}(\tau)=1$ ou $\theta\theta'$. Les deux cas sont encore conjugu\'es par $\delta$. Supposons donc $\tau_{G'}=\theta$. Le fixateur de $\hat{\alpha}_{0}$ ou $\hat{\alpha}_{1}$ est $\Gamma_{E'}$ o\`u $E'$ est l'extension quadratique de $F$ telle que $\Gamma_{Q/E'}=\{1,\tau\}$.  Le stabilisateur de $\hat{\alpha}_{n-1}$ ou $\hat{\alpha}_{n}$ est $\Gamma_{E''}$ ou $E''$ est l'extension quadratique de $F$ telle que $\Gamma_{Q/E''}=\{1,\tau\rho\}$. Les extensions $E'$ et $E''$ sont les deux extensions quadratiques ramifi\'ees de $F$. Supposons ${\cal O}=\{\hat{\alpha}_{m}\}$ pour un $m\in \{2,...,n-2\}$.  On voit alors que $G'_{SC}\simeq Spin_{E'/F}(2m)\times Spin_{E''/F}(2n-2m)$. En utilisant \ref{Dnpairpadiqueram} (1) et \ref{Dnimppadiqueram} (1) par r\'ecurrence, on a $FC^{st}(\mathfrak{g}'(F))\not=\{0\}$ si et seulement si $n-m=i^2$ et $m=j^2$ avec $i,j\in {\mathbb N}$ et $i,j$ impairs. Supposons ces conditions  v\'erifi\'ees. Alors l'espace $FC^{st}(\mathfrak{g}'(F))$ est une droite. Le groupe d'automorphismes ext\'erieurs  de ${\bf G}'$ est $\{1,\theta\theta'\}$. Il agit trivialement sur $FC^{st}(\mathfrak{g}'(F))$.   Si $i\not=j$,  on pose $ {\bf G}'_{i,j}={\bf G}'$ et $\tilde{FC}^{{\cal E}}_{i,j}=FC^{st}(\mathfrak{g}_{i,j}'(F))^{Out({\bf G}'_{i,j})}$.  Supposons  $i=j$. On a $n/2=i^2$ qui est impair. L'\'el\'ement $s_{sc}$ est le m\^eme que ci-dessus et on calcule $\tau_{G'}(s_{sc})s_{sc}^{-1}=1$, $\rho_{G'}(s_{sc})s_{sc}^{-1}=z$. La restriction de ce cocycle \`a $\Gamma_{E_{0}}$ n'est pas non ramifi\'e. Il  d\'etermine un caract\`ere $\xi\in \Xi$ qui est non trivial sur $\mathfrak{o}_{E_{0}}^{\times}/\mathfrak{o}_{E_{0}}^{\times,2}$. Donc $\xi\in \Xi_{0}$ puisque $n/2$ est impair et $(i,i,\xi)$ appartient \`a ${\cal Y}$. On pose $ {\bf G}'_{i,i,\xi}={\bf G}'$ et $FC^{{\cal E}}_{i,i,\xi}=FC^{st}(\mathfrak{g}_{i,i,\xi}'(F))^{Out({\bf G}'_{i,i,\xi})}$.   
 Maintenant ${\cal O}$ peut aussi \^etre \'egale \`a $\{\hat{\alpha}_{0},\hat{\alpha}_{1}\}$ ou \`a $\{\hat{\alpha}_{n-1},\hat{\alpha}_{n}\}$ (puisque $E'/F$ et $E''/F$ sont ramifi\'ees, ces cas ne sont plus exclus par le lemme \ref{centre}).    On s'aper\c{c}oit que ces cas sont similaires au pr\'ec\'edent, l'entier $j$, resp. $i$ \'etant alors \'egal \`a $1$ (ces cas \'etaient exclus auparavant puisqu'on avait $2\leq m\leq n-2$). On pose les m\^emes d\'efinitions que ci-dessus. 

Supposons que $\Gamma_{E_{G'}/F}\simeq {\mathbb Z}/4{\mathbb Z}$. On a $\rho_{G'}=1$ si $E_{G'}/F$ est non ramifi\'ee et $(\rho\tau^2)_{G'}=1$ si $E_{G'}/F$ est ramifi\'ee.   On doit avoir $\tau_{G'}^2=\theta\theta'$, ce qui impose $\tau_{G'}=\delta\theta$ ou $\tau_{G'}=\delta\theta'$. Ces deux cas sont conjugu\'es par $\delta$. Supposons donc que $\tau_{G'}=\delta\theta$.   Si l'orbite ${\cal O}$ a au moins deux \'el\'ements, le fixateur d'une racine dans ${\cal O}$ est de la forme $\Gamma_{L}$ ou $L$ est une extension non triviale de $F$ contenue dans $E_{G'}$. Alors $L$ contient $E_{0}$ et n'est pas totalement ramifi\'ee sur $F$. Le lemme \ref{centre} exclut ce cas. Il y a une seule orbite possible \`a un seul \'el\'ement, \`a savoir ${\cal O}=\{\hat{\alpha}_{n/2}\}$. Alors $G'_{SC}\simeq  Res_{E_{0}/F}(Spin_{E_{G'}/E_{0}}(n))$. On utilise (1) ci-dessous, \ref{Dnimppadiquenonram} (1), \ref{Dnpairpadiqueram} (1) et \ref{Dnimppadiqueram} (1)  par r\'ecurrence. L'espace $FC^{st}(\mathfrak{g}'(F))$ est non nul si et seulement si  les conditions suivantes sont v\'erifi\'ees: $n/2$ est de la forme $i^2$ pour un $i\in {\mathbb N}$, $E_{G'}/E_{0}$ est non ramifi\'ee si $i$ est pair et est ramifi\'ee si $i$ est impair. 
Supposons $n/2=i^2$ et supposons que $E_{G'}$ est l'extension ainsi d\'etermin\'ee. L'espace  $FC^{st}(\mathfrak{g}'(F))$ est une droite. Le groupe $Out({\bf G}')$ est \'egal \`a $\{1,\theta\theta'\}$ et on voit qu'il agit trivialement sur $FC^{st}(\mathfrak{g}'(F))$. Calculons $\xi_{{\bf G}'}$. L'\'el\'ement $s_{sc}$ est le m\^eme que ci-dessus. On obtient $\tau_{G'}(s_{sc})(s_{sc})^{-1}=z'z^{n/2}$, $\rho_{G'}(s_{sc})s_{sc}^{-1}=1$ si $E_{G'}/F$ est non ramifi\'ee et $\rho_{G'}(s_{sc})s_{sc}^{-1}=z$ si $E_{G'}/F$ est ramifi\'ee. Si $E_{G'}/F$ est non ramifi\'ee, c'est-\`a-dire si $n/2$ est pair, le cocycle est trivial sur le groupe d'inertie, donc d\'etermine un \'el\'ement $\xi\in \Xi$ qui est trivial sur $\mathfrak{o}_{E_{0}}^{\times}/\mathfrak{o}_{E_{0}}^{\times,2}$. On a $\xi\in \Xi_{0}$ puisque $n/2$ est pair.  Mais le cocycle est non cohomologue au cocycle trivial, donc $\xi\not={\bf 1}$. Si $E_{G'}/F$ est ramifi\'ee, c'est-\`a-dire si $n/2$ est impair, le cocycle est non trivial sur le groupe d'inertie, donc d\'etermine un \'el\'ement $\xi\in \Xi$ qui est non trivial sur $\mathfrak{o}_{E_{0}}^{\times}/\mathfrak{o}_{E_{0}}^{\times,2}$. On a $\xi\in \Xi_{0}$ puisque $n/2$ est impair. Mais on voit que le cocycle n'est pas cohomologue \`a celui que l'on a rencontr\'e plus haut. Dans les deux cas,  on a $(i,i,\xi)\in{\cal Y}$. On pose $ {\bf G}'_{i,i,\xi}={\cal G}'$ et $FC^{{\cal E}}_{i,i,\xi}=FC^{st}(\mathfrak{g}_{i,i,\xi}'(F))^{Out({\bf G}'_{i,i,\xi})}$.  

Si $\delta_{\square}(2n)=1$, on a d\'efini une droite $FC^{{\cal E}}_{i_{0},i_{0},\xi}$ pour tout \'el\'ement $\xi\in \Xi_{0}$ (qui intervient dans diff\'erents cas ci-dessus selon la parit\'e de $n/2$). Soit $(i,j)\in {\mathbb Y}$ avec $i>j$. Supposons  $j\not=0$. On a construit ci-dessus deux droites $\tilde{FC}^{{\cal E}}_{i,j}$ et $\tilde{FC}^{{\cal E}}_{j,i}$ (intervenant elles-aussi dans diff\'erents cas ci-dessus selon la parit\'e de $n/2$). On pose $FC^{{\cal E}}_{i,j}=\tilde{FC}^{{\cal E}}_{i,j}\oplus \tilde{FC}^{{\cal E}}_{j,i}$. Dans le cas o\`u il existe un \'el\'ement de ${\mathbb Y}$ de la forme $(i,0)$, on n'a rien associ\'e \`a cet \'el\'ement (dans les constructions ci-dessus, on avait toujours $ij\not=0$). On n'a pas trait\'e non plus la donn\'ee principale ${\bf G}$.

On d\'efinit une application $\phi:{\mathbb X}\to {\mathbb Y}$ par $\phi(k,h)=((k+h)/2,(k-h)/2)$. C'est une bijection qui se rel\`eve naturellement en une bijection $\varphi:{\cal X}\to {\cal Y}$. Il y a un \'el\'ement $(k,h)\in {\mathbb X}$ tel que $k=h$ si et seulement si $\delta_{\square}(n)=1$. Si cette condition est v\'erifi\'ee, on note $(k^{st},k^{st})$ cet \'el\'ement.  Il y a   un  \'el\'ement de ${\mathbb Y}$ de la forme $(i,0)$ si et seulement si $\delta_{\square}(n)=1$. Si cette condition est v\'erifi\'ee, l'\'el\'ement en question est \'egal \`a $(k^{st},0)$ (c'est l'\'el\'ement $y\in {\cal Y}$ auquel on n'a encore rien associ\'e). Si $\delta_{\square}(n)=0$, on pose ${\cal X}^{st}= {\cal Y}^{st}=\emptyset$. Si $\delta_{\square}(n)=1$,  on pose ${\cal X}^{st}= \{(k^{st},k^{st})\}$ et ${\cal Y}^{st}=\{(k^{st},0)\}$.  On voit que $\varphi({\cal X}^{st})={\cal Y}^{st}$. 

On remarque que, pour tout $x\in {\cal X}-{\cal X}^{st}$, on a $dim(FC_{x})=dim(FC^{{\cal E}}_{\varphi(x)})$. On peut alors appliquer l'habituel argument de dimension. On obtient que $FC^{st}(\mathfrak{g}(F))=\{0\}$ si ${\cal Y}^{st}=\emptyset$. Si ${\cal Y}^{st}\not=\emptyset$, $FC^{st}(\mathfrak{g}(F))$ est une droite et on la note $FC^{{\cal E}}_{k^{st},0}$. Cela ach\`eve la preuve de \ref{resultats}(2).

 On munit ${\mathbb X}$ de la relation $(k',h')\leq (k,h)$ si et seulement si $k'+h'\leq k+h$. On montre que c'est une relation d'ordre et on la rel\`eve en une relation de pr\'eordre sur ${\cal X}$. Le fait que certains espaces $FC^{{\cal E}}_{y}$ sont de dimension $2$ cr\'ee une petite difficult\'e. Pour la r\'esoudre, notons ${\cal Y}_{1}$, resp. ${\cal Y}_{2}$, le sous-ensemble des $y\in {\cal Y}$ tels que $dim(FC^{{\cal E}}_{y})=1$, resp. $dim(FC^{{\cal E}}_{y})=2$. L'ensemble ${\cal Y}_{2}$ est celui des couples $(i,j)\in {\mathbb N}^2$ tels que $i^2+j^2=n$ et $i>j>0$. On note $\tilde{{\cal Y}}_{2}$ l'ensemble des couples $(i,j)\in {\mathbb N}^2$ tels que $i^2+j^2=n$, $ij\not=0$ et $i\not=j$. Il y a une surjection \'evidente $\psi_{2}:\tilde{{\cal Y}}_{2}\to {\cal Y}_{2}$ qui envoie $(i,j)$ sur $(i,j)$ si $i>j$, sur $(j,i)$ si $i<j$. On pose $\tilde{{\cal Y}}={\cal Y}_{1}\cup \tilde{{\cal Y}}_{2}$ et on prolonge $\psi_{2}$ par l'identit\'e de ${\cal Y}_{1}$ en une surjection $\psi:\tilde{{\cal Y}}\to {\cal Y}$.  Pour unifier la notation, on pose $\tilde{FC}^{{\cal E}}_{y}=FC^{{\cal E}}_{y}$ pour $y\in {\cal Y}_{1}$. Dans les constructions pr\'ec\'edentes, on a  d\'efini une donn\'ee endoscopique ${\bf G}'_{\tilde{y}}$ et une droite $\tilde{FC}^{{\cal E}}_{\tilde{y}}$ pour tout $\tilde{y}\in \tilde{{\cal Y}}$. L'application $\tilde{y}\mapsto {\bf G}'_{\tilde{y}}$ est injective. Pour $y\in {\cal Y}_{2}$, on a $FC^{{\cal E}}_{y}=\oplus_{\tilde{y}\in \psi^{-1}(y)}\tilde{FC}^{{\cal E}}_{\tilde{y}}$. 
 Pour tout $\tilde{y}\in \tilde{{\cal Y}}$, on introduira en \ref{elementsYy} un \'el\'ement $Y_{\tilde{y}}\in \mathfrak{g}'_{\tilde{y},ell}(F)$ v\'erifiant les propri\'et\'es suivantes:
 
 (2) pour un \'el\'ement non nul $f'\in FC^{st}(\mathfrak{g}'_{\tilde{y}}(F))^{Out({\bf G}'_{\tilde{y}})}$, on a $S^{G'_{\tilde{y}}}(Y_{\tilde{y}},f')\not=0$;
 
 (3) soient $x\in {\cal X}$,  $f\in FC_{x}$ et $X$ un \'el\'ement de $\mathfrak{g}_{reg}(F)$ dont la classe de conjugaison stable correspond \`a celle de $Y_{\tilde{y}}$; supposons $I^G(X,f)\not=0$; alors $(x)\geq \varphi^{-1}((\psi(\tilde{y})))$.
 
 On utilise les notations de \ref{ingredients}. Pour $(y)\in \underline{{\cal Y}}$, l'entier $d_{(y)}=dim(FC^{\cal E}_{(y)})$ est le nombre d'\'el\'ements de l'ensemble $\psi^{-1}((y))$.   On associe \`a $(y)$ les familles $({\bf G}'_{\tilde{y}})_{\tilde{y}\in \psi^{-1}((y))}$ et $(Y_{\tilde{y}})_{\tilde{y}\in \psi^{-1}((y))}$.  Les propri\'et\'es (2) et (3) ci-dessus entra\^{\i}nent que  les hypoth\`eses de \ref{ingredients} sont v\'erifi\'ees (pour $\underline{{\cal Y}}^{\natural}=\underline{{\cal Y}}$). Cela entra\^{\i}ne
 
 $$transfert(FC_{(x)})=FC^{{\cal E}}_{\varphi((x))}$$
 pour tout $(x)\in \underline{{\cal X}}$. Quand la classe $(x)$ n'est pas r\'eduite \`a un seul \'el\'ement, on utilise omme en \ref{An-1quasidepram} l'action de $G_{AD}(F)/\pi(G(F))$ pour  raffiner cette \'egalit\'e et on obtient  \ref{resultats}(3).  
  L'assertion \ref{resultats}(4) se d\'eduit de \ref{resultats}(3) et de l'\'egalit\'e $FC^{st}(\mathfrak{g}(F)) =FC^{{\cal E}}_{k^{st},0}$ ci-dessus. Explicitons la cons\'equence de \ref{resultats}(4):
 
 (4) $dim(FC^{st}(\mathfrak{g}(F)))=\delta_{\square}(n)$.

\subsection{Forme int\'erieure du type $D_{n}$ quasi-d\'eploy\'e, $n$ pair, $E_{0}/F$ non ramifi\'ee}
On suppose que $G^*$ est du type pr\'ec\'edent et que $G$ en est la forme int\'erieure associ\'ee \`a l'unique caract\`ere non trivial de $Z(\hat{G}_{SC})^{\Gamma_{F}}=\{1,z\}$, ou encore \`a  l'image de  $\delta\in \hat{\Omega}\simeq N$ dans $N_{\Gamma_{F}^{nr}}$. Dans les tables de Tits, le groupe est de type $^4D_{n}$.

Si $\delta_{\square}(n)=0$, posons $ {\cal X} =\emptyset$. Si $\delta_{\square}(n)=1$, c'est-\`a-dire si $n=k^2$ pour un $k\in {\mathbb N}$ (forc\'ement pair), posons ${\cal X} =\{(k,k)\}$ et $\delta_{x}=1$ pour $x\in {\cal X}$. 

 On consid\`ere le diagramme ${\cal D}_{a}$ muni de l'action de $\Gamma_{{\mathbb F}_{q}}$ telle que le Frobenius $Fr$ agisse par $\delta\theta$. Cette action est triviale sur $\Gamma_{{\mathbb F}_{q^4}}$. L'ensemble $\underline{S}(G)$ est param\'etr\'e par les orbites de cette action.  Pour les sommets $s$ associ\'es \`a des  orbites qui ont au moins deux \'el\'ements,  le lemme \ref{orbites} montre que   $FC(\mathfrak{g}_{s}({\mathbb F}_{q}))=\{0\}$. Il reste l'unique sommet $s$ param\'etr\'e par l'orbite ${\cal O}=\{\hat{\alpha}_{n/2}\}$. Sur $\bar{{\mathbb F}}_{q}$, on a $G_{s}\simeq (Spin(n)\times Spin(n))/\{1,(z,z)\}$. L'action galoisienne est l'action d\'eploy\'ee tordue par l'action alg\'ebrique suivante: $Fr$ permute les deux copies de $Spin(n)$ et $Fr^2$ agit sur chaque facteur par l'automorphisme non trivial $\theta$ relatif \`a ce facteur. D'apr\`es \ref{Dnnondeppair}, l'espace $FC(\mathfrak{g}_{s}({\mathbb F}_{q}))$ est non nul si et seulement si $\delta_{\square}(n)=1$. Supposons cette condition v\'erifi\'ee, c'est-\`a-dire $n=k^2$ pour un entier $k\in {\mathbb N}$. L'espace $FC(\mathfrak{g}_{s}({\mathbb F}_{q}))$ est alors une droite munie d'un g\'en\'erateur $f_{N,\epsilon}$ avec $\epsilon(z)=1$. On note  cette droite $FC_{k,k}$. Cela d\'emontre l'assertion \ref{resultats}(1).
 
Si $\delta_{\square}(n)=0$, posons $ {\cal Y} =\emptyset$. Si $\delta_{\square}(n)=1$, c'est-\`a-dire si $n=k^2$ pour un $k\in {\mathbb N}$, posons ${\cal Y} =\{(k,0)\}$.

D'apr\`es le calcul ci-dessus de $FC(\mathfrak{g}(F))$ et  par \'egalit\'e des dimensions, $FC^{{\cal E}}(\mathfrak{g}(F))$ est nul si $\delta_{\square}(n)=0$ et est une droite si $\delta_{\square}(n)=1$. Il reste dans ce dernier cas \`a d\'eterminer cette droite. Mais alors $FC^{st}(\mathfrak{g}^*(F))$ est une droite  d'apr\`es le paragraphe pr\'ec\'edent. On pose alors ${\bf G}'_{k,0}={\bf G}$ et $FC^{{\cal E}}_{k,0}=FC^{st}(\mathfrak{g}^*(F))$. Cela d\'emontre \ref{resultats}(2). 

En d\'efinissant $\varphi$ comme l'unique bijection de ${\cal X}$ sur ${\cal Y}$, \ref{resultats}(3) est triviale. 

\subsection{Type $D_{n}$ quasi-d\'eploy\'e, $n$ impair, $E_{0}/F$ non ramifi\'ee}\label{Dnimppadiquenonram}
 On suppose que $G$ est quasi-d\'eploy\'e de type $D_{n}$, avec $n\geq4$, $n$ impair, et que $\Gamma_{F}$ agit sur le diagramme ${\cal D}$ de $G$ par l'action $\sigma\mapsto \sigma_{G}$ triviale sur $\Gamma_{E_{0}}$ et telle que $\sigma_{G}=\theta$ pour $\sigma\in \Gamma_{F}-\Gamma_{E_{0}}$. On note $\tau$ un \'el\'ement de Frobenius de $\Gamma_{F}$. L'action de $\tau$ sur $Z(G)$ est triviale si $\delta_{4}(q-1)=0$ et  envoie $z'$ sur $(z')^3$ si $\delta_{4}(q-1)=1$. Introduisons l'endomorphisme $\iota$ du groupe $E^{\times}\times F^{\times}$  d\'efini par $\iota(e,f)=(fe^2,e\tau(e))$. Un \'el\'ement de $T_{AD}(F)$ s'\'ecrit $\prod_{l=1,...,n}\check{\varpi}_{l}(x_{l})$, avec $x_{l}\in F^{\times}$ pour $l=1,...,n-2$, $x_{n-1},x_{n}\in E_{0}^{\times}$ et $x_{n}=\tau(x_{n-1})$. D\'efinissons un homomorphisme de $T_{AD}(F)$ dans $E^{\times}\times F^{\times}$ par $\prod_{l=1,...,n}\check{\varpi}_{l}(x_{l})\mapsto (x_{n},\prod_{l=1,...,n-2}x_{l}^l)$.  Il s'en  d\'eduit un isomorphisme $G_{AD}(F)/\pi(G(F))\to (E_{0}^{\times}\times F^{\times})/\iota(E_{0}^{\times}\times F^{\times})$. Il y a une suite exacte
$$1\to \mathfrak{o}_{E_{0}}^{\times}/(\mathfrak{o}_{F}^{\times}\mathfrak{o}_{E_{0}}^{\times,4})\to (E_{0}^{\times}\times F^{\times})/\iota(E_{0}^{\times}\times F^{\times})\to {\mathbb Z}/2{\mathbb Z}\to 0.$$
Le premier homomorphisme envoie $e\in \mathfrak{o}_{E_{0}}^{\times}$ sur $(e,1)$. Le second envoie $(e,f)$ sur $val_{F}(f)$ mod $2{\mathbb Z}$. L'image de $G_{AD}(F)_{0}$ est le groupe $\mathfrak{o}_{E_{0}}^{\times}/(\mathfrak{o}_{F}^{\times}\mathfrak{o}_{E_{0}}^{\times,4})$, qui a $2$ \'el\'ements si $\delta_{4}(q-1)=1$, $4$ \'el\'ements si $\delta_{4}(q-1)=0$. L'image de $SO_{E_{0}/F}(2n,F)$ est un sous-groupe d'indice $2$ de $G_{AD}(F)/\pi(G(F))$ qui s'envoie surjectivement sur ${\mathbb Z}/2{\mathbb Z}$ et dont l'intersection avec $\mathfrak{o}_{E_{0}}^{\times}/(\mathfrak{o}_{F}^{\times}\mathfrak{o}_{E_{0}}^{\times,4})$ est $\mathfrak{o}_{E_{0}}^{\times,2}/(\mathfrak{o}_{F}^{\times}\mathfrak{o}_{E_{0}}^{\times,4})$. Dans le cas o\`u $\delta_{4}(q-1)=0$, on note $\Xi_{0}$  l'ensemble des caract\`eres de $G_{AD}(F)/\pi(G(F))$ qui sont non triviaux sur $\mathfrak{o}_{E_{0}}^{\times,2}/(\mathfrak{o}_{F}^{\times}\mathfrak{o}_{E_{0}}^{\times,4})$, ou encore dont la restriction \`a $G_{AD}(F)_{0}/\pi(G(F))$ est d'ordre $4$.  

Si $\delta_{4}(q-1)=1$, on pose ${\cal X} =\emptyset$.

 Supposons $\delta_{4}(q-1)=0$. 
On note ${\mathbb X}$ l'ensemble des couples $(k,h)\in {\mathbb N}^2$ tels que $k(k+1)/2+h(h+1)/2=2n$, $k\geq h$ et $k(k+1)/2$ et $h(h+1)/2$ sont tous deux pairs.   On note ${\cal X} $ l'ensemble des triplets $(k,h,\xi)$ o\`u $(k,h)\in {\mathbb X}$ et $\xi\in \Xi_{0}$.   On pose $d_{x}=1$ pour tout $x\in {\cal X}$.

 L'ensemble $\underline{S}(G)$ s'envoie surjectivement sur l'ensemble des couples $(a,b)\in {\mathbb N}^2$ tels que $a+b=n$, $b\not=n$ et $b\not=1$. Les fibres ont un seul \'el\'ement sauf au-dessus de $(n,0)$ o\`u la fibre a deux \'el\'ements. L'action de $G_{AD}(F)$ fixe tous les sommets sauf ceux de cette fibre \`a deux \'el\'ements, qu'elle permute. Pour un sommet $s$ param\'etr\'e par $(a,b)$ avec $b\geq2$, on a $G_{s}\simeq (Spin_{{\mathbb F}_{q^2}/{\mathbb F}_{q}}(2a)\times Spin_{dep}(2b))/\{1,(z,z)\}$, avec la m\^eme convention  qu'en \ref{Dnpairpadiquenonram} si $a=1$.  Pour un sommet $s$ param\'etr\'e par $(n,0)$, on a $G_{s}\simeq Spin_{{\mathbb F}_{q^2}/{\mathbb F}_{q}}(2n)$.  Consid\'erons un sommet $s$ param\'etr\'e par $(a,b)$ avec $b\geq2$. Dans le cas $a=1$, on a $Z(G_{s})^0\not=\{1\}$ et $FC(\mathfrak{g}_{s}({\mathbb F}_{q}))=\{0\}$. Supposons $a\not=1$ donc $a\geq2$ puisque $a\not=0$.  Comme dans les paragraphes pr\'ec\'edents, on cherche les couples de fonctions $f_{N^{a},\epsilon^{a}}\times f_{N^b,\epsilon^b}$ telles que $\epsilon^{a}(z)=\epsilon^b(z)$. Supposons d'abord $\epsilon^{a}(z)=\epsilon^b(z)=1$. D'apr\`es \ref{Dnnondeppair} et \ref{Dnnondepimp}, on a $a$ pair et $\delta_{\square}(2a)=1$. D'apr\`es \ref{Dndeppair} et \ref{Dndepimp}, on a $b$ pair et $\delta_{\square}(2b)=1$. Ces conditions sont impossibles car l'\'egalit\'e $a+b=n$ entra\^{\i}ne que $a$ et $b$ sont de parit\'e distincte. Supposons maintenant $\epsilon^{a}(z)=\epsilon^b(z)=-1$.
 D'apr\`es \ref{Dnnondeppair} et \ref{Dnnondepimp}, on a $a$ impair, $\delta_{\triangle}(2a)=1$ et $\delta_{4}(q-1)=0$.  Ces conditions impliquent que $b$ est pair.  
 D'apr\`es \ref{Dndeppair}  on a $\delta_{\triangle}(2b)=1$. Supposons donc $\delta_{4}(q-1)=0$,  $\delta_{\triangle}(2a)=\delta_{\triangle}(2b)=1$, c'est-\`a-dire $2a=k'(k'+1)/2$ et $2b=2h'(h'+1)/2$. Il y a alors deux fonctions sur chaque composante de $G_{s}$ que l'on note $f_{N^{a},\epsilon^{'a}}$, $f_{N^{a},\epsilon^{''a}}$, $f_{N^{b},\epsilon^{'b}}$, $f_{N^{b},\epsilon^{''b}}$ conform\'ement \`a \ref{Dnnondepimp} et  \ref{Dndeppair}. L'image de $z'\in Z(G)$ dans $Z(G_{s})$ est $(z',z')$. Ce couple multiplie  chaque fonction $f_{N^{a},\epsilon^{a}}\times f_{N^b,\epsilon^b}$ par $\pm i$ et le signe n'est pas le m\^eme sur toutes les fonctions car $z'$ agit diff\'eremment sur $f_{N^{a},\epsilon^{'a}}$ et $f_{N^{a},\epsilon^{''a}}$, ainsi que sur  $f_{N^{b},\epsilon^{'b}}$ et $f_{N^{b},\epsilon^{'b}}$. Pour chacun des deux caract\`eres d'ordre $4$ de $G_{AD}(F)_{0}/\pi(G(F))$, on obtient deux fonctions se transformant selon ce caract\`ere. L'action de $G_{AD}(F)/\pi(G(F))$ tout entier est r\'ecup\'er\'ee par l'action naturelle d'un \'el\'ement $(x,y)\in O_{{\mathbb F}_{q^2}/{\mathbb F}_{q}}(2a,{\mathbb F}_{q})\times O_{dep}(2b,{\mathbb F}_{q})$ tel que $det(x)=det(y)=-1$. Or cette action permute nos fonctions, car l'action de $x$ permute $f_{N^{a},\epsilon^{'a}}$ et $f_{N^{a},\epsilon^{''a}}$ et celle de $y$ permute $f_{N^{b},\epsilon^{'b}}$ et $f_{N^{b},\epsilon^{''b}}$. Pour chaque caract\`ere $\xi\in \Xi$ dont la restriction \`a $G_{AD}(F)_{0}/\pi(G(F))$ est d'ordre $4$, c'est-\`a-dire pour tout $\xi\in \Xi_{0}$, il y a donc une combinaison lin\'eaire convenable de deux de nos fonctions qui se transforme selon $\xi$.  On note $FC_{k,h,\xi}$ la droite  de $FC(\mathfrak{g}(F))$ issue de cette fonction, o\`u   $(k,h)$ est l'unique \'el\'ement de $\{(k',h'),(h',k')\}$ tel que $k>h$ (on ne peut pas avoir $k=h$ puisque $a\not=b$). Remarquons que le couple $(k',h')$ est uniquement d\'etermin\'e par $(k,h)$ puisque $k'(k'+1)/4$ doit \^etre impair et $h'(h'+1)/4$ doit \^etre pair. Supposons maintenant que $s$ est param\'etr\'e par $(n,0)$. Il n'y a toujours pas de fonction $f_{N,\epsilon}$ avec $\epsilon(z)=1$ car $n$ est impair. Il y a une fonction $f_{N,\epsilon}$ avec $\epsilon(z)=-1$ si et seulement si $\delta_{4}(q-1)=0$ et $\delta_{\triangle}(2n)=1$, c'est-\`a-dire $2n=k(k+1)/2$ pour un $k\in {\mathbb N}$. Supposons ces conditions  v\'erifi\'ees. On a deux fonctions, qui se transforment selon les deux caract\`eres de $G_{AD}(F)_{0}/\pi(G(F))$ d'ordre $4$. Ici, ce groupe est le fixateur de $s$ dans $G_{AD}(F)/\pi(G(F))$. Selon \ref{actionsurFC}, pour tout caract\`ere $\xi\in \Xi$ prolongeant l'un de ces caract\`eres, c'est-\`a-dire pour tout $\xi\in \Xi_{0}$, on obtient un \'el\'ement de $FC(\mathfrak{g}(F))$ se transformant selon $\xi$. On note  $FC_{k,0,\xi}$ la droite port\'ee par cet \'el\'ement. Cela prouve \ref{resultats}(1). 
 
 Si $\delta_{4}(q-1)=1$, on pose ${\cal Y}=\emptyset$. Supposons $\delta_{4}(q-1)=0$. On note ${\mathbb Y}$ l'ensemble des couples $(i,j)$ avec $i,j\in {\mathbb N}$, $i$ est impair et  $2i^2+j(j+1)/2=n$. On note ${\cal Y} $ l'ensemble des triplets $(i,j,\xi)$ o\`u $(i,j)\in {\mathbb Y}$ et $\xi\in \Xi_{0}$.
 
 Dans le cas $\delta_{4}(q-1)=1$, on a d\'ej\`a vu que $FC(\mathfrak{g}(F))=\{0\}$, donc les assertions (2), (3) et (4) sont triviales. On suppose d\'esormais que $\delta_{4}(q-1)=0$.

 Consid\'erons un couple $(\sigma\mapsto \sigma_{G'},{\cal O})$. Remarquons que $\tau_{G'}\in \hat{\Omega}\theta$ et que tout \'el\'ement de cet ensemble est de carr\'e $1$, donc $\tau_{G'}^2=1$. 
 
 Consid\'erons le cas o\`u $ E_{G'}=E_{0}$. Le terme $\omega_{G'}(\tau)$ peut \^etre n'importe quel \'el\'ement de $\hat{\Omega}$ mais on voit que les cas $\omega_{G'}(\tau)=1$ et $\omega_{G'}(\tau)=\theta\theta'$ sont \'equivalents (conjugu\'es par $\delta\theta$) ainsi que les cas $\omega_{G'}(\tau)=\delta\theta$ et $\omega_{G'}(\tau)=\delta\theta'$. Supposons $\omega_{G'}(\tau)=1$ donc $\tau_{G'}=\theta$.  Si l'orbite ${\cal O}$ est r\'eduite \`a une racine $\hat{\alpha}_{0}$ ou $\hat{\alpha}_{1}$ (ces deux cas \'etant \'equivalents, conjugu\'es par $\theta\theta'$), la donn\'ee ${\bf G}'$ est la donn\'ee principale ${\bf G}$ et, \`a ce point, on ne peut rien dire de $FC^{st}(\mathfrak{g}'(F))$. Si ${\cal O}=\{\hat{\alpha}_{n-1},\hat{\alpha}_{n}\}$, le fixateur d'une de ces racines est $\Gamma_{E_{0}}$. Puisque $E_{0}/F$ est non ramifi\'ee, ce cas est exclu par le lemme \ref{centre}. Supposons ${\cal O}=\{\hat{\alpha}_{m}\}$ avec $m\in \{2,...,n-2\}$. On voit que $G'_{SC}\simeq Spin_{dep}(2m)\times Spin_{E_{0}/F}(2n-2m)$. On utilise \ref{Dndeppairpadique} (5) et \ref{Dnpairpadiquenonram} (4) par r\'ecurrence: on n'a $FC^{st}(\mathfrak{g}'(F))\not=\{0\}$ que si $m$ et $n-m$ sont tous deux des carr\'es pairs. C'est impossible puisque $n$ est impair. Supposons maintenant $\omega_{G'}(\tau)=\delta\theta$, donc $\sigma_{G'}(\tau)=\delta$. Puisque $n$ est impair, l'orbite ${\cal O}$ contient deux \'el\'ements et le stabilisateur de chacun d'eux est $\Gamma_{E_{0}}$. Ce cas est encore exclu par le lemme \ref{centre}.
 
 Consid\'erons le cas o\`u l'image de $\Gamma_{E_{G'}/E_{0}}$ dans $\hat{\Omega}$ est l'unique sous-groupe d'ordre $2$, \`a savoir $\{1,\theta\theta'\}$. On voit que $E_{G'}$ est forc\'ement l'extension biquadratique $Q$ de $F$. On note $\Gamma_{E_{G'}/E_{0}}=\{1,\rho\}$. De nouveau, \`a \'equivalence pr\`es, on a $\tau_{G'}=\theta$ ou $\tau_{G'}=\delta$. Supposons $\tau_{G'}=\theta$. Supposons ${\cal O}=\{\hat{\alpha}_{m}\}$ avec $m\in \{2,...,n-2\}$. On voit que $G'_{SC}\simeq Spin_{E'/F}(2m)\times Spin_{E''/F}(2n-2m)$, o\`u $E'$ et $E''$ sont les extensions quadratiques de $F$ telles que $\Gamma_{E_{G'}/E'}=\{1,\tau\}$, $\Gamma_{E_{G'}/E''}=\{1,\tau\rho\}$. Ce sont les deux extensions quadratiques ramifi\'ees de $F$. On utilise \ref{resultats}(4) par r\'ecurrence: on n'a $FC^{st}(\mathfrak{g}'(F))\not=\{0\}$ que si $m$ et $n-m$ sont tous deux des carr\'es impairs. C'est impossible puisque $n$ est impair. Supposons ${\cal O}=\{\hat{\alpha}_{0},\hat{\alpha}_{1}\}$ ou ${\cal O}=\{\hat{\alpha}_{n-1},\hat{\alpha}_{n}\}$. Le   lemme \ref{centre} ne s'applique pas. Toutefois, on a $G'_{SC}\simeq Spin_{E''/F}(2n-2)$ ou $Spin_{E'/F}(2n-2)$ avec les m\^emes extensions $E'$ et $E''$ que ci-dessus et on ne  peut avoir $FC^{st}(\mathfrak{g}'(F))\not=\{0\}$ que si $n-1$ est un carr\'e impair, ce qui est impossible. Supposons maintenant $\tau_{G'}=\delta$. On voit que l'orbite ${\cal O}$ a $2$ ou $4$ \'el\'ements et que le fixateur d'un de ces \'el\'ements est $\Gamma_{E_{0}}$ ou $\Gamma_{E_{G'}}$. En tout cas, cela est exclu par le lemme \ref{centre}.
 
 Il reste le cas o\`u l'image de $\Gamma_{E_{G'}/E_{0}}$ est $\hat{\Omega}$ tout entier. On a alors $\Gamma_{E_{G'}/E_{0}}\simeq {\mathbb Z}/4{\mathbb Z}$. Fixons un g\'en\'erateur $\rho$ de $\Gamma_{E_{G'}/E_{0}}$. On a forc\'ement $\rho_{G'}=\delta\theta$ ou $\rho_{G'}=\delta\theta'$. Puisque $\tau_{G'}\in \hat{\Omega}\theta$, on voit que $(\rho\tau)_{G'}=(\tau\rho^{-1})_{G'}$. D'apr\`es \ref{extensionsdiedrales} (3), et puisqu'on a suppos\'e $\delta_{4}(q-1)=0$, une telle extension $E_{G'}$  existe et  est unique. A \'equivalence pr\`es, il y a $4$ actions galoisiennes possibles, d\'efinies par  $\rho_{G'}=\delta\theta$ ou $\rho_{G'}=\delta\theta'$ et $\tau_{G'}=\theta$ ou $\tau_{G'}=\delta$. On note $\tau'$ et $\rho'$ les \'el\'ements de $\Gamma_{E_{G'}/F}$ tels que $\tau'_{G'}=\theta$ et $\rho'_{G'}=\delta\theta$. Selon les cas, le couple $(\tau',\rho')$ est l'un des quatre couples $(\tau,\rho)$, $(\tau,\rho^{-1})$, $(\tau\rho,\rho)$, $(\tau\rho^{-1},\rho^{-1})$. Supposons que ${\cal O}=\{\hat{\alpha}_{m}\}$ avec $m\in \{2,...,n-2\}$. Notons $E'$, resp. $K$,   la sous-extension de $E_{G'}$ telle que  $\Gamma_{E_{G'}/E'}=\{1,\tau',(\rho')^2,\tau'(\rho')^2\}$, resp. $\Gamma_{E_{G'}/K}=\{1,\tau'\}$. Les extensions $K/E'$ et $E'/F$ sont quadratiques ramifi\'ees.   Le fixateur de $\hat{\alpha}_{m}$ est $\Gamma_{E'}$.  D'apr\`es \ref{centre}, $FC^{st}(\mathfrak{g}'(F))=FC^{st}(\mathfrak{g}'_{SC}(F))$. On a
  $G'_{SC}\simeq Res_{E'/F}Spin_{K/E'}(2m)\times SU_{E'/F}(n-2m)$. On applique   \ref{Dnpairpadiqueram} (1)  , \ref{Dnimppadiqueram} (1) et \ref{An-1quasidepram} (5) par r\'ecurrence. On a 
 $FC^{st}(\mathfrak{g}'(F))\not=\{0\}$ si et seulement si $m=i^2$ avec $i\in {\mathbb N}$, $i$ impair et $n-2m=j(j+1)/2$ avec $j\in {\mathbb N}$.   Supposons ces conditions v\'erifi\'ees. Alors $FC^{st}(\mathfrak{g}'(F))$ est une droite. Le groupe d'automorphismes ext\'erieurs  $Out({\bf G}')$ est \'egal \`a $\{1,\theta\theta'\}$. On voit qu'il agit trivialement sur $FC^{st}(\mathfrak{g}'(F))$.  Calculons $\xi_{{\bf G}'}$. On peut choisir $s_{sc}$ comme en \ref{Dndeppairpadique} (2).
 On a $\theta(s_{sc})=s_{sc}$ et $\delta(s_{sc})s_{sc}^{-1}=(z')^n$. Donc $\tau'_{G'}(s_{sc})s_{sc}^{-1}=1$ et $\rho'_{G'}(s_{sc})s_{sc}^{-1}=(z')^n$. En rempla\c{c}ant $(\tau',\rho')$ par les quatre couples possibles, on voit que l'on obtient $4$ cocycles non cohomologues, dont la restriction \`a $\Gamma_{E_{0}}$ est un homomorphisme d'ordre $4$. Ces cocycles correspondent aux $4$\'el\'ements de $\Xi_{0}$. Donc, pour chaque $\xi\in \Xi_{0}$, il y a une et une seule de nos  donn\'ees ${\bf G}'$ telle que $\xi_{{\bf G}'}=\xi$. On la note ${\bf G}'_{i,j,\xi}$ et on pose $FC^{{\cal E}}_{i,j,\xi}=FC^{st}(\mathfrak{g}'_{i,j,\xi}(F))^{Out({\bf G}'_{i,j,\xi})}$. 
 On a $(i,j,\xi)\in {\cal Y}$.   Supposons maintenant que ${\cal O}=\{\hat{\alpha}_{0},\hat{\alpha}_{1},\hat{\alpha}_{n-1},\hat{\alpha}_{n}\}$. Le fixateur de $\hat{\alpha}_{0}$ est $\Gamma_{K}$ et $K/F$ est totalement ramifi\'ee.  Comme ci-dessus, on peut remplacer $G'$ par $G'_{SC}\simeq SU_{E'/F}(n-2)$. Ce cas  est similaire au pr\'ec\'edent, l'entier $i$ devenant $1$ (cas exclu auparavant puisqu'on avait $m\geq2$). 
 
  On a associ\'e une droite $FC_{y}^{{\cal E}}$ \`a tout \'el\'ement $y\in {\cal Y}$.   On n'a pas trait\'e la donn\'ee principale ${\bf G}$. 
  
   On d\'efinit une application $\phi:{\mathbb X}\to {\mathbb Y}$ par
 $$\phi(k,h)=\left\lbrace\begin{array}{cc}((k-h)/4,(k+h)/2),& \,\,si\,\,k \equiv h\,\,mod\,\,2{\mathbb Z},\\ ((k+h+1)/4,(k-h-1)/2),&\,\,si \,\,k \not\equiv h\,\,mod\,\,2{\mathbb Z}.\\   \end{array}\right.$$
Ce sont les m\^emes formules qui d\'efinissaient $\phi^-$ en \ref{Dndeppairpadique}. On v\'erifie que $\phi$ est une bijection. Elle se
  rel\`eve naturellement en une bijection $\varphi:{\cal X}\to {\cal Y}$. On pose ${\cal X}^{st}=\emptyset$. 
  
  On utilise  l'habituel argument de dimension. On obtient que la donn\'ee principale ne peut pas intervenir, c'est-\`a-dire que $FC^{st}(\mathfrak{g}(F))=\{0\}$. Cela  ach\`eve la preuve de \ref{resultats}(2) et d\'emontre en m\^eme temps \ref{resultats}(4).  On prouve \ref{resultats}(3) comme en \ref{Dndepimppadique}. 
   Explicitons la cons\'equence de \ref{resultats}(4):
 
 (1) $FC^{st}(\mathfrak{g}(F))=\{0\}$. 

 \subsection{Forme int\'erieure du type $D_{n}$ quasi-d\'eploy\'e, $n$ impair, $E_{0}/F$ non ramifi\'ee}
 On suppose que $G^*$ est comme dans le paragraphe pr\'ec\'edent et que $G$ en est la forme int\'erieure associ\'ee \`a l'unique caract\`ere non trivial de $Z(\hat{G}_{SC})^{\Gamma_{F}}=\{1,z\}$ ou encore \`a  l'image dans $N_{\Gamma_{F}^{nr}}$ de l'\'el\'ement $\delta\theta\in \hat{\Omega}\simeq N$.  Dans les tables de Tits, le groupe est de type $^2D''_{n}$.
 
 On pose ${\cal X} =\emptyset$.   On munit le diagramme ${\cal D}_{a}$ de l'action galoisienne triviale sur $\Gamma_{E_{0}}$ et telle que $\sigma\in \Gamma_{F}-\Gamma_{E_{0}}$ agisse par $\delta$. L'ensemble $\underline{S}(G)$ est param\'etr\'e par l'ensemble des orbites de cette action galoisienne. Or toute orbite a  deux \'el\'ements et le stabilisateur d'une racine est $\Gamma_{E_{0}}$.  Alors   $FC(\mathfrak{g}_{s}({\mathbb F}_{q}))=\{0\}$ d'apr\`es le lemme \ref{orbites}. Donc $FC(\mathfrak{g}(F))=\{0\}$. Cela d\'emontre  \ref{resultats} (1).
 
 En posant ${\cal Y}=\emptyset$, les assertions \ref{resultats}(2) et (3) sont triviales. 
 
 \subsection{Type $D_{n}$ quasi-d\'eploy\'e, $n$ pair, $E/F$ ramifi\'ee}\label{Dnpairpadiqueram}
On fixe une extension quadratique $E/F$ ramifi\'ee.    On suppose que $G$ est quasi-d\'eploy\'e de type $D_{n}$ avec $n\geq4$, $n$ pair, et que $\Gamma_{F}$ agit sur le diagramme ${\cal D}$ de $G$ par l'action $\sigma\mapsto \sigma_{G}$ triviale sur $\Gamma_{E}$ et telle que $\sigma_{G}=\theta$ pour $\sigma\in \Gamma_{F}-\Gamma_{E}$. Dans les tables de Tits, le groupe est de type $C-B_{n-1}$.
On fixe un \'el\'ement $\tau\in \Gamma_{F}-\Gamma_{E}$.  L'action de $\tau$ sur $Z(G)$ fixe $z$ et \'echange $z'$ et $z''$. 
 Un \'el\'ement de $T_{AD}(F)$ s'\'ecrit $\prod_{l=1,...,n}\check{\varpi}_{l}(x_{l})$ avec $x_{l}\in F^{\times}$ pour $l=1,...,n-2$, $x_{n-1},x_{n}\in E^{\times}$ et $x_{n}=\tau(x_{n-1})$. D\'efinissons un homomorphisme $T_{ad}(F)\to E^{\times}/E^{\times,2}$ par $\prod_{l=1,...,n}\check{\varpi}_{l}(x_{l})\mapsto \prod_{l=1,...,n}x_{l}^{l}$. De cet homomorphisme se d\'eduit un isomorphisme $G_{AD}(F)/\pi(G(F))\simeq E^{\times}/E^{\times,2}$. Les images de $SO_{E/F}(2n,F)$  et de $G_{AD}(F)_{0}$ sont les m\^emes, \`a savoir  $\mathfrak{o}_{E}^{\times}/\mathfrak{o}_{E}^{\times,2}\simeq\{\pm 1\}$. On note $\Xi_{0}$ l'ensemble des \'el\'ements de $\Xi$ dont la restriction \`a $G_{AD}(F)_{0}/\pi(G(F))$ est non triviale. 
 
 On note ${\mathbb X}$ l'ensemble des couples $(k,h)\in {\mathbb N}^2$ tels que $2n=k(k+1)/2+h(h+1)/2$, $k> h$ et $k(k+1)/2$ et $h(h+1)/2$ tous deux impairs. On pose ${\cal X} =\{(k,h,\xi); (k,h)\in {\mathbb X}, \xi\in \Xi_{0}\}$.

 L'ensemble $\underline{S}(G)$ s'envoie surjectivement sur l'ensemble des couples $(a,b)\in {\mathbb N}^2$ tels que $a+b=n-1$ et $a> b$ (puisque $n$ est pair, on ne peut pas avoir $a=b$). Les fibres ont deux \'el\'ements. L'action de $G_{AD}(F)$ pr\'eserve les fibres et permute les \'el\'ements de celles-ci. Pour un sommet $s$ param\'etr\'e par $(a,b)$, on a $G_{s}\simeq (Spin(2a+1)\times Spin(2b+1))/\{1,(z,z)\}$ si $b\not=0$, $G_{s}\simeq Spin(2a+1)$ si $b=0$. Supposons $b\not=0$. On cherche les fonctions $f_{N^{a},\epsilon^{a}}\times f_{N^b,\epsilon^b}$ telles que $\epsilon^{a}(z)=\epsilon^b(z)$.  On utilise \ref{Bn}. Si  $\epsilon^{a}(z)=\epsilon^b(z)=1$, on doit avoir $2a+1=k^2$ et $2b+1=h^2$ avec $k,h\in {\mathbb N}$ tous deux impairs. Alors $2n=2a+1+2b+1=k^2+h^2$ est congru \`a $2$ modulo $8{\mathbb Z}$, ce qui est impossible puisque $n$ est pair. Supposons $\epsilon^{a}(z)=\epsilon^b(z)=-1$. C'est possible si et seulement si $2a+1=k(k+1)/2$ et $2b+1=h(h+1)/2$ avec $k,h\in {\mathbb N}$. Supposons ces conditions v\'erifi\'ees. Il y a alors une unique fonction $f_{N^{a},\epsilon^{a}}\times f_{N^b,\epsilon^b}$. L'\'el\'ement $z\in Z(G)$ appartient \`a $Z(G)^{I_{F}}$ et s'envoie sur l'\'el\'ement $(1,z)=(z,1)$ de $G_{s}$ qui agit non trivialement sur notre fonction. Celle-ci se transforme donc selon le caract\`ere non trivial de $G_{AD}(F)_{0}/\pi(G(F))$. Conform\'ement \`a \ref{actionsurFC}, pour tout $\xi\in \Xi_{0}$, la fonction donne naissance \`a un \'el\'ement de $FC(\mathfrak{g}(F))$ qui se transforme selon le caract\`ere $\xi$ de $G_{AD}(F)/\pi(G(F))$. On note $FC_{k,h,\xi}$ la droite port\'ee par cet \'el\'ement. On obtient \ref{resultats}(1).
 
  On note ${\mathbb Y}$ l'ensemble des couples $(i,j)$ avec $i,j\in {\mathbb N}$,  $2i^2+j(j+1)/2=n$ et $i$ est impair. On pose ${\cal Y} =\{(i,j,\xi);(i,j)\in {\mathbb Y}, \xi\in \Xi_{0}\}$.

 Consid\'erons un couple $(\sigma\mapsto \sigma_{G'},{\cal O})\in {\cal E}_{ell}(G)$.  L'image de l'homomorphisme $\omega_{G'}:\Gamma_{E_{G'}/G}\to \hat{\Omega}$  est normalis\'ee par $\tau_{G'}$ qui appartient \`a $\hat{\Omega}\theta$. Un tel \'el\'ement permute $\delta$ et $\delta\theta\theta'$ et fixe $\theta\theta'$.  Les seules images possibles sont donc  $\{1\}$, $\{1,\theta\theta'\}$ ou $\hat{\Omega}$ tout entier. 
 
 Supposons $E_{G'}=E$. On a alors $\tau_{G'}^2=1$ donc $\omega_{G'}(\tau)=1$ ou $\omega_{G'}(\tau)=\theta\theta'$, c'est-\`a-dire $\tau_{G'}=\theta$ ou $\theta'$. Les deux cas sont \'equivalents (conjugu\'es par $\delta\in \hat{\Omega}$), on peut supposer $\tau_{G'}=\theta$. Si ${\cal O}=\{\hat{\alpha}_{0}\}$ ou ${\cal O}=\{\hat{\alpha}_{1}\}$ (ces deux cas sont \'equivalents), on a ${\bf G}'={\bf G}$ et on ne peut rien dire \`a pr\'esent de $FC^{st}(\mathfrak{g}(F))$. 
 Si ${\cal O}=\{\hat{\alpha}_{m}\}$ avec $m\in \{2,...,n-2\}$, on a $G'_{SC}\simeq Spin_{dep}(2m)\times Spin_{E/F}(2n-2m)$. En utilisant \ref{Dndeppairpadique} (5), \ref{Dndepimppadique} (1) et \ref{An-1quasidepram} (5) par r\'ecurrence, on n'a $FC^{st}(\mathfrak{g}'(F))\not=\{0\}$ que si $m=i^2$ et  $n-m=j^2$
avec $i,j\in {\mathbb N}$, $i$ pair et $j$ impair. C'est impossible puisque $n$ est pair.  Si ${\cal O}=\{\hat{\alpha}_{n-1},\hat{\alpha}_{n}\}$, le fixateur de chacune de ces racines est $\Gamma_{E}$.  D'apr\`es le lemme \ref{centre}, on peut remplacer $G'$ par $G'_{SC}$. On a $G'_{SC}\simeq Spin_{dep}(2n-2)$ et, comme ci-dessus, $FC^{st}(\mathfrak{g}'(F))=\{0\}$.

Supposons que $E_{G'}/F$ soit de degr\'e $4$ et que $\Gamma_{E_{G'}/F}$ soit isomorphe \`a $({\mathbb Z}/2{\mathbb Z})^2$. Alors $E_{G'}=Q$ est l'extension biquadratique de $F$ et l'extension $E_{G'}/E$ est non ramifi\'ee. On pose $\Gamma_{E_{G'}/E}=\{1,\rho\}$. On a forc\'ement $\rho_{G'}=\theta\theta'$. L'\'egalit\'e $\tau_{G'}^2=1$ force $\tau_{G'}=\theta$ ou $\tau_{G'}=\theta'$. Les deux cas sont \'equivalents (conjugu\'es par $\delta$). On suppose $\tau_{G'}=\theta$. 
Notons $E'$, resp. $E''$, l'extension quadratique de $F$ telles que $\Gamma_{E_{G'}/E'}=\{1,\tau\}$, resp. $\Gamma_{E_{G'}/E''}=\{1,\tau\rho\}$. Ce sont les deux extensions quadratiques de $F$ distinctes de $E$. L'une est non ramifi\'ee, l'autre est ramifi\'ee. Pour fixer la notation, supposons que $E'$ soit l'extension  non ramifi\'ee $E_{0}$.   Si ${\cal O}=\{\hat{\alpha}_{m}\}$ avec $m\in \{2,...,n-2\}$, on a $G'_{SC}\simeq Spin_{E_{0}/F}(2m)\times Spin_{E''/F}(2n-2m)$. Puisque $E_{0}/F$ est non ramifi\'ee et $E''/F$ est ramifi\'ee, on a encore   $FC^{st}(\mathfrak{g}'(F))\not=\{0\}$ seulement si $m=i^2$ et  $n-m=j^2$
avec $i,j\in {\mathbb N}$, $i$ pair et $j$ impair. C'est impossible puisque $n$ est pair. Si ${\cal O}=\{\hat{\alpha}_{0},\hat{\alpha}_{1}\}$, le fixateur d'une de ces racines est $\Gamma_{E_{0}}$ et le lemme \ref{centre} exclut ce cas. Si  ${\cal O}=\{\hat{\alpha}_{n-1},\hat{\alpha}_{n}\}$, on a $G'_{SC}\simeq Spin_{E_{0}/F}(2n-2)$ et, comme ci-dessus, $FC^{st}(\mathfrak{g}'(F))=\{0\}$.

Supposons que $E_{G'}/F$ soit cyclique de degr\'e $4$.  Puisque $E\subset E_{G'}$, $E_{G'}/F$ est totalement ramifi\'ee.  D'apr\`es \ref{extensionsbiquadratiques}(2) et (3), on a $\delta_{4}(q-1)=1$. Supposons cette condition v\'erifi\'ee et introduisons l'extension biquadratique $Q_{E}$ de $E$ et les notations de \ref{extensionsbiquadratiques}(6), en particulier les \'el\'ements $\rho'$ et $\rho''$ de $\Gamma_{Q_{E}/E}$.   L'extension $E_{G'}$ peut \^etre $K'$ ou $K''$. L'\'el\'ement $\tau_{G'}$ est d'ordre $4$ ce qui impose $\omega_{G'}(\tau)=\delta$ ou $\delta\theta\theta'$, c'est-\`a-dire $\tau_{G'}=\delta\theta$ ou $\delta\theta'$. Ces deux cas sont \'equivalents (conjugu\'es par $\delta$). Supposons $\tau_{G'}=\delta\theta$. Supposons ${\cal O}=\{\hat{\alpha}_{m},\hat{\alpha}_{n-m}\}$ avec $m\in \{2,...,n/2-1\}$. Le fixateur de chacune de ces racines est $\Gamma_{E}$.   D'apr\`es le lemme \ref{centre}, on  peut remplacer $G'$ par $G'_{SC}$ qui est isomorphe \`a $Res_{E/F}(Spin_{E_{G'}/E}(2m))\times SU_{E/F}(n-2m)$. L'extension $E_{G'}/E$ est ramifi\'ee. En utilisant (1) ci-dessous, \ref{Dnimppadiqueram} (1) et \ref{An-1quasidepram} (5) par r\'ecurrence, on a $FC^{st}(\mathfrak{g}'(F))\not=\{0\}$ si et seulement si $m=i^2$ avec $i\in {\mathbb N}$ et $i$ impair  et $n-2m=j(j+1)/2$ avec $j\in {\mathbb N}$.  Supposons ces conditions v\'erifi\'ees. Alors  $FC^{st}(\mathfrak{g}'(F))$ est une  droite. On a $Out({\bf G}')=\{1,\theta\theta'\}$ et on voit que ce groupe agit trivialement sur $FC^{st}(\mathfrak{g}'(F))$.  Calculons $\xi_{{\bf G}'}$. On peut choisir $s_{sc}$ comme en \ref{Dndeppairpadique}(2). 
  On a $\theta(s_{sc})=s_{sc}$ et $\delta(s_{sc})s_{sc}^{-1}=z'z^{n/2}$. On en d\'eduit $\tau_{G'}(s_{sc})s_{sc}^{-1}=z'z^{n/2}$, puis $\tau^2_{G'}(s_{sc})s_{sc}^{-1}=z$. Si $E_{G'}=K'$, on a  $\rho'_{G'}=1$ donc $\rho'_{G'}(s_{sc})s_{sc}^{-1}=1$. Si $E_{G'}=K''$, on a $\rho''_{G'}=1$ donc  $\rho''_{G'}(s_{sc})s_{sc}^{-1}=1$ mais $\rho''=\rho'\tau^2$ donc $\rho'_{G'}(s_{sc})s_{sc}^{-1}=z$. On voit que l'on obtient deux cocycles distincts dont la restriction \`a $I_{E}$ est non triviale. Ces cocycles correspondent aux deux \'el\'ement de $\Xi_{0}$. Donc, pour tout $\xi\in \Xi_{0}$, il y a une et une seule de nos donn\'ees (c'est-\`a-dire une et une seule des deux extensions $K'$ et $K''$ possibles) telle que $\xi_{{\bf G}'}=\xi$. On a $(i,j,\xi)\in {\cal Y}$. On pose ${\bf G}'_{i,j,\xi}={\bf G}'$ et $FC^{{\cal E}}_{i,j,\xi}=FC^{st}(\mathfrak{g}'_{i,j,\xi}(F))^{Out({\bf G}'_{i,j,\xi})}$.  Supposons maintenant ${\cal O}=\{\hat{\alpha}_{0},\hat{\alpha}_{1},\hat{\alpha}_{n-1},\hat{\alpha}_{n}\}$. Le fixateur d'une de ces racines est $\Gamma_{E_{G'}}$. Puisque $E_{G'}/F$ est totalement ramifi\'ee, on peut encore remplacer $G'$ par $G'_{SC}\simeq SU_{E/F}(n-2)$. Ce cas est similaire au pr\'ec\'edent, l'entier $i$ devenant $1$.

Supposons que l'image par $\sigma\mapsto \omega_{G'}(\sigma)$ de $\Gamma_{K/E}$ soit $\hat{\Omega}$ tout entier. Alors $E_{G'}$ est l'extension biquadratique $Q_{E}$ de $E$.   L'homomorphisme $\sigma\mapsto \sigma_{G'}$ est un isomorphisme de $\Gamma_{Q_{E}/F}$ sur $Aut(\hat{D}_{a})$. Ce dernier groupe n'est pas commutatif. D'apr\`es \ref{extensionsbiquadratiques}(6) et (7), on a $\delta_{4}(q-1)=0$. Supposons cette condition v\'erifi\'ee. On utilise les notations $\rho',\rho_{0},\rho''$ de \ref{extensionsbiquadratiques}(7) et on suppose que $\tau$ v\'erifie les conditions de cette assertion.    L'\'el\'ement $\rho_{0}$ est central dans $\Gamma_{Q_{E}/F}$ donc forc\'ement $\rho_{0,G'}=\theta\theta'$. On a alors $\rho'_{G'}=\delta$ ou $\delta\theta\theta'$. On a $\tau^2\in \Gamma_{Q_{E}}$ donc $\tau_{G'}^2=1$. Puisque  $\tau_{G'}\in \hat{\Omega}\theta$, on a forc\'ement $\tau_{G'}=\theta$ ou $\theta'$. Ces deux cas sont \'equivalents (conjugu\'es par $\delta$). On suppose $\tau_{G'}=\theta$.  Par contre, la valeur de $\rho'_{G'}$ ne change pas par \'equivalence, on a donc deux donn\'ees  ${\bf G}'$ correspondant \`a ces deux valeurs de $\rho'_{G'}$. Supposons ${\cal O}=\{\hat{\alpha}_{m},\hat{\alpha}_{n-m}\}$ avec $m\in \{2,...,n/2-1\}$. Le fixateur de $\hat{\alpha}_{m}$ dans $\Gamma_{Q_{E}/F}$ est le groupe $\{1,\tau,\rho_{0},\tau\rho_{0}\}$, c'est-\`a-dire $\Gamma_{Q_{E}/E'}$, o\`u $E'$ est l'extension quadratique ramifi\'ee de $F$ diff\'erente de $E$ (cf. \ref{extensionsbiquadratiques}(7)).  Notons $K$ l'extension de $E'$ telle que $\Gamma_{Q_{E}/K}=\{1,\tau\}$. L'extension $K/E'$ est ramifi\'ee (car l'extension non ramifi\'ee est fix\'ee par $\rho_{0}$ et non pas par $\tau$). Puisque $E'/F$ est ramifi\'ee, le lemme \ref{centre} nous permet de remplacer $G'$ par $G'_{SC}$ qui est isomorphe \`a $Res_{E'/F}(Spin_{K/E'}(2m))\times SU_{E'/F}(n-2m)$.  En utilisant  (1) ci-dessous, \ref{Dnimppadiqueram} (1) et \ref{An-1quasidepram} (5) par r\'ecurrence, on a $FC^{st}(\mathfrak{g}'(F))\not=\{0\}$ si et seulement si $m=i^2$ avec $i\in {\mathbb N}$ et $i$ impair   et $n-2m=j(j+1)/2$ avec $j\in {\mathbb N} $.  Supposons ces conditions v\'erifi\'ees. Alors  $FC^{st}(\mathfrak{g}'(F))$ est une  droite. On a $Out({\bf G}')=\{1,\theta\theta'\}$ et on voit que ce groupe agit trivialement sur $FC^{st}(\mathfrak{g}'(F))$.  Calculons $\xi_{{\bf G}'}$. L'\'el\'ement $s_{sc}$ est  le m\^eme qu'en \ref{Dndeppairpadique}(2) et on calcule $\tau_{G'}(s_{sc})s_{sc}^{-1}=1$, $\rho_{0}(s_{sc})s_{sc}^{-1}= z$, $\rho'_{G'}(s_{sc})s_{sc}^{-1}=z'z^{n/2}$ si $\rho_{G'}=\delta$, $\rho'_{G'}(s_{sc})s_{sc}^{-1}=z''z^{n/2}$ si $\rho_{G'}=\delta\theta\theta'$. On voit que l'on obtient deux cocycles distincts dont la restriction \`a $I_{E}$ est non triviale. Ces cocycles correspondent aux deux \'el\'ement de $\Xi_{0}$. Donc, pour tout $\xi\in \Xi_{0}$, il y a une et une seule de nos donn\'ees (c'est-\`a-dire une et une seule des deux  actions possibles de $\rho'_{G'}$) telle que $\xi_{{\bf G}'}=\xi$. On a $(i,j,\xi)\in {\cal Y}$. On pose ${\bf G}'_{i,j,\xi}={\bf G}'$ et $FC^{{\cal E}}_{i,j,\xi}=FC^{st}(\mathfrak{g}'_{i,j,\xi}(F))^{Out({\bf G}'_{i,j,\xi})}$.  
  Supposons maintenant ${\cal O}=\{\hat{\alpha}_{0},\hat{\alpha}_{1},\hat{\alpha}_{n-1},\hat{\alpha}_{n}\}$. Le fixateur d'une de ces racines est $\Gamma_{K}$. Puisque $K/F$ est totalement ramifi\'ee, on  peut encore  remplacer $G'$ par $G'_{SC}\simeq SU_{E/F}(n-2)$. Ce cas est similaire au pr\'ec\'edent, l'entier $i$ devenant $1$.
  
 On a associ\'e \`a  tout $y\in {\cal Y}$ une droite $FC^{{\cal E}}_{y}$ (avec deux constructions diff\'erentes selon la valeur de $\delta_{4}(q-1)$). On n'a pas trait\'e la donn\'ee principale ${\bf G}$. 

   On d\'efinit une bijection $\phi:{\mathbb X}\to {\mathbb Y}$ comme en \ref{Dndeppairpadique}. Elle se rel\`eve naturellement en une bijection $\varphi:{\cal X}\to {\cal Y}$. On pose ${\cal X}^{st}=\emptyset$. On peut alors achever la preuve de \ref{resultats}(2), (3) et (4) comme en 
    \ref{Dndepimppadique}. Explicitons la cons\'equence de \ref{resultats}(4):
    
    (1) $FC^{st}(\mathfrak{g}(F))=\{0\}$.

   \subsection{Forme int\'erieure du type $D_{n}$ quasi-d\'eploy\'e, $n$ pair, $E/F$ ramifi\'ee}
  On suppose que $G^*$ est comme ci-dessus. Le groupe $N$ est ici le groupe d'automorphismes  du diagramme de Dynkin local de type $C-B_{n-1}$. On a $N\simeq {\mathbb Z}/2{\mathbb Z}$.  On suppose que $G$ est la forme int\'erieure de $G^*$ param\'etr\'ee par l'\'el\'ement non trivial de ce groupe, ou encore par le caract\`ere non trivial de $Z(\hat{G}_{SC})^{\Gamma_{F}}=\{1,z\}$. 
  
  On pose ${\cal X} =\emptyset$.  
   L'ensemble $\underline{S}(G)$ est param\'etr\'e par l'ensemble des orbites de $N$ agissant sur le diagramme de Dynkin local. Parce que $n$ est pair, ces orbites ont deux \'el\'ements. D'apr\`es le lemme \ref{orbites}, on a $FC(\mathfrak{g}_{s}({\mathbb F}_{q}))=\{0\}$ pour tout $s\in \underline{S}(G)$. Cela entra\^{\i}ne $FC(\mathfrak{g}(F))=\{0\}$, d'o\`u la relation \ref{resultats}(1). 
   
   En posant ${\cal Y}=\emptyset$, les assertions \ref{resultats}(2) et (3) sont triviales.

  \subsection{Type $D_{n}$ quasi-d\'eploy\'e, $n$ impair, $E/F$ ramifi\'ee}\label{Dnimppadiqueram}
  On fixe une extension quadratique $E/F$ ramifi\'ee.    On suppose que $G$ est quasi-d\'eploy\'e de type $D_{n}$ avec $n\geq4$, $n$ impair, et que $\Gamma_{F}$ agit sur le diagramme ${\cal D}$ de $G$ par l'action $\sigma\mapsto \sigma_{G}$ triviale sur $\Gamma_{E}$ et telle que $\sigma_{G}=\theta$ pour $\sigma\in \Gamma_{F}-\Gamma_{E}$. Dans les tables de Tits, le groupe est de type $C-B_{n-1}$.
On fixe un \'el\'ement $\tau\in \Gamma_{F}-\Gamma_{E}$ qui agit trivialement sur l'extension quadratique non ramifi\'ee $E_{0}$ de $F$.  L'action de $\tau$ sur $Z(G)$  envoie $z'$ sur $(z')^{-1}$. Comme en \ref{Dnimppadiquenonram}, le groupe $G_{AD}(F)/\pi(G(F))$ est isomorphe \`a $(E^{\times}\times F^{\times})/\iota(E^{\times}\times F^{\times})$. Mais la structure de ce groupe n'est plus la m\^eme. On voit que l'homomorphisme
$$E^{\times}\times F^{\times}\to ({\mathbb Z}/2{\mathbb Z})\times F^{\times}/norme_{E/F}(E^{\times})$$
qui \`a $(e,f)$ associe le couple form\'e de $val_{E}(e)$ modulo $2{\mathbb Z}$ et de l'image de $f$ dans $F^{\times}/norme_{E/F}(E^{\times})$ est un isomorphisme. L'image de $G_{AD}(F)_{0}$ est le sous-groupe $\{0\}\times F^{\times}/norme_{E/F}(E^{\times})$. C'est aussi l'image de $SO_{E/F}(2n,F)$. Notons $\Xi^+$, resp. $\Xi^-$, l'ensemble des $\xi\in \Xi$ dont la restriction \`a $G_{AD}(F)_{0}/\pi(G(F))$ est triviale, resp. non triviale.

On note ${\mathbb X}^+$ l'ensemble des couples $(k,h)\in {\mathbb N}^2$ tels que $k^2+h^2=2n$, $k$ et $h$ sont impairs et $k\geq h$. Si $\delta_{\square}(n)=0$, on   note ${\cal X}^+$ l'ensemble des triplets $(k,h,\xi)$ avec $(k,h)\in {\mathbb X}^+$ et $\xi\in \Xi^+$. Si $\delta_{\square}(n)=1$, il y a dans ${\mathbb X}^+$ un unique couple $(k,h)$ avec $k=h$. On le note $(k^{st},k^{st})$. On  note ${\cal X}^+$ la r\'eunion de $\{(k^{st},k^{st})\}$ et de l'ensemble des  triplets $(k,h,\xi)$ avec $(k,h)\in {\mathbb X}^+$, $k>h$ et $\xi\in \Xi^+$. On note ${\mathbb X}^-$ l'ensemble des couples $(k,h)\in {\mathbb N}^2$ tels que $2n=k(k+1)/2+h(h+1)/2$, $k(k+1)/2$ et $h(h+1)/2$ sont impairs  et $k\geq h$. Si $\delta_{\triangle}(n)=0$, on note ${\cal X}^-$ l'ensemble des triplets $(k,h,\xi)$ avec $(k,h)\in {\mathbb X}^-$ et $\xi\in \Xi^-$. Si $\delta_{\triangle}(n)=1$, il y a dans ${\mathbb X}^-$ un unique couple $(k,h)$ avec $k=h$. On le note $(k_{0},k_{0})$. On note ${\cal X}^-$ la r\'eunion de $\{(k_{0},k_{0})\}$ et de l'ensemble des triplets $(k,h,\xi)$ avec $(k,h)\in {\mathbb X}^-$, $k>h$ et $\xi\in \Xi^-$. 
On pose ${\cal X}={\cal X}^+\sqcup {\cal X}^-$ et $d_{x}=1$ pour tout $x\in {\cal X}$.

 L'ensemble $\underline{S}(G)$ s'envoie surjectivement sur l'ensemble des couples $(a,b)\in {\mathbb N}^2$ tels que $a+b=n-1$ et $a\geq b$. Les fibres ont deux \'el\'ements, sauf au-dessus du couple $((n-1)/2,(n-1)/2)$ o\`u elle n'en a qu'un. L'action de $G_{AD}(F)$ pr\'eserve les fibres et permute les \'el\'ements de celles-ci. Pour un sommet $s$ param\'etr\'e par $(a,b)$, on a $G_{s}\simeq (Spin(2a+1)\times Spin(2b+1))/\{1,(z,z)\}$ si $b\not=0$, $G_{s}\simeq Spin(2a+1)$ si $b=0$. Supposons $b\not=0$ et $a\not=b$. On cherche les  couples de fonctions $(f_{N^{a},\epsilon^{a}}, f_{N^b,\epsilon^b})$ telles que $\epsilon^{a}(z)=\epsilon^b(z)$.  On utilise \ref{Bn}. Il y a un tel couple avec $\epsilon^{a}(z)=\epsilon^b(z)=1$ si et seulement si $2a+1=k^2$, $2b+1=h^2$ pour deux entiers $k,h$ forc\'ement impairs. Supposons ces conditions v\'erifi\'ees. Il y a une  unique fonction $f_{N^{a},\epsilon^{a}}\times f_{N^b,\epsilon^b}$. Puisque $\epsilon^{a}(z)=\epsilon^b(z)=1$, elle est invariante par l'action naturelle du groupe $SO(2a+1,{\mathbb F}_{q})\times SO(2b+1,{\mathbb F}_{q})$ et cette action est aussi celle de $G_{AD}(F)_{0}/\pi(G(F))$. Conform\'ement \`a \ref{actionsurFC}, pour tout $\xi\in \Xi^+$, la fonction d\'etermine un \'el\'ement de $FC(\mathfrak{g}(F))$ se transformant par $G_{AD}(F)/\pi(G(F))$ selon le caract\`ere $\xi$. On  a $(k,h,\xi)\in {\cal X}^+$. On note  $FC_{k,h,\xi}$ la droite port\'ee par l'\'el\'ement pr\'ec\'edent. Il y a un couple de fonctions $(f_{N^{a},\epsilon^{a}}, f_{N^b,\epsilon^b})$ telles que $\epsilon^{a}(z)=\epsilon^b(z)=-1$ si et seulement si $2a+1=k(k+1)/2$, $2b+1=h(h+1)/2$ pour deux entiers $k,h$ avec $k(k+1)/2$ et $h(h+1)/2$ impairs. Supposons ces conditions v\'erifi\'ees. Il y a une  unique fonction $f_{N^{a},\epsilon^{a}}\times f_{N^b,\epsilon^b}$. Puisque $\epsilon^{a}(z)=\epsilon^b(z)=-1$, elle se transforme selon le caract\`ere non trivial de  $G_{AD}(F)_{0}/\pi(G(F))$. Conform\'ement \`a \ref{actionsurFC}, pour tout $\xi\in \Xi^-$, la fonction d\'etermine un \'el\'ement de $FC(\mathfrak{g}(F))$ se transformant par $G_{AD}(F)/\pi(G(F))$ selon le caract\`ere $\xi$. On a $(k,h,\xi)\in {\cal X}^-$ et on note  $FC_{k,h,\xi}$ la droite port\'ee par l'\'el\'ement pr\'ec\'edent. Supposons maintenant $b=0$. On voit que le r\'esultat est le m\^eme, l'entier $h$ devenant $1$ dans les deux cas possibles. 
 Supposons enfin $a=b=(n-1)/2$. Le seul changement est que le sommet $s$ est conserv\'e par $G_{AD}(F)/\pi(G(F))$ tout entier, donc chacune de nos fonctions d\'etermine un seul \'el\'ement de $FC(\mathfrak{g}(F))$. Si $n $ est un carr\'e, c'est-\`a-dire $n=(k^{st})^2$, on a une fonction  $f_{N^{a},\epsilon^{a}}\times f_{N^b,\epsilon^b}$ avec $\epsilon^{a}(z)=\epsilon^b(z)=1$. Elle d\'etermine un \'el\'ement de $FC(\mathfrak{g}(F))$ qui se transforme selon un certain caract\`ere de $G_{AD}(F)/\pi(G(F))$ trivial sur $G_{AD}(F)_{0}/\pi(G(F))$. On note  $FC_{k^{st},k^{st}}$ la droite port\'ee par cet \'el\'ement (on a  $(k^{st},k^{st})\in {\cal X}^+$). Si $\delta_{\triangle}(n)=1$, c'est-\`a-dire $n=k_{0}(k_{0}+1)/2$, on a une fonction  $f_{N^{a},\epsilon^{a}}\times f_{N^b,\epsilon^b}$ avec $\epsilon^{a}(z)=\epsilon^b(z)=-1$. Elle d\'etermine un \'el\'ement de $FC(\mathfrak{g}(F))$ qui se transforme selon un certain caract\`ere de $G_{AD}(F)/\pi(G(F))$ non trivial sur $G_{AD}(F)_{0}/\pi(G(F))$. On note   $FC_{k_{0},k_{0}}$ la droite port\'ee par cet \'el\'ement (on a $(k_{0},k_{0})\in {\cal X}^-$).  Cette description d\'emontre \ref{resultats}(1). 

 On note ${\mathbb Y}^+$ l'ensemble des couples $(i,j)\in {\mathbb N}^2$ tels que $i^2+j^2=n$ et $i>j$. Remarquons que $i$ et $j$ sont de parit\'e distinctes. Si $\delta_{\square}(n)=0$, on note ${\cal Y}^+$ l'ensemble des triplets $(i,j,\xi)$ avec $(i,j)\in {\mathbb Y}^+$ et $\xi\in \Xi^+$. Si $\delta_{\square}(n)=1$, il y a dans ${\mathbb Y}^+$ un unique couple $(i,j)$ avec $j=0$, \`a savoir $(k^{st},0)$. On note ${\cal Y}^+$ la r\'eunion de $ \{(k^{st},0)\}$ et de l'ensemble des  triplets $(i,j,\xi)$ avec $(i,j)\in {\mathbb Y}^+$, $j\not=0$ et $\xi\in \Xi^+$. On note ${\mathbb Y}^-$ l'ensemble des couples $(i,j)$ avec $i,j\in {\mathbb N}$,   $2i^2+j(j+1)/2=n$ et $i$ est pair. Si $\delta_{\triangle}(n)=0$, on note ${\cal Y}^-$ l'ensemble des triplets $(i,j,\xi)$ avec $(i,j)\in {\mathbb Y}^-$ et $\xi\in \Xi^-$. Si $\delta_{\triangle}(n)=1$, il y a un couple $(i,j)\in {\mathbb Y}^-$ tel que $i=0$, \`a savoir le couple $(0,k_{0})$.  On note ${\cal Y}^-$ la r\'eunion de $\{(0,k_{0})\}$ et de l'ensemble des triplets $(i,j,\xi)$ avec $(i,j)\in {\mathbb Y}^-$, $i\not=0$  et $\xi\in \Xi^-$.

 Consid\'erons un couple $(\sigma\mapsto \sigma_{G'},{\cal O})\in {\cal E}_{ell}(G)$ et une orbite ${\cal O}$.   On a $\tau_{G'}\in \hat{\Omega}\theta$ et on voit que tout \'el\'ement de cet ensemble est de carr\'e $1$. Donc $\tau_{G'}^2=1$. 

Supposons $E_{G'}=E$. On voit qu'\`a \'equivalence pr\`es (conjugaison par $\delta\theta\in \hat{\Omega}$), on a $\tau_{G'}=\theta$ ou $\tau_{G'}=\delta$. Supposons d'abord $\tau_{G'}=\theta$. Si ${\cal O}=\{\hat{\alpha}_{0}\}$ ou ${\cal O}=\{\hat{\alpha}_{1}\}$ (ces deux cas sont \'equivalents, conjugu\'es par $\theta\theta'\in \hat{\Omega}$), la donn\'ee ${\bf G}'$ est la donn\'ee principale ${\bf G}$ et, \`a ce point, on ne peut rien dire de l'espace $FC^{st}(\mathfrak{g}(F))$. Supposons ${\cal O}=\{\hat{\alpha}_{m}\}$, avec $m\in \{2,...,n-2\}$.  Alors $G'_{SC}\simeq Spin_{dep}(2m)\times Spin_{E/F}(2n-2m)$. On utilise \ref{Dndeppairpadique} (5), \ref{Dndepimppadique} (1), \ref{Dnpairpadiqueram} (1) et (1) ci-dessous par r\'ecurrence. On a $FC^{st}(\mathfrak{g}'(F))\not=\{0\}$ si et seulement si $m=(i')^2$ et $n-m=(j')^2$ avec $i'$ pair et $j'$ impair. Supposons ces conditions v\'erifi\'ees. L'espace $FC^{st}(\mathfrak{g}'(F))$ est une droite. Le groupe $Out({\bf G}')$ est $\{1,\theta\theta'\}$, qui agit trivialement sur $FC^{st}(\mathfrak{g}'(F))$. On calcule
$$s_{sc}=\prod_{l=1,...,m}\hat{\hat{\alpha}}_{l}((-1)^l).$$
On voit que $\tau_{G'}(s_{sc})s_{sc}^{-1}=1$ donc $\xi_{{\bf G}'}={\bf 1}$. On note $(i,j)$ l'unique couple $(i',j')$ ou $(j',i')$ tel que $i\geq j$. Remarquons que $(i,j)$ d\'etermine $(i',j')$ puisque $i'$ est pair tandis que $j'$ est impair. Le triplet $(i,j,{\bf 1})$  appartient \`a ${\cal Y}^+$. On pose ${\bf G}'_{i,j,{\bf 1}}={\bf G}'$ et $FC^{{\cal E}}_{i,j,{\bf 1}}=FC^{st}(\mathfrak{g}'_{i,j,{\bf 1}}(F))^{Out({\bf G}'_{i,j,{\bf 1}})}$. Supposons ${\cal O}=\{\hat{\alpha}_{n-1},\hat{\alpha}_{n}\}$. Le fixateur de chacune de ces racines est $\Gamma_{E}$. Le lemme \ref{centre} permet de 
 remplacer $G'$ par $G'_{SC}$. Le r\'esultat est le m\^eme que ci-dessus, l'entier $j'$ devenant $1$. Supposons maintenant que $\tau_{G'}=\delta$. Supposons ${\cal O}=\{\hat{\alpha}_{m},\hat{\alpha}_{n-m}\}$ avec $m\in \{2,...,(n-1)/2\}$. De nouveau, on peut remplacer $G'$ par $G'_{SC}\simeq Res_{E/F}(Spin_{dep}(2m))\times SU_{E/F}(n-2m)$. On utilise \ref{Dndeppairpadique} (5), \ref{Dndepimppadique} (1) et \ref{An-1quasidepram} (5) par r\'ecurrence.
 On a $FC^{st}(\mathfrak{g}'(F))\not=\{0\}$ si et seulement si $m=i^2$ pour un entier $i$ pair et $n-2m=j(j+1)/2$ pour un entier $j\in {\mathbb N}$. Supposons ces conditions v\'erifi\'ees.  L'espace $FC^{st}(\mathfrak{g}'(F))$ est une droite. Le groupe $Out({\bf G}')$ est $\{1,\theta\theta'\}$, qui agit trivialement sur $FC^{st}(\mathfrak{g}'(F))$. Calculons $\xi_{{\bf G}'}$. L'\'el\'ement $s_{sc}$ est le m\^eme qu'en \ref{Dndeppairpadique}(2). 
  On a $\theta(s_{sc})=s_{sc}$ et $\delta(s_{sc})s_{sc}^{-1}=(z')^n$. D'o\`u $\tau_{G'}(s_{sc})s_{sc}^{-1}=(z')^n$ et $\sigma_{G'}(s_{sc})s_{sc}^{-1}=1$ pour $\sigma\in \Gamma_{E}$. Le cocycle  ainsi d\'efini est non trivial sur $I_{F}$ et d\'etermine un certain \'el\'ement $\xi\in \Xi^-$. Le triplet $(i,j,\xi)$ appartient \`a ${\cal Y}^-$. On pose ${\bf G}'_{i,j,\xi}={\bf G}'$  et $FC^{{\cal E}}_{i,j,\xi}=FC^{st}(\mathfrak{g}'_{i,j,\xi}(F))^{Out({\bf G}'_{i,j,\xi})}$. Supposons ${\cal O}=\{\hat{\alpha}_{0},\hat{\alpha}_{n}\}$ ou ${\cal O}=\{\hat{\alpha}_{1},\hat{\alpha}_{n-1}\}$. Ces deux cas sont \'equivalents (conjugu\'es par $\theta\theta'$). Le r\'esultat est le m\^eme que pr\'ec\'edemment, le couple $(i,j)$ devenant $(0,k_{0})$. On   pose simplement ${\bf G}'_{0,k_{0}}={\bf G}'$  et $FC^{{\cal E}}_{0,k_{0}}=FC^{st}(\mathfrak{g}'_{0,k_{0}}(F))^{Out({\bf G}'_{0,k_{0}})}$. 
 
 Supposons que $E_{G'}/E$ soit quadratique.   Puisque $\tau_{G'}^2=1$, $E_{G'}$ est forc\'ement l'extension biquadratique $Q$ de $F$. On pose $\Gamma_{E_{G'}/E}=\{1,\rho\}$ et on suppose que $\tau$ est trivial sur l'extension quadratique non ramifi\'ee $E_{0}/F$.  L'image de $\Gamma_{E_{G'}/E}$ dans $\hat{\Omega}$ par l'homomorphisme $\omega_{G'}$ est $\{1,\theta\theta'\}$, c'est-\`a-dire $\rho_{G'}=\theta\theta'$. Comme pr\'ec\'edemment, on peut avoir $\tau_{G'}=\theta$ ou $\tau_{G'}=\delta$.  Supposons d'abord $\tau_{G'}=\theta$. Supposons ${\cal O}=\{\hat{\alpha}_{m}\}$, avec $m\in \{2,...,n-2\}$.  Alors $G'_{SC}\simeq Spin_{E_{0}/F}(2m)\times Spin_{E'/F}(2n-2m)$, o\`u $E'$ est l'extension quadratique ramifi\'ee de $F$ diff\'erente de $E$. On a $FC^{st}(\mathfrak{g}'(F))\not=\{0\}$ si et seulement si $m=(i')^2$ et $n-m=(j')^2$ avec $i'$ pair et $j'$ impair. Supposons ces conditions v\'erifi\'ees. On note  encore $(i,j)$ l'unique couple $(i',j')$ ou $(j',i')$ tel que $i\geq j$.  Le calcul se poursuit comme plus haut. La diff\'erence est dans le caract\`ere $\xi_{{\bf G}'}$. On a cette fois $\tau_{G'}(s_{sc})s_{sc}^{-1}=1$ mais $\rho_{G'}(s_{sc})s_{sc}^{-1}=z$. Ce cocycle est non trivial mais trivial sur $I_{F}$. Il d\'etermine l'\'el\'ement non trivial $\xi$ de $\Xi^+$. Le triplet $(i,j,\xi)$ appartient \`a ${\cal Y}^+$. On pose ${\bf G}'_{i,j,\xi}={\bf G}'$  et $FC^{{\cal E}}_{i,j,\xi} =FC^{st}(\mathfrak{g}'_{i,j,\xi}(F))^{Out({\bf G}'_{i,j,\xi})}$. Supposons ${\cal O}=\{\hat{\alpha}_{n-1},\hat{\alpha}_{n}\}$. On voit que ce cas est similaire au pr\'ec\'edent, l'entier $j'$ devenant $1$. Supposons ${\cal O}=\{\hat{\alpha}_{0},\hat{\alpha}_{1}\}$. Alors le stabilisateur de $\hat{\alpha}_{0}$ est $\Gamma_{E_{0}}$. Puisque $E_{0}/F$ est non ramifi\'ee, ce cas est exclu par le lemme \ref{centre}. Supposons maintenant  que $\tau_{G'}=\delta$. Supposons ${\cal O}=\{\hat{\alpha}_{m},\hat{\alpha}_{n-m}\}$ avec $m\in \{2,...,(n-1)/2\}$. Le  fixateur de chacune de ces racines est $\Gamma_{E}$ et on peut remplacer $G'$ par $G'_{SC}\simeq Res_{E/F}(Spin_{Q/E}(2m))\times SU_{E/F}(n-2m)$. On a $FC^{st}(\mathfrak{g}'(F))\not=\{0\}$ si et seulement si $m=i^2$ pour un entier $i$ pair et $n-m=j(j+1)/2$ pour un entier $j\in {\mathbb N}$. Supposons ces conditions v\'erifi\'ees.    Le calcul se poursuit comme ci-dessus, le changement portant sur le caract\`ere $\xi_{{\bf G}'}$. On a $\tau_{G'}(s_{sc})s_{sc}^{-1}=(z')^n$ et $\rho_{G'}(s_{sc})s_{sc}^{-1}=z$.   Le cocycle  ainsi d\'efini est non trivial sur $I_{F}$  mais non cohomologue au cocycle obtenu plus haut. Il d\'etermine un  \'el\'ement $\xi\in \Xi^-$ qui est l'autre \'el\'ement que celui obtenu plus haut.  Le triplet $(i,j,\xi)$ appartient \`a ${\cal Y}^-$. On pose ${\bf G}'_{i,j,\xi}={\bf G}'$   et   $FC^{{\cal E}}_{i,j,\xi} =FC^{st}(\mathfrak{g}'_{i,j,\xi}(F))^{Out({\bf G}'_{i,j,\xi})}$. Supposons enfin ${\cal O}=\{\hat{\alpha}_{0},\hat{\alpha}_{1},\hat{\alpha}_{n-1},\hat{\alpha}_{n}\}$. Le fixateur d'une de ces racines est $\Gamma_{Q}$. Puisque $Q/F$ n'est pas totalement ramifi\'ee, ce cas est exclu par le lemme \ref{centre}. Remarquons que, contrairement au cas $E_{G'}=E$, on n'a pas obtenu ici de donn\'ee param\'etr\'ee par l'\'eventuel \'el\'ement $(i,j)\in {\mathbb Y}^-$ tel que $i=0$. 
 
 On a associ\'e une droite $FC^{{\cal E}}_{y}$ \`a tout \'el\'ement $y\in {\cal Y}$, sauf \`a l'\'el\'ement $(k^{st},0)$ dans le cas o\`u $\delta_{\square}(n)=1$. On n'a pas trait\'e le cas de la donn\'ee principale ${\bf G}$.  
 
 On d\'efinit des bijections $\phi^{\pm}:{\mathbb X}^{\pm}\to {\mathbb Y}^{\pm}$ comme en \ref{Dndeppairpadique}. Il s'en d\'eduit des bijections $\varphi^{\pm}:{\cal X}^{\pm}\to {\cal Y}^{\pm}$ que l'on r\'eunit en une bijection $\varphi:{\cal X}\to {\cal Y}$. Si $\delta_{\square}(n)=0$, on pose ${\cal X}^{st}=\emptyset$. Si $\delta_{\square}(n)=1$, on pose ${\cal X}^{st}=\{(k^{st},k^{st})\}$. On peut alors achever la preuve de \ref{resultats}(2), (3) et (4) de la m\^eme fa\c{c}on qu'en \ref{Dndeppairpadique}. En particulier, si $\delta_{\square}(n)=1$, on pose $FC^{{\cal E}}_{k^{st},0}=FC^{st}(\mathfrak{g}(F))$. Explicitons la cons\'equence de \ref{resultats}(4):
 
 (1) $dim(FC^{st}(\mathfrak{g}(F))=\delta_{\square}(n)$.

  \subsection{Forme int\'erieure du type $D_{n}$ quasi-d\'eploy\'e, $n$ impair, $E/F$ ramifi\'ee} 
  
   On suppose que $G^*$ est comme ci-dessus. Le groupe $N$ est ici le groupe d'automorphismes  du diagramme de Dynkin local de type $C-B_{n-1}$. On a $N\simeq {\mathbb Z}/2{\mathbb Z}$.  On suppose que $G$ est la forme int\'erieure de $G^*$ param\'etr\'ee par l'\'el\'ement non trivial de ce groupe, ou encore par le caract\`ere non trivial de $Z(\hat{G}_{SC})^{\Gamma_{F}}=\{1,z\}$. 
  
   Si $\delta_{\square}(n)=0$, on pose ${\cal X}^{+} =\emptyset$. Si $\delta_{\square}(n)=1$,  on \'ecrit $n=(k^{st})^2$. On pose ${\cal X}^{+}=\{(k^{st},k^{st})\}$. 
   Si $\delta_{\triangle}(n)=0$, on pose ${\cal X}^{-}= \emptyset$. Si $\delta_{\triangle}(n)=1$, on \'ecrit $n=k_{0}(k_{0}+1)/2$ et on pose ${\cal X}^-=\{(k_{0},k_{0})\}$. On pose ${\cal X}={\cal X}^+\sqcup {\cal X}^-$ et $d_{x}=1$ pour tout $x\in {\cal X}$.

   L'ensemble $\underline{S}(G)$ est param\'etr\'e par l'ensemble des orbites de $N$ agissant sur le diagramme de Dynkin local. Pour un sommet $s$ param\'etr\'e par une orbite \`a deux \'el\'ements,  on a  $FC(\mathfrak{g}_{s}({\mathbb F}_{q}))=\{0\}$ d'apr\`es le lemme \ref{orbites}. Parce que $n$ est impair, il y a un unique sommet $s$ param\'etr\'e par la racine centrale dans le diagramme. Pour celui-ci, on a $G_{s}\simeq (Spin(n)\times Spin(n)) /\{1,(z,z)\}$ sur $\bar{{\mathbb F}}_{q}$ avec une action galoisienne tordue par l'action alg\'ebrique de $\Gamma_{{\mathbb F}_{q}}$ qui est triviale sur  $\Gamma_{{\mathbb F}_{q^2}}$
et telle qu'un \'el\'ement de $\Gamma_{{\mathbb F}_{q}}-\Gamma_{{\mathbb F}_{q^2}}$ \'echange les deux copies de $Spin(n)$. Il y a une fonction $f_{N,\epsilon}$ avec $\epsilon(z)=1$ si et seulement si $\delta_{\square}(n)=1$ et une fonction $f_{N,\epsilon}$ avec $\epsilon(z)=-1$ si et seulement si $\delta_{\triangle}(n)=1$. Quand elles existent, ces fonctions donnent naissance \`a des \'el\'ements de $FC(\mathfrak{g}(F))$. On note $FC_{k^{st},k^{st}}$, resp. $FC_{k_{0},k_{0}}$, la droite port\'ee par la premi\`ere, resp. seconde, fonction. Comme dans le paragraphe pr\'ec\'edent, ces fonctions se transforment selon le caract\`ere trivial de $G_{AD}(F)_{0}/\pi(G(F))$ pour la premi\`ere, par le caract\`ere non trivial pour la seconde.  On obtient ainsi \ref{resultats}(1).

 Si $\delta_{\square}(n)=0$, on pose ${\cal Y}^{+} =\emptyset$. Si $\delta_{\square}(n)=1$,  on   pose ${\cal Y}^{+}=\{(k^{st},0)\}$. 
   Si $\delta_{\triangle}(n)=0$, on pose ${\cal Y}^{-}= \emptyset$. Si $\delta_{\triangle}(n)=1$, on pose    ${\cal Y}^-=\{(0,k_{0})\}$. On pose ${\cal Y}={\cal Y}^+\sqcup {\cal Y}^-$.  
   
La description des donn\'ees endoscopiques ${\bf G}'$ est la m\^eme que dans le paragraphe pr\'ec\'edent. Parce que le groupe $G$ est param\'etr\'e par le caract\`ere non trivial de $Z(\hat{G}_{SC})^{\Gamma_{F}}$, on voit que, quand le groupe $Out({\bf G}')$ contient    $\theta\theta'$, cet automorphisme agit maintenant par multiplication par $-1$ sur $FC^{st}(\mathfrak{g}'(F))$. Ces donn\'ees disparaissent donc. Il reste au plus deux donn\'ees dont le groupe d'automorphismes ext\'erieurs  ne contient pas $\theta\theta'$. Leurs groupes d'automorphismes ext\'erieurs  est d'ailleurs trivial. C'est la donn\'ee principale ${\bf G}$ et la donn\'ee que l'on avait param\'etr\'ee par  $(0,k_{0})$ (qui existe si et seulement si $\delta_{\triangle}(n)=1$). Pour la donn\'ee ${\bf G}$, d'apr\`es le r\'esultat du paragraphe pr\'ec\'edent,  l'espace $FC^{st}(\mathfrak{g}^*(F))$ est non nul si et seulement si $\delta_{\square}(n)=1$. Dans ce cas, cet espace est une droite. On    note $FC^{{\cal E}}_{k^{st},0}$ cette droite et on pose ${\bf G}'_{k^{st},0}={\bf G}$.   Pour l'autre donn\'ee, l'espace $FC^{st}(\mathfrak{g}(F))$ est non nul si et seulement si $\delta_{\triangle}(n)=1$. Dans ce cas, cet espace est une droite  que l'on note $FC^{{\cal E}}_{0,k_{0}}$ et on pose ${\bf G}'_{0,k_{0}}={\bf G}'$. Remarquons que, comme dans le paragraphe pr\'ec\'edent, $\xi_{{\bf G}'_{0,k_{0}}}$ est non trivial sur $G_{AD}(F)_{0}/\pi(G(F))$.  On a obtenu \ref{resultats}(2).

On note $\varphi:{\cal X}\to {\cal Y}$ la bijection \'evidente.  
L'assertion \ref{resultats}(3) est imm\'ediate puisque les droites $FC_{x}$ pour $x\in {\cal X}$ se distinguent par le caract\`ere de $G_{AD}(F)_{0}/\pi(G(F))$ par lequel  ce groupe agit sur elles et que les donn\'ees endoscopiques ${\bf G}'_{y}$ pour $y\in {\cal Y}$ se distinguent elles-aussi selon la restriction \`a ce groupe du caract\`ere $\xi_{{\bf G}'_{y}}$. 
  
\section{Descriptions explicites pour les groupes classiques}

\subsection{Type $B_{n}$ d\'eploy\'e}\label{Bnexplicite}
  
On suppose que $G$ est d\'eploy\'e de type $B_{n}$ avec $n\geq2$. On suppose aussi que $\delta_{2\triangle}(n)=1$. L'ensemble ${\cal X}^{st}$ de \ref{Bndeppadique} a un \'el\'ement que l'on note $(k,h)$. Rappelons que $k,h\in {\mathbb N}$, $2n+1=k^2+h^2$, $k$  est pair,  $h$ est impair et $\vert k-h\vert =1$. 

Fixons un espace $V$ sur $F$ de dimension $2n+1$, muni d'une base $(e_{i})_{i=1,...,2n+1}$ et de la forme quadratique $q$ d\'efinie par $q(v)=x_{n+1}^2+\sum_{i=1,...,2n+1, i\not=n+1}x_{i}x_{2n+2-i}$ pour tout $v=\sum_{i=1,...,2n+1}x_{i}e_{i}$. 
 Le groupe $G_{AD}$ s'identifie au  groupe sp\'ecial orthogonal de cet espace quadratique et $G$  \`a son rev\^etement simplement connexe.  

  Notons $R_{\geq0}\subset V$ le $\mathfrak{o}_{F}$-r\'eseau
engendr\'e par $e_{1},...,e_{2n+1-k^2/2},\varpi_{F}e_{2n+2-k^2/2},...,\varpi_{F}e_{2n+1}$. On pose $R_{\geq0}^*=\{v\in V; \forall v'\in R_{\geq0},\,\,q(v,v')\in \mathfrak{o}_{F}\}$. Alors $R_{\geq0}^*$ est le $\mathfrak{o}_{F}$-r\'eseau engendr\'e par $\varpi_{F}^{-1}e_{1},...,\varpi_{F}^{-1}e_{k^2/2},e_{1+k^{2}/2},...,e_{2n+1}$. Posons $I^h=\{k^2/2+1,...,2n+1-k^{2}/2\}$, $I^{k}_{+}=\{1,...,k^{2}/2\}$, $I_{-}^k=\{2n+2-k^{2}/2,...,2n+1\}$ et $I^k=I^k_{+}\cup I^k_{-}$. On pose $\underline{V}^h=R_{\geq0}/\mathfrak{p}_{F}R_{\geq0}^*$, $\underline{V}^{k}=R_{\geq0}^*/R_{\geq0}$. Pour $i\in I^h$, resp. $i\in I^k_{+}$, $i\in I^k_{-}$, notons $\underline{e}_{i}$ la r\'eduction dans $\underline{V}^h$, resp. $\underline{V}^k$, de $e_{i}$, resp. $\varpi_{F}^{-1}e_{i}$, $e_{i}$. L'espace $\underline{V}^h$ sur ${\mathbb F}_{q}$ est muni de la base $(\underline{e}_{i})_{i\in I^h}$ et de la r\'eduction $q^h$ de la forme $q$. L'espace $\underline{V}^k$ sur ${\mathbb F}_{q}$ est muni de la base $(\underline{e}_{i})_{i\in I^k}$ et de la r\'eduction $q^k$ de la forme $\varpi_{F}q$. Pour $\underline{v}=\sum_{i\in I^h}\underline{x}_{i}\underline{e}_{i}\in \underline{V}^h$, on a $q^h(\underline{v})=\underline{x}_{n+1}^2+\sum_{i\in I^h, i\not=n+1}\underline{x}_{i}\underline{x}_{2n+2-i}$. Pour $\underline{v}=\sum_{i\in I^k}\underline{x}_{i}\underline{e}_{i}\in \underline{V}^k$, on a $q^k(\underline{v})=\sum_{i\in I^k}\underline{x}_{i}\underline{x}_{2n+2-i}$. Pour $a=h$ ou $k$, 
 on note $\underline{G}^{a}$ le groupe spinoriel associ\'e \`a $(\underline{V}^{a},q^{a})$. On construit comme en \ref{Bn} ou \ref{Dndeppair} une  fonction $f_{N_{\square},\epsilon_{\square}}$ sur $\underline{\mathfrak{g}}^{a}({\mathbb F}_{q})$, que l'on note simplement $f^{a}$.

Notons $\mathfrak{k}$, resp. $\mathfrak{k}^{\perp}$, le sous-ensemble des $X\in \mathfrak{g}(F)$ tels que $X(R_{\geq0})\subset R_{\geq0}$, resp. $X(R_{\geq0}^*)\subset R_{\geq 0}$ et $X(R_{\geq 0})\subset \varpi_{F}R_{\geq 0}^*$.  On a l'\'egalit\'e
$$\mathfrak{k}/\mathfrak{k}^{\perp}=\underline{\mathfrak{g}}^k({\mathbb F}_{q})\oplus \underline{\mathfrak{g}}^h({\mathbb F}_{q}).$$
On rel\`eve la fonction $f^{k}\otimes f^h$ en une fonction sur $\mathfrak{k}$ et on l'\'etend par $0$ hors de $\mathfrak{k}$. On obtient une fonction $f_{k,h}\in FC(\mathfrak{g}(F))$ qui engendre la droite $FC_{k,h}$ d\'efinie en \ref{Bndeppadique}. 

  Posons $l=\frac{k+h-1}{2}$.  
  Soit $m\in \{1,..., l\}$. On fixe un \'el\'ement $\alpha_{m}\in \bar{F}^{\times}$ tel que $\alpha_{m}^{4m}=\varpi_{F}$. Comme en  \ref{lecasstable}, on fixe un ensemble de repr\'esentants ${\cal C}\subset \mathfrak{o}_{F}^{\times}$ du quotient $\mathfrak{o}_{F}^{\times}/\mathfrak{o}_{F}^{\times,2}$. On consid\`ere $F(\alpha_{m})$ comme un espace sur $F$ et, pour $\gamma\in {\cal C}$, on le munit de la forme quadratique $q_{m,\gamma}(v,v')=(4m)^{-1}\gamma trace_{F(\alpha_{m})/F}(\bar{v}v')$, o\`u $v\mapsto \bar{v} $ est la conjugaison galoisienne relative \`a l'extension $F(\alpha_{m})/F(\alpha_{m}^2)$. Soit $\underline{\gamma}=(\gamma_{m})_{m=1,...,l}\in {\cal C}^l$. On pose $\gamma_{0}=(-1)^l\prod_{m=1,...,l}\gamma_{m}$ si $k=h-1$, $\gamma_{0}=\varpi_{F}^{-1}(-1)^l\prod_{m=1,...,l}\gamma_{m}$ si $k=h+1$. On pose $V_{\underline{\gamma}}=F\oplus(\oplus_{m=1,...,l}F(\alpha_{m}))$. On munit $F$ de la forme quadratique $x\mapsto \gamma_{0}x^2$ et les $F(\alpha_{m})$ des formes $q_{m,\gamma_{m}}$. On munit $V_{\underline{\gamma}}$ de la somme orthogonale de ces formes quadratiques. On v\'erifie que cet espace quadratique est isomorphe \`a $(V,q)$ si et seulement si $\underline{\gamma}\in {\cal C}^m(V)$, o\`u ${\cal C}^m(V)$ est l'ensemble des $\underline{\gamma}=(\gamma_{m})_{m=1,...,l}\in {\cal C}^l$ tels que $\prod_{m=1,...,l}\gamma_{l}=(-1)^{(h-1)/2}$. On suppose cette condition v\'erifi\'ee et on fixe un isomorphisme d'espaces quadratiques $\iota_{\underline{\gamma}}:V_{\underline{\gamma}}\to V$. On note $X'_{\underline{\gamma}}$ l'endomorphisme de $V_{\underline{\gamma}}$ qui agit par multiplication par $\alpha_{m}$ sur chaque $F(\alpha_{m})$ et par $0$ sur $F$. On note $X_{\underline{\gamma}}$ l''endomorphisme de $V$ qui s'en d\'eduit par $\iota_{\underline{\gamma}}$. On v\'erifie que $X_{\underline{\gamma}}\in \mathfrak{g}_{ell}(F)$. Remarquons que

(1) il  existe une valeur propre $\alpha\in \bar{F}^{\times}$ de $X_{\underline{\gamma}}$ (consid\'er\'e comme un \'el\'ement de $End_{F}(V)$) telle que $val_{F}(\alpha)=\frac{1}{2(k+h-1)}$. 

Il est bien connu que les \'el\'ements $X_{\underline{\gamma}}$ sont tous stablement conjugu\'es, plus pr\'ecis\'ement que la famille $(X_{\underline{\gamma}})_{\underline{\gamma}\in {\cal C}^{l}(V)}$ est un ensemble de repr\'esentants des classes de conjugaison par $G_{AD}(F)$ dans leur classe de conjugaison stable. Ce n'est pas un ensemble de repr\'esentants des classes de conjugaison par $G(F)$ mais cela n'importe pas pour calculer les int\'egrales orbitales stables: on a
$$S^G(X_{\underline{\gamma}},f)=c\sum_{\underline{\gamma}\in {\cal C}^l(V)}I^{G_{AD}}(X_{\underline{\gamma}'},f)$$
pour tout $f\in C_{c}^{\infty}(\mathfrak{g}(F))$, avec une constante $c>0$ provenant de la comparaison entre les mesures sur $G(F)$ et sur $G_{AD}(F)$. Un calcul similaire \`a celui de \ref{lecasstable} prouve que

\begin{lem}{On a $S^G(X_{\underline{\gamma}},f_{k,h})\not=0$ pour tout $\underline{\gamma}\in {\cal C}^l(V)$.}\end{lem}

\subsection{Type $C_{n}$ d\'eploy\'e}\label{Cnexplicite}

On suppose que $G$ est d\'eploy\'e de type $C_{n}$ avec $n\geq2$. On suppose aussi que $\delta_{2\triangle}(n)=1$. L'ensemble ${\cal X}^{st}$ de \ref{Cndeppadique} a un \'el\'ement que l'on note $(k,k)$. Rappelons que $k\in {\mathbb N}$ et $n=k(k+1)$. On fixe un \'el\'ement non nul $f_{k,k}$ de la droite $FC_{k,k}$. On peut en donner une description analogue
 \`a celle de \ref{Bnexplicite}, nous la laissons au lecteur. 
 
 On fixe un espace $V$ sur $F$ de dimension $2n$, munie d'une forme symplectique $q$. Le groupe $G$ s'identifie au groupe symplectique de $(V,q)$.

Soit $m\in \{1,..., k\}$. On fixe un \'el\'ement $\alpha_{m}\in \bar{F}^{\times}$ tel que $\alpha_{m}^{4m}=\varpi_{F}$.  On consid\`ere $F(\alpha_{m})$ comme un espace sur $F$ et, pour $\gamma\in {\cal C}$, on le munit de la forme symplectique $q_{m,\gamma}(v,v')=(4m)^{-1}\gamma trace_{F(\alpha_{m})/F}(\alpha_{m}\bar{v}v')$, o\`u $v\mapsto \bar{v} $ est la conjugaison galoisienne relative \`a l'extension $F(\alpha_{m})/F(\alpha_{m}^2)$. Soit $\underline{\gamma}=(\gamma_{m})_{m=1,...,k}\in {\cal C}^k$.   On pose $V_{\underline{\gamma}}= \oplus_{m=1,...,k}F(\alpha_{m})$. On munit  les $F(\alpha_{m})$ des formes $q_{m,\gamma_{m}}$. On munit $V_{\underline{\gamma}}$ de la somme orthogonale de ces formes symplectiques. Cet espace symplectique est \'evidemment isomorphe \`a $(V,q)$. On fixe un isomorphisme d'espaces symplectiques $\iota_{\underline{\gamma}}:V_{\underline{\gamma}}\to V$. On note $X'_{\underline{\gamma}}$ l'endomorphisme de $V_{\underline{\gamma}}$ qui agit par multiplication par $\alpha_{m}$ sur chaque $F(\alpha_{m})$. On note $X_{\underline{\gamma}}$ l'endomorphisme de $V$ qui s'en d\'eduit par $\iota_{\underline{\gamma}}$. On v\'erifie que $X_{\underline{\gamma}}\in \mathfrak{g}_{ell}(F)$. Remarquons que

(1) il  existe une valeur propre $\alpha\in \bar{F}^{\times}$ de $X_{\underline{\gamma}}$ (consid\'er\'e comme un \'el\'ement de $End_{F}(V)$) telle que $val_{F}(\alpha)=\frac{1}{4k}$. 

Les \'el\'ements $X_{\underline{\gamma}}$ sont tous stablement conjugu\'es, plus pr\'ecis\'ement  la famille $(X_{\underline{\gamma}})_{\underline{\gamma}\in {\cal C}^{l}(V)}$ est un ensemble de repr\'esentants des classes de conjugaison par $G(F)$ dans leur classe de conjugaison stable.   Un calcul similaire \`a celui de \ref{lecasstable} prouve que

\begin{lem}{On a $S^G(X_{\underline{\gamma}},f_{k,k})\not=0$ pour tout $\underline{\gamma}\in {\cal C}^k$.}\end{lem}

\subsection{Type $D_{n}$}\label{Dnexplicite}

On suppose que $G$ est quasi-d\'eploy\'e de type $D_{n}$.  Dans les diff\'erents paragraphes concernant les groupes quasi-d\'eploy\'es de type $D_{n}$, on a d\'efini un ensemble ${\cal X}^{st}$ et on suppose qu'il est non vide. On a alors $\delta_{\square}(n)=1$, c'est-\`a-dire $n=k^2$ pour un entier $k\in {\mathbb N}$ et l'ensemble ${\cal X}^{st}$ est \'egal \`a $(k,k)$. De plus, on est forc\'ement dans l'une des situations suivantes:

(A) $n$ est pair et $G$ est d\'eploy\'e;

(B) $n$ est pair, $G$ n'est pas d\'eploy\'e et l'est  sur l'extension quadratique $E_{0}/F$ non ramifi\'ee;

(C) $n$ est impair, $G$ n'est pas d\'eploy\'e et l'est  sur une extension quadratique $E/F$ ramifi\'ee. 

Dans le cas (A), on pose $\lambda=1$. Dans les cas (B), resp. (C), on fixe un \'el\'ement $\lambda\in F^{\times}$ tel que $E_{0}=F(\sqrt{\lambda})$ et $val_{F}(\lambda)=0$, resp. $E=F(\sqrt{\lambda})$ et $val_{F}(\lambda)=-1$. On note $\underline{\lambda}$ la r\'eduction dans ${\mathbb F}_{q}$ de $\lambda$ dans les cas (A) ou (B), de $\varpi_{F}\lambda$ dans le cas (C).

On fixe un espace vectoriel $V$ sur $F$ de dimension $2n$, muni d'une base $(e_{i})_{i=1,...,2n}$. On d\'efinit une forme quadratique $q$ sur $V$ par l'\'egalit\'e suivante, pour $v=\sum_{i=1,...,2n}x_{i}v_{i}$:

dans les cas (A) ou (B), $q(v)=x_{n}^2-\lambda x_{n+1}^2+\sum_{i=1,...,2n,i\not=n,n+1}x_{i}x_{2n+1-i}$;

dans le cas (C), $q(v)=x_{n}^2-\lambda x_{2n}^2+\sum_{i=1,...,2n-1,i\not=n}x_{i}x_{2n-i}$. 

Posons $l^-=[n/2]$, $l^+=[(n+1)/2]$.   Notons $R_{\geq0}\subset V$ le $\mathfrak{o}_{F}$-r\'eseau
engendr\'e par $e_{1},...,e_{2n-l^+},\varpi_{F}e_{2n+1-l^+},...,\varpi_{F}e_{2n}$. On pose $R_{\geq0}^*=\{v\in V; \forall v'\in R_{\geq0},\,\,q(v,v')\in \mathfrak{o}_{F}\}$. Alors $R_{\geq0}^*$ est le $\mathfrak{o}_{F}$-r\'eseau engendr\'e par $\varpi_{F}^{-1}e_{1},...,\varpi_{F}^{-1}e_{l^-},e_{1+l^-},...,e_{2n}$. Posons $I'=\{l^-+1,...,2n-l^+\}$, $I''_{+}=\{1,...,l^-\}$, $I''_{-}=\{2n+1-l^+,...,2n\}$ et $I''=I''_{+}\cup I''_{-}$. On pose $\underline{V}'=R_{\geq0}/\mathfrak{p}_{F}R_{\geq0}^*$, $\underline{V}''=R_{\geq0}^*/R_{\geq0}$. Pour $i\in I'$, resp. $i\in I''_{+}$, $i\in I''_{-}$, notons $\underline{e}_{i}$ la r\'eduction dans $\underline{V}'$, resp. $\underline{V}''$, de $e_{i}$, resp. $\varpi_{F}^{-1}e_{i}$, $e_{i}$. L'espace $\underline{V}'$ sur ${\mathbb F}_{q}$ est muni de la base $(\underline{e}_{i})_{i\in I'}$ et de la r\'eduction $q'$ de la forme $q$. L'espace $\underline{V}''$ sur ${\mathbb F}_{q}$ est muni de la base $(\underline{e}_{i})_{i\in I''}$ et de la r\'eduction $q''$ de la forme $\varpi_{F}q$. Pour $\underline{v}=\sum_{i\in I'}\underline{x}_{i}\underline{e}_{i}\in \underline{V}'$, on a 

dans les cas (A) ou (B), $q'(\underline{v})=\underline{x}_{n}^2-\underline{\lambda}x_{n+1}^2+\sum_{i\in I', i\not=n,n+1}\underline{x}_{i}\underline{x}_{2n+1-i}$;

dans le cas (C), $q'(\underline{v})=\underline{x}_{n}^2 +\sum_{i\in I', i\not=n}\underline{x}_{i}\underline{x}_{2n-i}$.

\noindent Pour $\underline{v}=\sum_{i\in I''}\underline{x}_{i}\underline{e}_{i}\in \underline{V}''$, on a 

dans les cas (A) ou (B), $q''(\underline{v})= \sum_{i\in I''}\underline{x}_{i}\underline{x}_{2n+1-i}$;

dans le cas (C), $q''(\underline{v})=\underline{\lambda}\underline{x}_{2n}^2 +\sum_{i\in I'', i\not=2n}\underline{x}_{i}\underline{x}_{2n-i}$.

 \noindent Pour $a$ l'un des symboles $'$ ou $''$, on note $\underline{G}^{a}$ le groupe spinoriel associ\'e \`a $(\underline{V}^{a},q^{a})$. On construit selon les cas comme en \ref{Bn} ou \ref{Dndeppair} ou \ref{Dnnondeppair} une  fonction $f_{N_{\square},\epsilon_{\square}}$ sur $\underline{\mathfrak{g}}^{a}({\mathbb F}_{q})$, que l'on note simplement $f^{a}$.

Notons $\mathfrak{k}$, resp. $\mathfrak{k}^{\perp}$, le sous-ensemble des $X\in \mathfrak{g}(F)$ tels que $X(R_{\geq0})\subset R_{\geq0}$, resp. $X(R_{\geq0}^*)\subset R_{\geq 0}$ et $X(R_{\geq 0})\subset \varpi_{F}R_{\geq 0}^*$.  On a l'\'egalit\'e
$$\mathfrak{k}/\mathfrak{k}^{\perp}=\underline{\mathfrak{g}}''({\mathbb F}_{q})\oplus \underline{\mathfrak{g}}'({\mathbb F}_{q}).$$
On rel\`eve la fonction $f''\otimes f'$ en une fonction sur $\mathfrak{k}$ et on l'\'etend par $0$ hors de $\mathfrak{k}$. On obtient une fonction $f_{k,k}\in FC(\mathfrak{g}(F))$ qui engendre la droite $FC_{k,k}$ d\'efinie selon les cas en \ref{Dndeppairpadique}, \ref{Dnpairpadiquenonram}, \ref{Dnimppadiqueram}. 

On fixe une famille $\underline{\beta}=(\beta_{m})_{m=1,...,k}$ d'\'el\'ements de $\mathfrak{o}_{F}^{\times}$ telle que
 la r\'eduction dans ${\mathbb F}_{q}^{\times}$ de $\prod_{m=1,...,k}\beta_{m}$ appartient \`a $(-1)^k\underline{\lambda}{\mathbb F}_{q}^{\times,2}$.

Soit $m\in \{1,...,k\}$. On fixe un \'el\'ement   un \'el\'ement  $\alpha_{m}\in \bar{F}^{\times}$ tel que $\alpha_{m}^{4m-2}=\varpi_{F}\beta_{m}$. On consid\`ere $F(\alpha_{m})$ comme un espace sur $F$ et, pour $\gamma\in {\cal C}$, on le munit de la forme quadratique $q_{m,\gamma}(v,v')=(4m)^{-1}\gamma trace_{F(\alpha_{m})/F}(\bar{v}v')$, o\`u $v\mapsto \bar{v} $ est la conjugaison galoisienne relative \`a l'extension $F(\alpha_{m})/F(\alpha_{m}^2)$.

Soit $\underline{\gamma}=(\gamma_{m})_{m=1,...,k}\in {\cal C}^k$. On pose $V_{\underline{\gamma}}=\oplus_{m=1,...,l}F(\alpha_{m})$. On munit  les $F(\alpha_{m})$ des formes $q_{m,\gamma_{m}}$. On munit $V_{\underline{\gamma}}$ de la somme orthogonale de ces formes quadratiques. On v\'erifie que cet espace quadratique est isomorphe \`a $(V,q)$ si et seulement si $\underline{\gamma}\in {\cal C}^k(V)$, o\`u ${\cal C}^k(V)$ est l'ensemble des $\underline{\gamma}=(\gamma_{m})_{m=1,...,k}\in {\cal C}^k$ v\'erifiant l'\'egalit\'e suivante dans ${\cal C}\simeq \mathfrak{o}_{F}^{\times}/\mathfrak{o}_{F}^{\times,2}$:

dans les cas (A) ou (B), $\prod_{m=1,...,k}\gamma_{m}=\lambda(-1)^{k/2}$;

dans le cas (C), $\prod_{m=1,...,k}\gamma_{m}=(-1)^{(k-1)/2}$.

 \noindent On suppose $\underline{\gamma}\in {\cal C}^k(V)$  et on fixe un isomorphisme d'espaces quadratiques $\iota_{\underline{\gamma}}:V_{\underline{\gamma}}\to V$. On note $X'_{\underline{\gamma}}$ l'endomorphisme de $V_{\underline{\gamma}}$ qui agit par multiplication par $\alpha_{m}$ sur chaque $F(\alpha_{m})$ et par $0$ sur $F$. On note $X_{\underline{\gamma}}$ l'endomorphisme de $V$ qui s'en d\'eduit par $\iota_{\underline{\gamma}}$. On v\'erifie que $X_{\underline{\gamma}}\in \mathfrak{g}_{ell}(F)$. Remarquons que

(1) il  existe une valeur propre $\alpha\in \bar{F}^{\times}$ de $X_{\underline{\gamma}}$ (consid\'er\'e comme un \'el\'ement de $End_{F}(V)$) telle que $val_{F}(\alpha)=\frac{1}{4k-2}$. 

Les \'el\'ements $X_{\underline{\gamma}}$ ne sont pas forc\'ement  stablement conjugu\'es.  Ils sont inclus dans une seule  classe de conjugaison par le groupe orthogonal tout entier $O(2n,\bar{F})$. Mais cette classe se coupe en deux classes de conjugaison par $SO(2n,\bar{F})$. L'action d'un \'el\'ement de $O(2n,F)$ de d\'eterminant $-1$ \'echange ces deux classes. Quitte \`a remplacer chaque isomorphisme $\iota_{\underline{\gamma}}$ par l'action d'un tel \'el\'ement, on peut donc supposer que les \'el\'ements $X_{\underline{\gamma}}$  sont tous stablement conjugu\'es.
Alors la famille $(X_{\underline{\gamma}})_{\underline{\gamma}\in {\cal C}^{k}(V)}$ est un ensemble de repr\'esentants des classes de conjugaison par $G_{AD}(F)$ dans leur classe de conjugaison stable. Ce n'est pas un ensemble de repr\'esentants des classes de conjugaison par $G(F)$ mais, comme en \ref{Bnexplicite}, cela n'importe pas pour calculer les int\'egrales orbitales stables. Un calcul similaire \`a celui de \ref{lecasstable} prouve alors que

\begin{lem}{On a $S^G(X_{\underline{\gamma}},f_{k,k})\not=0$ pour tout $\underline{\gamma}\in {\cal C}^k(V)$.}\end{lem}

\subsection{Les \'el\'ements $Y_{y}$} \label{elementsYy}
Ici, le groupe $G$ est l'un des groupes trait\'es dans les paragraphes \ref{Bndeppadique} \`a \ref{Dnimppadiqueram}.  Dans plusieurs de ces paragraphes, on a affirm\'e l'existence d'\'el\'ements $Y_{y}$ pour $y\in {\cal Y}$.  

{\bf Remarque.} Il y a une variante en  \ref{Dnpairpadiquenonram}, o\`u les \'el\'ements ci-dessus sont remplac\'es par des  $Y_{\tilde{y}}$ pour $\tilde{y}\in \tilde{{\cal Y}}$. Cette variante ne cr\'ee pour ce paragraphe qu'une diff\'erence de notations . Pour simplifier la r\'edaction, nous r\'edigeons la preuve en n\'egligeant cette diff\'erence.  

\bigskip

Selon les cas, les ensembles ${\cal X}$ et ${\cal Y}$ \'etaient divis\'es en deux sous-ensembles se correspondant par la bijection $\varphi$. Dans ce cas on fixe de tels sous-ensembles ${\cal X}^{\star}$ et ${\cal Y}^{\star}$. Dans le cas o\`u on n'a pas divis\'e ${\cal X}$ et ${\cal Y}$ en de tels sous-ensembles, on pose ${\cal X}^{\star}={\cal X}$ et ${\cal Y}^{\star}={\cal Y}$ pour unifier la notation. 

A $y\in {\cal Y}^{\star}$ est associ\'ee une donn\'ee ${\bf G}'_{y}$. L'espace $FC^{st}(\mathfrak{g}_{y}'(F))$ est une droite. Les propri\'et\'es  exig\'ees de $Y_{y}$ sont les suivantes:

(1) c'est un \'el\'ement $G$-r\'egulier de $\mathfrak{g}'_{y,ell}(F) $;

(2) pour un \'el\'ement non nul $f'\in FC^{st}(\mathfrak{g}'_{y}(F))$, $S^{G'_{y}}(Y_{y},f')\not=0$;

(3) soit $x\in {\cal X}^{\star}$, $f\in FC_{x}$ et $X\in \mathfrak{g}(F)$ un \'el\'ement dont la classe de conjugaison stable correspond \`a celle de $Y_{y}$; supposons $I^G(X,f)\not=0$; alors $x\geq \varphi^{-1}(y)$. 

En examinant nos constructions, on voit que $G'_{y,SC}$ est toujours produit d'au plus $2$ groupes quasi-d\'eploy\'es de la forme suivante:

(a) $SU_{E/F}(m)$ pour une extension quadratique $E/F$ ramifi\'ee;

(b) $Spin(2m+1)$, $Sp(2m)$, $Spin_{dep}(2m)$ avec $m$ pair, $Spin_{E_{0}/F}$ avec $m$ pair, o\`u $E_{0}/F$ est l'extension  quadratique non ramifi\'ee, $Spin_{E/F}$ avec $m$ impair,  pour une extension quadratique $E/F$ ramifi\'ee;

(c) $Res_{K/F}(H)$ o\`u $K/F$ est une extension quadratique de $F$ et $H$ est du type (b) ci-dessus au changement pr\`es du corps de base $F$ qui est remplac\'e par $K$. 

On raffine de fa\c{c}on \'evidente les cas (b) et (c) en cas $(b,B_{m})$, $(b,C_{m})$, $(b,D_{m})$, $(c,B_{m})$, $(c,C_{m})$, $(c,D_{m})$.

Notons $G''$ un tel groupe. L'ensemble ${\cal X}^{''st}$ associ\'e \`a ce groupe est r\'eduit \`a un \'el\'ement que nous notons $(k'',h'')$. Fixons un \'el\'ement non nul de l'espace correspondant $FC_{k'',h''}$. Posons
$$d''=\left\lbrace\begin{array}{cc}2k''+2h''-1,& \rm{dans\, le\, cas\, (a)},\\ 2k''+2h''-2,& \rm{dans\, les\, cas\, }(b,B_{m}),(c,B_{m}), (b,D_{m}) \rm{\, et\, }(c,D_{m}),\\ 2k''+2h'',& \rm{dans\, les\, cas\, }(b,C_{m})\rm{ \,et\, }(c,C_{m}).\\ \end{array}\right.$$ 

 Dans les paragraphes \ref{lecasstable}, \ref{Bnexplicite}, \ref{Cnexplicite}, \ref{Dnexplicite}, on a d\'efini des \'el\'ement not\'es alors $X_{\underline{\gamma}}$ (dans le cas (C), il faut changer de corps de base).  Fixons  un tel \'el\'ement  que l'on  note $X''$. On a prouv\'e

(4) $S^{G''}(X'',f_{k'',h''})\not=0$.

On peut consid\'erer $X''$ comme un endomorphisme d'un espace vectoriel sur une extension $L$ de $F$: $L=E$ dans le cas (a), $L=F$ dans le cas (b), $L=K$ dans le cas (c). On a prouv\'e qu'il existait une valeur propre $\alpha''\in \bar{F}^{\times}$ de $X''$ telle que

(5) $val_{L}(\alpha)=\frac{1}{d''}$. 

Quand $G'_{y,SC}$ est r\'eduit \`a un seul groupe $G''$, on pose $Y_{y}=X''$ et $f'_{y}=f_{k'',h''}$. Quand $G'_{y,SC}=G''_{1}\times G''_{2}$, on affecte les notations ci-dessus d'indices $1$ et $2$ et on pose  $f'_{y}=f_{k''_{1},h''_{1}}\otimes f_{k''_{2},h''_{2}}$.  L'\'el\'ement $X''_{1}\oplus X''_{2}$ n'est pas forc\'ement $G$-r\'egulier. Mais, comme on l'a expliqu\'e en \ref{ingredients} et utilis\'e en \ref{preuve}, on peut remplacer $X''_{1}$ et $X''_{2}$ par des \'el\'ements assez voisins de sorte que les propri\'et\'es (4) et (5) restent v\'erifi\'ees et que $X''_{1}\oplus X''_{2}$ soit $G$-r\'egulier. On pose alors $Y_{y}=X''_{1}\oplus X''_{2}$.  On a

(6) $S^{G'_{y}}(Y_{y},f'_{y})\not=0$. 

Les rangs des groupes $G''$ intervenant sont inf\'erieurs  ou \'egaux \`a celui de $G$. Si $G$ n'est pas quasi-d\'eploy\'e  ou si ces rangs sont strictement inf\'erieurs   \`a celui de $G$, notre hypoth\`ese de r\'ecurrence \ref{resultats}(4) appliqu\'ee \`a ces groupes nous dit que $f'_{y}$ appartient \`a l'espace $FC^{st}(\mathfrak{g}'_{y}(F))$, qui est une droite. Alors (6) est \'equivalent \`a (2). On doit aussi traiter le cas o\`u $G$ est quasi-d\'eploy\'e et le  rang d'un des groupes $G''$ est \'egal \`a celui de $G$. Cette derni\`ere condition \'equivaut \`a ${\bf G}'_{y}={\bf G}$. Comme en \ref{preuve}, on remarque que, lorsqu'on utilise les \'el\'ements $Y_{y}$, on n'a pas encore d\'etermin\'e l'espace $FC^{st}(\mathfrak{g}(F))$ mais on a d\'ej\`a prouv\'e par l'habituel argument de dimensions qu'il \'etait de dimension au plus $1$. En fait, l'hypoth\`ese implicite qu'il existe un \'el\'ement $y\in {\cal Y}$ tel que ${\bf G}'_{y}={\bf G}$ implique que cette dimension est $1$. Alors, le m\^eme argument qu'en \ref{preuve} permet de d\'eduire (2) de (6). 
 
Remarquons que, dans chacun de ces paragraphes, $G$ est un groupe  vraiment classique et a donc une repr\'esentation naturelle dans un espace vectoriel sur $F$. On peut parler des valeurs propres d'un \'el\'ement $X\in \mathfrak{g}(F)$. 
L'\'el\'ement $\varphi^{-1}(y)$ se projette dans un ensemble ${\mathbb X}^{\star}$ qui est toujours un ensemble de couples $(k,h)\in {\mathbb N}^2$. Notons $(k_{y},h_{y})$ le couple associ\'e \`a $\varphi^{-1}(y)$. En  examinant tous les  cas un \`a un, on voit que la condition (5) entra\^{\i}ne

(7) soit $X\in \mathfrak{g}_{reg}(F)$ un \'el\'ement dont la classe de conjugaison stable correspond \`a celle de $Y_{y}$; alors il existe une valeur propre $\alpha\in \bar{F}^{\times}$ de $X$ telle que 
$val_{F}(\alpha)=\frac{1}{2k_{y}+2h_{y}-e(G)}$,

\noindent o\`u $e(G)\in {\mathbb N}$ est d\'efini par
$$e(G)= \left\lbrace\begin{array}{cc}2,&\rm{\,si\,}G\rm{\,est\,de\,type\,}B_{n}\rm{\, ou\, }D_{n},\\ 0,&\rm{\,si\,}G\rm{\,est\,de\,type\,}C_{n}.\\  \end{array}\right.$$

Soit $x\in {\cal X}^{st}$, notons $(k,h)$ l'\'el\'ement de ${\mathbb X}^{\star}$ qui lui est associ\'e. On peut construire une base de $FC_{x}$ par l'alg\`ebre lin\'eaire comme en \ref{descriptionexplicite}. Nous passons cette construction fastidieuse mais nous utilisons la cons\'equence, analogue au lemme  de ce paragraphe  \ref{descriptionexplicite}:

(8) soit $f\in FC_{x}$ et $X$ un \'el\'ement du support de $f$; soit $\alpha\in \bar{F}$ une valeur propre  de $X$; alors on a 
$val_{F}(\alpha)\geq \frac{1}{2k+2h-e(G)}$.

 Soient alors $x$, $f$ et $X$ comme en (3). L'hypoth\`ese $I^G(X,f)\not=0$ et la relation (7) impliquent qu'il existe un \'el\'ement $X'$ conjugu\'e \`a $X$ dans le support de $f$. La r\'eunion des relations (7) et (8) entra\^{\i}ne que $k+h\geq k_{y}+h_{y}$. Or cette relation \'equivaut \`a $x\geq \varphi^{-1}(y)$. Cela d\'emontre (3).
 
\subsection{Action d'automorphismes} \label{variance}
Notons $G$ l'un des groupes  not\'es $G''$ dans le paragraphe pr\'ec\'edent. Ils sont parfois munis d'automorphismes alg\'ebriques qui pr\'eservent  la droite $FC^{st}(\mathfrak{g}(F))$ donc y agissent par multiplication par un scalaire. On a utilis\'e dans certains cas la valeur de ce scalaire, il nous faut justifier ce calcul. Remarquons qu'\`a chaque fois que l'on a utilis\'e la valeur de ce scalaire, on pouvait appliquer au groupe $G$ nos hypoth\`eses de r\'ecurrence, donc on connaissait l'espace $FC^{st}(\mathfrak{g}(F))$. 
On a d\'ej\`a fait le calcul en \ref{actiondunautomorphisme} pour les groupes de type (a) de \ref{elementsYy}. En examinant   les paragraphes \ref{Bndeppadique} \`a \ref{Dnimppadiqueram}, on voit que les groupes $G$ qui nous importent sont des types suivants, o\`u $E/F$ est une extension quadratique et les groupes dans les parenth\`eses sont d\'efinis sur $E$:

(a) $Res_{E/F}(Spin(2n+1))$ avec $\delta_{2\triangle}(n)=1$;

(b) $Res_{E/F}(Spin_{dep}(2n))$ avec $n$ pair et $\delta_{\square}(n)=1$;

(c) $Res_{E/F}(Spin_{Q/E}(2n))$ avec $n$ pair et   $\delta_{\square}(n)=1$, o\`u on suppose $E/F$ ramifi\'ee et l'on rappelle que $Q/F$ est l'extension biquadratique de $F$;

(d) $Res_{E/F}(Spin_{K/E}(2n)$ avec $n$ impair et $\delta_{\square}(n)=1$, o\`u on suppose $\delta_{4}(q-1)=1$, $K/F$ est une extension galoisienne cyclique de degr\'e $4$ totalement ramifi\'ee et  contenant $E$.

Un seul des automorphismes de ces groupes nous int\'eresse, que l'on note $\Theta$. 

Traitons d'abord les cas (a), (b) et (c). Dans ces cas, le groupe entre parenth\`ese est obtenu par restriction des scalaires de $F$ \`a $E$ d'un groupe d\'efini sur $F$: le groupe est \'evident dans les cas (a) et (b); c'est $Spin_{E_{0}/F}(2n)$ dans le cas (c). Alors le groupe $G$ est naturellement muni d'une action alg\'ebrique du groupe $\Gamma_{E/F}$ et l'automorphisme $\Theta$ est l'action de l'\'el\'ement $\tau$ non trivial de $\Gamma_{E/F}$. Concr\`etement, r\'ealisons le groupe $G(F)=Spin(2n+1,E)$, resp. $G(F)=Spin_{dep}(2n,E)$, $G(F)=Spin_{Q/E}(2n,E)$, comme le groupe d'automorphismes d'un espace $(V,q)$ comme en \ref{Bnexplicite}, resp. \ref{Dnexplicite} cas (A) ou (B), cet espace \'etant maintenant un $E$-espace vectoriel. On munit $V$ d'une base comme dans ces paragraphes, en supposant dans le cas (c) que le terme $\lambda$ intervenant en \ref{Dnexplicite} cas (B) appartient \`a $\mathfrak{o}_{F}^{\times}$. Alors les \'el\'ements de $G(F)$ s'identifient \`a des matrices  \`a coefficients dans $E$. On v\'erifie que, pour $g\in G(F)$, $\Theta(g)$ est la matrice obtenue en appliquant $\tau$ aux coefficients de $g$.

Consid\'erons le cas (a) et reprenons les d\'efinitions de \ref{Bnexplicite}. On choisit pour uniformisante $\varpi_{E}$ l'\'element $\varpi_{F}$ si $E/F$ est non ramifi\'ee, la racine carr\'ee d'un \'el\'ement de $\varpi_{F}\mathfrak{o}_{F}^{\times}$ si $E/F$ est ramifi\'ee. 
Pour $i,j\in I^k$, posons $a_{i,j}=0$ si $i,j\in I^k_{+}$ ou $i,j\in I^k_{-}$, $a_{i,j}=1$ si $i\in I^k_{+}$ et $j\in I^k_{-}$ et $a_{i,j}=-1$ si $i\in I^k_{-}$ et $j\in I^k_{+}$. Pour $X=(x_{i,j})_{i,j=1,...,2n+1}\in \mathfrak{g}(F)$, notons $X^h$ la matrice extraite $(x_{i,j})_{i,j\in I^h}$ et $X^k$ la matrice $(\varpi_{E}^{a_{i,j}}x_{i,j})_{i,j\in I^k}$. Si $X\in \mathfrak{k}$, on v\'erifie que ces matrices sont \`a coefficients dans $\mathfrak{o}_{E}$ et que les images dans $\mathfrak{g}^{h}({\mathbb F}_{q})=\mathfrak{spin}(h^2,{\mathbb F}_{q_{E}})$, resp.$\mathfrak{g}^k({\mathbb F}_{q})=\mathfrak{spin}(k^2,{\mathbb F}_{q_{E}})$, sont les r\'eductions naturelles de $X^h$, resp. $X^k$. Cette description montre que $\Theta$ conserve $\mathfrak{k}$. De plus, \`a cause des termes $\varpi_{E}^{a_{i,j}}$ et parce que $\tau(\varpi_{E})=\varpi_{E}$ si $E/F$ n'est pas ramifi\'ee, $\tau(\varpi_{E})=-\varpi_{E}$ si $E/F$ est ramifi\'ee, on voit que $\Theta$ et se r\'eduit en les actions suivantes: 

si $E/F$ est non ramifi\'ee, l'action "matricielle" du Frobenius sur $\mathfrak{spin}(h^2,{\mathbb F}_{q^2})$ et $\mathfrak{spin}(k^2,{\mathbb F}_{q^2})$;

si $E/F$ est ramifi\'ee, l'action triviale sur $\mathfrak{spin}(h^2,{\mathbb F}_{q})$ et la conjugaison dans  $\mathfrak{spin}(k^2,{\mathbb F}_{q})$ par similitude de rapport $-1$ (celle qui agit par $-1$ sur les $\underline{e}_{i}$ pour $i\in I^k_{+}$ et par $1$ sur les $\underline{e}_{i}$ pour $i\in I_{-}^k$).

Dans le cas o\`u $E/F$ est non ramifi\'ee, on v\'erifie que les actions de Frobenius fixent les orbites nilpotentes  (pour le groupe sp\'ecial orthogonal) dans le support des  fonctions $f_{N_{\square},\epsilon_{\square}}$, donc fixent ces fonctions. On en d\'eduit que $\Theta(f_{k,h})=f_{k,h}$. Mais on sait que $f_{h,k}$ engendre la droite $FC^{st}(\mathfrak{g}(F))$. Donc

(1) si $E/F$ est non ramifi\'ee, $\Theta$ agit trivialement sur $FC^{st}(\mathfrak{g}(F))$.

Dans le cas o\`u $E/F$ est ramifi\'ee, consid\'erons un \'el\'ement  nilpotent  $N$ dans le support de la fonction $f_{N_{\square},\epsilon_{\square}}$ de $\mathfrak{spin}(k^2,{\mathbb F}_{q})$. Son orbite par le groupe sp\'ecial orthogonal est param\'etr\'ee par la partition $(2k-1,2k-3,...,1)$ et par une collection $ (\gamma_{j})_{j=2k-1,2k-3,...,1}\in ({\mathbb F}_{q}^{\times}/{\mathbb F}_{q}^{\times,2})^{k}$. Son image  $N'$ par la similitude de rapport $-1$ est param\'etr\'ee par la m\^eme partition et par la collection $(-\gamma_{j})_{j=2k-1,2k-3,...,1}$. En se r\'ef\'erant \`a \ref{Dndeppair}, on voit que $f_{k,h}(N')=sgn(-1)^{k/2}f_{k,h}(N)$. Donc $\Theta(f_{k,h})=sgn(-1)^{k/2}f_{k,h}$. D'o\`u

(2) si $E/F$ est ramifi\'ee, $\Theta$ agit sur $FC^{st}(\mathfrak{g}(F))$ par multiplication par $sgn(-1)^{k/2}$. 

Remarquons qu'en \'ecrivant $n=i(i+1)$ avec $i\in {\mathbb N}$, on a $i=inf(h,k)$ et $k/2=[(i+1)/2]$. C'est sous cette forme que nous avons utilis\'e (2) en \ref{Bndeppadique}.

Un calcul analogue vaut pour les cas (b) et (c). On voit que l'on obtient encore les relations (1) et (2).

 Le cas (d) est plus compliqu\'e. Notons $G_{0}$ le groupe $Spin_{dep}(2n)$ d\'efini sur $F$, muni d'un \'epinglage d\'efini sur $F$. Notons $\sigma\mapsto \sigma_{G_{0}}$ l'action galoisienne naturelle qui fixe l'\'epinglage et $\theta$ l'automorphisme ext\'erieur non trivial qui le fixe aussi. Fixons un g\'en\'erateur $\rho$ de $\Gamma_{K/F}$. Le groupe $G$ se d\'ecrit de la fa\c{c}on suivante. On a $G(\bar{F})=G_{0}(\bar{F})\times G_{0}(\bar{F})$ et il est muni de l'action galoisienne $\sigma\mapsto \sigma_{G}$ telle que 

pour $\sigma\in \Gamma_{K}$, $\sigma_{G}(g_{1},g_{2})=(\sigma_{G_{0}}(g_{1}),\sigma_{G_{0}}(g_{2}))$;

$\rho_{G}(g_{1},g_{2})=(\rho_{G_{0}}(g_{2}),\rho_{G_{0}}\theta(g_{1}))$. 

L'application $(g_{1},g_{2})\mapsto g_{1}$ identifie $G(F)$ \`a $Spin_{K/E}(2n,E)$. 
L'automorphisme $\Theta$  est d\'efini par $\Theta(g_{1},g_{2})= (g_{2},\theta(g_{1}))$. Il est d'ordre $4$.  Concr\`etement, r\'ealisons $G(F)=Spin_{K/E}(2n,E)$ comme le groupe d'automorphismes d'un espace $(V,q)$ comme en  \ref{Dnexplicite} cas (c),  cet espace \'etant maintenant un $E$-espace vectoriel. On munit $V$ d'une base comme dans ce paragraphe en supposant que $\lambda$ v\'erifie $\tau(\lambda)=-\lambda$.  Alors les \'el\'ements de $G(F)$ s'identifient \`a des matrices  \`a coefficients dans $E$. Pour une telle matrice $g$, notons $\tau_{mat}(g)$ la matrice obtenue en appliquant $\tau$ aux coefficients de $g$. Fixons une racine primitive $\zeta$ d'ordre $4$ de l'unit\'e dans $F^{\times}$ ( qui existe gr\^ace \`a l'hypoth\`ese $\delta_{4}(q-1)=1$). Notons $t$ la matrice qui multiplie l'\'el\'ement de base $e_{2n}$ par $\zeta$ et fixe les autres \'el\'ements de base. On v\'erifie que, pour $g\in G(F)$, $\Theta(g)=t^{-1}\tau_{mat}(g)t$ (plus exactement, on peut choisir les diverses identifications de sorte qu'il en soit ainsi).

Reprenons les notations de \ref{Dnexplicite} cas (C). On voit comme dans le cas (a) ci-dessus que $\Theta$ conserve $\mathfrak{k}$ et se r\'eduit en l'identit\'e de $\underline{\mathfrak{g}}'({\mathbb F}_{q})$ (parce que, $K/F$ \'etant totalement ramifi\'ee, l'action galoisienne de $\rho$ se r\'eduit en l'identit\'e de ${\mathbb F}_{q}$) et en la conjugaison dans $\underline{\mathfrak{g}}''(F)$ par une similitude de rapport $-1$. Il s'agit de la similitude qui multiplie $\underline{e}_{i}$ par $-1$ pour $i\in I''_{+}$, par  la r\'eduction $\boldsymbol{\zeta}$  de $\zeta$ dans ${\mathbb F}_{q}$ pour $i=2n$ et par $1$ pour $i\in I''_{-}-\{2n\}$. Or cette similitude  est le produit d'un \'el\'ement du groupe sp\'ecial orthogonal de $\underline{V}''$ et de la multiplication par $\boldsymbol{\zeta}$, dont l'action par conjugaison est triviale. Puisque la fonction $f''$ est invariante par conjugaison par le groupe sp\'ecial orthogonal, elle est aussi invariante par notre similitude. D'o\`u $\Theta(f_{k,k})=f_{k,k}$ et

(3) $\Theta$ agit trivialement sur $FC^{st}(\mathfrak{g}(F))$.

\section{Type $D_{4}$ trialitaire et types exceptionnels}

\subsection{Un lemme immobilier}\label{lemmeimmobilier}
Pour tout $x\in Imm(G_{AD})$ et tout r\'eel $r$, Moy et Prasad ont d\'efini le $\mathfrak{o}_{F}$-r\'eseau $\mathfrak{k}_{x,r}\subset   \mathfrak{g}(F)$. On pose $\mathfrak{k}_{x,r+}=\cup_{s>r}\mathfrak{k}_{x,s}$. En particulier $\mathfrak{k}_{x,0}=\mathfrak{k}_{{\cal F}}$ et $\mathfrak{k}_{x,0+}=\mathfrak{k}_{{\cal F}}^+$, o\`u ${\cal F}$ est la facette \`a laquelle appartient $x$. On pose
$\mathfrak{g}(F)_{r}=\cup_{x\in Imm(G_{AD})}\mathfrak{k}_{x,r}$. Pour $X\in \mathfrak{g}(F)$, on appelle profondeur de $X$ et on note $r(X)$ la borne sup\'erieure de l'ensemble des $r\in {\mathbb R}$ tels que $X\in \mathfrak{g}(F)_{r}$. Si $X$ est nilpotent, cette profondeur est $+\infty$. Sinon, cette profondeur est finie et est atteinte, c'est-\`a-dire que $X\in \mathfrak{g}(F)_{r(X)}$. Indiquons une d\'efinition \'equivalente de $r(X)$. D\'ecomposons $X$ en somme $X_{s}+X_{n}$ de ses parties semi-simple et nilpotente et fixons un sous-tore maximal $T$ de $G$ tel que $X_{s}\in \mathfrak{t}(F)$. Introduisons l'ensemble $\Sigma^T$ des racines de $T$ dans $G$ (sur la cl\^oture alg\'ebrique $\bar{F}$). Alors 

(1) $r(X)=inf_{\alpha\in \Sigma^T}val_{F}(\alpha(X_{s}))$. 

On utilisera la propri\'et\'e suivante. Soit $K$ une extension galoisienne finie de $F$ mod\'er\'ement ramifi\'ee. Notons $e$ l'indice de ramification de $K/F$. On a d\'ej\`a dit que $Imm_{F}(G_{AD})$ s'identifiait au sous-ensemble des points fixes de l'action de $\Gamma_{K/F} $ dans $Imm_{K}(G_{AD})$. A $x\in Imm_{F}(G_{AD})$ sont associ\'es des r\'eseaux $\mathfrak{k}_{x,r,F}\subset \mathfrak{g}(F)$ et $\mathfrak{k}_{x,r,K}\subset \mathfrak{g}(K)$. On a alors

(2) $\mathfrak{k}_{x,r,F}=\mathfrak{g}(F)\cap \mathfrak{k}_{x,er,K}$.

On aura besoin du lemme suivant.

\begin{lem} {Supposons $G$ d\'eploy\'e, soit $T$ un sous-tore maximal d\'eploy\'e de $G$ et soit $X\in \mathfrak{t}(F)$. Supposons $val_{F}(\beta(X))=0$ pour toute racine $\beta$ de $T$ dans $G$. Alors le sous-ensemble des $x\in Imm(G_{AD})$ tels que $X\in \mathfrak{k}_{x,0}$ est \'egal \`a l'appartement associ\'e \`a $T$.}\end{lem}

Preuve. On note $A(T)$ cet appartement.  Puisque $G$ est suppos\'e semi-simple et puisque $p$ est grand, pour tout $x\in Imm(G_{AD})$, $\mathfrak{k}_{x,0}$ n'est autre que l'image r\'eciproque dans $\mathfrak{g}(F)$  du r\'eseau $\mathfrak{k}_{AD,x,0}$. On ne perd rien \`a remplacer $G$ par $G_{AD}$ et \`a supposer, le temps de cette d\'emonstration, que $G$ est adjoint. Le tore $T$ \'etant d\'eploy\'e, a une structure naturelle sur $\mathfrak{o}_{F}$ et l'hypoth\`ese entra\^{\i}ne que $X\in \mathfrak{t}(\mathfrak{o}_{F})$. Or cet ensemble est inclus par construction dans $\mathfrak{k}_{x,0}$ pour tout $x\in A(T)$. Inversement, soit $x\in Imm(G_{AD})$ tel que $X\in \mathfrak{k}_{x,0}$. Fixons un point $y\in A(T)$. Comme on vient de le voir, on a aussi $X\in \mathfrak{k}_{y,0}$. Tra\c{c}ons la g\'eod\'esique $[x,y]$ reliant $x$ \`a $y$. 
Il est bien connu que $\mathfrak{k}_{x,0}\cap \mathfrak{k}_{y,0}\subset \mathfrak{k}_{z,0}$ pour tout $z\in [x,y]$. En particulier, $X\in \mathfrak{k}_{z,0}$. Notons $u$ le point le plus proche de $x$ dans l'ensemble ferm\'e $[x,y]\cap A(T)$. Si $u=x$, on a $x\in A(T)$ comme on le voulait. Supposons $u\not=x$.  Notons ${\cal F}\in Fac(G)$ la facette contenant $u$. Il existe une facette ${\cal F}'\in Fac(G)$ telle que tout point de $[x,u[$ suffisamment proche de $u$ appartienne \`a ${\cal F}'$. La facette ${\cal F}$ est contenue dans l'adh\'erence $\bar{{\cal F}}'$ de ${\cal F}'$. Comme on l'a dit ci-dessus, on a  $X\in \mathfrak{k}_{{\cal F}}\cap \mathfrak{k}_{{\cal F}'}$. 
Puisque $u\in A(T)$, le tore $T$ se r\'eduit en un sous-tore maximal $T_{{\cal F}}$ de $G_{{\cal F}_{u}}$. L'\'el\'ement $X\in \mathfrak{k}_{{\cal F},0}$ se r\'eduit en un \'el\'ement $X_{{\cal F}}\in \mathfrak{t}_{{\cal F}}({\mathbb F}_{q})$. Pour toute racine $\beta_{{\cal F}}$ de $T_{{\cal F}}$ dans $G_{{\cal F}}$, il existe une racine $\beta$ de $T$ dans $G$ de sorte que $\beta_{{\cal F}}(X_{{\cal F}})$ soit la r\'eduction dans ${\mathbb F}_{q}$ de $\beta(X)\in \mathfrak{o}_{F}$. L'hypoth\`ese entra\^{\i}ne donc que $X_{{\cal F}}$ est r\'egulier. Puisque ${\cal F}\subset \bar{{\cal F}}'$, il existe un sous-groupe parabolique $P_{{\cal F}'}\subset G_{{\cal F}}$ de sorte que l'image de $\mathfrak{k}_{{\cal F}'}\cap \mathfrak{k}_{{\cal F}}$  dans $\mathfrak{g}_{{\cal F}}({\mathbb F}_{q})$ soit \'egale \`a $\mathfrak{p}_{{\cal F}'}({\mathbb F}_{q})$ et on sait que ${\cal F}'\subset A(T)$ si et seulement si $P_{{\cal F}'}$ contient $T_{{\cal F}}$. Puisque $X\in  
\mathfrak{k}_{{\cal F}}\cap \mathfrak{k}_{{\cal F}'}$, on a $X_{{\cal F}}\in \mathfrak{p}_{{\cal F}'}({\mathbb F}_{q})$. Puisque $X_{{\cal F}}$ est r\'egulier, cela entra\^{\i}ne que $P_{{\cal F}'}$ contient $T_{{\cal F}}$. Donc ${\cal F}'\subset A(T)$ mais cela contredit la d\'efinition de $u$. Cette contradiction ach\`eve la preuve. $\square$

\subsection{Profondeur des \'el\'ements de $FC(\mathfrak{g}(F))$ }\label{profondeur}
On suppose que $G$ n'est pas de type $A_{n-1}$. En  \ref{alcoves}, on a introduit l'appartement ${\cal A}^{nr}$  de $Imm_{F^{nr}}(G_{AD})$ et l'alc\^ove $C^{nr}$ et on a d\'ecrit ces ensembles. Il y a une relation affine 
$$(1) \qquad \sum_{\alpha\in \Delta_{a}^{nr}}d(\alpha)\alpha^{aff}=1,$$
 cf. \ref{alcoves}(2). Fixons un sommet $s$ de $C^{nr}$ fix\'e par $\Gamma_{F}^{nr}$ donc correspondant \`a une racine $\alpha_{s}\in \Delta_{a}^{nr}$ \'egalement fix\'ee par ce groupe. Les coordonn\'ees de $s$ sont  $\alpha^{aff}(s)=0$ pour $\alpha\in \Delta_{a}^{nr}-\{\alpha_{s}\}$ et $\alpha_{s}^{aff}(s)=\frac{1}{d(\alpha_{s})}$.  
 
  Le groupe $G_{s}$ a un syst\`eme de racines sur $\bar{{\mathbb F}}_{q}$ dont une base s'identifie \`a $\Delta_{a}^{nr}-\{\alpha_{s}\}$. Consid\'erons une fonction $f\in fc(\mathfrak{g}_{s}({\mathbb F}_{q}))$. On a introduit en \ref{groupessurFq} une fonction $\tilde{f}$ \`a support dans le radical nilpotent $\mathfrak{u}_{P}({\mathbb F}_{q})$ de l'alg\`ebre de Lie d'un sous-groupe parabolique $P$ de $G_{s}$. Comme on l'a dit en \ref{lespaceFC}, on peut identifier ces deux fonctions \`a des fonctions sur $\mathfrak{g}(F)$, \`a support dans $\mathfrak{k}_{s}$. En poussant ces fonctions dans $I(\mathfrak{g}(F))$, on a alors l'\'egalit\'e
$$f=\vert G_{s}({\mathbb F}_{q})/P({\mathbb F}_{q})\vert \tilde{f}.$$
Le support de la fonction $\tilde{f}$ est contenu dans l'image r\'eciproque de $\mathfrak{u}_{P}({\mathbb F}_{q})$ dans $\mathfrak{k}_{s}$. Notons $\mathfrak{k}_{s}[\mathfrak{u}_{P}]$ cet ensemble. Le sous-groupe $P$ est standard donc associ\'e \`a un sous-ensemble $\Delta(U_{P})\subset \Delta_{a}^{nr}-\{\alpha_{s}\}$. C'est-\`a-dire que, pour $\beta\in \Delta_{a}^{nr}-\{\alpha_{s}\}$, un \'el\'ement de base $E_{\beta}$ de l'espace radiciel dans $\mathfrak{g}_{s}$ associ\'e \`a $\beta$ appartient \`a $\mathfrak{u}_{P}$ si et seulement si $\beta\in \Delta(U_{P})$. Cet ensemble $\Delta(U_{P})$ est invariant par l'action de $\Gamma_{{\mathbb F}_{q}}$ sur $\Delta_{a}^{nr}$.  
Posons 
$$r=\frac{1}{\sum_{\alpha\in \Delta(U_{P})\cup\{\alpha_{s}\}}d(\alpha)}.$$
   Introduisons le point $x\in {\cal A}^{nr}$ tel que $\alpha^{aff}(x)=0$ pour $\alpha\in \Delta_{a}^{nr}-(\Delta(U_{P})\cup\{\alpha_{s}\})$ et $\alpha^{aff}(x)=r$ pour $\alpha\in \Delta(U_{P})\cup\{\alpha_{s}\}$. Ce point est invariant par l'action galoisienne de ${\mathbb F}_{q}$ sur ${\cal A}^{nr}$ donc appartient \`a $Imm_{F}(G_{AD})$.
 
 \begin{lem}{On a l'\'egalit\'e $\mathfrak{k}_{s}[\mathfrak{u}_{P}]=\mathfrak{k}_{x,r}$. De plus, soient $y\in Imm_{F}(G_{AD})$ et $v\in {\mathbb R}$. Supposons $\mathfrak{k}_{s}[\mathfrak{u}_{P}]\subset \mathfrak{k}_{y,v}$. Alors $v\leq r$. }\end{lem}

 Preuve. La premi\`ere \'etape est d'\'etendre le corps de base $F$ en $F^{nr}$. Le $\mathfrak{o}_{F}$-r\'eseau $\mathfrak{k}_{s}[\mathfrak{u}_{P}]$ a un analogue sur $F^{nr}$, qui est un $\mathfrak{o}_{F^{nr}}$-r\'eseau $\mathfrak{k}_{s,F^{nr}}[\mathfrak{u}_{P}]$. Il est stable par l'action naturelle de $\Gamma_{F}^{nr}$ sur $\mathfrak{g}(F^{nr})$ et $\mathfrak{k}_{s}[\mathfrak{u}_{P}]$ est son sous-ensemble des points fixes. De m\^eme, pour $y\in Imm_{F}(G_{AD})$ et $v\in {\mathbb R}$, le r\'eseau $\mathfrak{k}_{y,v}$ a un analogue sur $F^{nr}$, qui est un $\mathfrak{o}_{F^{nr}}$-r\'eseau $\mathfrak{k}_{y,v,F^{nr}}$. Il est stable par l'action naturelle de $\Gamma_{F}^{nr}$ sur $\mathfrak{g}(F^{nr})$ et $\mathfrak{k}_{y,v}$ est son sous-ensemble des points fixes. La propri\'et\'e g\'en\'erale suivante est bien connue. Soit $V$ un espace de dimension finie sur $F$. Posons $V^{nr}=V\otimes_{F}F^{nr}$.  Le groupe $\Gamma_{F}^{nr}$ agit naturellement sur  cet espace. Soit $R$ un sous-$\mathfrak{o}_{F^{nr}}$-r\'eseau de $V^{nr}$. Supposons $R$ stable par l'action de $\Gamma_{F}^{nr}$, notons $R^{\Gamma_{F}^{nr}}$ son ensemble de points fixes. Alors $R^{\Gamma_{F}^{nr}}$ est un sous-$\mathfrak{o}_{F}$-r\'eseau de $V$ et $R$ est le sous-$\mathfrak{o}_{F^{nr}}$-module de $V^{nr}$ engendr\'e par $R^{\Gamma_{F}^{nr}}$. En utilisant cette propri\'et\'e, on voit que l'\'egalit\'e $\mathfrak{k}_{s}[\mathfrak{u}_{P}]=\mathfrak{k}_{y,v}$, resp. l'inclusion $\mathfrak{k}_{s}[\mathfrak{u}_{P}]\subset \mathfrak{k}_{y,v}$, est \'equivalente \`a l'\'egalit\'e $\mathfrak{k}_{s,F^{nr}}[\mathfrak{u}_{P}]=\mathfrak{k}_{y,v,F^{nr}}$, resp. l'inclusion $\mathfrak{k}_{s,F^{nr}}[\mathfrak{u}_{P}]\subset \mathfrak{k}_{y,v,F^{nr}}$. Alors le lemme r\'esulte de l'assertion plus forte:
 
 (2) on a l'\'egalit\'e $\mathfrak{k}_{s,F^{nr}}[\mathfrak{u}_{P}]=\mathfrak{k}_{x,r,F^{nr}}$; de plus, soient $y\in Imm_{F^{nr}}(G_{AD})$ et $v\in {\mathbb R}$; supposons $\mathfrak{k}_{s,F^{nr}}[\mathfrak{u}_{P}]\subset \mathfrak{k}_{y,v,F^{nr}}$;  alors $v\leq r$.
 
 Il convient de distinguer  deux cas selon que $G$ est d\'eploy\'e ou non sur $F^{nr}$. Nous ne traiterons que le second cas, le premier \'etant similaire et plus simple. On reprend les notations de \ref{orbites} et \ref{alcoves}.  On a compl\'et\'e la paire $(B,T)$ par un \'epinglage  $(E_{\beta})_{\beta\in \Delta}$. On note $E$ la plus petite extension de $F^{nr}$ sur laquelle $G$ est d\'eploy\'e. On suppose que $I_{F}/\Gamma_{E}=\Gamma_{E/F^{nr}}\simeq {\mathbb Z}/e{\mathbb Z}$ pour un entier $e\geq2$. On fixe un \'el\'ement $\rho\in I_{F}$ d'image $1$ dans ce groupe cyclique et on note $\theta$ l'automorphisme de $G$ tel que $\rho$ agisse (alg\'ebriquement) sur $G$ par   $\theta$. L'automorphisme $\theta$ est  d'ordre $e$ et  pr\'eserve $(B,T)$ et l'\'epinglage. On a compl\'et\'e l'\'epinglage en une base de Chevalley   de sorte que $\mathfrak{k}_{s^{nr},F^{nr}}$ soit le sous-ensemble des points fixes par $I_{F}$ dans le $\mathfrak{o}_{\bar{F}}$-r\'eseau engendr\'e par cette base. La base de Chevalley contient une famille $(E_{\beta})_{\beta\in \Sigma}$ et, parce que l'on suppose que $G$ n'est pas de type $A_{n-1}$, on peut supposer que $E_{\theta(\beta)}=\theta(E_{\beta})$ pour tout $\beta\in \Sigma$, cf. \cite{KS} 1.3.  L'ensemble $\Sigma^{nr}$ des restrictions \`a $X_{*}(T)^{I_{F}}$ des \'el\'ements de $\Sigma$ s'identifie \`a l'ensemble des orbites de l'action de $\theta$ dans $\Sigma$.  On a identifi\'e l'appartement  ${\cal A}^{nr}$ associ\'e \`a $T^{nr}$ dans $Imm_{F^{nr}}(G_{AD})$ \`a $X_{*}(T)^{I_{F}}\otimes _{{\mathbb Z}}{\mathbb R}$ de sorte que $s^{nr}$ s'identifie \`a $0$.  Soit $v\in {\mathbb R}$. On  pose $\mathfrak{t}(F^{nr})_{v}=\{X\in \mathfrak{t}(F^{nr}); \forall \beta\in \Sigma, val_{F}(\beta(X))\geq v\}$. Pour $\alpha\in \Sigma^{nr}$, notons $\mathfrak{u}(\alpha)_{v}$ l'ensemble des $X=\sum_{\beta\in \Sigma, \beta^{res}=\alpha}\lambda_{\beta}E_{\beta}$ avec des $\lambda_{\beta}\in E$, tels que $\lambda_{\theta(\beta)}=\rho(\lambda_{\beta})$ et $val_{F}(\lambda_{\beta})\geq v$ pour tout $\beta\in \Sigma$ tel que $\beta^{res}=\alpha$.

 Soient $y\in {\cal A}^{nr}$ et $v\in {\mathbb R}$ . Alors, par d\'efinition,
 
 (3)  le r\'eseau $ \mathfrak{k}_{y,v,F^{nr}}$ est la somme de $\mathfrak{t}(F^{nr})_{v}$ et des $\mathfrak{u}(\alpha)_{-\alpha(y)+v}$ sur les $\alpha\in \Sigma^{res}$.

 Du tore $T^{nr}$ se  d\'eduit un sous-tore maximal $T_{s}$ de $G_{s}$. Notons $\Sigma^{\star}(G_{s})$ le sous-ensemble des racines affines $\alpha^{aff}\in \Sigma^{aff}$ telles que $\alpha^{aff}(s)=0$. C'est l'ensemble des $\alpha[-\alpha(s)]$ quand $\alpha$ parcourt le sous-ensemble   $\Sigma^{nr}(G_{s})$ des $\alpha\in \Sigma^{nr}$ tels que $\alpha(s)\in \frac{1}{e(\alpha)}{\mathbb Z}$.  L'ensemble $\Sigma^{nr}(G_{s})$ est le syst\`eme de racines de $G_{s}$ relatif \`a $T_{s}$ et 
  $\Delta_{a}^{nr}-\{\alpha_{s}\}$ en  est une base.  On note $\alpha\mapsto \alpha^{\star}=\alpha[-\alpha(s)]$ la bijection de $\Sigma^{nr}(G_{s})$ sur $\Sigma^{\star}(G_{s})$. Pour $\alpha\in \Delta_{a}^{nr}-\{\alpha_{s}\}$, on a $\alpha^{\star}=\alpha^{aff}$.  L'espace $\mathfrak{g}_{s}(\bar{{\mathbb F}}_{q})$
s'identifie \`a
$$(4)\qquad \mathfrak{t}(F^{nr})_{0}/\mathfrak{t}(F^{nr})_{\frac{1}{e}}\oplus\oplus_{\alpha\in \Sigma^{nr}(G_{s})}\mathfrak{u}(\alpha)_{-\alpha(s)}/\mathfrak{u}(\alpha)_{-\alpha(s)+\frac{1}{e(\alpha)}}.$$
Remarquons que

(5) si $\Sigma^{nr}\not=\Sigma^{nr}(G_{s})$, alors $\alpha_{s}\not=\alpha_{0}$.

En effet, si $\alpha_{s}=\alpha_{0}$, alors $\alpha(s)=0$ pour tout $\alpha\in \Delta^{nr}$ donc aussi $\alpha(s)=0$ pour tout $\alpha\in \Sigma^{nr}$, a fortiori $\alpha(s)\in \frac{1}{e(\alpha)}{\mathbb Z}$ pour tout $\alpha\in \Sigma^{nr}$.

On a besoin des quelques propri\'et\'es suivantes. Soit $\alpha\in \Sigma^{nr}$, \'ecrivons $\alpha=\epsilon \sum_{\alpha'\in \Delta^{nr}}m(\alpha')\alpha'$ avec $\epsilon\in \{\pm 1\}$ et  des coefficients $m(\alpha')\in {\mathbb N}$. Alors

(6) $m(\alpha')e(\alpha)\leq d(\alpha')$ pour tout $\alpha'\in \Delta^{nr}$.

On peut supposer $\epsilon=1$. L'\'egalit\'e  $\alpha= \sum_{\alpha'\in \Delta^{nr}}m(\alpha')\alpha'$ se r\'ecrit
 $$\check{\alpha}=\sum_{\alpha'\in \Delta^{nr}}e(\alpha)m(\alpha')e(\alpha')^{-1}\check{\alpha}'.$$
  Puisque $\check{\alpha}_{0}=-\sum_{\alpha'\in \Delta^{nr}}d(\check{\alpha}')\check{\alpha}'$ est l'oppos\'ee de la plus grande coracine, on a alors $e(\alpha)m(\alpha')e(\alpha')^{-1}\leq d(\check{\alpha}')$ pour tout $\alpha'\in \Delta^{nr}$, d'o\`u $e(\alpha)m(\alpha')\leq e(\alpha')d(\check{\alpha}')=d(\alpha')$, d'o\`u (6). 
  
   Soit maintenant  $\alpha\in \Sigma^{nr}(G_{s})$. Ecrivons $\alpha=\epsilon\sum_{\alpha'\in \Delta_{a}^{nr}-\{\alpha_{s}\}}m(\alpha')\alpha'$ avec $\epsilon\in \{\pm 1\}$ et  des coefficients $m(\alpha')\in {\mathbb N}$. Alors
 
 (7) $m(\alpha')e(\alpha)\leq d(\alpha')$ pour tout $\alpha'\in \Delta_{a}^{nr}-\{\alpha_{s}\}$.
 
 On peut encore supposer $\epsilon=1$. Si $\alpha_{s}=\alpha_{0}$ ou si $\alpha_{0}\in \Delta_{a}^{nr}-\alpha_{s}$ et $m(\alpha_{0})=0$, la propri\'et\'e r\'esulte de (6). Supposons $\alpha_{0}\not=\alpha_{s}$ et $m(\alpha_{0})\geq1$. 
  Comme dans la preuve de (6), on r\'ecrit notre \'egalit\'e
 $$\check{\alpha}=\sum_{\alpha'\in \Delta_{a}^{nr}-\{\alpha_{s}\}}e(\alpha)m(\alpha')e(\alpha')^{-1}\check{\alpha}'.$$
 On remplace $\check{\alpha}_{0}$ par $-\sum_{\alpha'\in \Delta^{nr}}d(\check{\alpha}')\check{\alpha}'$. On obtient
 $$\check{\alpha}=-e(\alpha)m(\alpha_{0})e(\alpha_{0})^{-1}d(\check{\alpha}_{s})\check{\alpha}_{s}-\sum_{\alpha'\in \Delta^{nr}-\{\alpha_{s}\}} (e(\alpha)m(\alpha_{0})e(\alpha_{0})^{-1}d(\check{\alpha}')-e(\alpha)m(\alpha')e(\alpha')^{-1})\check{\alpha}'.$$
 C'est l'\'ecriture de $\check{\alpha}$ dans la base $\check{\Delta}^{nr}$. Le coefficient de $\check{\alpha}_{s}$ est strictement n\'egatif donc les autres sont n\'egatifs ou nuls. Pour le coefficient de $\check{\alpha}_{s}$, on a la m\^eme majoration que dans la preuve de (6):
 $$e(\alpha)m(\alpha_{0})e(\alpha_{0})^{-1}d(\check{\alpha}_{s})\leq d(\check{\alpha}_{s}),$$
 d'o\`u
 $$(8) \qquad e(\alpha)m(\alpha_{0})e(\alpha_{0})^{-1}\leq 1.$$
 Puisque $d(\check{\alpha}_{0})=1$, on a $d(\alpha_{0})=e(\alpha_{0})$ et l'in\'egalit\'e pr\'ec\'edente d\'emontre (7) pour $\alpha'=\alpha_{0}$. Pour $\alpha'\in \Delta^{nr}-\{\alpha_{s}\}$, la
 n\'egativit\'e du coefficient de $\check{\alpha}'$ dans l'expression de $\check{\alpha}$ entra\^{\i}ne
 $$e(\alpha)m(\alpha')e(\alpha')^{-1}\leq e(\alpha)m(\alpha_{0})e(\alpha_{0})^{-1}d(\check{\alpha}'),$$
 puis, gr\^ace \`a (8),
 $$e(\alpha)m(\alpha')e(\alpha')^{-1}\leq d(\check{\alpha}').$$
 Cela \'equivaut \`a (7) pour $\alpha'$. Cela ach\`eve la preuve de cette relation.

         On note $\Sigma^{nr}(U_{P})$, resp. $\Sigma^{nr}(\bar{P})$, le sous-ensemble des $\alpha\in \Sigma^{nr}(G_{s})$ qui interviennent, resp. n'interviennent pas, dans le radical nilpotent $\mathfrak{u}_{P}$ de $\mathfrak{p}$. Pour $\alpha\in \Sigma^{nr}$, posons
 $$r(\alpha)=\left\lbrace\begin{array}{cc}\frac{[-e(\alpha)\alpha(s)]+1}{e(\alpha)},&\text{\,si\,}\alpha\not\in \Sigma^{nr}(G_{s}),\\ -\alpha(s),&\text{\, si\,}\alpha\in \Sigma^{nr}(U_{P}),\\ -\alpha(s)+\frac{1}{e(\alpha)},&\text{\, si\, }\alpha\in \Sigma^{nr}(\bar{P}).\\ \end{array}\right.$$
 On voit que
 
 (9) le r\'eseau $\mathfrak{k}_{s,F^{nr}}[\mathfrak{u}_{P}]$ est somme de $\mathfrak{t}(F^{nr})_{\frac{1}{e}}$, des $\mathfrak{u}(\alpha)_{r(\alpha)}$ pour tout $\alpha\in \Sigma^{nr}$.

 En comparant (3) et (9), on voit que, pour prouver l'\'egalit\'e $\mathfrak{k}_{s,F^{nr}}[\mathfrak{u}_{P}]=\mathfrak{k}_{x,r,F^{nr}}$, il suffit de prouver les in\'egalit\'es
 
  (10) $0<r\leq \frac{1}{e}$;
 
 (11) pour tout $\alpha\in \Sigma^{nr}$, $r(\alpha)-\frac{1}{e(\alpha)}<-\alpha(x)+r\leq r(\alpha)$.
 
 En reprenant la d\'efinition des nombres $d(\alpha)$ et en inspectant tous les syst\`emes de racines irr\'eductibles dont les racines n'ont pas toutes m\^eme longueur, on voit que $d(\alpha)$ est un multiple de $e$ pour tout $\alpha\in \Delta_{a}^{nr}$. L'in\'egalit\'e (10) r\'esulte alors de la d\'efinition de $r$.
 
 Soit $\alpha\in \Sigma^{nr}(G_{s})$.   Ecrivons $\alpha=\epsilon\sum_{\alpha'\in \Delta_{a}^{nr}-\{\alpha_{s}\}}m(\alpha')\alpha'$ avec $\epsilon\in \{\pm 1\}$ et  des coefficients $m(\alpha')\in {\mathbb N}$. Posons $c(\alpha)=\sum_{\alpha'\in \Delta(U_{P})}m(\alpha')$. On a $\alpha(x)-\alpha(s)=\alpha^{\star}(x)=\epsilon c(\alpha)r$. On voit que la relation (11) pour $\alpha$ \'equivaut \`a
 $$ \left\lbrace\begin{array}{cc}r\leq  c(\alpha)r<r+\frac{1}{e(\alpha)},&\rm{\, si\,}\epsilon=1\rm{ \, et\, }c(\alpha)\geq1,\\ r-\frac{1}{e(\alpha)}\leq \epsilon c(\alpha)r<r,&\rm{\, sinon.}\\ \end{array}\right.$$
 La premi\`ere in\'egalit\'e de la premi\`ere ligne ainsi que la seconde in\'egalit\'e de la seconde ligne sont triviales. La seconde in\'egalit\'e de la premi\`ere ligne est moins forte que la premi\`ere in\'egalit\'e de la seconde ligne appliqu\'ee \`a $-\alpha$. Il suffit de d\'emontrer cette premi\`ere in\'egalit\'e de la seconde ligne. Si $c(\alpha)=0$, elle r\'esulte de (10). Supposons $c(\alpha)\geq1$. Les hypoth\`eses de la seconde ligne impliquent alors $\epsilon=-1$ et l'in\'egalit\'e \'equivaut \`a 
 
 (12) $c(\alpha)r\leq \frac{1}{e(\alpha)}-r$. 
 
 \noindent La relation (7) nous dit que  $e(\alpha)c(\alpha)\leq \frac{1}{r}-d(\alpha_{s})$. D'o\`u
 
 (13) $c(\alpha)r\leq \frac{1}{e(\alpha)}-r\frac{d(\alpha_{s})}{e(\alpha)}$.
 
 \noindent  Comme on l'a dit ci-dessus, on a $d(\alpha_{s})\geq e\geq e(\alpha)$. Alors (13) entra\^{\i}ne (12), ce qui prouve (11) dans le cas o\`u $\alpha\in \Sigma^{nr}(G_{s})$. 
 
 Soit $\alpha\in \Sigma^{nr}-\Sigma^{nr}(G_{s})$. On a $\alpha_{s}\not=\alpha_{0}$ d'apr\`es (5). Ecrivons   $\alpha=\epsilon \sum_{\alpha'\in \Delta^{nr}}m(\alpha')\alpha'$ avec $\epsilon\in \{\pm 1\}$ et  des coefficients $m(\alpha')\in {\mathbb N}$. On a $\alpha(s)=\epsilon \frac{m(\alpha_{s})}{d(\alpha_{s})}$. D'apr\`es (5), on a $0\leq \frac{m(\alpha_{s})}{d(\alpha_{s})}\leq \frac{1}{e(\alpha)}$, c'est-\`a-dire $0\leq \vert \alpha(s)\vert \leq \frac{1}{e(\alpha)}$. Mais l'hypoth\`ese $\alpha\not\in \Sigma^{nr}(G_{s})$ signifie que $\alpha(s)\not\in \frac{1}{e(\alpha)}{\mathbb Z}$, donc $0<  \vert \alpha(s)\vert < \frac{1}{e(\alpha)}$. On calcule alors $r(\alpha)=0$ si $\epsilon=1$, $r(\alpha)=\frac{1}{e(\alpha)}$ si $\epsilon=-1$. Posons $c(\alpha)=m(\alpha_{s})+\sum_{\alpha'\in \Delta^{nr}\cap \Delta(U_{P})}m(\alpha')$. Alors $\alpha(x)=\epsilon c(\alpha)r$. On voit que la relation (11) pour $\alpha$ \'equivaut \`a
 $$ \left\lbrace\begin{array}{cc}r\leq  c(\alpha)r<r+\frac{1}{e(\alpha)},&\rm{\, si\,}\epsilon=1 \\ -r<  c(\alpha)r\leq \frac{1}{e(\alpha)}-r,&\rm{\, si\,}\epsilon=-1.\\ \end{array}\right.$$
 Il suffit de d\'emontrer que
 $$(14) \qquad r\leq c(\alpha)r\leq \frac{1}{e(\alpha)}-r.$$
 On a vu que $m(\alpha_{s})\geq1$ donc $c(\alpha)\geq1$. Cela entra\^{\i}ne la premi\`ere in\'egalit\'e. La relation (6) nous dit que $e(\alpha)c(\alpha)\leq \frac{1}{r}-d(\alpha_{0})$. D'o\`u
 
 (15) $c(\alpha)r\leq \frac{1}{e(\alpha)}-r\frac{d(\alpha_{0})}{e(\alpha)}$.
 
 \noindent  Comme on l'a dit ci-dessus, on a $d(\alpha_{0})\geq e\geq e(\alpha)$ (en fait $d(\alpha_{0})=e$). Alors (15) entra\^{\i}ne la deuxi\`eme in\'egalit\'e de  (14), ce qui ach\`eve de prouver (11). 
 
 On a ainsi prouv\'e la premi\`ere assertion de (2).

   Notons $W^{G_{s}}$ le groupe de Weyl de $G_{s}$ relatif \`a $T_{s}$. Il s'identifie au quotient $(K_{s}\cap Norm_{G(F^{nr})}(T^{nr}))/T(\mathfrak{o}_{\bar{F}})^{I_{F}}$. Pour $w\in W^{G_{s}}$, fixons un rel\`evement $w_{G}$ de $w$ dans $K_{s}\cap Norm_{G(F^{nr})}(T^{nr})$. L'\'el\'ement $w$ agit naturellement 
  sur $\Sigma^{nr}(G_{s})$. L'\'el\'ement $w_{G}$ agit naturellement sur   ${\cal A}^{nr}$ et il s'en d\'eduit une action sur $\Sigma^{aff}$. On a la relation  suivante entre ces deux actions: 
  
  (16) pour $\alpha\in \Sigma^{nr}(G_{s})$ et $w\in W^{G_{s}}$, $(w(\alpha))^{\star}=w_{G}(\alpha^{\star})$. 
  
 \noindent Le groupe $K_{s}\cap Norm_{G(F^{nr})}(T^{nr})$ agit par conjugaison sur $\mathfrak{g}(F^{nr})$ en conservant $\mathfrak{k}_{s}$. L'action du sous-groupe $T(\mathfrak{o}_{\bar{F}})^{I_{F}}$ conserve $\mathfrak{k}_{s,F^{nr}}[\mathfrak{u}_{P}]$. Ainsi, pour $w\in W^{G_{s}}$, le r\'eseau  $w_{G}(\mathfrak{k}_{s,F^{nr}}[\mathfrak{u}_{P}])$ est bien d\'efini. 
Notons $M$ le Levi standard de $P$. Son ensemble de racines a pour base $\Delta_{a}^{nr}-(\Delta(U_{P})\cup\{\alpha_{s}\})$. Notons  $W^M\subset W^{G_{s}}$ le groupe de Weyl de $M$ relatif \`a $T_{s}$.  On voit \`a l'aide de (16) que $w_{G}$ conserve le r\'eseau $\mathfrak{k}_{s,F^{nr}}[\mathfrak{u}_{P}]$ pour tout $w\in W^M$. 

 Soient $y\in {\cal A}^{nr}$ et $v\in {\mathbb R}$. Supposons  
 
 (17) $\mathfrak{k}_{s,F^{nr}}[\mathfrak{u}_{P}]\subset \mathfrak{k}_{y,v,F^{nr}}$. 
 
 \noindent D'apr\`es ce que l'on vient de voir, on a aussi $\mathfrak{k}_{s,F^{nr}}[\mathfrak{u}_{P}]\subset \mathfrak{k}_{w_{G}(y),v,F^{nr}}$ pour tout $w\in W^M$. En rempla\c{c}ant $y$ par $w_{G}(y)$ pour un $w\in W^M$ bien choisi, on peut supposer $\alpha^{aff}(y)\geq0$ pour tout $\alpha\in \Delta_{a}^{nr}-(\Delta(U_{P})\cup\{\alpha_{s}\})$. Soit $\alpha\in \Delta(U_{P})\cup\{\alpha_{s}\}$. En comparant (3) et (9), l'inclusion (17) implique $-\alpha(y)+v\leq r(\alpha)$. Pour $\alpha\in \Delta(U_{P})$, on a $r(\alpha)=-\alpha(s)$ et l'in\'egalit\'e pr\'ec\'edente \'equivaut \`a   $\alpha^{aff}(y)\geq v$. 
 Consid\'erons la racine $\alpha_{s}$. Par construction, on a $d(\alpha_{s})\geq e(\alpha_{s})$. Supposons d'abord $d(\alpha_{s})> e(\alpha_{s})$, ce qui implique $\alpha_{s}\not=\alpha_{0}$ d'apr\`es la d\'efinition de $d(\alpha_{s})$. Alors $\alpha_{s}\not\in \Sigma^{nr}(G_{s})$ et $r(\alpha_{s})=0$. L'in\'egalit\'e $-\alpha_{s}(y)+v\leq r(\alpha_{s})$ \'equivaut encore \`a $\alpha_{s}^{aff}(y)\geq v$. Supposons maintenant $d(\alpha_{s})=e(\alpha_{s})$.  Alors $\alpha_{s}\in \Sigma^{nr}(G_{s})$. La relation (1) entra\^{\i}ne que $\alpha_{s}$ est combinaison lin\'eaire \`a coefficients strictement n\'egatifs des \'el\'ements de $\Delta_{a}^{nr}-\{\alpha_{s}\}$ donc $\alpha_{s}\in \Sigma^{nr}(\bar{P})$. Alors $r(\alpha_{s})=-\alpha_{s}(s)+\frac{1}{e(\alpha_{s})}=-\alpha_{s}(s)+\frac{1}{d(\alpha_{s})}=-\alpha_{s}(s)+\alpha_{s}^{aff}(s)$. L'in\'egalit\'e $-\alpha_{s}(y)+v\leq r(\alpha_{s})$ \'equivaut \`a $\alpha_{s}(y)-\alpha_{s}(s)+\alpha_{s}^{aff}(s)\geq v$. Mais $\alpha_{s}(y)-\alpha_{s}(s)=\alpha^{aff}_{s}(y)-\alpha^{aff}_{s}(s)$ et on obtient encore $\alpha_{s}^{aff}(y)\geq v$. On a obtenu
 
 $\alpha^{aff}(y)\geq0$ pour tout $\alpha\in \Delta_{a}^{nr}-(\Delta(U)\cup\{\alpha_{s}\})$;
 
 $\alpha^{aff}(y)\geq v$ pour tout $\alpha\in \Delta(U)\cup \{\alpha_{s}\}$. 
 
 \noindent L'\'egalit\'e (1) et la d\'efinition de $r$ entra\^{\i}nent alors $v\leq r$. C'est la deuxi\`eme assertion de (2) avec la restriction que l'on a suppos\'e $y\in {\cal A}^{nr}$. 
 
 Levons cette restriction. Soit $y\in Imm_{F^{nr}}(G_{AD})$ et $v\in {\mathbb R}$. Supposons l'inclusion (17) v\'erifi\'ee.  Introduisons la g\'eod\'esique $[x,y]$ joignant $x$ \`a $y$ dans $Imm_{F^{nr}}(G_{AD})$. Elle est incluse dans un appartement qui est un espace euclidien dont on note la distance $\vert .\vert $.  Pour $t\in [0,1]$, on note $x_{t}\in [x,y]$ le point tel que  $\vert x_{t}-x\vert =t\vert y-x\vert $ et on pose $v_{t}=(1-t)v+tr$. On v\'erifie sur la formule (3) que
 $\mathfrak{k}_{x_{t}, v_{t},F^{nr}}$ contient $\mathfrak{k}_{y,v,F^{nr}}\cap \mathfrak{k}_{x,r,F^{nr}}$. D'apr\`es (17) et la premi\`ere assertion de (2) d\'ej\`a prouv\'ee, on a donc  $\mathfrak{k}_{s,F^{nr}}[\mathfrak{u}_{P}]\subset \mathfrak{k}_{x_{t},v_{t},F^{nr}}$. Pour un \'el\'ement $t\in ]0,1]$ assez petit, $x_{t}$ appartient \`a une facette ${\cal F}$ dont l'adh\'erence contient $x$. On sait qu'alors, il existe un \'el\'ement $k\in K_{x,F^{nr}}^0$ tel que $k{\cal F}\subset {\cal A}^{nr}$. Fixons un tel $k$. D'apr\`es la d\'efinition de $x$, $K_{x,F^{nr}}^0$ est l'image r\'eciproque dans $K_{s}^0$ de $P(\bar{{\mathbb F}}_{q})$. L'assertion (9) montre que la conjugaison par un tel $k$ conserve $\mathfrak{k}_{s,F^{nr}}[\mathfrak{u}_{P}]$. Donc $\mathfrak{k}_{s,F^{nr}}[\mathfrak{u}_{P}]\subset \mathfrak{k}_{k(x_{t}),v_{t},F^{nr}}$. On applique ce que l'on a d\'ej\`a prouv\'e \`a $k(x_{t})\in {\cal A}^{nr}$ et $v_{t}$. On en d\'eduit $v_{t}\leq r$. Puisque $t>0$, cela entraine $v\leq r$. Cela ach\`eve la preuve de (2) et du lemme. $\square$

 \subsection{Type $D_{4}$ trialitaire}\label{D4trialitairepadique}
  
  On suppose $G$ quasi-d\'eploy\'e de type $D_{4}$. On rappelle que le diagramme ${\cal D}$ de ce groupe a un groupe d'automorphismes isomorphe \`a $\mathfrak{S}_{3}$ et on a introduit des g\'en\'erateurs $\theta$ et $\theta_{3}$ de ce diagramme en \ref{D4trialitaire}. Le diagramme du groupe dual $\hat{G}$ est le m\^eme mais, pour plus de pr\'ecision, on le note $\hat{{\cal D}}$ et on note $\hat{\theta}$ et $\hat{\theta}_{3}$ les g\'en\'erateurs de son groupe d'automorphismes. On suppose que l'homorphisme $\sigma\mapsto \sigma_{G}$ de $\Gamma_{F}$ dans $Aut(\hat{\cal D})$ a pour image un sous-groupe d'ordre $\geq 3$. A l'aide de \ref{extensionsdiedrales}, on voit qu'il y a trois cas possibles:

  (A) soit $E/F$ l'extension non ramifi\'ee de degr\'e $3$;  $\Gamma_{E}$ agit trivialement sur $\hat{{\cal D}}$ et il existe un \'el\'ement de Frobenius $\rho\in \Gamma_{F}$ tel que $\rho_{G}= \hat{\theta}_{3}
$;

(B) on suppose que $\delta_{3}(q-1)=1$; soit $E$ une extension galoisienne ramifi\'ee de $F$ de degr\'e $3$;      $\Gamma_{E}$ agit trivialement sur $\hat{{\cal D}}$ et il existe  un \'el\'ement $\rho\in \Gamma_{F}-\Gamma_{E}$ tel que  $\rho_{G}=\hat{\theta}_{3}
$;

  (C) on suppose que $\delta_{3}(q+1)=1$; on note $E$ l'unique extension galoisienne de $F$ telle que $\Gamma_{E/F}\simeq \mathfrak{S}_{3}$ (cf. \ref{extensionsdiedrales}(3)); elle contient l'extension quadratique $E_{0}$ non ramifi\'ee de $F$; $\Gamma_{E}$ agit trivialement sur $\hat{{\cal D}}$ et il existe  un \'el\'ement $\rho\in \Gamma_{E_{0}}-\Gamma_{E}$ et un \'el\'ement de Frobenius $\tau\in \Gamma_{F}$ tels que $\rho_{G}=\hat{\theta}_{3}
$ et  $\tau_{G}=\hat{\theta}$.

Dans le cas (A), on pose ${\cal X}=\{2\}$. Dans les cas (B) ou (C), on pose ${\cal X}=\{0,134\}$. On pose $d_{x}=1$ pour tout $x\in {\cal X}$. 

 On se place d'abord dans le cas (A). Le groupe est de type $^3D_{4}$ dans les tables de Tits. L'index local (c'est-\`a-dire le diagramme ${\cal D}_{a}^{nr}$ de \ref{orbites}) est le diagramme de Dynkin affine  de type $D_{4}$ avec l'action galoisienne suivante: elle est triviale sur $\Gamma_{{\mathbb F}_{q^3}}$ et le Frobenius agit par $\theta_{3}$. Les \'el\'ements de $\underline{S}(G)$ sont en bijection avec les orbites de cette action galoisienne. L'action de $G_{AD}(F)/\pi(G(F))$  sur $\underline{S}(G)$ est triviale. L'orbite $\{\alpha_{1},\alpha_{3},\alpha_{4}\}$ est exclue par le lemme  \ref{orbites}. Pour le sommet $s$ associ\'e \`a l'orbite $\{\alpha_{0}\}$, on voit que $G_{s}$ est isog\`ene au groupe trialitaire de type $D_{4}$ sur ${\mathbb F}_{q}$, donc $FC(\mathfrak{g}_{s}({\mathbb F}_{q}))=\{0\}$ d'apr\`es \ref{D4trialitaire}. Consid\'erons  le sommet $s$ param\'etr\'e par l'orbite $\{\alpha_{2}\}$. On voit que $G_{s}\simeq (Res_{{\mathbb F}_{q^3}/{\mathbb F}_{q}}(SL(2))\times SL(2))/diag(\{\pm 1\})$, le groupe $\{\pm 1\}$ s'identifiant \'evidement aux centres de chacun des deux facteurs.  D'apr\`es \ref{An-1dep}, l'espace $FC(\mathfrak{g}_{s}({\mathbb F}_{q}))$ est une droite. Elle  donne naissance  \`a une droite dans   $FC(\mathfrak{g}(F))$ que l'on note $FC_{2}$. Cela d\'emontre \ref{resultats}(1) dans le cas (A).
  
  Consid\'erons le cas (B). Le groupe est de type $G_{2}^1$ dans les tables de Tits.  Explicitons les descriptions de \ref{alcoves} et \ref{profondeur}. On a ${\cal A}^{nr}=X_{*}(T)^{I_{F}}\otimes_{{\mathbb Z}}{\mathbb R}$. L'ensemble $\Sigma^{nr}$ est celui des $\beta^{res}$ pour $\beta\in \Sigma$. Puisque nous notons $\beta$ les \'el\'ements de $\Sigma$, on note l'ensemble $\Delta=\{\beta_{i};i=1,...,4\}$. La base $\Delta^{nr}$ est \'egale \`a $\{\alpha_{2},\alpha_{134}\}$, o\`u $\alpha_{2}=\beta_{2}^{res}$ et $\alpha_{134}=\beta_{1}^{res}=\beta_{3}^{res}=\beta_{4}^{res}$. Le sous-ensemble de racines positives dans $\Sigma^{nr}$ est $\{\alpha_{2},\alpha_{134},\alpha_{134}+\alpha_{2},2\alpha_{134}+\alpha_{2},3\alpha_{134}+\alpha_{2},3\alpha_{134}+2\alpha_{2}\}$. On a $\alpha_{0}=-\alpha_{2}-2\alpha_{134}$, $e(\alpha_{0})=e(\alpha_{134})=3$ et  $e(\alpha_{2})=1$. La relation (2) de \ref{alcoves} est
  $$ (1) \qquad 3\alpha_{0}+3\alpha_{2}+6\alpha_{134}=0.$$
  C'est-\`a-dire $d(\alpha_{0})=d(\alpha_{2})=3$, $d(\alpha_{134})=6$. On a $\alpha^{aff}=\frac{1}{3}+\alpha_{0}$. 
   L'ensemble $\underline{S}(G)$ s'identifie \`a celui des sommets de l'alc\^ove, donc \`a  l'ensemble $\Delta_{a}^{nr}=\{\alpha_{2},\alpha_{134},\alpha_{0}\}$. Introduisons un tore   $\underline{T}$ d\'efini  sur ${\mathbb F}_{q}$ tel que $X_{*}(\underline{T})=X_{*}(T)^{I_{F}}$, cet ensemble \'etant  muni de l'action triviale de $\Gamma_{{\mathbb F}_{q}}$.
    On le munit du m\^eme ensemble de racines que ci-dessus. Pour le sommet $s$ associ\'e \`a  $\alpha\in \Delta_{a}^{nr}$, le groupe $G_{s}$ est d\'eploy\'e, il a pour tore maximal   $\underline{T}$  et pour ensemble de racines le sous-ensemble  $\Sigma^{nr}(G_{s})$ de $\Sigma^{nr}$  d\'ecrit en \ref{profondeur}.  
    Pour la racine $\alpha_{0}$,  on voit que $G_{s}$ est de type $G_{2}$. Alors  $FC(\mathfrak{g}_{s}({\mathbb F}_{q}))$ est une droite d'apr\`es \ref{G2}. Il s'en d\'eduit une droite de    $FC(\mathfrak{g}(F))$ que l'on note $FC_{0}$.  Pour la racine $\alpha_{1,3,4}$, $G_{s}\simeq (SL(2)\times SL(2))/diag(\{\pm 1\})$. L'espace $FC(\mathfrak{g}_{s}({\mathbb F}_{q}))$ est une droite. Il s'en d\'eduit une droite de    $FC(\mathfrak{g}(F))$ que l'on note $FC_{134}$.    Pour la racine $\alpha_{2}$, on trouve $G_{s}=PGL(3)$ et $FC(\mathfrak{g}_{s}({\mathbb F}_{q}))=\{0\}$. Cela d\'emontre \ref{resultats}(1) dans le cas (B).
  
  Le calcul est le m\^eme dans le cas (C) car l'action de $\Gamma_{{\mathbb F}_{q}}$ sur ${\cal D}_{a}^{nr}$ (qui se d\'eduit de l'action de $\theta$) est triviale. Cela ach\`eve la preuve de \ref{resultats}(1).

On pose   ${\cal Y}={\cal X}$, ${\cal X}^{st}={\cal X}$ et on note $\varphi$ l'identit\'e de ${\cal X}$ sur ${\cal Y}$. Montrons que les assertions (2), (3) et (4) de \ref{resultats} sont v\'erifi\'ees. Il suffit de prouver que

(2) $FC(\mathfrak{g}(F))=FC^{st}(\mathfrak{g}(F))$.

En effet, on pose alors $FC^{{\cal E}}_{y}=FC_{y}$ pour tout $y\in {\cal Y}={\cal X}$ et ces assertions (2), (3) et (4) de \ref{resultats} deviennent triviales. 

Prouvons (1). On a d\'ecrit le groupe $\hat{\Omega} \simeq ({\mathbb Z}/2{\mathbb Z})^2$ en \ref{Dndeppairpadique} . Le groupe $Aut(\hat{{\cal D}}_{a})$ est isomorphe \`a $\mathfrak{S}_{4}=\Omega\rtimes \mathfrak{S}_{3}$.    On consid\`ere un  couple $(\sigma\mapsto \sigma_{G'},{\cal O})\in {\cal E}_{ell}(G)$. Notons $\hat{\Omega}_{G'}$ l'image de  $\Gamma_{E}$ dans $\hat{\Omega}$ par l'homomorphisme $\sigma\mapsto \omega_{G'}(\sigma)$.  C'est un sous-groupe de $\hat{\Omega}$ invariant par l'image de $\Gamma_{E/F}$ dans $\mathfrak{S}_{3}$ par l'application $\sigma\mapsto \sigma_{G}$. Cette image contient toujours $\hat{\theta}_{3}$. L'action par conjugaison de $\hat{\theta}_{3}$ dans $\hat{\Omega}$ permute cycliquement les $3$ \'el\'ements non triviaux de $\hat{\Omega}$. Donc $\hat{\Omega}_{G'}$ est \'egal \`a $\{1\}$ ou a $\hat{\Omega}$ tout entier. Si $\hat{\Omega}_{G'}=\hat{\Omega}$, on voit que le groupe $\Gamma_{E_{G'}/F}$ dans les cas (A) ou (B), resp. $\Gamma_{E_{G'}/E_{0}}$ dans le cas (C), a une structure qui est exclue par   \ref{extensions12}. Donc $\hat{\Omega}_{G'}=\{1\}$, c'est-\`a-dire $E_{G'}=E$. Parce que l'action par conjugaison de $\hat{\theta}_{3}$ dans $\hat{\Omega}$ permute cycliquement les $3$ \'el\'ements non triviaux de $\hat{\Omega}$, on voit que tout \'el\'ement de $\hat{\Omega}\hat{\theta}_{3}$ est conjugu\'e \`a $\hat{\theta}_{3}$ par un \'el\'ement de $\hat{\Omega}$. A \'equivalence pr\`es, on peut donc supposer $\omega_{G'}(\rho)=1$ et $\rho_{G'}=\rho_{G}=\hat{\theta}_{3}$. 
Dans les cas (A) ou (B), cela implique $\sigma_{G'}=\sigma_{G}$ pour tout $\sigma\in \Gamma_{F}$. Dans le cas (C), cela implique $\sigma_{G'}=\sigma_{G}$ pour tout $\sigma\in \Gamma_{E_{0}}$.  Dans ce cas (C), puisque $\tau^{-1}\rho\tau \in \Gamma_{E_{0}}$, on a $\tau_{G'}^{-1}\rho_{G'}\tau_{G'}=\tau_{G}^{-1}\rho_{G}\tau_{G}$, c'est-\`a-dire $\omega_{G'}(\tau)^{-1}\hat{\theta}_{3}\omega_{G'}(\tau)=\hat{\theta}_{3}$. D'o\`u $\omega_{G'}(\tau)=1$ et, de nouveau,  $\sigma_{G'}=\sigma_{G}$ pour tout $\sigma\in \Gamma_{F}$. 
    Il y a trois orbites possibles: ${\cal O}=\{\hat{\alpha}_{0}\}$, ${\cal O}=\{\hat{\alpha}_{2}\}$, ${\cal O}=\{\hat{\alpha}_{1},\hat{\alpha}_{3},\hat{\alpha}_{4}\}$. Dans le premier cas, on a ${\bf G}'={\bf G}$. Dans les deux autres cas, on voit que $G'_{SC}$ contient un facteur $SL(2)$ dans le cas ${\cal O}=\{\hat{\alpha}_{2}\}$, $SL(3)$ dans le cas ${\cal O}=\{\hat{\alpha}_{1},\hat{\alpha}_{3},\hat{\alpha}_{4}\}$. D'apr\`es \ref{An-1deppadique} (1), on a $FC^{st}(\mathfrak{g}'(F))=\{0\}$.  Cela d\'emontre que, si ${\bf G}'\not={\bf G}$, alors $FC^{st}(\mathfrak{g}'(F))^{Out({\bf G}')}=\{0\}$
  L'assertion (2) en r\'esulte. 
  
  Explicitons la cons\'equence de \ref{resultats}(4):
  
 $$ (3) \qquad dim(FC^{st}(\mathfrak{g}(F))=\left\lbrace\begin{array}{cc}1,&\rm{\,dans\,le\, cas\,}(A),\\ 2,&\rm{\,dans\,les\, cas\,}(B)\rm{\,et\,}(C).\\ \end{array}\right.$$
  
    \subsection{Type $D_{4}$ trialitaire, s\'eparation des \'el\'ements de $FC(\mathfrak{g}(F))$}\label{D4trialitaireseparation}
    On se place dans le cas (B) du paragraphe pr\'ec\'edent. On fixe des \'el\'ements non nuls $f_{0}\in FC_{0}$ et $f_{134}\in FC_{134}$. On va introduire un \'el\'ement  r\'egulier $X\in \mathfrak{g}_{ell}(F)$ et prouver
    
    (1) $S^G(X,f_{0})=0$, $S^G(X,f_{134})\not=0$. 
    
    Explicitons les constructions de \ref{alcoves} et \ref{profondeur} pour les deux sommets $s_{0}$ et $s_{134}$ associ\'es aux racines $\alpha_{0}, \alpha_{134}\in \Delta_{a}^{nr}$. On a attach\'e \`a chacun de ces sommets divers objets que l'on affecte d'un indice $0$ ou $134$. 
    
    On a $\alpha_{2}(s_{0})=\alpha_{134}(s_{0})=0$, $\Delta(U_{P_{0}})=\{\alpha_{2}\}$ et, en utilisant \ref{D4trialitairepadique}(1), on trouve que $r_{0}= \frac{1}{6}$, $\alpha_{2}(x_{0})=\frac{1}{6}$, $\alpha_{134}(x_{0})=0$.  On a $r_{0}(\alpha)=0$ pour $\alpha=\alpha_{2},\alpha_{2}+\alpha_{134}, \alpha_{2}+2\alpha_{134},\alpha_{2}+3\alpha_{134}, 2\alpha_{2}+3\alpha_{134}$, $r_{0}(\alpha)=\frac{1}{3}$ pour $\alpha=\pm \alpha_{134}, -\alpha_{2}-\alpha_{134},-\alpha_{2}-2\alpha_{134}$ et $r_{0}(\alpha)=1$ pour $\alpha=-\alpha_{2},-\alpha_{2}-3\alpha_{134},-2\alpha_{2}-3\alpha_{134}$.  Le r\'eseau  $\mathfrak{k}_{x_{0},r_{0},F^{nr}}$ est d\'ecrit par \ref{profondeur} (9) et ces valeurs de $r_{0}(\alpha)$. On voit que
       l'image de ce r\'eseau dans $\mathfrak{g}_{s_{0}}({\mathbb F}_{q})$ est 
  le sous-espace  des invariants par $\Gamma_{{\mathbb F}_{q}}$ dans 
  
  (2) $\oplus_{\alpha\in \Sigma^{nr}(U_{P_{0}})}\mathfrak{u}(\alpha)_{0}/\mathfrak{u}(\alpha)_{\frac{1}{e(\alpha)}}$.  
  
  \noindent On peut supposer que $f_{0}$ est la fonction $\tilde{f}$ de \ref{profondeur} associ\'ee \`a $s_{0}$. Le lemme de ce paragraphe entra\^{\i}ne que
    
    (3) le support de $f_{0}$ est contenu dans $\mathfrak{g}(F)_{\frac{1}{6}}$.
    
    On a $\alpha_{2}(s_{134})=0$, $\alpha_{134}(s_{134})=\frac{1}{6}$, $\Delta_{U_{P_{134}}}=\{\alpha_{0},\alpha_{2}\}$. On trouve que $r_{134}=\frac{1}{12}$, $\alpha_{2}(x_{134})=\alpha_{134}(x_{134})=\frac{1}{12}$.  On a $r_{134}(\alpha)=0$ pour toute racine positive $\alpha\in \Sigma^{nr}$, $r_{134}(\alpha)=1$ pour $\alpha=\alpha_{2},\alpha_{2}+3\alpha_{134},2\alpha_{2}+3\alpha_{134}$ et $r_{134}(\alpha)=\frac{1}{3}$ pour $\alpha=\alpha_{134},\alpha_{2}+\alpha_{134},\alpha_{2}+2\alpha_{134}$.  Le r\'eseau  $\mathfrak{k}_{x_{134},r_{134},F^{nr}}$ est d\'ecrit par \ref{profondeur} (9) et ces valeurs de $r_{134}(\alpha)$. On voit que
       l'image de ce r\'eseau dans $\mathfrak{g}_{s_{134}}({\mathbb F}_{q})$ est 
  le sous-espace  des invariants par $\Gamma_{{\mathbb F}_{q}}$ dans 
  
  (4) $\mathfrak{u}(\alpha_{2})_{0}/\mathfrak{u}(\alpha_{2})_{1}\oplus \mathfrak{u}(\alpha_{0})_{\frac{1}{3}}/\mathfrak{u}(\alpha_{0})_{\frac{2}{3}}$. 

    On peut supposer que $f_{134}$ est la fonction $\tilde{f}$ de \ref{profondeur} associ\'ee \`a $s_{134}$. Le lemme de ce paragraphe  entra\^{\i}ne que
   le support de $f_{134}$ est contenu dans $\mathfrak{g}(F)\cap \mathfrak{k}_{x_{134},r_{134},F^{nr}}$.      
    On a fix\'e une base de Chevalley de $\mathfrak{g}(\bar{F})$, qui contient un ensemble $(E_{\beta})_{\beta\in \Sigma}$. On a suppos\'e que $\theta_{3}(E_{\beta})=E_{\theta_{3}(\beta)}$ pour tout $\beta\in \Sigma$. Fixons une uniformisante $\varpi_{E}$ de $E$ telle que $\varpi_{E}^3\in F^{\times}$, notons $\zeta\in \boldsymbol{\zeta}_{3}(\bar{F})$ la racine telle que $\rho(\varpi_{E})=\zeta\varpi_{E}$. Posons
$$X=(\sum_{l=1,...,4}E_{l})+\varpi_{E}(\zeta^{-1}E_{-123}+ E_{-234}+\zeta E_{-124}),$$
o\`u on a not\'e $E_{l}=E_{\beta_{l}}$ et, par exemple, $E_{-123}=E_{-\beta_{1}-\beta_{2}-\beta_{3}}$. On voit que $X\in \mathfrak{g}(F)\cap  \mathfrak{k}_{x_{134},r_{134},F^{nr}}$ et que la r\'eduction de $X$ dans 
 $\mathfrak{g}_{s_{134}}({\mathbb F}_{q})$ appartient \`a $\tilde{\mathfrak{g}}_{s_{134},2}$. Donc $f_{134}(X)\not=0$. 
 
 Montrons que
 
 (5) $x_{134}$ est l'unique point $y\in Imm(G_{AD})$ tel que $X\in \mathfrak{k}_{y,\frac{1}{12}}$.
 
 Plongeons $\mathfrak{g}(\bar{F})$ dans $\mathfrak{gl}(8,\bar{F})$. On peut choisir le plongement de telle sorte que l'image de $X$ soit la matrice
 $$\left\lbrace\begin{array}{cccccccc}0&1&0&0&0&0&0&0\\0&0&1&0&0&0&0&0\\0&0&0&1&1&0&0&0\\ \varpi_{E}\zeta^{-1}&0&0&0&0&-1&0&0\\ \varpi_{E}\zeta&0&0&0&0&-1&0&0\\0&\varpi_{E}&0&0&0&0&-1&0\\0&0&-\varpi_{E}&0&0&0&0&-1\\0&0&0&-\varpi_{E}\zeta&-\varpi_{E}\zeta^{-1}&0&0&0\\ \end{array}\right\rbrace$$
 On calcule le polyn\^ome caract\'eristique de $X$. On obtient
$$T^8+6\varpi_{E} T^4-3\varpi_{E}^2.$$
L'alg\`ebre $\mathfrak{g}(E)$ est l'alg\`ebre d\'eploy\'ee du groupe $D_{4}$ et un \'el\'ement ayant un tel polyn\^ome caract\'eristique est elliptique r\'egulier, c'est-\`a-dire que $X$ est r\'egulier et elliptique dans $\mathfrak{g}(E)$, a fortiori dans $\mathfrak{g}(F)$. Fixons une extension galoisienne finie $K$ de $F$, contenant $E$ et telle que le tore $T_{X}$ commutant de $X$ soit d\'eploy\'e. Elle contient un \'el\'ement $u$ tel que $val_{F}(u)=1/12$. Le calcul ci-dessus du polyn\^ome caract\'eristique de $X$ entra\^{\i}ne que $u^{-1}X$ appartient \`a $\mathfrak{t}_{X}(\mathfrak{o}_{K})$ et que $val_{K}(\beta(u^{-1}X))=0$ pour toute racine $\beta$ de $T_{X}$. Soit $y\in Imm(G_{AD})$, supposons $X\in \mathfrak{k}_{y,1/12}$. Plongeons $Imm(G_{AD})$ dans l'immeuble $Imm_{K}(G_{AD},K)$.
  Alors $u^{-1}X\in \mathfrak{k}_{y,0,K}=\mathfrak{k}_{y,K}$. D'apr\`es le lemme \ref{lemmeimmobilier}, $y$ appartient \`a l'appartement $S_{K}(T_{X})$ de $Imm_{K}(G_{AD})$ associ\'e au tore $T_{X}$. Donc $y\in S_{K}(T_{X})\cap Imm(G_{AD})$. Cet ensemble est l'ensemble des points fixes par l'action galoisienne dans $S_{K}(T_{X})$ et est isomorphe \`a l'immeuble \'etendu de $T_{X}$ sur $F$. Puisque $T_{X}$ est elliptique, cet immeuble est r\'eduit \`a un point et $y$ est ce point. Cela d\'emontre (5). 
  
   Calculons $I^G(X,f_{134})$. On a 
$$I^G(X,f_{134})=\int_{G(F)}f_{134}(g^{-1}Xg)\,dg.$$
Pour $g$ tel que $f_{134}(g^{-1}Xg)\not=0$, $g^{-1}Xg$ appartient au support de $f_{134}$, donc \`a $ \mathfrak{k}_{x_{134},1/12}$.  Donc $X\in \mathfrak{k}_{gx_{134},1/12}$.  D'apr\`es (5),  cela entra\^{\i}ne $gx_{134}=x_{134}$, donc $g\in K_{x_{134}}^0$. Mais ce groupe conserve $f_{134}$. Donc $I^G(X,f_{134})$ est le produit d'une mesure et de $f_{134}(X)$, qui est non nul.  Donc $I^G(X,f_{134})\not=0$. Puisque $X$ est elliptique et que $f_{134}$ appartient \`a $I_{cusp}^{st}(\mathfrak{g}(F))$, $S^G(X,f_{134})$ est un multiple non nul de $I^G(X,f_{134})$ et est non nul. Cela d\'emontre la deuxi\`eme relation de (1). 

Enfin, puisque $val_{F}(\beta(X))=1/12
$ pour toute racine $\beta$ de $T_{X}$, on a

(6) $r(X)=\frac{1}{12}$.

\noindent D'apr\`es (3), aucun \'el\'ement stablement conjugu\'e \`a $X$ n'appartient au support de $f_{0}$. Donc 
  $S^G(X,f_{0})=0$. Cela d\'emontre (1).

  Pla\c{c}ons-nous maintenant dans le cas (C) du paragraphe pr\'ec\'edent. On va introduire un \'el\'ement r\'egulier $X\in \mathfrak{g}_{ell}(F)$ qui a les m\^emes propri\'et\'es que ci-dessus.  Notons $K$ le sous-corps de $E$ tel que $\Gamma_{K}$ soit engendr\'e par $\Gamma_{E}$ et $\tau$. L'extension $E/K$ est quadratique et la d\'efinition de $\tau $ implique que $E$ est la compos\'ee des extensions $K$ et $E_{0}$ de $F$. Il r\'esulte de \ref{extensionsdiedrales}(2) que $K$ contient une uniformisante $\varpi_{K}$ telle que $\varpi_{K}^3\in F^{\times}$. On fixe un tel \'el\'ement, on note $\zeta\in \boldsymbol{\zeta}_{3}(\bar{F})$ la racine telle que $\rho(\varpi_{K})=\zeta\varpi_{K}$. Posons
$$X=(\sum_{l=1,...,4}E_{l})+\varpi_{K}(\zeta^{-1}E_{-123}+ E_{-234}+\zeta E_{-124}).$$
Montrons que $X\in \mathfrak{g}(F)$. D'apr\`es la d\'efinition d'une base de Chevalley, on  a $E_{-123}=\epsilon[E_{-1},[E_{-2},E_{-3}]]$ pour un $\epsilon\in \{\pm 1\}$. Puisque $E_{-234}=\theta_{3}(E_{-123})$ et $E_{-124}=\theta_{3}^2(E_{-123})$, on a $E_{-234}=\epsilon[E_{-3},[E_{-2},E_{-4}]]$, $E_{-124}=\epsilon[E_{-4},[E_{-2},E_{-1}]$. L'\'el\'ement $\theta$ fixe $E_{-1}$, $E_{-2}$ et permute $E_{-3}$ et $E_{-4}$. On constate que $\theta$ fixe $E_{-234}$ et permute $E_{-123}$ et $E_{-124}$. Les conditions pour que $X$ appartienne \`a $\mathfrak{g}(F)$ sont donc

$\rho(\zeta^{-1}\varpi_{K})=\varpi_{K}$, $\rho(\varpi_{K})=\zeta\varpi_{K}$;

$\tau(\varpi_{K})=\varpi_{K}$, $\tau(\zeta^{-1}\varpi_{K})=\zeta\varpi_{K}$.

Elles r\'esultent de la d\'efinition de $K$ et $\zeta$, en tenant compte des relations $\rho(\zeta)=\zeta$ (car $\zeta\in E_{0}$ et $\rho$ fixe tout \'el\'ement de ce corps) et $\tau(\zeta)=\zeta^{-1}$ (car on a suppos\'e $\delta_{3}(q+1)=1$ donc $\zeta\not\in F^{\times}$ et $\tau$ ne fixe pas $\zeta$). 
  
  Ensuite, la m\^eme preuve que ci-desus s'applique: les r\'eseaux sont les m\^emes que pr\'ec\'edemment car l'action de $\Gamma_{F}^{nr}$ sur $\Sigma^{nr}$ est triviale, tous les \'el\'ements de cet ensemble \'etant fixes par $\theta$. 
  
  \subsection{Type $D_{4}$ trialitaire, action d'un automorphisme}\label{D4trialitaireautomorphisme}
  On se place dans le cas (B) du paragraphe \ref{D4trialitairepadique}. Le groupe $G$ a un automorphisme $\theta$. Cet automorphisme agit trivialement sur ${\cal A}^{nr}$ par d\'efinition de cet ensemble, donc fixe les sommets $s_{0}$ et $s_{134}$. Pour $\alpha\in \Sigma^{nr}$, on voit que $\theta$ agit trivialement sur $\mathfrak{u}(\alpha)_{0}/\mathfrak{u}(\alpha)_{\frac{1}{e(\alpha)}}$. Il r\'esulte alors de \ref{D4trialitaireseparation}(2) que $\theta$ fixe tout \'el\'ement de $FC_{0}$. Rappelons que l'on a suppos\'e $\delta_{3}(q-1)=1$ donc $\boldsymbol{\zeta}_{3}({\mathbb F}_{q})$ est d'ordre $3$. Pour $\alpha\in \Sigma^{nr}$ tel que $e(\alpha)=3$, on voit que $\theta$ agit sur $\mathfrak{u}(\alpha)_{\frac{1}{3}}/\mathfrak{u}(\alpha)_{\frac{2}{3}}$ par multiplication par une racine primitive de l'unit\'e d'ordre $3$ dans ${\mathbb F}_{q}^{\times}$ (on le voit par exemple en r\'eduisant l'\'el\'ement $X$ introduit dans le paragraphe pr\'ec\'edent). La relation \ref{D4trialitaireseparation}(4) nous dit comment agit $\theta$ sur le support de la fonction $f_{134}$ associ\'ee \`a $s_{134}$. On se rappelle que cette fonction est produit tensoriel de deux fonctions vivant sur les deux composantes $\mathfrak{sl}(2)({\mathbb F}_{q})$ de $G_{s_{134}}({\mathbb F}_{q})$. Mais 
  la multiplication par un carr\'e ne change pas l'orbite d'un \'el\'ement nilpotent de $\mathfrak{sl}(2)({\mathbb F}_{q})$ (et une racine de l'unit\'e d'ordre $3$ est un carr\'e). Il en r\'esulte que $\theta$ fixe tout \'el\'ement de $FC_{134}$. En r\'esum\'e, $\theta$ fixe tout \'el\'ement de $FC(\mathfrak{g}(F))$. 
      
 \subsection{Le type $E_{6}$ d\'eploy\'e}\label{E6deppadique}
   On suppose que $G$ est d\'eploy\'e de type $E_{6}$. On d\'efinit un homomorphisme $T_{ad}(F)\to F^{\times}/F^{\times,3}$ par $\prod_{l=1,...,6}\varpi_{l}(x_{l})\mapsto x_{1}x_{3}^2x_{5}x_{6}^2$. De cet homomorphisme se d\'eduit un isomorphisme $G_{AD}(F)/\pi(G(F))\simeq T_{ad}(F)/\pi(T(F))\simeq F^{\times}/F^{\times,3}$. L'image de $G_{AD}(F)_{0}$ est $\mathfrak{o}_{F}^{\times}/\mathfrak{o}_{F}^{\times,3}$.
   
   Si $\delta_{3}(q-1)=0$, on pose ${\cal X}=\emptyset$. 
   
   Supposons $\delta_{3}(q-1)=1$. On note $\Xi^{ram}$ l'ensemble des \'el\'ements de $\Xi$ dont la restriction \`a $G_{AD}(F)_{0}$ est non triviale. On  note ${\cal X}=\{4\}\cup\{(016,\xi);\xi\in \Xi^{ram}\}\cup\{(235,\xi);\xi\in \Xi^{ram}\}$.  On pose $d_{x}=2$ si $x=4$ et $d_{x}=1$ pour tout $x\in {\cal X}$, $x\not=4$.

   Le diagramme ${\cal D}$ a un automorphisme $\theta$ que l'on a d\'ecrit en \ref{E6dep}. Le diagramme ${\cal D}_{a}$ a un groupe d'automorphismes isomorphe \`a $\mathfrak{S}_{3}$. On note $\theta_{3}$ l'automorphisme d'ordre $3$ qui permute cycliquement $\alpha_{0},\alpha_{1},\alpha_{6}$ ainsi que $\alpha_{2},\alpha_{3},\alpha_{5}$ et fixe $\alpha_{4}$. L'action de $G_{AD}(F)$ sur le diagramme se fait par le groupe d'ordre $3$ engendr\'e par $\theta_{3}$. L'ensemble $\underline{S}(G)$ s'envoie surjectivement sur celui des orbites de cette action, qui s'identifie \`a  l'ensemble de repr\'esentants $ \{\alpha_{0},\alpha_{2},\alpha_{4}\}$. Les fibres ont trois \'el\'ements au-dessus de $\alpha_{0}$ et $\alpha_{2}$ et un seul au-dessus de $\alpha_{4}$. Pour un sommet $s$ param\'etr\'e par $\alpha_{0}$, le groupe $G_{s}$ est simplement connexe et d\'eploy\'e de type $E_{6}$. On a $FC(\mathfrak{g}_{s}({\mathbb F}_{q}))\not=\{0\}$ si et seulement si $\delta_{3}(q-1)=1$.  Supposons cette condition v\'erifi\'ee. D'apr\`es \ref{E6dep}, $FC(\mathfrak{g}_{s}({\mathbb F}_{q}))$ est de dimension $2$, engendr\'e par deux  fonctions $f_{N,\epsilon}$ pour lesquelles les restrictions de $\epsilon$ au centre de $G_{s}$ sont les deux caract\`eres non triviaux de ce centre. Le centre de $G$ s'envoie surjectivement sur celui de $G_{s}$. Il en r\'esulte que chacune des fonctions $f_{N,\epsilon}$ se transforme par un caract\`ere non trivial de $G_{AD}(F)_{0}$ et que les caract\`eres en question sont diff\'erents pour les deux fonctions. D'apr\`es \ref{actionsurFC}, chaque fonction donne naissance \`a trois \'el\'ements de $FC(\mathfrak{g}(F))$ qui se transforment selon les caract\`eres de $G_{AD}(F)/\pi(G(F))$ prolongeant le caract\`ere de $G_{AD}(F)_{0}$ attach\'e \`a la fonction en question. En r\'eunissant ces deux ensembles \`a trois \'el\'ements, on obtient pour tout $\xi\in \Xi^{nr}$ un \'el\'ement de $FC(\mathfrak{g}(F))$ se transformant selon $\xi$. On note  $FC_{016,\xi}$ la droite port\'ee par cet \'el\'ement. Pour un sommet $s$ param\'etr\'e par $\alpha_{2}$, on a $G_{s}=(SL(2)\times SL(6))/diag(\{\pm 1\})$, le groupe $\{\pm 1\}$ s'identifiant au centre de $SL(2)$ et au sous-groupe d'ordre $2$ de celui de $SL(6)$.  D'apr\`es \ref{An-1dep}, $FC(\mathfrak{g}_{s}({\mathbb F}_{q}))\not=\{0\}$ si et seulement si $\delta_{3}(q-1)=1$.  Supposons cette condition v\'erifi\'ee. Alors $FC(\mathfrak{g}_{s}({\mathbb F}_{q}))$ est de dimension $2$, engendr\'e par deux  fonctions $f_{N_{2},\epsilon_{2}}\otimes f_{N_{6},\epsilon_{6}}$, o\`u $f_{N_{2},\epsilon_{2}}$ est une unique fonction sur $SL(2)$ et $f_{N_{6},\epsilon_{6}}$ est une fonction sur $SL(6)$, $\epsilon_{6}$ \'etant l'un des deux caract\`eres d'ordre $6$ du centre de $SL(6)$. Le centre de $G$ s'identifie au sous-groupe d'ordre $3$ de celui de $SL(6)$. Le r\'esultat est alors le m\^eme que pour la racine $\alpha_{0}$: de nos fonctions  est issue une famille d'\'el\'ements   de $FC(\mathfrak{g}(F))$ param\'etr\'ee par $\Xi^{ram}$, de sorte que l'\'el\'ement param\'etr\'e par $\xi\in \Xi^{ram}$ se tranforme selon ce caract\`ere par le groupe 
  $G_{AD}(F)/\pi(G(F))$. On note $FC_{235,\xi}$ la droite port\'ee par cet \'el\'ement.   Pour un sommet $s$ param\'etr\'e par $\alpha_{4}$, on a $G_{s}= (SL(3)\times SL(3)\times SL(3))/diag(\zeta_{3}(\bar{{\mathbb F}}_{q}))$, le groupe $\zeta_{3}(\bar{{\mathbb F}}_{q})$ s'identifiant aux centres de chaque facteur. Pr\'ecis\'ement, indexons les trois facteurs par $0,1,6$, le facteur index\'e par $l$ contenant la coracine $\check{\alpha}_{l}$. 
Un \'el\'ement $\zeta\in \boldsymbol{\zeta}_{3}(\bar{{\mathbb F}}_{q})$ s'identifie aux \'el\'ements centraux   $\iota_{0}(\zeta)=\check{\alpha}_{0}(\zeta)\check{\alpha}_{2}(\zeta^2)$, $\iota_{1}(\zeta)=\check{\alpha}_{1}(\zeta)\check{\alpha}_{3}(\zeta^2)$, $\iota_{6}(\zeta)=\check{\alpha}_{6}(\zeta)\check{\alpha}_{5}(\zeta^2)$. Ces formules permettent d'identifier les centres de chacun des trois facteurs $SL(3)$. D'apr\`es \ref{An-1dep}, on a $FC(\mathfrak{g}_{s}({\mathbb F}_{q}))\not=\{0\}$ si et seulement si $\delta_{3}(q-1)=1$.  Supposons cette condition v\'erifi\'ee.  Alors, sur chaque facteur $SL(3)$, on a deux fonctions $f_{N_{l},\epsilon_{l}}$ o\`u $\epsilon_{l}$ est l'un des deux caract\`eres d'ordre $3$ du centre de $SL(3)$. Un produit tensoriel $f_{N_{0},\epsilon_{0}}\otimes f_{N_{1},\epsilon_{1}}\times f_{N_{6},\epsilon_{6}}$ se quotiente par  $diag(\zeta_{3}(\bar{{\mathbb F}}_{q}))$ si et seulement si $\epsilon_{0}\epsilon_{1}\epsilon_{6}=1$ (en identifiant les centres de chaque facteur), ce qui \'equivaut \`a $\epsilon_{0}=\epsilon_{1}=\epsilon_{6}$. Il reste deux fonctions de la forme pr\'ec\'edente. Le centre de $G$ s'envoie surjectivement sur le groupe $\{ \iota_{0}(\zeta_{0})\iota_{1}(\zeta_{1})\iota_{6}(\zeta_{6});\zeta_{0},\zeta_{1},\zeta_{6}\in \zeta(\bar{{\mathbb F}}_{q}), \zeta_{0}\zeta_{1}\zeta_{6}=1\}/diag(\zeta_{3}(\bar{{\mathbb F}}_{q}))$. Puisque  $\epsilon_{0}=\epsilon_{1}=\epsilon_{6}$, le caract\`ere $\epsilon_{0}\otimes \epsilon_{1}\otimes \epsilon_{6}$ est trivial sur ce groupe. Donc chaque fonction est  invariante par $G_{AD}(F)_{0}$. L'action du groupe $G_{AD}(F)$ tout entier est r\'ecup\'er\'e par la permutation cyclique des facteurs qui fixe \'evidemment chacune des fonctions. On obtient donc deux \'el\'ements de $FC(\mathfrak{g}(F))$ invariants par l'action de $G_{AD}(F)/\pi(G(F))$. On note $FC_{4}$ l'espace engendr\'e par ces deux fonctions. 
 On a obtenu  \ref{resultats}(1).
 
 Dans le cas o\`u $\delta_{3}(q-1)=0$, on pose ${\cal Y}=\emptyset$. Puisqu'on a d\'ej\`a prouv\'e que $FC(\mathfrak{g}(F))=\{0\}$, les assertions (2), (3) et (4) de \ref{resultats}(1) sont triviales.
 
 On suppose d\'esormais $\delta_{3}(q-1)=1$. On pose ${\cal Y}=\{0\}\cup\{(016,0,\xi);\xi\in \Xi^{ram}\}\cup\{(016,134,\xi);\xi\in \Xi^{ram}\}$.

    On consid\`ere un \'el\'ement $(\sigma\mapsto \sigma_{G'},{\cal O})$ de ${\cal E}_{ell}(G)$.  Le groupe $\hat{\Omega}$ est d'ordre $3$ et est engendr\'e par l'automorphisme $\hat{\theta}_{3}$ similaire au $\theta_{3}$ ci-dessus (le diagramme $\hat{{\cal D}}_{a}$ est le m\^eme que ${\cal D}_{a}$). Puisque $G$ est d\'eploy\'e, on a $\sigma_{G'}=\omega_{G'}(\sigma)$ pour tout $\sigma$ et $\sigma\mapsto \omega_{G'}(\sigma)$ est un homomorphisme injectif de $\Gamma_{E_{G'}/F}$ dans $\hat{\Omega}$.

       Supposons d'abord $E_{G'}=F$, donc l'action $\sigma\mapsto \sigma_{G'}$ est triviale. A conjugaison par $\hat{\Omega}$ pr\`es, l'orbite ${\cal O}$ peut \^etre \'egale \`a $\{\hat{\alpha}_{0}\}$, $\{\hat{\alpha}_{2}\}$ ou $\{\hat{\alpha}_{4}\}$. Dans le premier cas ${\bf G}'={\bf G}$ et, \`a ce point, on ne peut rien dire de $FC^{st}(\mathfrak{g}(F))$. Dans les deux autres cas, on voit que $G'_{SC}$ contient un facteur $SL(2)$ ou $SL(3)$ donc $FC^{st}(\mathfrak{g}'(F))=\{0\}$ d'apr\`es \ref{An-1deppadique}(1). 
   
   Supposons maintenant que $E_{G'}$ soit l'extension non ramifi\'ee de degr\'e $3$ de $F$. Le lemme \ref{centre} exclut les orbites ${\cal O}$ form\'ees de trois \'el\'ements. Il reste l'orbite ${\cal O}=\{\hat{\alpha}_{4}\}$ pour laquelle on peut remplacer $G'$ par $G'_{SC}$. On a $G'_{SC}\simeq Res_{E_{G'}/F}(SL(3))$ donc $FC^{st}(\mathfrak{g}'(F))=\{0\}$ d'apr\`es \ref{An-1deppadique}(1). 
   
   Supposons enfin que $E_{G'}$ soit une extension cyclique d'ordre $3$ et ramifi\'ee de $F$.  Une telle extension existe puisque  $\delta_{3}(q-1)=1$.   Il y a  $3$ extensions possibles. Fixons-en une et choisissons un g\'en\'erateur $\tau$ de $\Gamma_{E_{G'}/F}$. Il y a deux actions $\sigma\mapsto \sigma_{G'}$ possibles, l'une telle que $\tau_{G'}=\hat{\theta}_{3}$, l'autre telle que $\tau_{G'}=\hat{\theta}_{3}^2$. Si ${\cal O}=\{\hat{\alpha}_{4}\}$, on a comme ci-dessus  $G'_{SC}\simeq Res_{E_{G'}/F}(SL(3))$ donc $FC^{st}(\mathfrak{g}'(F))=\{0\}$. Si ${\cal O}=\{\hat{\alpha}_{2},\hat{\alpha}_{3},\hat{\alpha}_{5}\}$, on a $G'_{SC}\simeq Res_{E_{G'}/F}(SL(2))\times SL(2)$ donc  $FC^{st}(\mathfrak{g}'(F))=\{0\}$. Supposons ${\cal O}=\{\hat{\alpha}_{0},\hat{\alpha}_{1},\hat{\alpha}_{6}\}$. D'apr\`es le lemme \ref{centre}, on peut remplacer $G'$ par $G'_{SC}\simeq Spin_{E_{G'}/F}(8)$, o\`u on d\'esigne ainsi la forme trialitaire de $Spin(8)$ associ\'ee \`a l'extension $E_{G'}/F$. Alors $FC^{st}(\mathfrak{g}'(F))$ est de dimension $2$ d'apr\`es \ref{D4trialitairepadique}. Plus pr\'ecis\'ement, $FC^{st}(\mathfrak{g}'(F))$ est somme de deux droites not\'ees $FC_{0}$ et $FC_{134}$ dans ce paragraphe, notons-les $FC_{0}(\mathfrak{g}'(F))$ et $FC_{134}(\mathfrak{g}'(F))$. 
  Le groupe $Out({\bf G}')$ est $\hat{\Omega}$ tout entier. L'action de $\hat{\theta}_{3}$ sur $G'_{SC}$ est l'automorphisme $\theta_{3}$ de ce groupe. Il  agit trivialement sur $FC^{st}(\mathfrak{g}'(F))$ d'apr\`es \ref{D4trialitaireautomorphisme}. Il nous reste \`a d\'eterminer $\xi_{{\bf G}'}$.  Rappelons (cf. \ref{description}) que l'\'el\'ement $s$ de la donn\'ee ${\bf G}'$ d\'epend du choix de racines de l'unit\'e dans ${\mathbb C}^{\times}$, ici d'une racine primitive d'ordre $3$ que l'on note $j$.   L'\'el\'ement $s$ de la donn\'ee ${\bf G}'$ v\'erifie $\hat{\alpha}_{1}(s)=\hat{\alpha}_{6}(s)=j$ et $\hat{\alpha}_{i}(s)=1$ pour $i=2,...,5$.   On peut choisir $s_{sc}=\check{\hat{\alpha}}_{1}(j^2)\check{\hat{\alpha}}_{2}(j^2)\check{\hat{\alpha}}_{4}(j)\check{\hat{\alpha}}_{6}(j^2)$. Alors $\hat{\theta}_{3}(s_{sc})s_{sc}^{-1}=z$ o\`u $z=\check{\hat{\alpha}}_{1}(j^2)\check{\hat{\alpha}}_{3}(j)\check{\hat{\alpha}}_{5}(j^2)\check{\hat{\alpha}}_{6}(j)\in Z(\hat{G}_{SC})$. L'action galoisienne sur $Z(\hat{G}_{SC})$ \'etant triviale, les cocycles s'identifient \`a des caract\`eres de $\Gamma_{F}$ dans ce groupe d'ordre $3$. Pour $E_{G'}$ fix\'ee, les deux actions $\sigma\mapsto \sigma_{G'}$ possibles d\'eterminent les deux caract\`eres d'ordre $3$ triviaux sur $\Gamma_{E_{G'}}$. Quand on fait varier $E_{G'}$ on obtient les $6$ caract\`eres d'ordre $3$ de $\Gamma_{F}$ qui sont ramifi\'es. Ils correspondent aux $6$ \'elements de $\Xi^{ram}$. Donc, pour tout $\xi\in \Xi^{ram}$, il y a une et une seule de nos donn\'ee ${\bf G}'$ telle que $\xi_{{\bf G}'}=\xi$. On la note ${\bf G}'_{016,\xi}$ cette donn\'ee et on pose $FC^{{\cal E}}_{016,0,\xi}=FC_{0}(\mathfrak{g}'_{016,\xi}(F))$, $FC^{{\cal E}}_{016,134,\xi}=FC_{134}(\mathfrak{g}'_{016,\xi}(F))$. 
  
  On pose ${\mathbb X}=\{4,016,235\}$, ${\mathbb Y}=\{0,(016,0),(016,134)\}$ et on d\'efinit une bijection $\phi:{\mathbb X}\to {\mathbb Y}$ par $\phi(4)=0$, $\phi(016)=(016,0)$, $\phi(235)=(016,134)$. Il y a des surjections \'evidentes ${\cal X}\to {\mathbb X}$ et ${\cal Y}\to {\mathbb Y}$ et la bijection $\phi$ se rel\`eve naturellement (en ajoutant des $\xi\in \Xi^{ram}$ dans la d\'efinition) en une bijection $\varphi:{\cal X}\to {\cal Y}$. On pose ${\cal X}^{st}=\{4\}$. 
  
  On a associ\'e ci-dessus une droite $FC^{{\cal E}}_{y}$ \`a tout \'el\'ement $y\in {\cal Y}$ diff\'erent de $0$. On n'a pas trait\'e la donn\'ee principale ${\bf G}$. A l'aide de la bijection $\varphi$, l'argument habituel de dimension montre que $FC^{st}(\mathfrak{g}(F))$ est de dimension $2$ et on pose $FC^{{\cal E}}_{0}=FC^{st}(\mathfrak{g}(F))$. On a ainsi compl\'et\'e la description \ref{resultats}(3). 
  
  L'action de $G_{AD}(F)/\pi(G(F))$ suffit \`a d\'emontrer \ref{resultats}(4). En effet, on a forc\'ement
  $$FC^{st}(\mathfrak{g}(F))\subset FC(\mathfrak{g}(F))\cap I_{cusp,{\bf 1}}(\mathfrak{g}(F))=FC_{4}(\mathfrak{g}(F)).$$
  Les espaces extr\^emes \'etant tous deux de dimension $2$, ils sont \'egaux.

Posons ${\mathbb X}^{\star}=\{016,235\}$, ${\mathbb Y}^{\star}=\{(016,0),(016,134)\}$. On munit ${\mathbb X}^{\star}$ de l'ordre $016<235$. Par la surjection ${\cal X}\to {\mathbb X}$, il se rel\`eve en un pr\'eordre sur  l'image r\'eciproque ${\cal X}^{\star}$ de ${\mathbb X}^{\star}$ et on applique les constructions de \ref{ingredients}. L'ensemble $\underline{{\cal Y}}^{\star}$ s'identifie \`a ${\mathbb Y}^{\star}$ dont le plus petit \'el\'ement $(y_{min})$ est $(016,0)$. On pose $\underline{{\cal Y}}^{\sharp}=\underline{{\cal Y}}^{\star}-\{(y_{min})\}$. Il a un unique \'el\'ement que l'on  note $(y)$ (on peut l'identifier \`a l'\'el\'ement $(016,134)$ de ${\mathbb Y}^{\star}$).   On a $d_{(y)}=\vert \Xi^{ram}\vert =6$. On associe \`a $(y)$ la collection $({\bf G}'_{016,\xi})_{\xi\in \Xi^{ram}}$.   Pour chaque groupe $G'_{016,\xi}$, on a construit en \ref{D4trialitaireseparation} un \'el\'ement $X\in \mathfrak{g}'_{016,\xi,reg}(F)$. Un calcul fastidieux permet de prouver que $X$ est $G$-r\'egulier. Nous n'avons pas besoin de ce calcul: par l'argument de \ref{ingredients}, on peut remplacer $X$ par un \'el\'ement assez voisin de sorte  les propri\'et\'es (1) et (6) de \ref{D4trialitaireseparation} soient encore v\'erifi\'ees et que $X$ soit $G$-r\'egulier. 
  On fixe un tel \'el\'ement que l'on note
  $Y_{\xi}$. On associe \`a $(y)$ la collection $(Y_{\xi})_{\xi\in \Xi^{ram}}$. 

Montrons que les conditions de \ref{ingredients} sont satisfaites. La relation (1) de ce paragraphe est v\'erifi\'ee par simple compatibilit\'e du transfert avec les actions de $G_{AD}(F)/\pi(G(F))$. La condition (2) est satisfaite d'apr\`es notre description des espaces en question. Les conditions (3) et (4) r\'esultent de la condition (1) de \ref{D4trialitaireseparation} et du fait que, pour un \'el\'ement $\xi\in \Xi^{ram}$, un \'el\'ement $y'\in {\cal Y}^{\star}$ et une fonction $f'_{y'}\in FC^{{\cal E}}_{y'}$, les d\'efinitions entra\^{\i}nent que $S^{G'_{016,\xi}}(Y_{\xi},f'_{y'})$ ne peut \^etre non nul que si  l'\'el\'ement $\xi'$ qui figure dans $y'$ est \'egal \`a $\xi$. Il reste \`a prouver la condition (5). Celle-ci se r\'ecrit sous la forme suivante:

 (1) soient $x=(016,\xi')\in {\cal X}$, $f\in FC_{x}$, $\xi\in \Xi^{ram}$ et $X$ un \'el\'ement de $\mathfrak{g}_{reg}(F)$ dont la classe de conjugaison stable correspond \`a celle de $Y_{\xi}$; alors $I^G(X,f)=0$. 
   
 D'apr\`es   \ref{D4trialitaireseparation} (6), on a $r(X)=\frac{1}{12}$. Il suffit donc de prouver que le support de $f$ est contenu dans $\mathfrak{g}(F)_{r}$ pour un $r>r(X)$. La fonction $f$ est par d\'efinition combinaison lin\'eaire de trois fonctions issues des sommets associ\'es aux trois racines $\alpha_{0}$, $\alpha_{1}$, $\alpha_{6}$. Ces fonctions sont conjugu\'ees par l'action de $G_{AD}(F)/\pi(G(F))$ et leurs supports sont aussi conjugu\'es. Il suffit de consid\'erer l'une d'elles, par exemple celle associ\'ee au sommet $s_{0}$ attach\'e  \`a la racine $\alpha_{0}$.  C'est une fonction $\tilde{f}$ du type consid\'er\'e en \ref{profondeur}. Le syst\`eme de racines $\Sigma^{nr}$ de  $G_{s_{0}}$ est $\Sigma$ tout entier. D'apr\`es \ref{E6dep}, on a $\Delta(U_{P_{0}})=\{\alpha_{1},\alpha_{4},\alpha_{6}\}$.  On a $d(\alpha_{0})=d(\alpha_{1})=d(\alpha_{6})=1$ et $d(\alpha_{4})=3$. Alors le nombre $r$ de \ref{profondeur} est \'egal \`a $\frac{1}{6}$, donc $\tilde{f}$ est \`a support dans $\mathfrak{g}(F)_{\frac{1}{6}}$. Cela d\'emontre (1). 
 
 Alors le lemme de \ref{ingredients} nous dit que $transfert(FC_{(x)})=FC^{{\cal E}}_{\varphi((x))}$ pour tout $(x)\in \underline{{\cal X}}^{\star}$. Comme d'habitude, l'action de $G_{AD}(F)/\pi(G(F))$ permet de raffiner cette \'egalit\'e en $transfert(FC_{x})=FC^{{\cal E}}_{\varphi(x)}$ pour tout $x\in {\cal X}^{\star}$. Cela prouve \ref{resultats}(3). 
 
 Explicitons la cons\'equence de \ref{resultats}(4):
 
 (2) $dim(FC^{st}(\mathfrak{g}(F)))=2\delta_{3}(q-1)$. 
   
 \subsection{Forme int\'erieure du type $E_{6}$ d\'eploy\'e}\label{E6nonquasidep}
On suppose que $G^*$ est du type pr\'ec\'edent et que $G$ en est une forme int\'erieure non d\'eploy\'ee. Dans les tables de Tits, le groupe est de type $^3E_{6}$. 

Si $\delta_{3}(q-1)=0$, on pose ${\cal X}= \emptyset$. Si $\delta_{3}(q-1)=1$, on pose ${\cal X} =\{ 4\}$ et $d_{4}=2$.
  
Le diagramme ${\cal D}^{nr}_{a}$ est le diagramme de Dynkin affine d'un syst\`eme de racines de type  $E_{6}$ sur lequel le Frobenius $Fr\in \Gamma_{{\mathbb F}_{q}}$ agit par $\theta_{3}$. L'ensemble $\underline{S}(G)$ est param\'etr\'e par l'ensemble des orbites de  cette action. Pour un sommet param\'etr\'e par une orbite \`a $3$ \'el\'ements, le lemme \ref{orbites} implique que $FC(\mathfrak{g}_{s}({\mathbb F}_{q}))=\{0\}$. Il reste le sommet $s$ param\'etr\'e par l'orbite r\'eduite \`a  $\alpha_{4}$. Sur $\bar{{\mathbb F}}_{q}$, on a $G_{s}\simeq (SL(3)\times SL(3)\times SL(3))/diag(\zeta_{3}(\bar{{\mathbb F}}_{q}))$. On a d\'ecrit pr\'ecis\'ement ce groupe dans le paragraphe pr\'ec\'edent en indexant les trois facteurs $SL(3)$ par $0,1,6$. En notant $ Fr_{SL(3)}$ l'action  usuelle du Frobenius sur $SL(3)$, celle sur $G_{s}$ est $(g_{0},g_{1},g_{6})\mapsto (Fr_{SL(3)}(g_{6}),Fr_{SL(3)}(g_{0}),Fr_{SL(3)}(g_{1}))$. Il y a deux faisceaux-caract\`eres cuspidaux sur chaque facteur $SL(3)$, associ\'es \`a des couples $(N,\epsilon)$ o\`u $\epsilon$ est l'un des deux caract\`eres de $Z(SL(3))$ d'ordre $3$. Consid\'erons un produit de trois tels faisceaux sur chacune des composantes, d'o\`u trois caract\`eres $\epsilon_{0},\epsilon_{1},\epsilon_{6}$. Pour que le produit soit invariant par $diag(\zeta_{3}(\bar{{\mathbb F}}_{q}))$, on doit avoir $\epsilon_{0}=\epsilon_{1}=\epsilon_{6}$. Pour que le produit soit invariant par l'action galoisienne, la formule ci-dessus entra\^ine qu'il faut que ce caract\`ere commun $\epsilon_{0}=\epsilon_{1}=\epsilon_{6}$ soit invariant par l'action galoisienne naturelle, c'est-\`a-dire que $\delta_{3}(q-1)=1$. Si cette condition est v\'erifi\'ee, on r\'ecup\`ere deux faisceaux-caract\`eres, d'o\`u deux \'el\'ements de $FC(\mathfrak{g}_{s}({\mathbb F}_{q}))$. Il s'en d\'eduit deux \'el\'ements de $FC(\mathfrak{g}(F))$ et on note $FC_{4}(\mathfrak{g}(F))$ le plan qu'ils engendrent. On  obtient ainsi l'assertion \ref{resultats}(1). 

Si $\delta_{3}(q-1)=0$, on pose ${\cal Y}=\emptyset$. Puisqu'on vient de prouver que $FC(\mathfrak{g}(F))=\{0\}$, les assertions (2) et (3) de \ref{resultats} sont tautologiques. Supposons $\delta_{3}(q-1)=1$. On pose ${\cal Y}=\{0\}$. On a vu dans le paragraphe pr\'ec\'edent que $FC^{st}(\mathfrak{g}^*(F))$ etait de dimension $2$. On note cet espace $FC^{{\cal E}}_{0}$. L'argument de comparaison des dimensions montre que  $FC^{st}(\mathfrak{g}'(F))^{Out({\bf G}')}=\{0\}$ pour toute donn\'ee endoscopique ${\bf G}'\not={\bf G}$. Les assertions (2) et (3) de \ref{resultats}  s'en d\'eduisent (o\`u $\varphi$ est l'unique bijection de ${\cal X}$ sur ${\cal Y}$).

\subsection{Type $E_{6}$ quasi-d\'eploy\'e, $E/F$ non ramifi\'ee}
On note $E_{0}$ l'extension quadratique non ramifi\'ee de $F$. On fixe $\tau\in \Gamma_{F}-\Gamma_{E_{0}}$.  On suppose que $G$ est quasi-d\'eploy\'e de type $E_{6}$ et que $\Gamma_{F}$ agit sur ${\cal D}$ de la fa\c{c}on suivante: l'action de $\Gamma_{E_{0}}$ est triviale et  $\tau$ agit par l'automorphisme $\theta$. Dans les tables de Tits, le groupe est de type $^2E_{6}$. Un \'el\'ement de $T_{ad}(F)$ s'\'ecrit $\prod_{l=1,...,6}\check{\varpi}_{l}(x_{l})$ avec $x_{2},x_{4}\in F^{\times}$, $x_{1},x_{3}\in E_{0}^{\times}$ et $x_{6}=\tau(x_{1})$, $x_{5}=\tau(x_{3})$. Notons $E_{0}^1$ le groupe des \'el\'ements de $E_{0}^{\times}$ de norme $1$ et $(E_{0}^1)^3=\{z^3; z\in E_{0}^1\}$. On d\'efinit un homomorphisme $T_{ad}(F)\to E_{0}^1/(E_{0}^1)^3$ par $\prod_{l=1,...,6}\check{\varpi}_{l}(x_{l})\mapsto \frac{x_{1}x_{3}^2}{\tau(x_{1}x_{3}^2)}$. De cet homomorphisme est issu un isomorphisme $G_{AD}(F)/\pi(G(F))\simeq T_{ad}(F)/\pi(T(F))\simeq E_{0}^1/(E_{0}^1)^3$. On a $G_{AD}(F)_{0}=G_{AD}(F)$. Remarquons que $E_{0}^1/(E_{0}^1)^3$ est $\{1\}$ si $\delta_{3}(q-1)=1$ et est d'ordre $3$ si $\delta_{3}(q-1)=0$. Dans ce dernier cas, on note $\Xi_{\not={\bf 1}}$ l'ensemble des deux \'el\'ements non triviaux de $\Xi$.

Si $\delta_{3}(q-1)=1$, on pose   $ {\cal X}=\{4\}$ et $d_{4}=2$. Si $\delta_{3}(q-1)=0$, on pose ${\cal X}=\{(0,\xi);\xi\in \Xi_{\not={\bf 1}}\}\cup \{(2,\xi);\xi\in \Xi_{\not={\bf 1}}\}$ et $d_{x}=1$ pour tout $x\in {\cal X}$.

  L'ensemble $\underline{S}(G)$ est en bijection avec les orbites de l'action galoisienne sur ${\cal D}_{a}$. Pour un sommet $s$ param\'etr\'e par une orbite \`a deux \'el\'ements, le lemme  \ref{orbites} dit que $FC(\mathfrak{g}_{s}({\mathbb F}_{q}))=\{0\}$. Pour le sommet $s$ param\'etr\'e par $\alpha_{0}$, le groupe $G_{s}$ est simplement connexe de type $E_{6}$, avec action du Frobenius par $\theta$. D'apr\`es \ref{E6nondep}, $FC(\mathfrak{g}_{s}({\mathbb F}_{q}))$ est nul si $\delta_{3}(q-1)=1$, de dimension $2$ si $\delta_{3}(q-1)=0$.  Dans ce dernier cas,  on voit comme en \ref{E6deppadique} que les deux g\'en\'erateurs  naturels de cet espace se transforment par le groupe $G_{AD}(F)/\pi(G(F))$ selon les deux \'el\'ements non triviaux de $\Xi$. Pour chaque $\xi\in \Xi_{\not={\bf 1}}$, il  s'en d\'eduit un \'el\'ement de $FC(\mathfrak{g}(F))$ et on note $FC_{0,\xi}$ la droite qu'il engendre.  Pour le sommet $s$ param\'etr\'e par $\alpha_{2}$, on a $G_{s}=(SL(2)\times SU_{{\mathbb F}_{q^2}/{\mathbb F}_{q}}(6))/\{\pm 1\}$. D'apr\`es \ref{An-1nondep}, $FC(\mathfrak{g}_{s}({\mathbb F}_{q}))$ est nul si $\delta_{3}(q-1)=1$, de dimension $2$ si $\delta_{3}(q-1)=0$.  Dans ce dernier cas, on voit comme en \ref{E6deppadique} que les deux g\'en\'erateurs naturels de l'espace se transforment par le groupe $G_{AD}(F)/\pi(G(F))$ selon les deux \'el\'ements non triviaux de $\Xi$. On en d\'eduit comme ci-dessus des droites dans $FC(\mathfrak{g}(F))$ que l'on note $FC_{2,\xi}$ pour $\xi\in \Xi_{\not={\bf 1}}$. 
   Pour le sommet $s$ param\'etr\'e par $\alpha_{4}$, on a $G_{s}=(SL(3)\times Res_{{\mathbb F}_{q^2}/{\mathbb F}_{q}}(SL(3)))/\boldsymbol{\zeta}_{3}(\bar{{\mathbb F}}_{q})$, le plongement de $\boldsymbol{\zeta}_{3}(\bar{{\mathbb F}}_{q})$ \'etant le m\^eme qu'en  \ref{E6deppadique}. On voit comme dans ce paragraphe que $FC(\mathfrak{g}_{s}({\mathbb F}_{q}))$ est nul si $\delta_{3}(q-1)=0$, de dimension $2$ si $\delta_{3}(q-1)=1$, et que, dans ce dernier cas, les deux g\'en\'erateurs naturels  sont invariants par $G_{AD}(F)/\pi(G(F))$. On en d\'eduit deux \'el\'ements de $FC(\mathfrak{g}(F))$ et on note $FC_{4}$ le plan qu'ils engendrent.   Cela d\'emontre \ref{resultats}(1). 
 
 Si $\delta_{3}(q-1)=1$, on pose   $ {\cal Y}=\{0\}$. Si $\delta_{3}(q-1)=0$, on pose ${\cal Y}=\{(016,0,\xi;\xi\in \Xi_{\not={\bf 1}}\}\cup \{(016,134,\xi;\xi\in \Xi_{\not={\bf 1}}\}$. 
 
 Consid\'erons un couple $(\sigma\mapsto \sigma_{G'},{\cal O})\in {\cal E}_{ell}(G)$.   
 
Supposons d'abord $E_{G'}=E_{0}$. On a $\tau_{G'}\in \hat{\Omega}\tau_{G}=\hat{\Omega}\theta$. Or tout \'el\'ement de $\hat{\Omega}\theta$ est conjugu\'e \`a $\theta$ par un \'el\'ement de $\hat{\Omega}$ donc, \`a \'equivalence pr\`es, on peut supposer $\tau_{G'}=\theta$. Le cas d'une orbite ${\cal O}$ \`a deux \'el\'ements est exclu par le lemme \ref{centre}. Si ${\cal O}=\{\hat{\alpha}_{4}\}$, resp. ${\cal O}=\{\hat{\alpha}_{2}\}$, le groupe $G'_{SC}$ contient un facteur $SL(3)$, resp. $SL(2)$, donc $FC^{st}(\mathfrak{g}'(F))=\{0\}$. Si ${\cal O}=\{\hat{\alpha}_{0}\}$, on a ${\bf G}'={\bf G}$ et, \`a ce point, on ne peut rien dire de $FC^{st}(\mathfrak{g}(F))$. 

Supposons maintenant que $E_{G'}$ soit une extension de $E_{0}$ de degr\'e $3$. La structure de $\Gamma_{E_{G'}/F}$ est alors celle \'etudi\'ee en \ref{extensionsdiedrales}(3): une telle extension $E_{G'}$ existe si et seulement si $\delta_{3}(q-1)=0$ et, si cette \'egalit\'e est v\'erifi\'ee, cette extension  est unique. Supposons $\delta_{3}(q-1)=0$ et consid\'erons cette extension $E_{G'}/F$. Elle est ramifi\'ee.  On a d\'ej\`a fix\'e $\tau$ (on a forc\'ement $\tau^2\in \Gamma_{E_{G'}}$) et on fixe un g\'en\'erateur $\rho$ de $\Gamma_{E_{G'}/E_{0}}$. Comme ci-cessus, \`a \'equivalence pr\`es, on peut supposer $\tau_{G'}=\theta$. Il y a deux possibilit\'es pour $\rho_{G'}$: $\rho_{G'}=\theta_{3}$, $\rho_{G'}=\theta_{3}^2$, qui ne sont pas \'equivalentes. Supposons $\rho_{G'}=\theta_{3}$. Notons $K$ la sous-extension de $E_{G'}$ telle que $\Gamma_{E_{G'}/K}=\{1,\tau\}$. L'extension $E_{G'}/K$ est non ramifi\'ee et l'extension $K/F$ est ramifi\'ee et non galoisienne. Remarquons que $\Gamma_{E_{G'}/K}$ est le fixateur de $\hat{\alpha}_{0}$ et $\hat{\alpha}_{2}$ pour l'action galoisienne $\sigma\mapsto \sigma_{G'}$. Si ${\cal O}=\{\hat{\alpha}_{4}\}$, $G'_{SC}\simeq Res_{K/F}(SL(3))$ et $FC^{st}(\mathfrak{g}'(F))=\{0\}$. Si ${\cal O}=\{\hat{\alpha}_{2},\hat{\alpha}_{3},\hat{\alpha}_{5}\}$, $G'_{SC}\simeq SL(2)\times Res_{K/F}(SL(2))$ et encore $FC^{st}(\mathfrak{g}'(F))=\{0\}$. Supposons ${\cal O}=\{\hat{\alpha}_{0},\hat{\alpha}_{1},\hat{\alpha}_{6}\}$. D'apr\`es le lemme \ref{centre}, on peut remplacer $G'$ par $G'_{SC} =Spin_{E_{G'}/F}(8)$, c'est-\`a-dire le groupe de type (C) \'etudi\'e en \ref{D4trialitairepadique}.  On a vu dans ce paragraphe que l'espace   $FC^{st}(\mathfrak{g}'(F))$ est de dimension $2$ et est  somme de deux droites not\'ees alors $FC_{0}$ et $FC_{134}$ et que nous notons ici $FC_{0}(\mathfrak{g}'(F))$ et $FC_{134}(\mathfrak{g}'(F))$.    Le groupe $Out({\bf G}')$ est contenu dans celui des $\omega\in \hat{\Omega}$ qui commutent \`a $\sigma_{G'}$ pour tout $\sigma$, donc est r\'eduit \`a $1$. Calculons $\xi_{{\bf G}'}$. L'\'el\'ement $s$ figurant dans  la donn\'ee  v\'erifie $\hat{\alpha}_{1}(s)=\hat{\alpha}_{6}(s)=j$ et $\hat{\alpha}_{i}(s)=1$ pour $i=2,...,5$,  o\`u $j$ est une racine cubique de $1$ dans ${\mathbb C}^{\times}$.   On peut choisir $s_{sc}=\check{\hat{\alpha}}_{1}(j^2)\check{\hat{\alpha}}_{2}(j^2)\check{\hat{\alpha}}_{4}(j)\check{\hat{\alpha}}_{6}(j^2)$. On a $\theta(s_{sc})=s_{sc}$ et
 $\hat{\theta}_{3}(s_{sc})s_{sc}^{-1}=z$ o\`u $z=\check{\hat{\alpha}}_{1}(j^2)\check{\hat{\alpha}}_{3}(j)\check{\hat{\alpha}}_{5}(j^2)\check{\hat{\alpha}}_{6}(j)\in Z(\hat{G}_{SC})$.
  D'o\`u $\rho_{G'}(s_{sc})s_{sc}^{-1}=z$. Cela signifie que $\xi_{{\bf G}'}$ est un caract\`ere d'ordre $3$ de $G_{AD}(F)/\pi(G(F))$.  
  Supposons maintenant $\rho_{G'}=\theta_{3}^2$. Les r\'esultats sont les m\^emes, la seule diff\'erence \'etant le calcul de $\xi_{{\bf G}'}$: on a cette fois $\delta(s_{sc})s_{sc}^{-1}=z^2$. Le caract\`ere $\xi_{{\bf G}'}$ est encore un \'el\'ement de $\Xi$ d'ordre $3$ mais il est diff\'erent du caract\`ere pr\'ec\'edent. Pour les deux \'el\'ements  de $\Xi_{\not={\bf 1}}$, on a donc une unique donn\'ee ${\bf G}'$ comme ci-dessus telle que $\xi_{{\bf G}'}=\xi$. On la note ${\bf G}'_{016,\xi}$ et on pose $FC^{{\cal E}}_{016,0,\xi}=FC_{0}(\mathfrak{g}_{016,\xi}'(F))$, $FC^{{\cal E}}_{016,134,\xi}=FC_{134}(\mathfrak{g}_{016,\xi}'(F))$. 
  
  Supposons $\delta_{3}(q-1)=1$. On pose ${\cal X}^{st}={\cal X}$. On n'a construit ci-dessus aucune donn\'ee endoscopique elliptique ${\bf G}'$ telle que $FC^{st}(\mathfrak{g}'(F))^{Out({\bf G}')}\not=\{0\}$ mais on n'a pas trait\'e la donn\'ee principale ${\bf G}$. L'argument de comparaison des dimensions montre que $FC(\mathfrak{g}(F))=FC^{st}(\mathfrak{g}(F))$. On pose $FC^{{\cal E}}_{0}=FC^{st}(\mathfrak{g}(F))$ et les assertions (2), (3) et (4) de \ref{resultats} sont imm\'ediates ($\varphi$ \'etant l'unique bijection de ${\cal X}$ sur ${\cal Y}$). 
  
  Supposons $\delta_{3}(q-1)=0$.  On pose ${\mathbb X}=\{016,235\}$, ${\mathbb Y}=\{(016,0),(016,134)\}$ 
  et on d\'efinit une bijection $\phi:{\mathbb X}\to {\mathbb Y}$ par $\phi(016)=(016,0)$, $\phi(235)=(016,134)$. Elle se rel\`eve en une bijection $\varphi:{\cal X}\to {\cal Y}$. On pose ${\cal X}^{st}=\emptyset$. On a  associ\'e ci-dessus une droite $FC^{{\cal E}}_{y}$ \`a tout \'el\'ement de ${\cal Y}$. L'argument de comparaison des dimensions entra\^{\i}ne que $FC^{st}(\mathfrak{g}(F))=\{0\}$. On a ainsi achev\'e la preuve de \ref{resultats}(2) et prouv\'e \ref{resultats}(4). Remarquons que $Imm_{F^{nr}}(G_{AD})$ est le m\^eme qu'en \ref{E6deppadique}. On prouve alors  \ref{resultats}(3) de la m\^eme fa\c{c}on que dans ce paragraphe.
  
  Explicitons la cons\'equence de \ref{resultats}(4):
  
  (1) $dim(FC^{st}(\mathfrak{g}(F))=2\delta_{3}(q-1)$.
  
 \subsection{Type $E_{6}$ quasi-d\'eploy\'e, $E/F$ ramifi\'ee}\label{E6quasideppadiqueram}
On fixe une extension quadratique $E/F$ ramifi\'ee. On fixe $\tau\in \Gamma_{F}-\Gamma_{E}$.  On suppose que $G$ est quasi-d\'eploy\'e de type $E_{6}$ et que $\Gamma_{F}$ agit sur ${\cal D}$ de la fa\c{c}on suivante: l'action de $\Gamma_{E}$ est triviale et  $\tau$ agit par l'automorphisme $\theta$. Dans les tables de Tits, le groupe est de type $F_{4}^{I}$.

Si $\delta_{3}(q-1)=0$, on pose ${\cal X}=\{0\}$ et $d_{0}=1$. Si $\delta_{3}(q-1)=1$, on pose ${\cal X}=\{0,35\}$, $d_{0}=1$ et $d_{0,35}=2$. 

Explicitons les descriptions de \ref{alcoves} et \ref{profondeur}. On a ${\cal A}^{nr}=X_{*}(T)^{I_{F}}\otimes_{{\mathbb Z}}{\mathbb R}$. L'ensemble $\Sigma^{nr}$ est celui des $\beta^{res}$ pour $\beta\in \Sigma$. Puisqu'on note $\beta$ les \'el\'ements de $\Sigma$, on note $\Delta=\{\beta_{i};i=1,...,6\}$. La base $\Delta^{nr}$ est \'egale \`a $\{\alpha_{2},\alpha_{4},\alpha_{35},\alpha_{16}\}$, o\`u $\alpha_{2}=\beta_{2}^{res}$, $\alpha_{4}=\beta_{4}^{res}$, $\alpha_{35}=\beta_{3}^{res}=\beta_{5}^{res}$, $\alpha_{16}=\beta_{1}^{res}=\beta_{6}^{res}$. L'ensemble de racines $\Sigma^{nr}$ est de type $F_{4}$. On a $\alpha_{0}=-\alpha_{2}-2\alpha_{4}-3\alpha_{35}-2\alpha_{16}$ (rappelons que $\alpha_{0}$ n'est pas l'oppos\'ee de la plus grande racine, c'est la racine telle que $-\check{\alpha}_{0}$ est la plus grande coracine). On a $e(\alpha_{2})=e(\alpha_{4})=1$, $e(\alpha_{0})=e(\alpha_{35})=e(\alpha_{16})=2$. La relation (2) de \ref{alcoves} est

(1) $2\alpha_{0}+2\alpha_{2}+4\alpha_{4}+6\alpha_{35}+4\alpha_{16}=0$. 

C'est-\`a-dire $d(\alpha_{0})=d(\alpha_{2})=2$, $d(\alpha_{4})=d(\alpha_{16})=4$ et $d(\alpha_{35})=6$. On a $\alpha_{0}^{aff}=\frac{1}{2}+\alpha_{0}$.  L'ensemble $\underline{S}(G)$ s'identifie \`a celui des sommets de l'alc\^ove, donc \`a l'ensemble $\Delta_{a}^{nr}=\{\alpha_{2},\alpha_{4},\alpha_{35},\alpha_{16},\alpha_{0}\}$. Introduisons le tore $\underline{T}$ d\'efini et d\'eploy\'e sur ${\mathbb F}_{q}$ tel que $X_{*}(\underline{T})=X_{*}(T)^{I_{F}}$. On le munit de l'ensemble de racines $\Sigma^{nr}$. Pour le sommet $s$ associ\'e \`a $\alpha\in \Delta_{a}^{nr}$, le groupe $G_{s}$ est d\'eploy\'e, il a pour tore maximal $\underline{T}$ et pour ensemble de racines le sous-ensemble $\Sigma^{nr}(G_{s})$ d\'ecrit en \ref{profondeur}. On note $s_{0}$, $s_{35}$ etc... le sommet attach\'e \`a $\alpha_{0}$, $\alpha_{35}$ etc... On voit que $G_{s_{0}}$ est de type $F_{4}$, $G_{s_{2}}=Sp(8)/\{\pm 1\}$, $G_{s_{\alpha_{4}}}=(SL(2)\times SL(4))/\boldsymbol{\zeta}_{4}(\bar{{\mathbb F}}_{q})$, o\`u $\zeta\in  \boldsymbol{\zeta}_{4}(\bar{{\mathbb F}}_{q})$ s'envoie sur les \'el\'ements centraux \'evidents $\zeta^2$ de $SL(2)$ et $\zeta$ de $SL(4)$, $G_{s_{\alpha_{3,5}}}=(SL(3)\times SL(3))/diag(\boldsymbol{\zeta}_{3}(\bar{{\mathbb F}}_{q}))$,  $G_{s_{\alpha_{1,6}}}=(SL(2)\times Spin(7))/diag(\{\pm 1\})$. D'apr\`es  \ref{Cn}, \ref{An-1dep} et \ref{Bn}, on a $FC(\mathfrak{g}_{s}({\mathbb F}_{q}))=\{0\}$ pour $s=s_{2},s_{4},s_{16}$. D'apr\`es \ref{F4}, $FC(\mathfrak{g}_{s_{0}}({\mathbb F}_{q}))$ est une droite. Il s'en d\'eduit une droite dans $FC(\mathfrak{g}(F))$ que l'on note $FC_{0}$. D'apr\`es \ref{An-1dep}, $FC(\mathfrak{g}_{s_{35}}({\mathbb F}_{q}))$ est  de dimension $2$ si $\delta_{3}(q-1)=1$ et est nul sinon. Si $\delta_{3}(q-1)=1$, on note $FC_{35}$ le sous-espace de dimension $2$ de $FC(\mathfrak{g}(F))$ issu de $FC(\mathfrak{g}_{s_{35}}({\mathbb F}_{q}))$. Cela d\'emontre \ref{resultats}(1). 

On pose ${\cal Y}={\cal X}$, ${\cal X}^{st}={\cal X}$ et on note $\varphi:{\cal X}\to {\cal Y}$ l'identit\'e. Prouvons que

(2) $FC(\mathfrak{g}(F))=FC^{st}(\mathfrak{g}(F))$.

En admettant cela, on pose $FC^{{\cal E}}_{y}=FC_{y}$ pour tout $y\in {\cal Y}={\cal X}$ et les assertions (2), (3) et (4) de \ref{resultats} sont triviales. 

Prouvons (2). Consid\'erons un couple 
  $(\sigma\mapsto \sigma_{G'},{\cal O})\in {\cal E}_{ell}(G)$.   Le groupe $\hat{\Omega}$ est le m\^eme qu'en \ref{E6deppadique}, il est d'ordre $3$ et engendr\'e par $\theta_{3}$.  Rappelons que $\omega_{G'}$ est un homomorphisme injectif de $\Gamma_{E_{G'}/E}$ dans $\hat{\Omega}$.   Si l'image de cet homomorphisme est $\hat{\Omega}$ tout entier, le groupe $\Gamma_{E_{G'}/F}$ a une structure qui est interdite par \ref{extensionsdiedrales}(4). Donc l'image est r\'eduite \`a $\{1\}$ et on a $E_{G'}=E$. Tout \'el\'ement de $\hat{\Omega}\hat{\theta}$ est conjugu\'e \`a $\hat{\theta}$ par un \'el\'ement de $\hat{\Omega}$. A \'equivalence pr\`es, on peut donc supposer $\omega_{G'}(\tau)=1$, c'est-\`a-dire $\sigma_{G'}=\sigma_{G}$ pour tout $\sigma\in \Gamma_{F}$. Si ${\cal O}=\{\hat{\alpha}_{0}\}$, on a ${\bf G}'={\bf G}$ et, comme toujours, on ne sait encore rien de $FC^{st}(\mathfrak{g}(F))$. Si ${\cal O}=\{\hat{\alpha}_{2}\}$, $\{\hat{\alpha}_{4}\}$, $\{\hat{\alpha}_{3},\hat{\alpha}_{5}\}$, on voit que $G'_{SC}$ contient un facteur $SL(2)$, resp. $SL(3)$, $SL(4)$, donc $FC^{st}(\mathfrak{g}'(F))=\{0\}$ d'apr\`es \ref{An-1deppadique}. Si ${\cal O}=\{\hat{\alpha}_{1},\hat{\alpha}_{6}\}$, le groupe $G'_{SC}$ contient un facteur $Spin_{E/F}(10)$ et encore $FC^{st}(\mathfrak{g}'(F))=\{0\}$ d'apr\`es \ref{Dnimppadiqueram} (1). La seule donn\'ee possible est donc ${\bf G}'={\bf G}$, ce qui d\'emontre (2).
  
  Explicitons la cons\'equence de \ref{resultats}(4):
  
  (3) $dim(FC^{st}(\mathfrak{g}(F))=1+2\delta_{3}(q-1)$.

 \subsection{Type $E_{6}$ quasi-d\'eploy\'e, $E/F$ ramifi\'ee, s\'eparation des \'el\'ements de $FC(\mathfrak{g}(F))$}\label{E6quasidepseparation}
    Le groupe $G$ est celui du  paragraphe pr\'ec\'edent et on suppose $\delta_{3}(q-1)=1$.   On 
    va introduire deux \'el\'ements $X_{i}\in \mathfrak{g}_{ell}(F)$ pour $i=1,2$ tels que
    
    (1) les formes lin\'eaires $f\mapsto S^G(X_{i},f)$ pour $i=1,2$ sont nulles sur $FC_{0}$ et se restreignent en une base du dual de $FC_{35}$. 
    
      Explicitons les constructions de \ref{alcoves} et \ref{profondeur} pour les deux sommets $s_{0}$ et $s_{35}$ associ\'es aux racines $\alpha_{0}, \alpha_{35}\in \Delta_{a}^{nr}$. On a attach\'e \`a chacun de ces sommets divers objets que l'on affecte d'un indice $0$ ou $35$. 
      
      On a $\alpha(s_{0})=0$ pour tout $\alpha\in \Delta^{nr}$, 
     $\Delta(U_{P_{0}})=\{\alpha_{4}\}$ (d'apr\`es \ref{F4}) et, en utilisant \ref{E6quasideppadiqueram}(1), on trouve que $r_{0}= \frac{1}{6}$. Soit $f_{0}$ un \'el\'ement non nul de $FC_{0}$.  On peut supposer que $f_{0}$ est la fonction $\tilde{f}$ de \ref{profondeur} associ\'ee \`a $s_{0}$. Le lemme de ce paragraphe entra\^{\i}ne que
    
    (2) le support de $f_{0}$ est contenu dans $\mathfrak{g}(F)_{\frac{1}{6}}$.
    
    On a $\alpha_{2}(s_{35})=\alpha_{4}(s_{35})=\alpha_{16}(s_{35})=0$ et $\alpha_{35}(s_{35})=\frac{1}{6}$. On a $\Delta_{U_{P_{35}}}=\Delta_{a}^{nr}-\{\alpha_{35}\}$ d'apr\`es \ref{An-1dep}. On trouve que $r_{35}=\frac{1}{18}$. On a $\alpha(x_{35})=\frac{1}{18}$ pour tout $\alpha\in \Delta^{nr}$.  On calcule $r_{35}(\alpha)=0$ pour toute racine positive $\alpha\in \Sigma^{nr}$ et $r_{35}(\alpha)=\frac{1}{e(\alpha)}$ pour toute racine n\'egative. Rappelons que, d'apr\`es \ref{profondeur}(9), $\mathfrak{k}_{x_{35},r_{35}}=\mathfrak{k}_{s_{35}}[U_{P_{35}}]$ est le sous-espace des points fixes par $\Gamma_{F}^{nr}$ dans

    (3)  $\mathfrak{k}_{x_{35},r_{35},F^{nr}}=\mathfrak{t}(F^{nr})_{\frac{1}{2}}\oplus\oplus_{\alpha\in \Sigma^{nr}}\mathfrak{u}(\alpha)_{r(\alpha)}$.  
    
    Rappelons que $\mathfrak{g}_{s_{35}}\simeq \mathfrak{sl}(3)\oplus \mathfrak{sl}(3)$. On a d\'efini en \ref{groupessurFq} un sous-espace $\mathfrak{g}_{s_{35},2}$. On peut l'identifier \`a la somme des deux sous-espaces de $\mathfrak{sl}(3)$ engendr\'es par les deux \'el\'ements de l'\'epinglage standard de cette alg\`ebre. Les deux sous-espaces en question sont les images dans  $\mathfrak{g}_{s_{35}}({\mathbb F}_{q})$ des sous-espaces des points fixes par $\Gamma_{F}^{nr}$ dans $\mathfrak{u}(\alpha_{2})_{0}\oplus \mathfrak{u}(\alpha_{4})_{0}$, resp. $\mathfrak{u}(\alpha_{16})_{0}\oplus \mathfrak{u}(\alpha_{0})_{\frac{1}{2}}$.  L'espace $FC_{35}(\mathfrak{g}(F))$ a deux g\'en\'erateurs naturels qui sont issus d'\'el\'ements de $fc(\mathfrak{g}_{s_{35}}({\mathbb F}_{q}))$. Notons $\tilde{f}_{i}$, pour $i=1,2$ ces deux \'el\'ements.  On a plus pr\'ecis\'ement $\tilde{f}_{i}=\tilde{f}_{N',\epsilon_{i}}\otimes \tilde{f}_{N'',\epsilon_{i}}$. Les deux facteurs vivent sur $\mathfrak{sl}(3)$; $N'$ et $N''$ sont des \'el\'ements fix\'es de l'orbite ouverte dans chacun des espaces de dimension $2$ que l'on vient de d\'ecrire, ils ne d\'ependent pas de $i$; on a identifi\'e les centres des deux facteurs $SL(3)$ et  les $\epsilon_{i}$ pour $i=1,2$ d\'ecrivent les deux caract\`eres non triviaux de ce centre.   On a fix\'e une base de Chevalley de $\mathfrak{g}(\bar{F})$, qui contient un ensemble $(E_{\beta})_{\beta\in \Sigma}$. On a suppos\'e que $\theta(E_{\beta})=E_{\theta(\beta)}$ pour tout $\beta\in \Sigma$. Fixons une uniformisante $\varpi_{E}$ de $E$ telle que  $\tau(\varpi_{E})=-\varpi_{E}$. On peut supposer que $N'$ est l'image de $E_{2}+E_{4}$ dans la premi\`ere copie de $\mathfrak{sl}(3)$, o\`u on a pos\'e $E_{l}=E_{\beta_{l}}$, et que $N''$ est l'image de $E_{1}+E_{6}+\varpi_{E}(E_{-123^24^256}-E_{-1234^25^26})$ dans la seconde copie, o\`u on a pos\'e par exemple, $E_{-123^24^256}=E_{-\alpha_{1}-\alpha_{2}-2\alpha_{3}-2\alpha_{4}-\alpha_{5}-\alpha_{6}}$. 
  Fixons une uniformisante $\varpi_{E}$ de $E$ telle que  $\tau(\varpi_{E})=-\varpi_{E}$. Soit $v$ une racine de l'unit\'e dans $F^{\times}$ d'ordre premier \`a $p$.  Posons
$$(4) \qquad X(v)=vE_{2}+E_{4}+(E_{3}+E_{5})+(E_{1}+E_{6})+\varpi_{E}(E_{-123^24^256}-E_{-1234^25^26}).$$
 Gr\^ace \`a (3), on voit que $X(v)$ appartient \`a $\mathfrak{k}_{x_{35},r_{35}}$. On voit aussi  que son image $\underline{X}(v)$ dans $\mathfrak{g}_{s_{35}}({\mathbb F}_{q})=\mathfrak{sl}(3)({\mathbb F}_{q})\oplus \mathfrak{sl}(3)({\mathbb F}_{q})$ est l'image du couple $(vE_{2}+E_{4},E_{1}+E_{6}+\varpi_{E}(E_{-123^24^256}-E_{-1234^25^26}))$.  La deuxi\`eme composante est $N''$. La premi\`ere est conjugu\'ee \`a $N'$ par un \'el\'ement de $GL(3)({\mathbb F}_{q})$ de d\'eterminant $v$ (en notant encore $v$ l'image de ce terme dans ${\mathbb F}_{q}^{\times}$). Pour $i=1,2$, $\epsilon_{i}$ s'identifie \`a un caract\`ere de ${\mathbb F}_{q}^{\times}/{\mathbb F}_{q}^{\times,3}$ et on obtient  
 
 (5) $\tilde{f}_{i}(\underline{X}(v))=\epsilon_{i}(v)$.

    Montrons que
    
    (6) $x_{35}$ est l'unique point $x\in Imm_{F^{nr}}(G_{AD})$ tel que $X(v)\in \mathfrak{k}_{x,r_{35},F^{nr}}$.
    
    Soit un tel point $x$. Alors, pour tout $y$ dans la g\'eod\'esique joignant $x$ \`a $x_{35}$, on a $X(v)\in \mathfrak{k}_{y,r_{35},F^{nr}}$. Le point $x_{35}$ est \`a l'int\'erieur de l'alc\^ove $C^{nr}$. 
   Si $x\not=x_{35}$,  on peut donc choisir un point $y$ dans la g\'eod\'esique tel que  $y\not=x_{35}$ et $y\in C^{nr}$. Ecrivons $X(v)=\sum_{\alpha\in \Sigma^{nr}}X(v,\alpha)$, avec $X(v,\alpha)\in \mathfrak{u}(\alpha)$. La condition $X(v)\in \mathfrak{k}_{y,r_{35},F^{nr}}$ \'equivaut \`a $X(v,\alpha)\in \mathfrak{u}(\alpha)_{-\alpha(y)+r_{35}}$ pour tout $\alpha$. On voit sur la d\'efinition (4) que cela entra\^{\i}ne $-\alpha(y)+r_{35}\leq 0$ pour $\alpha\in \Delta^{nr}$ et $-\alpha_{0}(y)+r_{35}\leq \frac{1}{2}$. En posant $z=y-x_{35}$, cela \'equivaut \`a $\alpha(z)\geq 0$ pour tout $\alpha\in \Delta_{a}^{nr}$. La relation (1) de \ref{E6quasideppadiqueram} entra\^{\i}ne alors $z=0$ contrairement \`a l'hypoth\`ese. Cette contradiction prouve (6).

    Montrons que
    
    (7) $X(v)$ est semi-simple r\'egulier elliptique.
    
    On d\'ecompose $X(v)$ en somme de sa partie semi-simple et de sa partie nilpotente:   $X(v)=X(v)_{s}+X(v)_{n}$. La profondeur de $X(v)_{s}$ est la m\^eme que $X(v)$, donc $r(X(v)_{s})\geq r_{35}$. Introduisons le commutant connexe $G'=G_{X(v)_{s}}$, qui devient un Levi de $G$ sur $\bar{F}$. Notons $Imm_{F^{nr}}(G')$ son immeuble \'elargi. On affecte d'un  $'$ les objets relatifs \`a cet immeuble.  Cet immeuble se plonge dans $Imm_{F^{nr}}(G_{AD})$. Modulo ce plongement, on a $\mathfrak{k}'_{x,r,F^{nr}}\subset \mathfrak{k}_{x,r,F^{nr}}$  pour tout point $x\in Imm_{F^{nr}}(G')$ et tout $r\in {\mathbb R}$.
      Si $X(v)_{n}\not=0$, le groupe $G'$ contient un tore maximal $T'$ d\'efini  et maximalement d\'eploy\'e sur $F^{nr}$ et  tel qu'il existe $t\in T'(F^{nr})$ de sorte que $ad(t)^m(X(v)_{n})$ tend vers $0$ quand $m$ tend vers l'infini. Soit ${\cal A}'$ l'appartement de $Imm_{F^{nr}}(G')$ associ\'e \`a ce tore. Puisque $X(v)_{s}$ est central dans $G'$ et $r(X(v)_{s})\geq r_{35}$, on  a $X(v)_{s}\in \mathfrak{k}'_{x,r_{35},F^{nr}}$ pour tout $x\in {\cal A}'$. Il est clair qu'il existe une infinit\'e de points $x\in {\cal A}'$ tels que l'on ait aussi $X(v)_{n}\in \mathfrak{k}'_{x,r_{35},F^{nr}}$ (\`a partir d'un point $x$ quelconque, le point $t^{-m}x$ v\'erifie cette condition pour $m$ assez grand). Alors $X(v)\in \mathfrak{k}'_{x,r_{35},F^{nr}}$ donc aussi \`a $\mathfrak{k}_{x,r_{35},F^{nr}}$ pour une infinit\'e de $x$, ce qui est exclu par (6). Cela prouve que $X(v)_{n}=0$ et $X(v)$ est semi-simple. Par le m\^eme argument, $X(v)$ appartient \`a $\mathfrak{k}'_{x,r_{35},F^{nr}}$ pour tout point $x\in Imm_{F^{nr}}(G')$ et on conclut comme pr\'ec\'edemment que cet immeuble est r\'eduit \`a un point. Cela entra\^{\i}ne que tout sous-tore de $G'$ qui est d\'eploy\'e sur $F^{nr}$ est r\'eduit \`a $\{1\}$. Un tel groupe $G'$ est forc\'ement un tore anisotrope sur $F^{nr}$. 
   Un \'el\'ement de $\mathfrak{g}(F^{nr})$ dont le commutant v\'erifie ces conditions est r\'egulier elliptique. Cela d\'emontre (7). 
    
    Pour tout point $x\in Imm(G_{AD})$ et tout $r\in {\mathbb R}$ tel que $r> r_{35}$, on a $X(v)\not\in \mathfrak{k}_{x,r}$: sinon on a a fortiori $X(v)\in \mathfrak{k}_{x,r_{35}}$ donc $x=x_{35}$ d'apr\`es (6), or un calcul analogue \`a celui de (6) montre que $X(v)\not\in \mathfrak{k}_{x_{35},r}$. Donc 
    
    (8) $r(X(v))=r_{35}=\frac{1}{18}$. 
    
    On a aussi  $r(X)=\frac{1}{18}$ pour tout \'el\'ement $X$ stablement conjugu\'e \`a $X(v)$.  Alors l'assertion (2) implique
    
    (9) $S^G(f_{0},X(v))=0$. 
    
      A l'aide de (5) et (6), on voit comme en \ref{D4trialitaireseparation} que
   
    $S^G(\tilde{f}_{i})=c\epsilon_{i}(v)$ pour $i=1,2$,
    
   \noindent o\`u $c$ est une constante non nulle. On choisit un \'el\'ement $v$ tel que $\epsilon_{1}(v)\not=\epsilon_{2}(v)$. On pose $X_{1}=X(1)$, $X_{2}=X(v)$. Alors la deuxi\`eme condition de (1) est v\'erifi\'ee.  La premi\`ere l'est d'apr\`es (9). Cela prouve (1). 
   
   \subsection{Type $E_{6}$ quasi-d\'eploy\'e, $E/F$ ramifi\'ee, action d'un automorphisme}\label{E6automorphisme}
   Le groupe $G$ est comme en \ref{E6quasideppadiqueram}. Il  a un automorphisme $\theta$. Celui-ci agit trivialement sur ${\cal A}^{nr}$ par d\'efinition de cet ensemble. Il fixe donc les sommets $s_{0}$ et $s_{35}$ et donc aussi les r\'eseaux $\mathfrak{k}_{s_{0}}$ et $\mathfrak{k}_{s_{35}}$.
   
    On voit que $\mathfrak{k}_{s_{0},F^{nr}}$ est le $\mathfrak{o}_{F^{nr}}$-module engendr\'e par les $E_{\beta}$ pour $\beta\in \Sigma$ tel que $\theta(\beta)=\beta$ et les $E_{\beta}+E_{\theta(\beta)}$ pour $\beta\in \Sigma$ tel que $\theta(\beta)\not=\beta$. Alors $\theta$ agit trivialement sur la r\'eduction $\mathfrak{g}_{s_{0}}({\mathbb F}_{q})$, donc aussi sur $FC_{0}$.
    
  Supposons $\delta_{3}(q-1)=1$. On calcule $\theta(X(v))$ avec les notations du paragraphe pr\'ec\'edent. On obtient
  $$\theta(X(v))=vE_{2}+E_{4}+(E_{3}+E_{5})+(E_{1}+E_{6})-\varpi_{E}(E_{-123^24^256}-E_{-1234^25^26}),$$
  le dernier terme a chang\'e de signe. Donc $\theta$ n'agit pas trivialement sur $\mathfrak{g}_{s_{35}}({\mathbb F}_{q})$. Toutefois, la r\'eduction $red(\theta(X(v)))$  de $\theta(X(v))$ dans cet espace, lequel est isomorphe \`a $\mathfrak{sl}(3)({\mathbb F}_{q})\oplus \mathfrak{sl}(3)({\mathbb F}_{q})$,  a la m\^eme premi\`ere composante que la r\'eduction $red(X(v))$ de $X(v)$ et une seconde composante conjugu\'ee \`a celle de $red(X(v))$ par un \'el\'ement de $GL(3,{\mathbb F}_{q})$ de d\'eterminant $-1$. Or $-1$ est un cube de ${\mathbb F}_{q}^{\times}$. Cela entra\^{\i}ne que les fonctions $\tilde{f}_{i}$, pour $i=1,2$, prennent la m\^eme valeur sur $red(X(v))$ et  $red(\theta(X(v)))$, o\`u encore que   $\theta(\tilde{f}_{i})(X(v))=\tilde{f}_{i}(X(v))$ pour $i=1,2$. Puisque les \'evaluations sur $X(v)$ s\'eparent les deux fonctions $\tilde{f}_{1}$ et $\tilde{f}_{2}$ quand $v$ d\'ecrit les racines de l'unit\'e de $F^{\times}$ d'ordre premier \`a $p$,  on a $\theta(\tilde{f}_{i})=\tilde{f}_{i}$ pour $i=1,2$. Donc $\theta$ agit trivialement sur $FC_{35}$. 

\subsection{Type $E_{7}$ d\'eploy\'e}\label{E7padique}
On suppose $G$ d\'eploy\'e de type $E_{7}$. On d\'efinit un homomorphisme $T_{ad}(F)\to F^{\times}/F^{\times,2}$ par $\prod_{l=1,...,7}\check{\varpi}_{l}(x_{l})\mapsto x_{2}x_{5}x_{7}$. De cet homomorphisme est issu un isomorphisme

\noindent $G_{AD}(F)/\pi(G(F))\simeq T_{ad}(F)/\pi(T(F))\simeq F^{\times}/F^{\times,2}$. L'image de $G_{AD}(F)_{0}$ est $\mathfrak{o}_{F}^{\times}/\mathfrak{o}_{F}^{\times,2}$. On note $\Xi^{ram}$ le sous-ensemble \`a deux \'el\'ements de $\Xi$ form\'e des \'el\'ements non triviaux sur $G_{AD}(F)_{0}/\pi(G(F))$. 

On pose  ${\cal X}_{0}= \{(0,\xi); \xi\in \Xi^{ram}\}$. Si $\delta_{3}(q-1)=0$, on pose $ {\cal X}_{3}=\emptyset$;  si $\delta_{3}(q-1)=1$, on pose  ${\cal X}_{3}=\{(3,\xi); \xi\in \Xi^{ram}\}$. Si $\delta_{4}(q-1)=0$, on pose $ {\cal X}_{4}= \emptyset$; si $\delta_{4}(q-1)=1$, on pose  $ {\cal X}_{4}= \{ 4\}$. On pose ${\cal  X}={\cal X}_{0}\cup {\cal  X}_{3}\cup {\cal X}_{4}$, $d_{x}=1$ pour $x\in {\cal X}_{0}$ et $d_{x}=2$ pour $x\in {\cal X}_{3}\cup {\cal X}_{4}$.

Le groupe d'automorphismes du diagramme ${\cal D}_{a}$ est le groupe a deux \'el\'ements $\{1,\omega\}$, o\`u $\omega$ permute les sommets $\alpha_{0}$ et $\alpha_{7}$, $\alpha_{1}$ et $\alpha_{6}$, $\alpha_{3}$ et $\alpha_{5}$ et fixe les deux sommets $\alpha_{2}$ et $\alpha_{4}$. Le groupe $G_{AD}(F)$ agit sur le diagramme ${\cal D}_{a}$ par ce groupe d'automorphismes. L'ensemble $\underline{S}(G)$ s'envoie surjectivement sur celui des orbites de cette action, la fibre au-dessus d'une orbite ayant pour nombre d'\'el\'ements le nombre d'\'el\'ements de l'orbite. Pour un sommet $s$ param\'etr\'e par l'orbite de $\alpha_{0}$, le groupe $G_{s}$ est simplement connexe de type $E_{7}$. D'apr\`es \ref{E7}, l'espace $FC(\mathfrak{g}_{s}({\mathbb F}_{q}))$ est une droite port\'ee par un g\'en\'erateur  $f_{N,\epsilon}$. Le centre de $G$ s'identifie \`a celui de $G_{s}$ et $\epsilon$ est non trivial sur ce centre. La fonction se transforme donc par le caract\`ere non trivial de  $G_{AD}(F)_{0}/\pi(G(F))$. Conform\'ement \`a \ref{actionsurFC}, il se d\'eduit de cette fonction deux \'el\'ements de $FC(\mathfrak{g}(F))$ qui se transforment selon les deux caract\`eres de $G_{AD}(F)/\pi(G(F))$ prolongeant ce caract\`ere non trivial, c'est-\`a-dire selon les deux \'el\'ements $\xi\in \Xi^{ram}$. Pour un tel caract\`ere $\xi\in \Xi^{ram}$, on note $FC_{0,\xi}$ la droite port\'ee par la fonction en question. 
 Pour un sommet $s$ param\'etr\'e par l'orbite de $\alpha_{1}$, on a $G_{s,SC}=SL(2)\times Spin_{dep}(12)$ et $FC(\mathfrak{g}_{s}({\mathbb F}_{q}))=\{0\}$ d'apr\`es \ref{Dndeppair}. Pour un sommet param\'etr\'e par l'orbite de $\alpha_{3}$, on a $G_{s}=(SL(3)\times SL(6))/\boldsymbol{\zeta}_{3}(\bar{{\mathbb F}}_{q})$. Pr\'ecis\'ement $\zeta\in \boldsymbol{\zeta}_{3}(\bar{{\mathbb F}}_{q})$ s'envoie sur l'\'el\'ement $\check{\alpha}_{0}(\zeta)\check{\alpha}_{1}(\zeta^2)$ du centre de $SL(3)$ et sur l'\'el\'ement $\check{\alpha}_{2}(\zeta^2)\check{\alpha}_{4}(\zeta)\check{\alpha}_{6}(\zeta^2)\check{\alpha}_{7}(\zeta)$ du centre de $SL(6)$. Si $\delta_{3}(q-1)=0$, on a $FC(\mathfrak{g}_{s}({\mathbb F}_{q}))=\{0\}$ d'apr\`es \ref{An-1dep}. Supposons $\delta_{3}(q-1)=1$. Sur chacune des composantes $SL(3)$ et $SL(6)$, il y a deux fonctions du type $f_{N,\epsilon}$, d'o\`u $4$ fonctions produits tensoriels. Mais deux seulement se factorisent par $\boldsymbol{\zeta}_{3}(\bar{{\mathbb F}}_{q})$. Le centre de $G$ s'envoie sur le sous-groupe d'ordre $2$ du centre de $SL(6)$ et le caract\`ere $\epsilon$ de chacune de nos fonctions est non trivial sur ce groupe. Donc les deux fonctions se transforment selon le caract\`ere non trivial de $G_{AD}(F)_{0}/\pi(G(F))$.  Il se d\'eduit de nos deux fonctions $4$ \'el\'ements de $FC(\mathfrak{g}(F))$. Pour chaque caract\`ere $\xi\in \Xi$   prolongeant ce caract\`ere non trivial, c'est-\`a-dire  pour chaque $\xi\in \Xi^{ram}$, on obtient deux \'el\'ements de $FC(\mathfrak{g}(F))$ se transformant selon $\xi$. On note $FC_{3,\xi}$ le plan engendr\'e par ces deux fonctions. 
  Pour un sommet $s$ param\'etr\'e par $\alpha_{4}$, on a $G_{s}=(SL(4)\times SL(4)\times SL(2))/\boldsymbol{\zeta}_{4}(\bar{{\mathbb F}}_{q})$. Pr\'ecis\'ement, $\zeta\in \boldsymbol{\zeta}_{4}(\bar{{\mathbb F}}_{q})$ s'envoie sur l'\'el\'ement $\check{\alpha}_{0}(\zeta)\check{\alpha}_{1}(\zeta^2)\check{\alpha}_{3}(\zeta^3)$ du centre de la premi\`ere composante, sur l'\'el\'ement $\check{\alpha}_{7}(\zeta)\check{\alpha}_{6}(\zeta^2)\check{\alpha}_{5}(\zeta^3)$ du centre  de la seconde et sur l'\'el\'ement $\check{\alpha}_{2}(\zeta^2)$ du centre de $SL(2)$. Si $\delta_{4}(q-1)=0$, on a $FC(\mathfrak{g}_{s}({\mathbb F}_{q}))=\{0\}$ d'apr\`es \ref{An-1dep}. Supposons $\delta_{4}(q-1)=1$. Il y a deux fonctions du type $f_{N,\epsilon}$ sur chaque composante $SL(4)$ et une fonction de ce type sur $SL(2)$. D'o\`u $4$ fonctions produits tensoriels. Mais deux seulement se factorisent par $\boldsymbol{\zeta}_{4}({\mathbb F}_{q})$. L'\'el\'ement non trivial du centre de $G$ s'envoie sur le produit de l'\'el\'ement trivial de la premi\`ere composante $SL(4)$, de l'\'el\'ement d'ordre $2$ du centre de la seconde composante et de l'\'element non trivial du centre de $SL(2)$. On voit que le produit des caract\`eres $\epsilon$ vaut $1$ sur cet \'el\'ement donc nos fonctions sont invariantes par $G_{AD}(F)_{0}/\pi(G(F))$. Il se d\'eduit de nos deux fonctions deux \'el\'ements de $FC(\mathfrak{g}(F))$. On note $FC_{4}$ le plan qu'elles engendrent.    Pour le sommet $s$ param\'etr\'e par $\alpha_{2}$, on a $G_{s}=SL(8)/\{\pm 1\}$ et $FC(\mathfrak{g}_{s}({\mathbb F}_{q}))=\{0\}$ d'apr\`es \ref{An-1dep}. Cela d\'emontre \ref{resultats}(1).

 On pose  ${\cal Y}_{0}=\{(07,0,\xi); \xi\in \Xi^{ram}\}$. Si $\delta_{3}(q-1)=0$, on pose $ {\cal Y}_{3} =\emptyset$;  si $\delta_{3}(q-1)=1$, on pose  ${\cal Y}_{3}=\{( 07,35,\xi); \xi\in \Xi^{ram}\}$.  On pose  ${\cal Y}_{4}={\cal X}_{4}$. On pose ${\cal  Y}={\cal  Y}_{0}\cup {\cal  Y}_{3}\cup {\cal  Y}_{4}$.

  Consid\'erons un couple $(\sigma\mapsto \sigma_{G'},{\cal O})\in {\cal E}_{ell}(G)$. Le groupe $\hat{\Omega}$ est \'egal \`a $\{1,\hat{\omega}\}$, o\`u $\hat{\omega}$ est  similaire au $\omega$ ci-dessus (les diagrammes ${\cal D}_{a}$ et $\hat{{\cal D}}_{a}$ sont les m\^emes).  
  
Supposons que l'action $\sigma\mapsto \sigma_{G'}$ soit triviale. Si ${\cal O}=\{\hat{\alpha}_{0}\}$ ou $\{\hat{\alpha}_{7}\}$ (ces deux cas sont conjugu\'es), on a ${\bf G}'={\bf G}$ et, \`a ce point, on ne peut rien dire de $FC^{st}(\mathfrak{g}(F))$. Si ${\cal O}=\{\hat{\alpha}_{m}\}$, avec $m\in \{1,...,6\}$, on voit que le groupe $G'_{SC}$ contient un facteur $SL(l)$ avec $l\geq2$. Donc $FC^{st}(\mathfrak{g}'(F))=\{0\}$ d'apr\`es \ref{An-1deppadique}. 

Supposons que $E_{G'}/F$ soit quadratique et fixons un \'el\'ement $\tau\in \Gamma_{F}-\Gamma_{E_{G'}}$. Alors $\tau_{G'}=\hat{\omega}$. Si ${\cal O}=\{\hat{\alpha}_{4}\}$, on a $G'_{SC}=Res_{E_{G'}/F}(SL(4))\times SL(2)$ et $FC^{st}(\mathfrak{g}'(F))=\{0\}$. Si ${\cal O}=\{\hat{\alpha}_{2}\}$, on a $G'_{SC}=SU_{E_{G'}/F}(8)$ et encore $FC^{st}(\mathfrak{g}'(F))=\{0\}$ d'apr\`es \ref{An-1quasidepnonram} et \ref{An-1quasidepram} (5). Il ne reste que des orbites ${\cal O}$ a deux \'el\'ements. Le lemme \ref{centre} exclut le cas o\`u $E_{G'}/F$ est non ramifi\'ee. Supposons donc $E_{G'}/F$ ramifi\'ee. On peut alors remplacer le groupe $G'$ par $G'_{SC}$.  Si ${\cal O}=\{\hat{\alpha}_{1},\hat{\alpha}_{6}\}$ ou ${\cal O}=\{\hat{\alpha}_{3},\hat{\alpha}_{5}\}$, on voit que $G'_{SC}$ contient un facteur $Res_{E_{G'}/F}(SL(2))$ ou $Res_{E_{G'}/F}(SL(3))$ et $FC^{st}(\mathfrak{g}'(F))=\{0\}$ d'apr\`es \ref{An-1deppadique}. Si ${\cal O}=\{\hat{\alpha}_{0},\hat{\alpha}_{7}\}$, le groupe $G'_{SC}$ est de type $E_{6}$ avec action galoisienne non triviale par $\Gamma_{E_{G'}/F}$. D'apr\`es \ref{E6quasideppadiqueram}, $FC^{st}(\mathfrak{g}'(F))$ est de dimension $1+2\delta_{3}(q-1)$. Pr\'ecis\'ement, il est somme d'une droite not\'ee $FC_{0}$ dans \ref{E6quasideppadiqueram} et, si $\delta_{3}(q-1)=1$, d'un plan not\'e $FC_{35}$. Nous notons ici ces espaces $FC_{0}(\mathfrak{g}'(F))$ et $FC_{35}(\mathfrak{g}'(F))$. 
 Le groupe d'automorphismes ext\'erieurs  de ${\bf G}'$ est $\hat{\Omega}$ et  ce groupe agit trivialement sur $FC^{st}(\mathfrak{g}'(F))$ d'apr\`es \ref{E6automorphisme}. Calculons le caract\`ere $\xi_{{\bf G}'}$. L'\'el\'ement $s$ de la donn\'ee ${\bf G}'$ v\'erifie $\hat{\alpha}_{l}(s)=1$ pour $l=1,...,6$ et $\hat{\alpha}_{7}(s)=-1$. On peut  choisir
$$s_{sc}=\check{\hat{\alpha}}_{1}(-1)\check{\hat{\alpha}}_{2}(i^{-1})\check{\hat{\alpha}}_{4}(-1)\check{\hat{\alpha}}_{5}(i)\check{\hat{\alpha}}_{7}(i^{-1}).$$
Alors $\hat{\omega}(s_{sc})s_{sc}^{-1}$
est l'\'element non trivial du centre de $G$ donc  $\xi_{{\bf G}'}$ est le caract\`ere d'ordre $2$ de $\Gamma_{F}$ qui se factorise par $\Gamma_{E_{G'}/F}$. Quand $E_{G'}/F$ d\'ecrit les deux extensions quadratiques ramifi\'ees de $F$, $\xi_{{\bf G}'}$ d\'ecrit les deux \'el\'ements de $\Xi$ qui sont non triviaux sur le sous-groupe $\mathfrak{o}_{F}^{\times}/\mathfrak{o}_{F}^{\times,2}$ de $F^{\times}/F^{\times,2} $. Autrement dit les deux \'el\'ements de $\Xi^{ram}$. Pour un caract\`ere $\xi\in \Xi^{ram}$, on note ${\bf G}'_{07,\xi}$ la donn\'ee ${\bf G}'$ telle que $\xi_{{\bf G}'}=\xi$. On pose $FC^{{\cal E}}_{07,0,\xi}=FC_{0}(\mathfrak{g}'_{07,\xi}(F))$ et, si $\delta_{3}(q-1)=0$, $FC^{{\cal E}}_{07,35,\xi}=FC_{35}(\mathfrak{g}'_{07,\xi}(F))$.

 On note  $\varphi:{\cal X}\to {\cal Y}$ la bijection \'evidente qui envoie ${\cal X}_{i}$ sur ${\cal Y}_{i}$ pour $i=0,3,4$.   On pose ${\cal X}^{st}= {\cal X}_{4}$. 
 
 Nous avons associ\'e un espace $FC^{{\cal E}}_{y}$ \`a tout \'el\'ement de ${\cal Y}-{\cal Y}_{4}$. On n'a pas trait\'e la donn\'ee principale ${\bf G}$. Par l'argument habituel de comparaison des dimensions, on voit que $dim(FC^{st}(\mathfrak{g}(F))=dim(FC_{4})=2\delta_{4}(q-1)$. Si $\delta_{4}(q-1)=0$, la donn\'ee principale ne contribue pas \`a l'espace $FC^{{\cal E}}(\mathfrak{g}(F))$ et on a achev\'e de prouver \ref{resultats}(2), et en m\^eme temps \ref{resultats}(4). 
 Supposons $\delta_{4}(q-1)=1$. Alors $FC^{st}(\mathfrak{g}(F))\subset I_{cusp}^{st}(\mathfrak{g}(F))\cap FC(\mathfrak{g}(F))$. Or, pour $x\in {\cal X}_{0}\cup {\cal X}_{3}$, $FC_{x}$ est contenu dans un espace $I_{cusp,\xi}(\mathfrak{g}(F))$ pour un caract\`ere $\xi$ non trivial. Donc $I_{cusp}^{st}(\mathfrak{g}(F))\cap FC(\mathfrak{g}(F))$ est inclus dans l'orthogonal de ces espaces et est donc contenu dans $FC_{4}$. On obtient $FC^{st}(\mathfrak{g}(F))\subset FC_{4}$. Ces espaces \'etant de dimension $2$, ils sont \'egaux. On compl\`ete la description \ref{resultats}(2) en posant $FC^{{\cal E}}_{4}=FC_{4}$. On obtient en  m\^eme temps \ref{resultats}(4).

Il reste \`a d\'emontrer \ref{resultats}(3). Si $\delta_{3}(q-1)=0$, c'est \'evident:  on a forc\'ement 

\noindent $transfert(\oplus_{x\in {\cal X}_{0}}FC_{x})=\oplus_{y\in {\cal Y}_{0}}FC^{{\cal E}}_{y}$ d'apr\`es ce qui pr\'ec\`ede et les espaces int\'erieurs se s\'eparent par les  caract\`eres $\xi\in \Xi^{ram}$ par lesquels agit le groupe $G_{AD}(F)/\pi(G(F))$. Supposons $\delta_{3}(q-1)=1$. Posons ${\cal  X}^{\star}={\cal X}_{0}\cup {\cal X}_{3}$, ${\cal  Y}^{\star}={\cal Y}_{0}\cup {\cal Y}_{3}$. On munit ${\cal X}^{\star}$ de la relation de pr\'eordre d\'efini ainsi: pour $x\in {\cal X}_{i}$ et $x'\in {\cal X}_{i'}$, avec $i,i'\in \{0,3\}$, $x\leq x'$ si et seulement si $i\leq i'$.  On applique les constructions de \ref{ingredients}. L'ensemble $\underline{{\cal Y}}^{\star}$ s'identifie \`a $ \{0,3\}$ dont le plus petit \'el\'ement $(y_{min})$ est $0$. On pose $\underline{{\cal Y}}^{\sharp}=\underline{{\cal Y}}^{\star}-\{(y_{min})\}$. Il a un unique \'el\'ement que l'on  note $(y)$ et qui n'est autre que ${\cal Y}_{3}$.   On a $d_{(y)}=4$. Pour r\'econcilier nos notations avec celles de \ref{ingredients}, on identifie l'ensemble $\{1,...,4\}$ avec $I=\{(\xi,i); \xi\in \Xi^{ram},i=1,2\}$.  Pour $(\xi,i)\in I$, on pose ${\bf G}'_{\xi,i}={\bf G}'_{07,\xi}$. 
  Pour chaque groupe $G'_{07,\xi}$, on a construit en \ref{E6quasidepseparation} deux \'el\'ements $X_{1},X_{2}\in \mathfrak{g}'_{07,\xi,reg}(F)$. Il n'est pas clair qu'ils soient $G$-r\'eguliers. Mais, comme en \ref{ingredients}, on peut les remplacer par des \'el\'ements tr\`es voisins qui le sont. On fixe de tels \'el\'ements que l'on note $Y_{\xi,1}$ et $Y_{\xi,2}$.   On associe \`a $(y)$ les collections $(G'_{\xi,i})_{(\xi,i)\in I}$ et $(Y_{\xi,i})_{(\xi,i)\in I}$. 

Montrons que les conditions de \ref{ingredients} sont satisfaites. La relation (1) de ce paragraphe est v\'erifi\'ee comme plus haut par simple compatibilit\'e du transfert avec les actions de $G_{AD}(F)/\pi(G(F))$. La condition (2) est satisfaite d'apr\`es notre description des espaces en question. Les conditions (3) et (4) r\'esultent de la condition (1) de \ref{E6quasidepseparation} et du fait que, pour un \'el\'ement $(\xi,i)\in I$, un \'el\'ement $y'\in {\cal Y}^{\star}$ et une fonction $f'_{y'}\in FC^{{\cal E}}_{y'}$, les d\'efinitions entra\^{\i}nent que $S^{G'_{07,\xi}}(Y_{\xi,i},f'_{y'})$ ne peut \^etre non nul que si  l'\'el\'ement $\xi'$ qui figure dans $y'$ est \'egal \`a $\xi$. Il reste \`a prouver la condition (5). Celle-ci se r\'ecrit sous la forme suivante:

 (1) soient $x=(0,\xi')\in {\cal X}$, $f\in FC_{x}$, $(\xi,i)\in I$ et $X$ un \'el\'ement de $\mathfrak{g}_{reg}(F)$ dont la classe de conjugaison stable correspond \`a celle de $Y_{\xi,i}$; alors $I^G(X,f)=0$. 
   
 D'apr\`es   \ref{E6quasidepseparation} (6), on a $r(X)=\frac{1}{18}$. Il suffit donc de prouver que le support de $f$ est contenu dans $\mathfrak{g}(F)_{r}$ pour un $r>r(X)$. La fonction $f$ est par d\'efinition combinaison lin\'eaire de deux  fonctions issues des sommets associ\'es aux deux racines $\alpha_{0}$ et $\alpha_{7}$. Ces fonctions sont conjugu\'ees par l'action de $G_{AD}(F)/\pi(G(F))$ et leurs supports sont aussi conjugu\'es. Il suffit de consid\'erer l'une d'elles, par exemple celle associ\'ee au sommet $s_{0}$ attach\'e  \`a la racine $\alpha_{0}$.  C'est une fonction $\tilde{f}$ du type consid\'er\'e en \ref{profondeur}. Le syst\`eme de racines $\Sigma^{nr}$ de  $G_{s_{0}}$ est $\Sigma$ tout entier. D'apr\`es \ref{E7}, on a $\Delta(U_{P_{0}})=\{\alpha_{4},\alpha_{7}\}$.  On a $d(\alpha_{0})=d(\alpha_{7})=1$ et $d(\alpha_{4})=4$. Alors le nombre $r$ de \ref{profondeur} est \'egal \`a $\frac{1}{6}$, donc $\tilde{f}$ est \`a support dans $\mathfrak{g}(F)_{\frac{1}{6}}$. Cela d\'emontre (1). 
 
 Alors le lemme de \ref{ingredients} nous dit que $transfert(FC_{(x)})=FC^{{\cal E}}_{\varphi((x))}$ pour tout $(x)\in \underline{{\cal X}}^{\star}$. Comme d'habitude, l'action de $G_{AD}(F)/\pi(G(F))$ permet de raffiner cette \'egalit\'e en $transfert(FC_{x})=FC^{{\cal E}}_{\varphi(x)}$ pour tout $x\in {\cal X}^{\star}$. Cela prouve \ref{resultats}(3). 
 
 Explicitons la cons\'equence de \ref{resultats}(4):
 
 (2) $dim(FC^{st}(\mathfrak{g}(F)))=2\delta_{4}(q-1)$. 
 
 \subsection{Forme int\'erieure du type $E_{7}$}
On suppose que $G^*$ est de type $E_{7}$ et que $G$ en est la forme int\'erieure non d\'eploy\'ee. 

Si $\delta_{4}(q-1)=0$, on pose ${\cal X} =\emptyset$. Si $\delta_{4}(q-1)=1$, on pose ${\cal X} =\{ 4\}$ et $d_{4}=2$. 

Le diagramme  ${\cal D}_{a}^{nr}$ est le diagramme de Dynkin affine d'un syst\`eme de racines  de type  $E_{7}$ sur lequel le Frobenius $Fr$ agit par l'automorphisme $\omega$ non trivial. L'ensemble $\underline{S}(G)$ est param\'etr\'e par l'ensemble des orbites de cette action. Pour un sommet $s$ param\'etr\'e par une orbite \`a $2$ \'el\'ements,  le lemme \ref{orbites} implique que $FC(\mathfrak{g}_{s}({\mathbb F}_{q}))=\{0\}$. Pour le sommet $s$ param\'etr\'e par $\{\alpha_{2}\}$, le groupe $G_{s}$ est le m\^eme que dans le paragraphe pr\'ec\'edent, avec une action galoisienne diff\'erente: on voit que $G_{s}=SU_{{\mathbb F}_{q^2}/{\mathbb F}_{q}}(8)/\{\pm 1\}$, d'o\`u $FC(\mathfrak{g}_{s}({\mathbb F}_{q}))=\{0\}$ d'apr\`es \ref{An-1nondep}. Pour le sommet $s$ param\'etr\'e par $\alpha_{4}$, on trouve de m\^eme $G_{s}=(Res_{{\mathbb F}_{q^2}/{\mathbb F}_{q}}(SL(4))\times SL(2))/\boldsymbol{\zeta}_{4}(\bar{{\mathbb F}}_{q})$. Sur $\bar{{\mathbb F}}_{q}$, on a $G_{s}=(SL(4)\times SL(4)\times SL(2))/\boldsymbol{\zeta}_{4}(\bar{{\mathbb F}}_{q})$. Notons $Fr_{dep}$ l'action "naturelle" du Frobenius sur chacun des facteurs. L'action du Frobenius sur $G_{s}$ est alors $(x,y,z)\mapsto (Fr_{dep}(y),Fr_{dep}(x),Fr_{dep}(z))$. Il  y a sur chacun des facteurs $SL(4)$ deux faisceaux-caract\`eres cuspidaux, param\'etr\'es par les deux caract\`eres d'ordre $4$ de $\boldsymbol{\zeta}_{4}(\bar{{\mathbb F}}_{q})$. Notons ces caract\`eres $\epsilon'$ et $\epsilon''$ et, pour $\epsilon$ l'un de ces caract\`eres,  notons ${\cal F}_{N_{4},\epsilon}$  le faisceau-caract\`ere associ\'e. Il y a sur $SL(2)$ un unique faisceau-caract\`ere, not\'e ${\cal F}_{N_{2}, \mu}$, o\`u $\mu$ est le caract\`ere non trivial de $\{\pm 1\}$. L'image par le Frobenius du faisceau ${\cal F}_{N_{4},\epsilon_{1}}\otimes {\cal F}_{N_{4},\epsilon_{2}}\otimes {\cal F}_{N_{2},\mu}$ est ${\cal F}_{N_{4},\epsilon_{2}\circ Fr}  \otimes {\cal F}_{N_{4},\epsilon_{1}\circ Fr}\otimes{\cal F}_{N_{2},\mu}$, o\`u ici, $Fr$ est l'action du Frobenius sur $  \bar{{\mathbb F}}_{q}$. Le faisceau-caract\`ere est invariant par le Frobenius si et seulement si $\epsilon_{1}=\epsilon_{2}\circ Fr$ et $\epsilon_{2}=\epsilon_{1}\circ Fr$. Ces deux \'egalit\'es sont \'equivalentes car $Fr^2$ agit trivialement sur $\boldsymbol{\zeta}_{4}(\bar{{\mathbb F}}_{q})$. Supposons donc    $\epsilon_{1}=\epsilon_{2}\circ Fr$. Pour que le faisceau-caract\`ere se descende en un tel faisceau sur $G_{s}$, il faut et il suffit que $\epsilon_{1}\times \epsilon_{2}\times  \mu$ soit trivial sur $(\zeta,\zeta,\zeta^2)$ pour tout $\zeta\in   \boldsymbol{\zeta}_{4}(\bar{{\mathbb F}}_{q})$, c'est-\`a-dire $\epsilon_{1}(\zeta Fr(\zeta))=\mu(\zeta^2)$. Cette \'egalit\'e ne peut \^etre v\'erifi\'ee que si $Fr$ est trivial sur $\boldsymbol{\zeta}_{4}(\bar{{\mathbb F}}_{q})$, c'est-\`a-dire si $\delta_{4}(q-1)=1$. Si $\delta_{4}(q-1)=1$, elle est v\'erifi\'ee pour les deux caract\`eres $\epsilon_{1}=\epsilon',\epsilon''$. On obtient dans ce cas deux faisceaux-caract\`eres cuspidaux sur $G_{s}$, donc $FC(\mathfrak{g}_{s}({\mathbb F}_{q}))$ est de dimension $2$.  De cet espace se d\'eduit un sous-espace de dimension $2$ de $FC(\mathfrak{g}(F))$ que l'on note  $FC_{4}$. Cela d\'emontre \ref{resultats}(1). 

On pose ${\cal Y}={\cal X}$ et on note $\varphi:{\cal X}\to {\cal Y}$ l'identit\'e. On montre comme en \ref{E6nonquasidep} que l'espace $FC^{{\cal E}}(\mathfrak{g}(F))$ est r\'eduit \`a $FC^{st}(\mathfrak{g}^*(F))$. Si $\delta_{4}(q-1)=1$, on pose $FC^{{\cal E}}_{4}=FC^{st}(\mathfrak{g}^*(F))$   et les assertions (2) et (3) de \ref{resultats} sont triviales.

\subsection{Le type $E_{8}$}
On suppose $G$ d\'eploy\'e de type $E_{8}$.
 
 On pose ${\cal X}_{0}=\{0,1,8\}$. Si $\delta_{3}(q-1)=0$, resp. $\delta_{4}(q-1)=0$, resp. $\delta_{5}(q-1)=0$, on pose ${\cal X}_{3}=\emptyset$, resp. ${\cal X}_{4}=\emptyset$, resp. ${\cal X}_{5}=\emptyset$. Si $\delta_{3}(q-1)=1$, resp. $\delta_{4}(q-1)=1$, resp. $\delta_{5}(q-1)=1$, on pose ${\cal X}_{3}=\{4,7\}$, resp. ${\cal X}_{4}=\{6\}$, resp. ${\cal X}_{5}=\{5\}$. On pose ${\cal X}={\cal X}_{0}\cup {\cal X}_{3}\cup{\cal X}_{4}\cup{\cal X}_{5}$, $d_{x}=1$ pour $x\in {\cal X}_{0}$, $d_{x}=2$ pour $x\in {\cal X}_{3}\cup {\cal X}_{4}$, $d_{x}=4$ pour $x\in {\cal X}_{5}$.

On a $G_{AD}=G$ donc l'action de $G_{AD}(F)$ sur $\underline{S}(G)$ est triviale. Cet ensemble   est en bijection avec l'ensemble des racines du diagramme ${\cal D}_{a}$.  On calcule ais\'ement les divers groupes $G_{s}$ en se rappelant l'\'egalit\'e 
$$\check{\alpha}_{0}+2\check{\alpha}_{1}+3\check{\alpha}_{2}+4\check{\alpha}_{3}+6\check{\alpha}_{4}+5\check{\alpha}_{5}+4\check{\alpha}_{6}+3\check{\alpha}_{7}+2\check{\alpha}_{8}=0.$$ 
 Pour le sommet $s$ param\'etr\'e par $\alpha_{0}$, $G_{s}$ est de type $E_{8}$ donc $dim(FC(\mathfrak{g}_{s}({\mathbb F}_{q})))=1$ d'apr\`es \ref{E8}. Il s'en d\'eduit une droite dans $FC(\mathfrak{g}(F))$ que l'on note $FC_{0}$. Pour le sommet $s$ param\'etr\'e par $\alpha_{1}$,  $G_{s}\simeq Spin_{dep}(16)/\{1,z'\}$, avec les notations de \ref{Dndeppair}.  Donc $dim(FC(\mathfrak{g}_{s}({\mathbb F}_{q})))=1$ d'apr\`es \ref{Dndeppair}.  Il s'en d\'eduit une droite dans $FC(\mathfrak{g}(F))$ que l'on note $FC_{1}$. Pour le sommet $s$ param\'etr\'e par $\alpha_{2}$, $G_{s}=SL(9)/\boldsymbol{\zeta}_{3}(\bar{{\mathbb F}}_{q})$, donc  $dim(FC(\mathfrak{g}_{s}({\mathbb F}_{q})))=0$ d'apr\`es \ref{An-1dep}. Pour le sommet $s$ param\'etr\'e par $\alpha_{3}$, $G_{s}=(SL(8)\times SL(2))/\boldsymbol{\zeta}_{4}(\bar{{\mathbb F}}_{q})$, un \'el\'ement $\zeta\in \boldsymbol{\zeta}_{4}(\bar{{\mathbb F}}_{q})$ d'ordre $4$ s'envoyant sur un \'el\'ement central d'ordre $4$ de $SL(8)$ et sur l'\'el\'ement central d'ordre $2$ de $SL(2)$. Donc $dim(FC(\mathfrak{g}_{s}({\mathbb F}_{q})))=0$ d'apr\`es \ref{An-1dep}. Pour le sommet $s$ param\'etr\'e par $\alpha_{4}$, on a $G_{s}=(SL(2)\times SL(3)\times SL(6))/\boldsymbol{\zeta}_{6}(\bar{{\mathbb F}}_{q})$, un \'el\'ement $\zeta\in \boldsymbol{\zeta}_{6}(\bar{{\mathbb F}}_{q})$ d'ordre $6$ s'envoyant sur un \'el\'ement central d'ordre $6$ de $SL(6)$, un \'el\'ement central d'ordre $3$ de $SL(3)$ et l'\'el\'ement central d'ordre $2$ de $SL(2)$. Si $\delta_{3}(q-1)=0$, les facteurs $SL(6)$ et $SL(3)$ ne portent pas de fonctions de type $f_{N,\epsilon}$ donc $dim(FC(\mathfrak{g}_{s}({\mathbb F}_{q})))=0$. Supposons $\delta_{3}(q-1)=1$. Chaque facteur $SL(6)$ et $SL(3)$ porte deux telles fonctions et le facteur $SL(2) $ en porte une. On obtient $4$ telles fonctions produits tensoriels mais seulement deux se factorisent par  $\boldsymbol{\zeta}_{6}(\bar{{\mathbb F}}_{q})$. D'o\`u $dim(FC(\mathfrak{g}_{s}({\mathbb F}_{q})))=2$.  Il s'en d\'eduit un plan dans $FC(\mathfrak{g}(F))$ que l'on note $FC_{4}$. Pour le sommet $s$ param\'etr\'e par $\alpha_{5}$, on a $G_{s}=(SL(5)\times SL(5))/\boldsymbol{\zeta}_{5}(\bar{{\mathbb F}}_{q})$, un \'el\'ement $\zeta\in \boldsymbol{\zeta}_{5}(\bar{{\mathbb F}}_{q})$ d'ordre $5$ s'envoyant sur le produit de deux \'el\'ements centraux d'ordre $5$ de chaque facteur $SL(5)$. Le groupe $SL(5)$ porte des fonctions $f_{N,\epsilon}$ si et seulement si $\delta_{5}(q-1)=1$. Supposons cette condition v\'erifi\'ee. Le groupe $SL(5)$   porte alors $4$ telles fonctions. On obtient $16$ fonctions produits tensoriels mais seulement $4$ d'entre elles se factorisent par $\boldsymbol{\zeta}_{5}(\bar{{\mathbb F}}_{q})$. D'o\`u $dim(FC(\mathfrak{g}_{s}({\mathbb F}_{q})))=4$. Il s'en d\'eduit un sous-espace de dimension $4$  de$FC(\mathfrak{g}(F))$ que l'on note $FC_{5}$. Pour le sommet $s$ param\'etr\'e par $\alpha_{6}$, on a $G_{s}=(Spin_{dep}(10)\times SL(4))/\boldsymbol{\zeta}_{4}(\bar{{\mathbb F}}_{q})$, un  \'el\'ement d'ordre $4$ de $\boldsymbol{\zeta}_{4}(\bar{{\mathbb F}}_{q})$ s'envoyant sur le produit de deux \'el\'ements centraux  d'ordre $4$ de $SL(4)$ et $Spin_{dep}(10)$. Si $\delta_{4}(q-1)=0$, les facteurs $SL(4)$ et $Spin_{dep}(10)$ ne portent pas de fonctions $f_{N,\epsilon}$. Supposons $\delta_{4}(q-1)=1$. Chacun de ces facteurs porte $2$ telles fonctions. On obtient $4$ fonctions produits tensoriels mais seulement $2$ d'entre elles se factorisent par $\boldsymbol{\zeta}_{4}(\bar{{\mathbb F}}_{q})$. D'o\`u $dim(FC(\mathfrak{g}_{s}({\mathbb F}_{q})))=2$. Il s'en d\'eduit un plan  dans $FC(\mathfrak{g}(F))$ que l'on note $FC_{6}$. Pour le sommet $s$ param\'etr\'e par $\alpha_{7}$, on  a $G_{s}=(G_{E_{6},SC}\times SL(3))/\boldsymbol{\zeta}_{3}(\bar{{\mathbb F}}_{q})$. On a not\'e $G_{E_{6},SC}$ le groupe d\'eploy\'e simplement connexe de type $E_{6}$.  Un  \'el\'ement d'ordre $3$ de $\boldsymbol{\zeta}_{3}(\bar{{\mathbb F}}_{q})$ s'envoie sur le produit de deux \'el\'ements centraux  d'ordre $3$ de $SL(3)$ et $G_{E_{6},SC}$. Par un m\^eme calcul que ci-dessus et gr\^ace \`a \ref{E6dep}, $dim(FC(\mathfrak{g}_{s}({\mathbb F}_{q})))=2\delta_{3}(q-1)$. Si $\delta_{3}(q-1)=1$, il s'en d\'eduit un plan  dans $FC(\mathfrak{g}(F))$ que l'on note $FC_{7}$. Pour le sommet $s$ param\'etr\'e par $\alpha_{8}$, on a $G_{s}=(G_{E_{7},SC})\times SL(2))/\{\pm 1\}$. On a not\'e $G_{E_{7},SC}$ le groupe  simplement connexe de type $E_{7}$. 
Le groupe $\{\pm 1\}$ s'envoyant diagonalement dans le produit des centres des deux facteurs. Gr\^ace \`a \ref{E7}, on  obtient $dim(FC(\mathfrak{g}_{s}({\mathbb F}_{q})))=1$.  Il s'en d\'eduit une droite dans $FC(\mathfrak{g}(F))$ que l'on note $FC_{8}$.  Cela d\'emontre  \ref{resultats} (1). 

On pose ${\cal Y}={\cal X}$, ${\cal X}^{st}={\cal X}$ et on note $\varphi:{\cal X}\to {\cal Y}$ l'identit\'e. Pour prouver les assertions (2), (3) et (4) de \ref{resultats}, il suffit de prouver que

(1) $FC(\mathfrak{g}(F))=FC^{st}(\mathfrak{g}(F))$.

On pose alors $FC^{{\cal E}}_{y}=FC_{y}$ pour tout $y\in {\cal Y}={\cal X}$ et les assertions \`a d\'emontrer sont triviales. 

Prouvons (1).  Consid\'erons un couple $(\sigma\mapsto \sigma_{G'},{\cal O})\in {\cal E}_{ell}(G)$. Le groupe $\hat{\Omega}$ est trivial et l'action galoisienne sur $G$ aussi, donc l'action $\sigma\mapsto \sigma_{G'}$ est triviale. Pour ${\cal O}=\{\hat{\alpha}_{0}\}$, ${\bf G}'={\bf G}$. Pour ${\cal O}=\{\hat{\alpha}_{1}\}$, $G'_{SC}=Spin_{dep}(16)$ donc $FC^{st}(\mathfrak{g}'(F))=\{0\}$ d'apr\`es \ref{Dndeppairpadique} (5). Pour ${\cal O}=\{\hat{\alpha}_{i}\}$ avec $i\in \{2,...,8\}$, $G'_{SC}$ contient un facteur $SL(m)$ avec $m\geq2$. Donc encore $FC^{st}(\mathfrak{g}'(F))=\{0\}$ d'apr\`es \ref{An-1deppadique} (1).  Donc $FC^{{\cal E}}(\mathfrak{g}(F))$ est r\'eduit \`a $FC^{st}(\mathfrak{g}(F))$, ce qui prouve (1).

 Explicitons la cons\'equence de \ref{resultats}(4):
  
(2)   $dim(FC^{st}(\mathfrak{g}(F))=3+4\delta_{3}(q-1)+2\delta_{4}(q-1)+4\delta_{5}(q-1)$.

\subsection{Le type $F_{4}$}\label{F4padique}
    On suppose $G$ de type $F_{4}$. Donc $G$ est d\'eploy\'e et $Z(G)=\{1\}$. Le groupe dual $\hat{G}$ est aussi de type $F_{4}$. On utilise les notations de Bourbaki pour les diagrammes des deux groupes mais mais on prendra garde que la correspondance naturelle entre les racines des deux groupes n'est pas compatible avec la num\'erotation de ces racines. Cela vient du fait que cette correspondance naturelle provient d'un isomorphisme $X^*(T)\simeq X_{*}(\hat{T})$ et \'echange en fait racines et coracines. Or une racine est longue si et seulement sa coracine est courte. Avec les notations de Bourbaki, la correspondance est donc $\alpha_{i}\mapsto \hat{\alpha}_{5-i}$ pour $i=1,...,4$.

 On pose ${\cal X}_{0}=\{0,1,4\}$. Si $\delta_{3}(q-1)=0$, resp. $\delta_{4}(q-1)=0$, on pose ${\cal X}_{3}=\emptyset$, resp. ${\cal X}_{4}=\emptyset$. Si $\delta_{3}(q-1)=1$, resp. $\delta_{4}(q-1)=1$, on pose ${\cal X}_{3}=\{2\}$, resp. ${\cal X}_{4}=\{3\}$. On pose ${\cal X}={\cal X}_{0}\cup {\cal X}_{3}\cup {\cal X}_{4}$, $d_{x}=1$ pour $x\in {\cal X}_{0}$, $d_{x}=2$ pour $x\in {\cal X}_{3}\cup {\cal X}_{4}$. 
   
   On a $G_{AD}=G$ donc l'action de $G_{AD}(F)$ sur $S(G)$ est triviale. L'ensemble $\underline{S}(G)$ est en bijection avec l'ensemble des racines du diagramme ${\cal D}_{a}$. On calcule les groupes $G_{s}$ en se rappelant la relation
   $$\check{\alpha}_{0}+2\check{\alpha}_{1}+3\check{\alpha}_{2}+2\check{\alpha}_{3}+\check{\alpha}_{4}=0.$$
     Pour le sommet $s$ correspondant \`a $\alpha_{0}$, le groupe $G_{s}$ est de type $F_{4}$ donc $FC(\mathfrak{g}_{s}({\mathbb F}_{q}))$ est une droite d'apr\`es \ref{F4}. Il s'en d\'eduit une droite dans $FC(\mathfrak{g}(F))$ que l'on note $FC_{0}$.  Pour le sommet $s$ correspondant \`a $\alpha_{1}$, le groupe $G_{s}$ est isomorphe \`a $(SL(2)\times Sp(6))/\{\pm 1\}$, o\`u $\{\pm 1\}$ s'identifie aux centres des deux composantes. Alors $FC(\mathfrak{g}_{s}({\mathbb F}_{q}))$ est une droite d'apr\`es \ref{Cn}. Il s'en d\'eduit une droite dans $FC(\mathfrak{g}(F))$ que l'on note $FC_{1}$. 
   Pour le sommet $s$ correspondant \`a $\alpha_{2}$, le groupe $G_{s}$ est isomorphe \`a $(SL(3)\times SL(3))/\boldsymbol{\zeta}_{3}(\bar{{\mathbb F}}_{q})$, o\`u $\boldsymbol{\zeta}_{3}(\bar{{\mathbb F}}_{q})$ s'identifie aux centres des deux composantes. Si $\delta_{3}(q-1)=0$,  on a $FC(\mathfrak{g}_{s}({\mathbb F}_{q}))=\{0\}$.  Supposons $\delta_{3}(q-1)=1$. Il y a deux fonctions du type $f_{N,\epsilon}$ sur chaque composante $SL(3)$ donc $4$ fonctions produits tensoriels mais seulement deux d'entre elles  se quotientent par $\boldsymbol{\zeta}_{3}(\bar{{\mathbb F}}_{q})$. D'o\`u $dim(FC(\mathfrak{g}_{s}({\mathbb F}_{q})))=2 $.   Il s'en d\'eduit un plan dans $FC(\mathfrak{g}(F))$ que l'on note $FC_{2}$. Pour le sommet $s$ correspondant \`a $\alpha_{3}$, le groupe $G_{s}$ est isomorphe \`a $(SL(4)\times SL(2))/\{\pm 1\}$, le groupe $\{\pm 1\}$ s'identifiant au centre de $SL(2)$ et au sous-groupe d'ordre $2$  du centre de $SL(4)$. L'espace $\mathfrak{sl}(2)({\mathbb F}_{q})$ porte une fonction du type $f_{N,\epsilon}$. Si $\delta_{4}(q-1)=0$, l'espace $\mathfrak{sl}(4)({\mathbb F}_{q})$ n'en porte pas. Supposons $\delta_{4}(q-1)=1$. Alors cet espace porte $2$ telles fonctions   et le produit tensoriel  de chacune d'elles avec la fonction sur $\mathfrak{sl}(2)({\mathbb F}_{q})$ se quotiente par $\{\pm 1\}$. Donc $dim(FC(\mathfrak{g}_{s}({\mathbb F}_{q})))=2$.  Il s'en d\'eduit un plan dans $FC(\mathfrak{g}(F))$ que l'on note $FC_{3}$. Pour le sommet $s$ correspondant \`a $\alpha_{4}$, le groupe $G_{s}$ est isomorphe \`a $Spin_{dep}(9)$, donc $FC(\mathfrak{g}_{s}({\mathbb F}_{q}))$ est une droite d'apr\`es \ref{Bn}.  Il s'en d\'eduit une droite dans $FC(\mathfrak{g}(F))$ que l'on note $FC_{4}$.Cela ach\`eve la preuve de  \ref{resultats}(1).
   
   On pose ${\cal Y}={\cal X}$, ${\cal X}^{st}={\cal X}$ et on note $\varphi:{\cal X}\to {\cal Y}$ l'identit\'e. Pour prouver les assertions (2), (3) et (4) de \ref{resultats}, il suffit de prouver que

(1) $FC(\mathfrak{g}(F))=FC^{st}(\mathfrak{g}(F))$.

On pose alors $FC^{{\cal E}}_{y}=FC_{y}$ pour tout $y\in {\cal Y}={\cal X}$ et les assertions \`a d\'emontrer sont triviales. 

Prouvons (1).  Consid\'erons un couple $(\sigma\mapsto \sigma_{G'},{\cal O})\in {\cal E}_{ell}(G)$. Le groupe $\hat{\Omega}$ est trivial et l'action galoisienne sur $G$ aussi, donc l'action $\sigma\mapsto \sigma_{G'}$ est triviale.  Pour ${\cal O}=\{\hat{\alpha}_{0}\}$, on a ${\bf G}'={\bf G}$.  Pour ${\cal O}=\{\hat{\alpha}_{1}\}$, resp. $\{\hat{\alpha}_{2}\}$, $\{\hat{\alpha}_{3}\}$, le groupe $G'_{SC}$ contient un facteur $SL(2)$, resp. $SL(3)$, $SL(4)$, donc $FC^{st}(\mathfrak{g}'(F))=\{0\}$ d'apr\`es \ref{An-1deppadique}. Pour ${\cal O}=\{\hat{\alpha}_{4}\}$, le groupe $G'_{SC}$ est $Sp(8)$ et encore $FC^{st}(\mathfrak{g}'(F))=\{0\}$ d'apr\`es \ref{Dndeppairpadique}. Cela d\'emontre que $FC^{{\cal E}}(\mathfrak{g}(F))=FC^{st}(\mathfrak{g}(F))$, d'o\`u (1).  
   
    Explicitons la cons\'equence de \ref{resultats}(4):
      
  (2) $dim(FC^{st}(\mathfrak{g}(F)))=3+2\delta_{3}(q-1)+2\delta_{4}(q-1)$.

     \subsection{Le type $G_{2}$}
  On suppose $G$ de type $G_{2}$. Donc $G$ est d\'eploy\'e et $Z(G)=\{1\}$. Le groupe dual $\hat{G}$ est aussi de type $G_{2}$. On utilise les notations de Bourbaki pour les diagrammes des deux groupes mais, comme en \ref{F4padique},  on prendra garde que la correspondance naturelle entre les racines des deux groupes n'est pas compatible avec la num\'erotation de ces racines. 
  
  On pose ${\cal X}_{0}=\{0,2\}$. Si $\delta_{3}(q-1)=0$, on pose ${\cal X}_{3}=\emptyset$. Si $\delta_{3}(q-1)=1$, on pose ${\cal X}_{3}=\{1\}$. On pose ${\cal X}={\cal X}_{0}\cup{\cal X}_{3}$, $d_{x}=1$ pour $x\in {\cal X}_{0}$ et $d_{x}=2$ pour $x\in {\cal X}_{3}$.
      
   On a $G_{AD}=G$ donc l'action de $G_{AD}(F)$ sur $S(G)$ est triviale. L'ensemble $\underline{S}(G)$ est en bijection avec l'ensemble des racines du diagramme ${\cal D}_{a}$. Pour le sommet $s$ correspondant \`a $\alpha_{0}$, le groupe $G_{s}$ est de type $G_{2}$ donc $FC(\mathfrak{g}_{s}({\mathbb F}_{q}))$ est une droite d'apr\`es \ref{G2}. Il s'en d\'eduit une droite dans $FC(\mathfrak{g}(F))$ que l'on note $FC_{0}$.  Pour le sommet $s$ correspondant \`a $\alpha_{1}$, on a $G_{s}\simeq SL(3)$. D'o\`u $dim(FC(\mathfrak{g}_{s}({\mathbb F}_{q})))=2\delta_{3}(q-1)$.  Si $\delta_{3}(q-1)=1$, il s'en d\'eduit un plan dans $FC(\mathfrak{g}(F))$ que l'on note $FC_{1}$. Pour le sommet $s$ correspondant \`a $\alpha_{2}$, on a $G_{s}\simeq (SL(2)\times SL(2))/\{\pm 1\}$, donc $FC(\mathfrak{g}_{s}({\mathbb F}_{q}))$ est une droite. Il s'en d\'eduit une droite dans $FC(\mathfrak{g}(F))$ que l'on note $FC_{2}$. Cela d\'emontre \ref{resultats}(1). 
   
      On pose ${\cal Y}={\cal X}$, ${\cal X}^{st}={\cal X}$ et on note $\varphi:{\cal X}\to {\cal Y}$ l'identit\'e. Pour prouver les assertions (2), (3) et (4) de \ref{resultats}, il suffit de prouver que

(1) $FC(\mathfrak{g}(F))=FC^{st}(\mathfrak{g}(F))$.

On pose alors $FC^{{\cal E}}_{y}=FC_{y}$ pour tout $y\in {\cal Y}={\cal X}$ et les assertions \`a d\'emontrer sont triviales. 
   
       Consid\'erons un couple  $(\sigma\mapsto \sigma_{G'},{\cal O})\in {\cal E}_{ell}(G)$.    Ici, le groupe $\hat{\Omega}$ est  trivial et l'action $\sigma\mapsto \sigma_{G}$ aussi.  Donc l'action $\sigma\mapsto \sigma_{G'}$ est triviale. 
   Pour ${\cal O}=\{\hat{\alpha}_{0}\}$, on a ${\bf G}'={\bf G}$. Pour ${\cal O}=\{\hat{\alpha}_{1}\}$, resp. $\{\hat{\alpha}_{2}\}$, on voit que $G'_{SC}$  est $SL(3)$, resp. $SL(2)\times SL(2)$. D'apr\`es  \ref{An-1deppadique},  $FC^{st}(\mathfrak{g}'(F))=\{0\}$. Cela d\'emontre que $FC^{{\cal E}}(\mathfrak{g}(F))=FC^{st}(\mathfrak{g}(F))$, d'o\`u (1).

 Explicitons la cons\'equence de \ref{resultats}(4):
  
  (2)   $dim(FC^{st}(\mathfrak{g}(F)))=2(1+\delta_{3}(q-1))$.

\bigskip

jean-loup.waldspurger@imj-prg.fr

Institut de math\'ematiques de Jussieu-Paris rive gauche

4 place Jussieu

Bo\^{\i}te courrier 247

75252 Paris cedex 05

\end{document}